\documentclass{article} 
\usepackage{pdfpages}
\usepackage{geometry}
\usepackage{amsmath}  
\usepackage{graphicx} 
\usepackage{float}
\usepackage{color}
\usepackage{amsmath}
\usepackage{amssymb}
\usepackage{amsthm}

\newtheorem{theorem}{Theorem}

\theoremstyle{definition}
\newtheorem{definition}{Definition}[section]

\theoremstyle{lemma}
\newtheorem{lemma}[theorem]{Lemma}

\theoremstyle{remark}

\theoremstyle{boldremark}
\newtheorem*{boldremark}{Remark} 

\newtheorem{corollary}[theorem]{Corollary}

\newcommand{\maketitletwo}[2][]{\begin{center}
        \Large{\textbf{Solving Moving Sofa Problem Using Calculus of Variations}
            
            } 
        \vspace{5pt}
        
        \normalsize{Zhipeng Deng

        }       
        \vspace{15pt}
        
\end{center}}
\begin{document}
    \maketitletwo[5]

    \begin{abstract}
In 1966, Leo Moser introduced the "moving sofa problem," which seeks to determine the largest area of a shape that can be maneuvered through a 90-degree hallway of unit-width. This problem remains unsolved and open yet. In this paper, we employ calculus of variations method to solve this problem. Assuming the trajectories and envelopes are convex, the sofa's area is formulated as an integral functional on a set of parametric equations for curves. The final shape is determined by solving the Euler-Lagrange equations. Utilizing numerical methods, we obtain the non-trivial area of 2.2195316, consistent with the previously well-known Gerver's constant since 1992. We prove that both the results of Gerver's sofa and Romik's car satisfy the Euler-Lagrange equations for the necessary condition of maximal area. We also explore additional cases and asymmetric conditions, and discuss other variant problems.
\end{abstract}

\tableofcontents

\section{Introduction}

In 1966, Leo Moser first published the moving sofa problem \cite{Moser1966} \cite{wiki}. It remains an unsolved open problem in mathematics asking "what is the largest area of a shape that can be maneuvered through an L-shaped hallway of unit width?" 

In 1968, John Hammersley established a lower bound for the sofa area with $\frac{\pi }{2} +\frac{2}{\pi}$ \cite{Hammersley1968}. Then the previous book \cite{Croft2012} summarized the optimization of area results ranging from 2.2074 to 2.215649 prior to 1992. In 1992, Joseph Gerver at Rutgers University further improved the lower bound for the sofa constant to approximately 2.21953 \cite{Gerver1992}. He used geometry-based methods and Balanced Polygon condition to construct a shape. The boundary consisted of analytic pieces including lines, circular arcs, involutes of circles, and involutes of involutes of circles. Although this shape was considered optimal or near-optimal, no proof of it being the largest area was provided. Since then, no larger sofa area has been discovered \cite{mathworld}. In 2018, Romik further explored this problem using differential equations and explicit formulas, also investigating variant car problems \cite{Romik2018}. Despite some numerical calculation results \cite{Tyrrell2015} \cite{Batsch2022}, the moving sofa problem remains open and unsolved in geometry \cite{Prieto2021} \cite{Croft2012}. 

As for the upper bound of the sofa area, some researchers have made attempts as well. Hammersley proposed an upper bound of $2\sqrt{2}$ \cite{Hammersley1968}. Kallus and Romik established a new upper bound in 2018, capping the sofa constant at 2.37 \cite{Kallus2018}. 

The contribution of this paper is that we use calculus of variations method \cite{Elsgolc2012}. The sofa shape is defined as the envelope curve formed by the walls it navigates. The sofa’s area is formulated as a functional on a set of parametric equations for curves. The final shape is determined by solving the Euler-Lagrange equations. We achieve an area of 2.21953, matching Gerver’s constant from 1992 \cite{Gerver1992}. This verifies that Gerver's constructed non-trivial sofa shape \cite{Gerver1992} aligns with our results. This is the first study using the calculus of variations to solve the moving sofa problem, proving the necessity of maximum area rather than relying on geometric construction.

\section{Proof of concept}

Figure 1 illustrates the proof of concept in this paper. Our success lies in representing the shape and area corresponding to the moving path as an integral functional, then using the calculus of variations to prove optimality. Previous researchers did not recognize the power and inspiration of this approach. It is crucial because this is the authentic method to prove a maximum value of the functional, especially in an optimization problem in a space of continuous functions.

\begin{figure}[H]
  \centering
  \includegraphics[width=12cm]{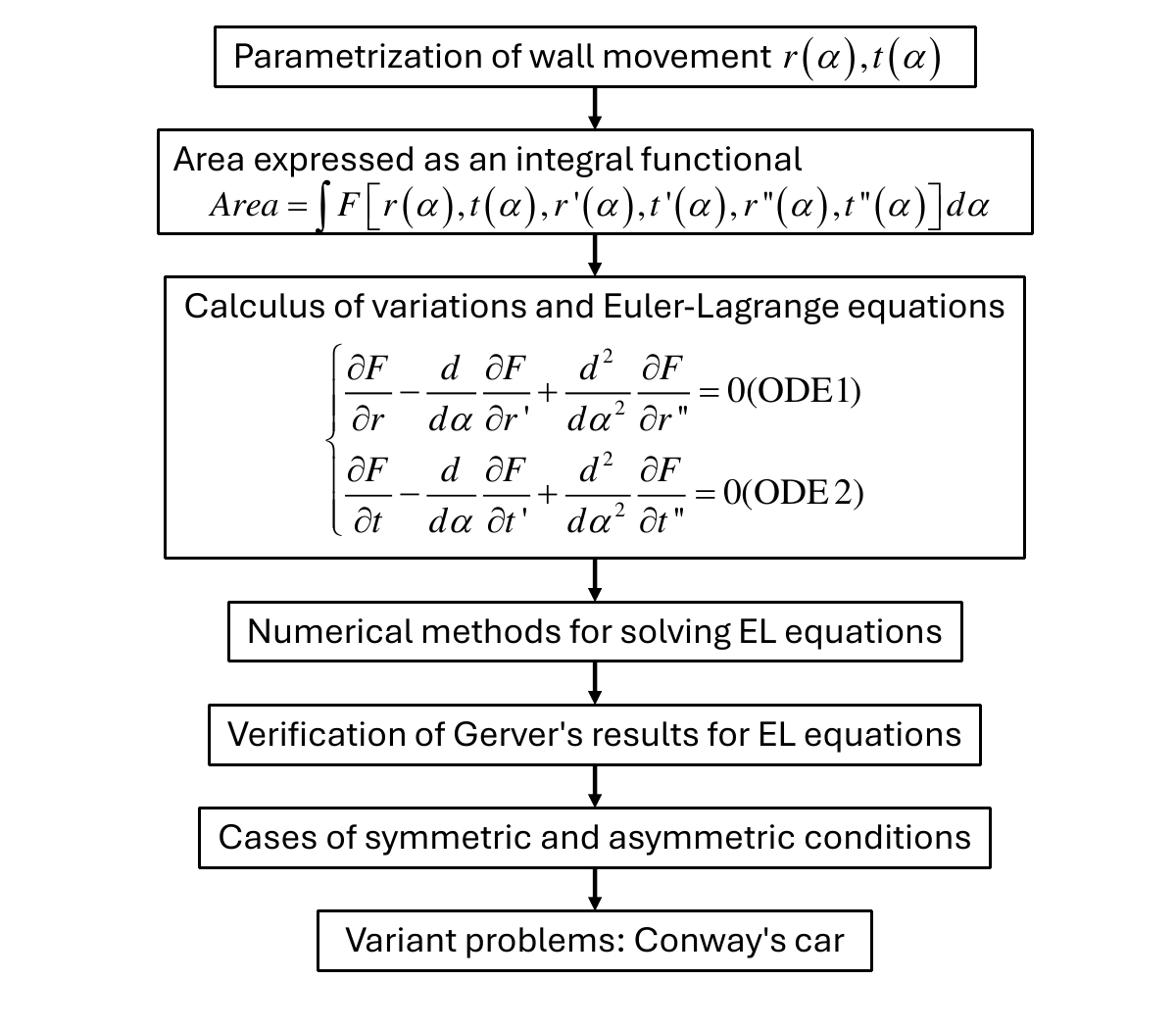}
  \caption{Proof of concept}
\end{figure}

\section{Notations and definitions}

\begin{definition}[coordinate system and path]
To formulate the problem, we first create a 2D coordinate system as shown in Figure 2. We assume the sofa shape rotates and moves clockwise, equivalent to the hallway wall rotating counterclockwise. Thus, lines AB, AC, ApBp, and ApCp rotate counterclockwise.

\begin{figure}[H]
  \centering
  \includegraphics[width=12cm]{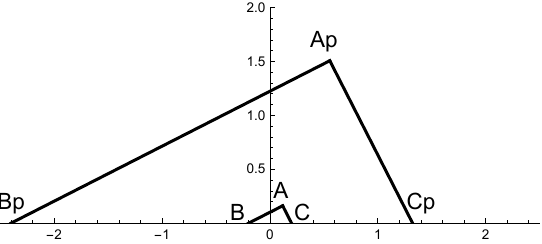}
  \caption{Points and lines in 2D coordinate system}
\end{figure}

Moving sofa in a hallway with unit width, we assume the hallway has a lower edge along the X-axis $y=0$, and the upper edge along the line $y=1$. 

We employ parametric equations with two functions $r(\alpha)$ and $t(\alpha)$ to define a family of curves to be determined, where $\alpha $ represents the rotation angle with 
$0\le \alpha\le \pi$. These curves are assumed to be differentiable and well-behaved \cite{Romik2018}.

We will show that the area of the sofa can be regarded as a functional of a set of parametric curve equations.

\end{definition}

\begin{definition}[points]
The coordinates of the points in Figure 2 are defined as follows:

\begin{equation*}
\left \{ A_{x},A_{y} \right \} =\left \{ r\left ( \alpha  \right )\cos\left ( \alpha  \right ),t\left ( \alpha  \right )\sin\left ( \alpha  \right ) \right \}
\end{equation*}

\begin{equation*}
\left \{ Ap_{x},Ap_{y} \right \} =\left\{\sqrt{2} \cos \left(\frac{\alpha }{2}+\frac{\pi }{4}\right)+\cos (\alpha ) r(\alpha ),\sqrt{2} \sin \left(\frac{\alpha }{2}+\frac{\pi
   }{4}\right)+\sin (\alpha ) t(\alpha )\right\}
\end{equation*}

\begin{equation*}
\left \{ B_{x},B_{y} \right \} =\left\{\cos (\alpha ) r(\alpha )-2 \cos ^2\left(\frac{\alpha }{2}\right) t(\alpha ),0\right\}
\end{equation*}

\begin{equation*}
\left \{ Bp_{x},Bp_{y} \right \} =\left\{-\csc \left(\frac{\alpha }{2}\right)+\cos (\alpha ) r(\alpha )-2 \cos ^2\left(\frac{\alpha }{2}\right)t(\alpha ),0\right\}
\end{equation*}

\begin{equation*}
\left \{ C_{x},C_{y} \right \} =\left\{\cos (\alpha ) r(\alpha )+\sin (\alpha ) \tan \left(\frac{\alpha }{2}\right) t(\alpha ),0\right\}
\end{equation*}

\begin{equation*}
\left \{ Cp_{x},Cp_{y} \right \} =\left\{\sec \left(\frac{\alpha }{2}\right)+\cos (\alpha ) r(\alpha )+\sin (\alpha ) \tan \left(\frac{\alpha }{2}\right) t(\alpha ),0\right\}
\end{equation*}

Point $A$ is the interior corner. Points $B$ and $C$ are the intersections of the interior walls with the X-axis. Point $Ap$ is the corner of the exterior wall. Points $Bp$ and $Cp$ are the intersections of the exterior walls with the X-axis.

\end{definition}

\begin{boldremark}
If $r\left ( \alpha  \right ) =0,t\left ( \alpha  \right ) =0 \left ( 0\le \alpha \le \pi  \right ) $ , it yields a semicircle that can easily pass through the hallway by rotating 90 degrees.
\end{boldremark}

\begin{boldremark}
If $r\left ( \alpha  \right ) =2/\pi,t\left ( \alpha  \right ) =2/\pi \left ( 0\le \alpha \le \pi  \right ) $ , it yields the Hammersley's sofa with an area of $\pi /2+2/\pi $.
\end{boldremark}

\begin{boldremark}
Some special points for the hallway
\begin{equation*}
A_{x}\bigg|_{\alpha =0}=r\left ( 0 \right ) 
\end{equation*}
\begin{equation*}
A_{x}\bigg|_{\alpha =\pi }=-r\left ( \pi \right ) 
\end{equation*}
\begin{equation*}
B_{x}\bigg|_{\alpha =0 }=r\left ( 0 \right ) -2t\left ( 0 \right ) 
\end{equation*}
\begin{equation*}
B_{x}\bigg|_{\alpha =\pi }=-r\left ( \pi \right ) 
\end{equation*}
\begin{equation*}
Ap_{x}\bigg|_{\alpha =0 }=1+r\left ( 0 \right ) 
\end{equation*}
\begin{equation*}
Ap_{x}\bigg|_{\alpha =\pi  }=-1-r\left ( \pi \right ) 
\end{equation*}
\end{boldremark}

\begin{definition}[hallway]
We define the hallway as two unit-width strips combined. For an intermediate state shown in Figure 2, the region inside the hallway can be expressed as follows
\begin{equation*}
\begin{split}
L\left ( \alpha  \right ) =&\left \{ (x, y) \in \mathbb{R}^2 \mid  y \leq \tan \frac{\alpha }{2}\left ( x-Ap_{x} \right ) +Ap_{y} \right \}\\
&\bigcap\left \{ (x, y) \in \mathbb{R}^2 \mid  y \leq -\cot \frac{\alpha }{2}\left ( x-Ap_{x} \right ) +Ap_{y} \right \}\\
&\bigcap\left \{ (x, y) \in \mathbb{R}^2 \mid \tan \frac{\alpha }{2}\left ( x-A_{x} \right ) +A_{y}\leq y  \right \}\\
&\bigcap\left \{ (x, y) \in \mathbb{R}^2 \mid -\cot \frac{\alpha }{2}\left ( x-A_{x} \right ) +A_{y}\leq y  \right \}\\
\end{split} 
\end{equation*}
\end{definition}

\begin{corollary}
The shape of the sofa that can move through the hallway must satisfy
\begin{equation*}
S\subseteq \bigcap_{0\le \alpha \le \pi }^{} L\left ( \alpha  \right ) 
\end{equation*}
\end{corollary}

\begin{boldremark}
Before hallway movement, 
\begin{equation*}
L\left ( 0 \right ) =\left \{ (x, y) \in \mathbb{R}^2 \mid x\le Ap_{x}\bigg|_{\alpha =0},0\le y\le 1 \right \}\bigcup \left \{ (x, y) \in \mathbb{R}^2 \mid A_{x}\bigg|_{\alpha =0} \le x\le Ap_{x}\bigg|_{\alpha =0},y\le 1 \right \} 
\end{equation*}
\end{boldremark}

\begin{boldremark}
After hallway movement, 
\begin{equation*}
L\left ( \pi  \right ) =\left \{ (x, y) \in \mathbb{R}^2 \mid Ap_{x}\bigg|_{\alpha =\pi }\le x ,0\le y\le 1 \right \}\bigcup \left \{ (x, y) \in \mathbb{R}^2 \mid Ap_{x}\bigg|_{\alpha =\pi } \le x\le A_{x}\bigg|_{\alpha =0},y\le 1 \right \} 
\end{equation*}
\end{boldremark}

\begin{definition}[parametric equations of envelope ]
From Figure 2 and Definition 3.2, we derive the explicit horizontal and vertical coordinates of the parametric equations of envelope of the hallway walls by $L_{AB}=\frac{y-A_{y}}{B_{y}-A_{y}}- \frac{x-A_{x}}{B_{x}-A_{x}}=0$ and $\frac{\partial L_{AB}}{\partial \alpha }=0$, $L_{AC}=0$ and $\frac{\partial L_{AC}}{\partial \alpha }=0$, $L_{ApBp}=0$ and $\frac{\partial L_{ApBp}}{\partial \alpha }=0$, $L_{ApCp}=0$ and $\frac{\partial L_{ApCp}}{\partial \alpha }=0$ as

\begin{equation*}
\begin{split}
E_{ABx} =&\frac{1}{2} (-1+2 \cos (\alpha )+\cos (2 \alpha )) r(\alpha )-2 \cos ^2\left(\frac{\alpha }{2}\right) \cos (\alpha ) t(\alpha )+\sin (\alpha )
   \left(\cos (\alpha ) r'(\alpha )\right.\\
   &\left.-(1+\cos (\alpha )) t'(\alpha )\right)
\end{split}
\end{equation*}
\begin{equation*}
E_{ABy} =-2 \sin ^2\left(\frac{a}{2}\right) \left(r(\alpha ) \sin (\alpha )-\sin (\alpha ) t(\alpha )-\cos (\alpha ) r'(\alpha )+t'(\alpha )+\cos (\alpha
   ) t'(\alpha )\right)
\end{equation*}
\begin{equation*}
\begin{split}
E_{ACx} =&r(\alpha ) \left(\cos (\alpha )+\sin ^2(\alpha )\right)-2 \cos (\alpha ) \sin ^2\left(\frac{\alpha }{2}\right) t(\alpha )+\sin (\alpha )
   \left(-\cos (\alpha ) r'(\alpha )\right.\\
   &\left. +(-1+\cos (\alpha )) t'(\alpha )\right)
\end{split}
\end{equation*}
\begin{equation*}
E_{ACy} =2 \cos ^2\left(\frac{\alpha }{2}\right) \left(-r(\alpha ) \sin (\alpha )+\sin (\alpha ) t(\alpha )+\cos (\alpha ) r'(\alpha )+t'(\alpha )-\cos
   (\alpha ) t'(\alpha )\right)
\end{equation*}
\begin{equation*}
\begin{split}
E_{ApBpx} =&\frac{1}{2} \left[(-1+2 \cos (\alpha )+\cos (2 \alpha )) r(\alpha )-2 \sin \left(\frac{\alpha }{2}\right)-4 \cos ^2\left(\frac{\alpha }{2}\right)\cos (\alpha ) t(\alpha )\right.\\
   & \left. +\sin (2 \alpha ) r'(\alpha )-2 \sin (\alpha ) t'(\alpha )-\sin (2 \alpha ) t'(\alpha )\right]
\end{split}
\end{equation*}
\begin{equation*}
\begin{split}
E_{ApBpy} =&\frac{1}{2} \left[2 \cos \left(\frac{\alpha }{2}\right)+r(\alpha ) (-2 \sin (\alpha )+\sin (2 \alpha ))-2 (-1+\cos (\alpha )) \sin (\alpha )t(\alpha )\right.\\
   &\left. -r'(\alpha )+2 \cos (\alpha ) r'(\alpha )-\cos (2 a) r'(\alpha )-t'(\alpha )+\cos (2 a) t'(\alpha )\right]
\end{split}
\end{equation*}
\begin{equation*}
\begin{split}
E_{ApCpx} =&\frac{1}{2} \left[2 \cos \left(\frac{\alpha }{2}\right)+2 r(\alpha ) \left(\cos (\alpha )+\sin ^2(\alpha )\right)+(1-2 \cos (\alpha )+\cos (2 \alpha )) t(\alpha ) \right.\\
  &\left. -\sin (2 a) r'(\alpha )-2 \sin (\alpha ) t'(\alpha )+\sin (2 a) t'(\alpha )\right]
\end{split}
\end{equation*}
\begin{equation*}
\begin{split}
E_{ApCpy} =&\frac{1}{2} \left[2 \sin \left(\frac{\alpha }{2}\right)-2 (1+\cos (\alpha )) r(\alpha ) \sin (\alpha )+(2 \sin (\alpha )+\sin (2 \alpha ))t(\alpha )\right.\\
  &\left. +r'(\alpha )+2 \cos (\alpha ) r'(\alpha )+\cos (2 a) r'(\alpha )+t'(\alpha )-\cos (2 \alpha ) t'(\alpha )\right]
\end{split}
\end{equation*}

These envelopes constitute the possible boundaries of the moving sofa.
\end{definition}

\begin{boldremark}
Some special points for the envelopes
\begin{equation*}
E_{ApCpx}\bigg|_{\alpha =0}=1+r\left ( 0 \right ) 
\end{equation*}
\begin{equation*}
E_{ApCpx}\bigg|_{\alpha =\pi }=-r\left ( \pi \right )  + 2 t\left ( \pi \right )
\end{equation*}
\begin{equation*}
E_{ApCpy}\bigg|_{\alpha =0}=2r'\left ( 0 \right ) 
\end{equation*}
\begin{equation*}
E_{ApCpy}\bigg|_{\alpha =\pi}=1
\end{equation*}
\begin{equation*}
E_{ApBpx}\bigg|_{\alpha =0}=r\left ( 0 \right )  - 2 t\left ( 0 \right )
\end{equation*}
\begin{equation*}
E_{ApBpx}\bigg|_{\alpha =\pi }=-1 - r\left ( \pi \right ) 
\end{equation*}
\begin{equation*}
E_{ApBpy}\bigg|_{\alpha =0}=1
\end{equation*}
\begin{equation*}
E_{ApBpy}\bigg|_{\alpha =\pi}=-2 r'\left ( \pi  \right ) 
\end{equation*}
\end{boldremark}

\begin{definition}[shape of sofa]

\begin{figure}[H]
  \centering
  \includegraphics[width=14cm]{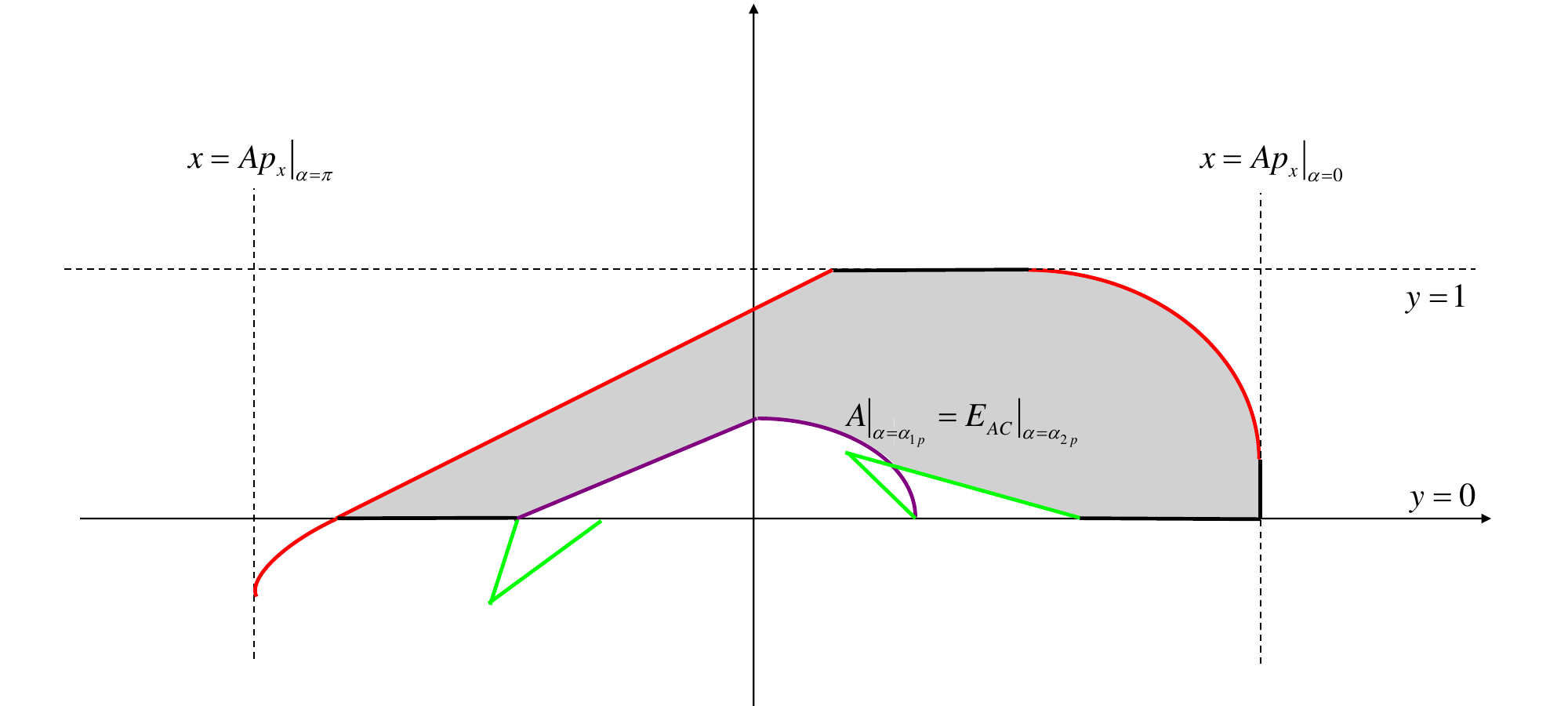}
  \caption{Shape of sofa with trajectories and envelopes.}
\end{figure}

Considering the envelopes as shown in Figure 3 and assuming the envelopes and trajectory are convex and do not self-intersect, the sofa's shape must lie above the line $y=0$, the trajectory of point A, the envelopes of AB and AC, and be below line $y=1$, envelopes of ApBp and ApCp, and between line $x=Ap_{x}\bigg|_{\alpha =0}$ and $x=Ap_{x}\bigg|_{\alpha =\pi }$ 

So that the maximum sofa region is
\begin{equation}
\begin{split}
S_{max}=&\bigcap_{0\le \alpha \le \pi }^{} \left \{ (x, y) \in \mathbb{R}^2 \mid  x=A_{x},y\ge A_{y} \right \}\\
&\bigcap \bigcap_{0\le \alpha \le \pi }^{} \left \{ (x, y) \in \mathbb{R}^2 \mid  x=E_{ABx},y\ge E_{ABy} \right \}\\
&\bigcap \bigcap_{0\le \alpha \le \pi }^{} \left \{ (x, y) \in \mathbb{R}^2 \mid  x=E_{ACx},y\ge E_{ACy} \right \}\\
&\bigcap \bigcap_{0\le \alpha \le \pi }^{} \left \{ (x, y) \in \mathbb{R}^2 \mid  x=E_{ApBpx},y\le E_{ApBpy} \right \}\\
&\bigcap \bigcap_{0\le \alpha \le \pi }^{} \left \{ (x, y) \in \mathbb{R}^2 \mid  x=E_{ApCpx},y\le E_{ApCpy} \right \}\\
&\bigcap \left \{ (x, y) \in \mathbb{R}^2 \mid  0\le y\le 1 \right \}\\
&\bigcap \left \{ (x, y) \in \mathbb{R}^2 \mid  x\ge Ap_{x}\bigg|_{\alpha =0} \right \}\\
&\bigcap \left \{ (x, y) \in \mathbb{R}^2 \mid  x\ge Ap_{x}\bigg|_{\alpha =\pi} \right \}\\
\end{split} 
\end{equation}

Please note some parts above may not exist or need to be piecewise.

The sofa shape 
\begin{equation*}
S \subseteq S_{max}
\end{equation*}

\begin{corollary}
To maximize the area, the shape should not be a subset on the left side but the intersection $S_{max}$ on the right side.
\end{corollary}

\end{definition}

\begin{definition}[path to area mapping]
So the following mapping associates the area of the sofa shape with the moving paths and parametric functions.
\begin{equation}
\left \{ r(\alpha ),t\left ( \alpha  \right )  \right \} \mapsto Area
\end{equation}
\end{definition}

\section{Calculus of variations for type of Gerver’s sofa}
\begin{definition}[intersection points]
If there is an intersection of the trajectory of A and the envelope of AC, let the intersection point be p. Let the parameter at point p on trajectory of A be $\alpha _{1p}\left(\mathrm{assume}\ \alpha_{1p}\le \pi/2  \right )  $. And the parameter on envelope of AC is $\alpha _{2p} \left(\mathrm{assume}\ \pi /2\le \alpha_{2p}\le \pi  \right ) $, as shown in Figure 3.

So 
\begin{equation}
\left \{ A_{x},A_{y} \right \} \bigg|_{\alpha =\alpha_{1p}} =\left \{ E_{ACx},E_{ACy} \right \} \bigg|_{\alpha =\alpha_{2p}}
\end{equation}

\end{definition}

\begin{definition}[area]
To calculate the area, we integrate the area under the curve using parametric equations that 
\begin{equation*}
\int y\left ( \alpha  \right ) x'\left ( \alpha  \right ) d\alpha 
\end{equation*}

Based on equation (1), it is the enclosed area by the envelopes of ApCp and ApBp, line $y=1$, the trajectory of point A, and the envelopes of AC and AB, and line $y=0$.

So that
\begin{equation}
\begin{split}
Area=& -r\left ( 0 \right ) +2t\left ( 0 \right ) -r\left ( \pi  \right ) +2t\left ( \pi  \right )\\
&+\int_{\pi }^{0} \left ( E_{ApBpy} \frac{dE_{ApBpx}}{d\alpha } \right ) d\alpha \\
&+\int_{\pi }^{0} \left ( E_{ApCpy} \frac{dE_{ApCpx}}{d\alpha } \right ) d\alpha \\
&-\int_{\pi -\alpha _{1p}}^{\alpha _{1p}}\left ( A_{y} \frac{dA_{x}}{d\alpha } \right ) d\alpha \\
&-\int_{0}^{\pi -\alpha _{2p}}\left ( E_{ABy} \frac{dE_{ABx}}{d\alpha } \right ) d\alpha \\
&-\int_{\alpha _{2p}}^{\pi}\left ( E_{ACy} \frac{dE_{ACx}}{d\alpha } \right ) d\alpha \\
\end{split}
\end{equation}

\begin{boldremark}
In order for the above equation (4) to hold, we also need to assume that 
\begin{equation}
E_{ApCpy}\bigg|_{\alpha =0  } =2r'\left ( 0 \right ) \ge 0
\end{equation}

\begin{equation}
E_{ApBpy}\bigg|_{\alpha =\pi  } =-2r'\left ( \pi \right ) \ge 0
\end{equation}

So that the envelopes of ApCp and ApBp lie above the X-axis, and there is no need to add absolute values when calculating the area.

We also need to assume that 
\begin{equation}
E_{ApCpx}\bigg|_{\alpha =\pi  }=-r\left ( \pi  \right ) +2t\left ( \pi  \right ) \ge 0
\end{equation}

\begin{equation}
E_{ApBpx}\bigg|_{\alpha =0  } =r\left ( 0  \right ) -2t\left ( 0  \right ) \le  0
\end{equation}
\end{boldremark}

It is also assumed here that the parametric envelopes and trajectory are convex and do not self-intersect, ensuring the hallway walls are always on one side of the envelopes and do not interact with them. This ensures the integral equation for calculating the area is valid. These assumptions can be verified after obtaining the results.
\end{definition}

\begin{definition}[area with symmetry to $\alpha=\pi /2$]
If assuming and taking advantage of symmetry to $\alpha=\pi /2$ for $r\left ( \alpha  \right )$ and $t\left ( \alpha  \right ) $ with $r\left ( \alpha  \right ) =r\left ( \pi -\alpha  \right ),t\left ( \alpha  \right ) =t\left ( \pi -\alpha  \right )$. Then area in equation (4) yields
\begin{equation}
\begin{split}
Area=&  -2r\left ( 0 \right ) +4t\left ( 0 \right )  \\
&+2\int_{0}^{\pi/2} \left ( -E_{ApBpy} \frac{dE_{ApBpx}}{d\alpha }-E_{ApCpy} \frac{dE_{ApCpx}}{d\alpha } \right ) d\alpha \\
&+2\int_{\alpha _{1p}}^{\pi /2}\left ( A_{y} \frac{dA_{x}}{d\alpha } \right ) d\alpha \\
&+2\int_{0}^{\pi -\alpha _{2p}}\left ( -E_{ABy} \frac{dE_{ABx}}{d\alpha } \right ) d\alpha \\
\end{split}
\end{equation}
\end{definition}

\begin{definition}[extended functions]

Let
\begin{equation*}
F_r=\left\{\begin{matrix} 
  -r\left ( 0 \right ),0\le \alpha \le 1  \\  
  0, 1<\alpha \le \pi /2
\end{matrix}\right. 
\end{equation*}

\begin{equation*}
F_t=\left\{\begin{matrix} 
  2t\left ( 0 \right ),0\le \alpha \le 1  \\  
  0, 1<\alpha \le \pi /2
\end{matrix}\right. 
\end{equation*}

\begin{equation*}
F_{ApBp}=-E_{ApBpy}\frac{d E_{ApBpx}}{d \alpha },0\le \alpha \le \pi /2
\end{equation*}

\begin{equation*}
F_{ApCp}=-E_{ApCpy}\frac{d E_{ApCpx}}{d \alpha },0\le \alpha \le \pi /2
\end{equation*}

\begin{equation*}
F_{A}=\begin{cases} 
  A_{y}\dfrac{d A_{x}}{d \alpha },\alpha _{1p}\le \alpha \le \pi /2 \\  
  0,\alpha < \alpha _{1p}
\end{cases} 
\end{equation*}

\begin{equation*}
F_{AB}=\begin{cases} 
  -E_{ABy}\dfrac{d E_{ABx}}{d \alpha },0\le \alpha \le \pi -\alpha_{2p}\\  
  0,\pi -\alpha_{2p} < \alpha \le \pi/2
\end{cases} 
\end{equation*}

We have six extended functions $F_r, F_t, F_{ApBp}, F_{ApCp}, F_{A}, F_{AB}$. 

So that 
\begin{equation}
Area=2\int_{0}^{\pi /2} \left ( F_r+F_t+F_{ApBp}+F_{ApCp}+F_A+F_{AB} \right ) d\alpha 
\end{equation}

Let
\begin{equation*}
F=F_r+F_t+F_{ApBp}+F_{ApCp}+F_A+F_{AB}
\end{equation*}
\begin{equation*}
F=F\left [ r\left ( \alpha  \right ) , t\left ( \alpha  \right ),r'\left ( \alpha  \right ) , t'\left ( \alpha  \right ),r''\left ( \alpha  \right ) , t''\left ( \alpha  \right )\right ] 
\end{equation*}
Thus
\begin{equation}
Area=2\int_{0}^{\pi /2} F\left [ r\left ( \alpha  \right ) , t\left ( \alpha  \right ),r'\left ( \alpha  \right ) , t'\left ( \alpha  \right ),r''\left ( \alpha  \right ) , t''\left ( \alpha  \right )\right ]d\alpha  
\end{equation}
The area is an integral functional of $r\left ( \alpha  \right ) , t\left ( \alpha  \right ),r'\left ( \alpha  \right ) , t'\left ( \alpha  \right ),r''\left ( \alpha  \right ) , t''\left ( \alpha  \right )$ as the mapping in equation (2).
\end{definition}

\begin{definition}[moving sofa problem-to maximize area]
The moving sofa problem aims to find the maximum value of the following integral functional

\begin{equation}
\text{maximize} \int_{0}^{\pi/2} Fd\alpha 
\end{equation}
\end{definition}

\begin{theorem}
The necessary condition for the sofa area to reach the maximum value is that the functional derivative is zero, which is given by the Euler-Lagrange equations.
\end{theorem}

\begin{proof}
It is well-known in textbooks of Calculus of Variations\cite{Elsgolc2012} \cite{Gelfand2000}.
\end{proof}

\begin{definition}[EL equations]
The necessary condition to maximize the area is that the functional derivative equals zero, which involves solving the associated EL equations below.
\begin{equation}
\left ( ODE1 \right )\frac{\partial F}{\partial r}-\frac{d }{d \alpha }\frac{\partial F}{\partial r'}+\frac{d^2 }{d \alpha^2 }\frac{\partial F}{\partial r''}=0  
\end{equation}

\begin{equation}
\left ( ODE2 \right )\frac{\partial F}{\partial t}-\frac{d }{d \alpha }\frac{\partial F}{\partial t'}+\frac{d^2 }{d \alpha^2 }\frac{\partial F}{\partial t''}=0  
\end{equation}

With the explicit form of functional $F$ in Definition 4.4, we obtain two ODEs (See Appendix 1 for details) as 
\begin{equation}
ODE1=\begin{cases}
  &-\dfrac{1}{2} \cos (\alpha ) \left(-\cos \left(\dfrac{\alpha }{2}\right)-2 (-3 \cos (\alpha )+\cos (2 \alpha )) r(\alpha )+\sin \left(\dfrac{\alpha
   }{2}\right)+2 (-3 \cos (\alpha )\right.\\
   &\left.+\cos (2 \alpha )) t(\alpha )\right)+12 \sin (\alpha ) r'(\alpha )-3 \sin (2 \alpha ) r'(\alpha )-6 \sin (\alpha )
   t'(\alpha )+3 \sin (2 \alpha ) t'(\alpha )\\
   &\left.+r''(\alpha )-6 \cos (\alpha ) r''(\alpha )+\cos (2 \alpha ) r''(\alpha )+t''(\alpha )-\cos (2 \alpha) t''(\alpha )\right), 0<\alpha \leq \alpha_{1p} \\   \\
   
  &-\dfrac{1}{2} \cos (\alpha ) \left(-\cos \left(\dfrac{\alpha }{2}\right)-2 (-3 \cos (\alpha )+\cos (2 \alpha )) r(\alpha )+\sin \left(\dfrac{\alpha
   }{2}\right)+2 (-2 \cos (\alpha )\right.\\
   &\left.+\cos (2 \alpha )) t(\alpha )\right)+12 \sin (\alpha ) r'(\alpha )-3 \sin (2 \alpha ) r'(\alpha )-4 \sin (\alpha )
   t'(\alpha )+3 \sin (2 \alpha ) t'(\alpha )\\
   &\left.+r''(\alpha )-6 \cos (\alpha ) r''(\alpha )+\cos (2 \alpha ) r''(\alpha )+t''(\alpha )-\cos (2 \alpha
   ) t''(\alpha )\right), \alpha_{1p}<\alpha \leq \pi- \alpha_{2p}\\ \\
   
  &-\dfrac{1}{2} \cos (\alpha ) \left(-\cos \left(\dfrac{\alpha }{2}\right)+4 \cos (\alpha ) r(\alpha )+\sin \left(\dfrac{\alpha }{2}\right)-2 \cos
   (\alpha ) t(\alpha )+8 \sin (\alpha ) r'(\alpha )\right)\\ 
   &\left.-2 \sin (\alpha ) t'(\alpha )-4 \cos (\alpha ) r''(\alpha )\right) , \pi- \alpha_{2p}<\alpha \leq \pi/2\\
\end{cases}
\end{equation}

\begin{equation}
ODE2=\begin{cases}
  &\dfrac{1}{2} \sin (\alpha ) \left(\cos \left(\dfrac{\alpha }{2}\right)+\sin \left(\dfrac{\alpha }{2}\right)+2 r(\alpha ) (3 \sin (\alpha )+\sin (2
   \alpha ))-2 (3 \sin (\alpha )+\sin (2 \alpha )) t(\alpha )\right.\\
   &\left.+r'(\alpha )-6 \cos (\alpha ) r'(\alpha )-3 \cos (2 \alpha ) r'(\alpha )+t'(\alpha
   )+12 \cos (\alpha ) t'(\alpha )+3 \cos (2 \alpha ) t'(\alpha )\right.\\
   &\left.-\sin (2 \alpha ) r''(\alpha )+6 \sin (\alpha ) t''(\alpha )+\sin (2 \alpha )
   t''(\alpha )\right), 0<\alpha \leq \alpha_{1p} \\  \\
   
  &\dfrac{1}{2} \sin (\alpha ) \left(\cos \left(\dfrac{\alpha }{2}\right)+\sin \left(\dfrac{\alpha }{2}\right)+4 (1+\cos (\alpha )) r(\alpha ) \sin
   (\alpha )-2 (3 \sin (\alpha )+\sin (2 \alpha )) t(\alpha )\right.\\
   &\left.+r'(\alpha )-4 \cos (\alpha ) r'(\alpha )-3 \cos (2 \alpha ) r'(\alpha )+t'(\alpha
   )+12 \cos (\alpha ) t'(\alpha )+3 \cos (2 \alpha ) t'(\alpha )\right.\\
   &\left.-\sin (2 \alpha ) r''(\alpha )+6 \sin (\alpha ) t''(\alpha )+\sin (2 \alpha )
   t''(\alpha )\right) , \alpha_{1p}<\alpha \leq \pi- \alpha_{2p}\\ \\
   
  &\dfrac{1}{2} \sin (\alpha ) \left(\cos \left(\dfrac{\alpha }{2}\right)+\sin \left(\dfrac{\alpha }{2}\right)+2 r(\alpha ) \sin (\alpha )-4 \sin
   (\alpha ) t(\alpha )-2 \cos (\alpha ) r'(\alpha )\right.\\
   &\left.+8 \cos (\alpha ) t'(\alpha )+4 \sin (\alpha ) t''(\alpha )\right) , \pi- \alpha_{2p}<\alpha \leq \pi/2\\
\end{cases}
\end{equation}

These are systems of second-order ODEs. 
\end{definition}

\begin{definition}[boundary conditions]

To solve these ODEs, four boundary conditions are needed.

With symmetry in Definition 4.3, we have
\begin{equation}
r'\left ( \frac{\pi }{2}  \right ) =0
\end{equation}
\begin{equation}
t'\left ( \frac{\pi }{2}  \right ) =0
\end{equation}

For the variational problems with movable boundaries $r(0)$ and $t(0)$, we also have the transversality conditions based on boundary terms from first variation
\begin{equation*}
\left ( \frac{\partial F}{\partial r'} -\frac{d}{d\alpha } \frac{\partial F}{\partial r''} \right ) \bigg|_{\alpha =0}=0
\end{equation*}

\begin{equation*}
\left ( \frac{\partial F}{\partial t'} -\frac{d}{d\alpha } \frac{\partial F}{\partial t''} \right ) \bigg|_{\alpha =0}=0
\end{equation*}

Thus we also obtain additional boundary condition
\begin{equation}
-\frac{1}{2} +2r\left(0\right)-2t\left (0\right)-2r''\left (0\right)=0
\end{equation}

In a previous paper \cite{Song2016}, the author also obtained the above boundary condition.

Thus we obtain three additional boundary conditions. The last condition needed can be $t(0)$, and we numerically search for the optimal value of $t(0)$ to maximize the area (See Appendix 2 for details). 

\end{definition}

\section{Results of numerical calculations}
We use finite-difference methods to discretize and solve the above ODEs in equations (15) and (16) with boundary conditions in equations (17)-(19) and constraint in equation (3) (See Appendix 3 for details). The interval $\left [ 0,\frac{\pi }{2}  \right ] $ is divided into uniform grid with 2000 points using a fourth-order difference scheme. Five iterations are performed to numerically obtain the values of $\alpha_{1p}$ and $\alpha_{2p}$ such that the ODE and additional equations are satisfied. We search for the parameter values and boundary conditions that meet the equations and continuity conditions. 
\begin{boldremark}
It is noted that a more refined search or equation solving should yield more precise parameters and results.
\end{boldremark}

The obtained sofa curves are shown in Figure 4, consisting of five curves and three line segments. The obtained area of the sofa
\begin{equation*}
Area\approx 2.2195316...
\end{equation*}

The difference between Gerver's sofa constant is less than \( 4 \times 10^{-8} \). 

\begin{figure}[H]
  \centering
  \includegraphics[width=12cm]{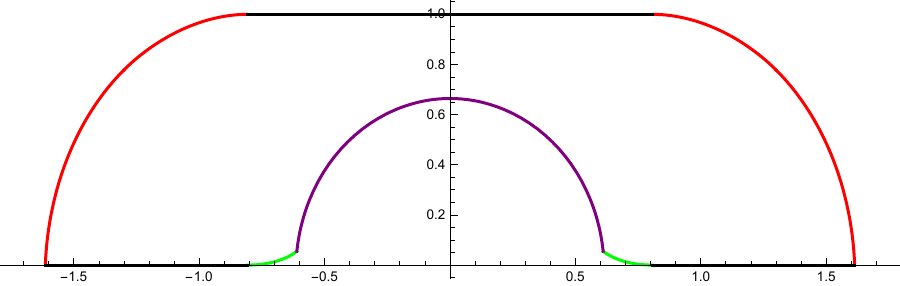}
  \caption{Obtained sofa shape with 5 curves and 3 lines.}
\end{figure}

Figure 5 shows the obtained functions $r\left ( \alpha  \right ), t\left ( \alpha  \right ) $ and $r'\left ( \alpha  \right ), t'\left ( \alpha  \right ) $.

\begin{figure}[H]
  \centering
  \includegraphics[width=12cm]{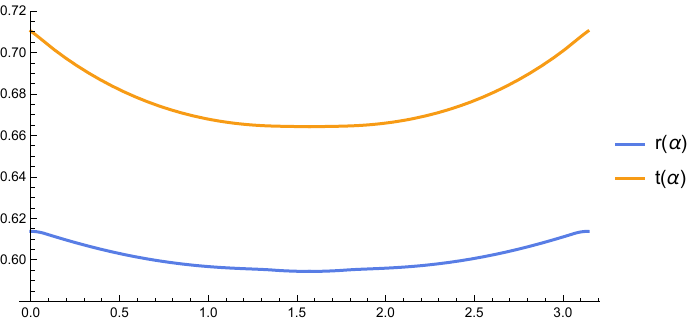}
  \includegraphics[width=12cm]{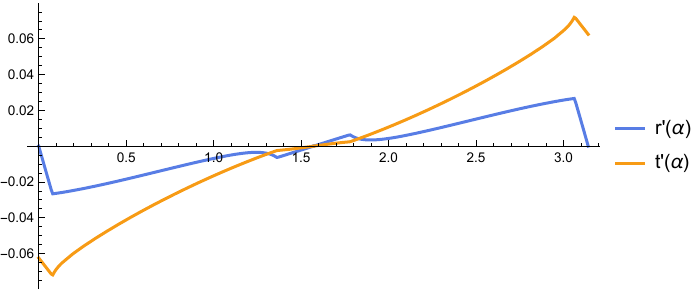}
  \caption{Obtained functions $r\left ( \alpha  \right ), t\left ( \alpha  \right ), r'\left ( \alpha  \right ), t'\left ( \alpha  \right ) $.}
\end{figure}

\section{Verification of Gerver’s results }
In this section, we will demonstrate that Gerver's previously constructed results coincide with those we obtain by solving the EL equations. We will specifically verify the trajectory results of point A, as the verification for other trajectories and envelopes follows a similar process.
\begin{definition}[recall Gerver’s results]
Here are the results from Gerver's work \cite{Gerver1992}. Gerver discovered the "balanced polygons" and provided conditions for the shape with optimal area. The boundary and area of the sofa are given by constants derived by solving implicit equations
\begin{equation*}
A(cos\theta -cos\phi )-2Bsin\phi+(\theta-\phi-1)cos\theta-sin\theta+cos\phi+sin\phi=0
\end{equation*}
\begin{equation*}
A(3sin\theta +sin\phi )-2Bcos\phi+3(\theta-\phi-1)sin\theta+3cos\theta-sin\phi+cos\phi=0
\end{equation*}
\begin{equation*}
Acos\phi-(sin\phi+\frac{1}{2}-\frac{1}{2}cos\phi+B sin\phi )=0
\end{equation*}
\begin{equation*}
\left ( A+\frac{1}{2}\pi -\phi -\theta   \right ) -\left [ B-\frac{1}{2}\left ( \theta -\phi  \right )(1+A) -\frac{1}{4} \left ( \theta -\phi  \right )^2  \right]=0
\end{equation*}

Where
\begin{equation*}A=0.094426560843653...\end{equation*}
\begin{equation*}B=1.399203727333547...\end{equation*}
\begin{equation*}\phi =0.039177364790084...\end{equation*}
\begin{equation*}\theta =0.681301509382725...\end{equation*}

\end{definition}

However, Gerver’s results are implicit.

\begin{definition}[recall Romik’s explicit and analytic formulas]
We will use explicit and analytic formulas, as well as symbols in the previous paper \cite{Romik2018} and demonstration \cite{Romikweb}. Recall the parametric equation of the trajectory of point A by Romik:

\begin{equation*}
\begin{array}{l} 
  \mathbf{x} \left ( t \right )=
\begin{cases} 
  \mathbf{x}_{1} \left ( t \right ),\text{if } 0\le t< \varphi \\ 
  \mathbf{x}_{2} \left ( t \right ),\text{if } \varphi\le t< \theta  \\ 
  \mathbf{x}_{3} \left ( t \right ),\text{if } \theta\le t\le \pi /2-\theta \\ 
  \mathbf{x}_{4} \left ( t \right ),\text{if } \pi /2-\theta< t\le \pi /2-\varphi  \\ 
  \mathbf{x}_{5} \left ( t \right ),\text{if } \pi /2-\varphi< t\le \pi /2
\end{cases}    
\end{array} 
\end{equation*}
where
\begin{equation*}
\begin{split}
&\mathbf{x}_{1} \left ( t \right )=\left\{\kappa _{1,1}-\sin (t) \left(a_1 \sin (t)-a_2 \cos (t)-\frac{1}{2}\right)+\cos (t) \left(a_2 \sin (t)+a_1 \cos (t)-1\right),\right.\\
&\left.\kappa_{1,2}+\cos (t) \left(a_1 \sin (t)-a_2 \cos (t)-\frac{1}{2}\right)+\sin (t) \left(a_2 \sin (t)+a_1 \cos (t)-1\right)\right\}\\
\end{split}
\end{equation*}

\begin{equation*}
\begin{split}
&\mathbf{x}_{2} \left ( t \right )=\left\{\kappa _{2,1}+\left(b_1 t+b_2-\frac{t^2}{4}\right) \cos (t)+\left(-b_1+\frac{t}{2}-1\right) (-\sin (t)),\right.\\
&\left.\kappa _{2,2}+\left(b_1 t+b_2-\frac{t^2}{4}\right) \sin (t)+\left(-b_1+\frac{t}{2}-1\right) \cos (t)\right\}\\
\end{split}
\end{equation*}

\begin{equation*}
\begin{split}
&\mathbf{x}_{3} \left ( t \right )=\left\{\kappa _{3,1}-\left(c_2+t\right) \sin (t)+\left(c_1-t\right) \cos (t),\right.\\
&\left.\kappa _{3,2}+\left(c_1-t\right) \sin (t)+\left(c_2+t\right) \cos(t)\right\}\\
\end{split}
\end{equation*}

\begin{equation*}
\begin{split}
&\mathbf{x}_{4} \left ( t \right )=\left\{\kappa _{4,1}-\left(d_1 t+d_2-\frac{t^2}{4}\right) \sin (t)+\left(d_1-\frac{t}{2}-1\right) \cos (t),\right.\\
&\left.\kappa _{4,2}+\left(d_1 t+d_2-\frac{t^2}{4}\right) \cos (t)+\left(d_1-\frac{t}{2}-1\right) \sin (t)\right\}\\
\end{split}
\end{equation*}

\begin{equation*}
\begin{split}
&\mathbf{x}_{5} \left ( t \right )=\left\{\kappa _{5,1}-\sin (t) \left(e_1 \sin (t)-e_2 \cos (t)-1\right)+\cos (t) \left(e_2 \sin (t)+e_1 \cos (t)-\frac{1}{2}\right),\right.\\
&\left.\kappa_{5,2}+\cos (t) \left(e_1 \sin (t)-e_2 \cos (t)-1\right)+\sin (t) \left(e_2 \sin (t)+e_1 \cos (t)-\frac{1}{2}\right)\right\}\\
\end{split}
\end{equation*}
\end{definition}

To align with the coordinate system defined in Figure 2, we convert $t = \alpha/2$ and translate the curves by $(-(1 - e_{1} + \kappa_{5, 1})/2, 0) \approx (0.6137632, 0)$ to achieve symmetry about the Y-axis.

Thus, we obtain the parametric equation of the trajectory of point A with $0 \le \alpha \le \pi$ as shown in Figure 2 as

\begin{equation*}
\begin{cases} 
  \mathbf{x}_{1} \left ( \alpha/2 \right )-\left \{ \dfrac{1}{2} \left(1-e_1+\kappa _{5,1}\right),0 \right \},\text{if } 0\le \alpha< 2\varphi \\ 
  \mathbf{x}_{2} \left ( \alpha/2 \right )-\left \{ \dfrac{1}{2} \left(1-e_1+\kappa _{5,1}\right),0 \right \},\text{if } 2\varphi\le \alpha< 2\theta  \\ 
  \mathbf{x}_{3} \left ( \alpha/2 \right )-\left \{ \dfrac{1}{2} \left(1-e_1+\kappa _{5,1}\right),0 \right \},\text{if } 2\theta\le \alpha\le \pi-2\theta \\ 
  \mathbf{x}_{4} \left ( \alpha/2 \right )-\left \{ \dfrac{1}{2} \left(1-e_1+\kappa _{5,1}\right),0 \right \},\text{if } \pi -2\theta< \alpha\le \pi -2\varphi  \\ 
  \mathbf{x}_{5} \left ( \alpha/2 \right )-\left \{ \dfrac{1}{2} \left(1-e_1+\kappa _{5,1}\right),0 \right \},\text{if } \pi -2\varphi< \alpha\le \pi 
\end{cases}
\end{equation*}

\begin{corollary}
One explicit formula is
\begin{equation*}
\begin{split}
&\mathbf{x}_{1} \left ( \alpha/2 \right )-\left \{ \frac{1}{2} \left(1-e_1+\kappa _{5,1}\right),0 \right \}=\\
&\left\{\frac{1}{2} \left(2 \kappa _{1,1}-\kappa _{5,1}+2 a_2 \sin (\alpha )+2 a_1 \cos (\alpha )+\sin \left(\frac{\alpha }{2}\right)-2 \cos
   \left(\frac{\alpha }{2}\right)+e_1-1\right),\right.\\
&\left.\kappa _{1,2}+a_1 \sin (\alpha )-a_2 \cos (\alpha )-\sin \left(\frac{\alpha }{2}\right)-\frac{1}{2}
   \cos \left(\frac{\alpha }{2}\right)\right\}
\end{split}
\end{equation*}
\end{corollary}

\begin{theorem}
Gerver’s and Romik’s sofa results satisfy the EL equations, which are necessary conditions for maximum area.
\end{theorem}

\begin{proof}
To prove that these formulas satisfy the EL equations ODE1 and ODE2 in equations (15) and (16), we employ substitution and symbolic calculations (see Appendix 4 for details). For instance, for $0<\alpha<\alpha_{1p}$
\begin{equation*}
r\left ( \alpha \right )=\frac{1}{2} \left(2 \kappa _{1,1}-\kappa _{5,1}+2 a_2 \sin (\alpha )+2 a_1 \cos (\alpha )+\sin \left(\frac{\alpha }{2}\right)-2 \cos
   \left(\frac{\alpha }{2}\right)+e_1-1\right)/\cos(\alpha )
\end{equation*}
\begin{equation*}
t\left ( \alpha \right )=\left [ \kappa _{1,2}+a_1 \sin (\alpha )-a_2 \cos (\alpha )-\sin \left(\frac{\alpha }{2}\right)-\frac{1}{2}
   \cos \left(\frac{\alpha }{2}\right) \right ] /\sin(\alpha )
\end{equation*}
satisfy ODE1 in equation (15) as
\begin{equation*}
\begin{split}
&-\cos \left(\frac{\alpha }{2}\right)-2 (-3 \cos (\alpha )+\cos (2 \alpha )) r(\alpha )+\sin \left(\frac{\alpha
   }{2}\right)+2 (-3 \cos (\alpha )+\cos (2 \alpha )) t(\alpha )\\
   &+12 \sin (\alpha ) r'(\alpha )-3 \sin (2 \alpha ) r'(\alpha )-6 \sin (\alpha )
   t'(\alpha )+3 \sin (2 \alpha ) t'(\alpha )+r''(\alpha )-6 \cos (\alpha ) r''(\alpha )\\
   &+\cos (2 \alpha ) r''(\alpha )+t''(\alpha )-\cos (2 \alpha
   ) t''(\alpha )=0
\end{split}
\end{equation*}
and ODE2 in equation (16) as
\begin{equation*}
\begin{split}
&\cos \left(\frac{\alpha }{2}\right)+\sin \left(\frac{\alpha }{2}\right)+2 r(\alpha ) (3 \sin (\alpha )+\sin (2
   \alpha ))-2 (3 \sin (\alpha )+\sin (2 \alpha )) t(\alpha )+r'(\alpha )\\
&-6 \cos (\alpha ) r'(\alpha )-3 \cos (2 \alpha ) r'(\alpha )+t'(\alpha
   )+12 \cos (\alpha ) t'(\alpha )+3 \cos (2 \alpha ) t'(\alpha )-\sin (2 \alpha ) r''(\alpha )\\
   &+6 \sin (\alpha ) t''(\alpha )+\sin (2 \alpha )
   t''(\alpha )=0
\end{split}
\end{equation*}

Similarly, the remaining $\mathbf{x}_{2},\mathbf{x}_{3}$ can also be verified by symbolic calculations (see Appendix 4 for details). Then $\mathbf{x}_{4},\mathbf{x}_{5}$ also satisfy since they and $\mathbf{x}_{2},\mathbf{x}_{1}$ are symmetrical.
\end{proof}

Fortunately, the ODEs corresponding to the Euler-Lagrange equations possess analytical solutions, which are the results obtained by Gerver and Romik. Therefore, Gerver's results satisfy the EL equations, meeting the necessary condition for maximum area.

\section{Symmetric conditions}

In the previous section, we discussed the scenario for intersection point of trajectory of point A and the envelope of AC in Figure 3 with equation (3). Two possible scenarios for intersection points are illustrated in Figure 6. We have discussed the second scenario in the previous section as shown in Figure 6(b).

\begin{figure}[H]
  \centering
  \includegraphics[width=12cm]{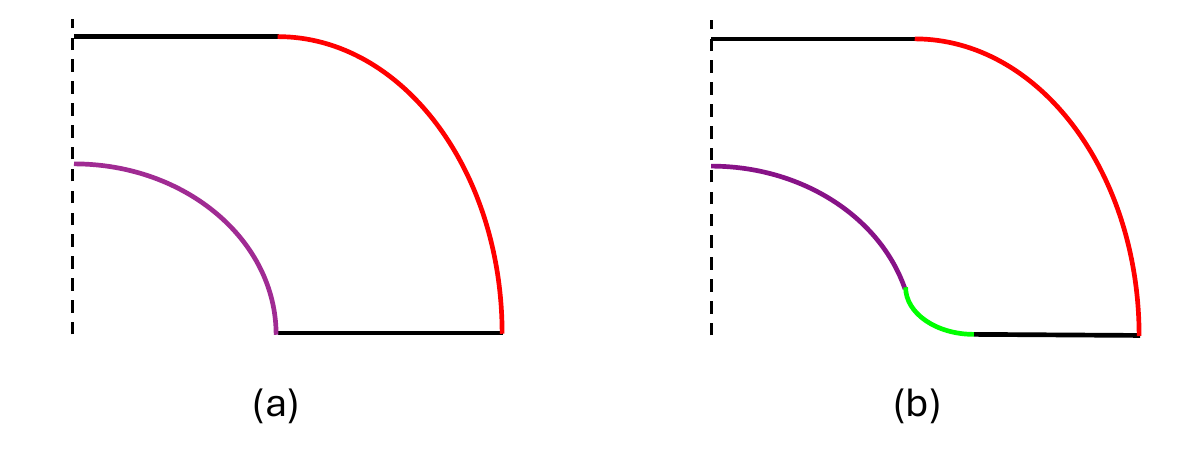}
  \caption{Two possible scenarios for intersection points of trajectory of A and the envelope of AC. The dashed line is the Y-axis of symmetry.}
\end{figure}

\subsection{EL equations and boundary conditions}

\begin{itemize}
\item \textbf{Scenario in Figure 6(a)}
\end{itemize}

\begin{figure}[H]
  \centering
  \includegraphics[width=12cm]{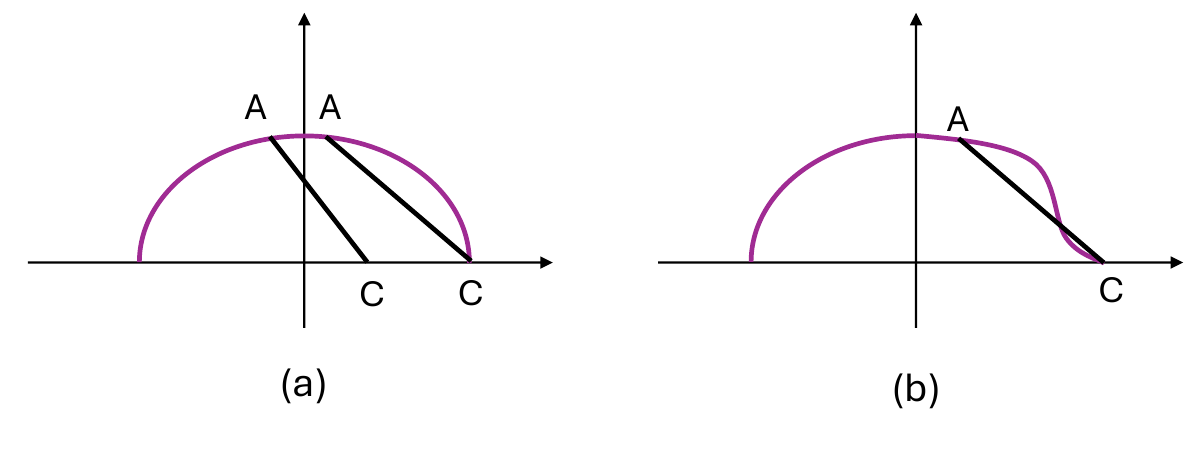}
  \caption{Convexity of trajectory of A and intersection of trajectory of A and the envelope of AC}
\end{figure}

\begin{lemma}
If the trajectory of point A is convex, and
\begin{equation*}
\forall 0\le \alpha \le \pi ,C_{x}\le A_{x}\bigg|_{\alpha =0}
\end{equation*}
then there is no intersection of the trajectory of point A and the envelope of AC.

Similarly, if
\begin{equation*}
\forall 0\le \alpha \le \pi , A_{x}\bigg|_{\alpha =\pi }\le B_{x}
\end{equation*}
It leads to no intersection of the trajectory of point A and the envelope of AB.

If the two above inequalities do not hold, or trajectory of point A is non-convex as shown in Figure 7(b), then there can be intersections of the trajectory of point A with envelope of AB and AC.
\end{lemma}
\begin{proof}
It follows from the definition of convexity.
\end{proof}

If there is no intersection point of the trajectory of A and the envelope of AC as in Figure 6(a), the area calculation would differ from equation (4).

\begin{definition}[area with symmetry to $\alpha =\pi /2$]
Assuming symmetry to $\pi /2$, and making similar assumptions as in equations (5)-(8), the area can be calculated as 
\begin{equation*}
\begin{split}
Area=&  -2r\left ( 0 \right ) +4t\left ( 0 \right )  \\
&+2\int_{0}^{\pi/2} \left ( -E_{ApBpy} \dfrac{dE_{ApBpx}}{d\alpha }-E_{ApCpy} \dfrac{dE_{ApCpx}}{d\alpha } \right ) d\alpha \\
&+2\int_{0}^{\pi /2}\left ( A_{y} \frac{dA_{x}}{d\alpha } \right ) d\alpha \\
\end{split}
\end{equation*}

For this condition, there is no intersection point of the trajectory of A and the envelope of AC. With above Lemma 3, we have the constraints
\begin{equation*}
r\left ( \alpha  \right ) \cos \alpha +t\left ( \alpha  \right ) \sin \alpha\tan \frac{\alpha }{2} \le r\left ( 0 \right ) 
\end{equation*}
\begin{equation*}
-r\left ( \pi  \right )\le r\left ( \alpha  \right ) \cos \alpha -2 \cos ^2\left(\frac{\alpha }{2}\right)t\left ( \alpha  \right ) 
\end{equation*}

Note that $r\left ( \pi  \right )=r\left ( 0  \right )$.
\end{definition}

\begin{definition}[moving sofa problem-to maximize area]
The problem then becomes a variational problem with constraints.
\begin{equation*}
\begin{split}
\text{maximize}&\int_{0}^{\pi/2} \left ( -E_{ApBpy} \frac{dE_{ApBpx}}{d\alpha }-E_{ApCpy} \frac{dE_{ApCpx}}{d\alpha }+A_{y} \frac{dA_{x}}{d\alpha } \right ) d\alpha -r\left ( 0 \right ) +2t\left ( 0 \right ) \\
\text{s.t.} &r\left ( \alpha  \right ) \cos \alpha +t\left ( \alpha  \right ) \sin \alpha\tan \frac{\alpha }{2} \le r\left ( 0 \right ) \\
&-r\left ( 0  \right )\le r\left ( \alpha  \right ) \cos \alpha -2 \cos ^2\left(\frac{\alpha }{2}\right)t\left ( \alpha  \right ) 
\end{split}
\end{equation*}

We can use a Lagrange multiplier so that 
\begin{equation*}
\begin{split}
\text{maximize} &\int_{0}^{\pi/2} \left ( -E_{ApBpy} \frac{dE_{ApBpx}}{d\alpha }-E_{ApCpy} \frac{dE_{ApCpx}}{d\alpha }+A_{y} \frac{dA_{x}}{d\alpha } \right ) d\alpha -r\left ( 0 \right ) +2t\left ( 0 \right )\\
&+\int_{0}^{\pi/2} \lambda_1\left ( \alpha  \right )\left [ r\left ( \alpha  \right ) \cos \alpha +t\left ( \alpha  \right ) \sin \alpha\tan \frac{\alpha }{2} - r\left ( 0 \right ) \right ] d\alpha\\
&+\int_{0}^{\pi/2} \lambda_2\left ( \alpha  \right )\left [ -r\left ( 0  \right )- r\left ( \alpha  \right ) \cos \alpha +2\cos^2 \left ( \frac{\alpha }{2}  \right ) t\left ( \alpha  \right ) \right ] d\alpha  
\end{split}
\end{equation*}
\end{definition}

\begin{definition}[EL equations and boundary conditions]
Similar to equations (15) and (16), we derive the associated EL equations (see Appendix 5 for details) as

\begin{equation}
\begin{split}
ODE1=&-\cos \left(\dfrac{\alpha }{2}\right)+4 \cos (\alpha ) r(\alpha )+\sin \left(\dfrac{\alpha }{2}\right)-2 \cos (\alpha ) t(\alpha )-2
   \lambda_{1}(\alpha )+2 \lambda_{2}(\alpha )+8 \sin (\alpha ) r'(\alpha )\\
   &-2 \sin (\alpha ) t'(\alpha )-4 \cos (\alpha )   r''(\alpha )
   \end{split}
\end{equation}

\begin{equation}
\begin{split}
ODE2=&\cos \left(\dfrac{\alpha }{2}\right)-\cos \left(\dfrac{3 \alpha }{2}\right)+\sin \left(\dfrac{\alpha }{2}\right)+4 r(\alpha ) \sin ^2(\alpha )+\sin
   \left(\frac{3 \alpha }{2}\right)-8 \sin ^2(\alpha ) t(\alpha )\\
   &+4 \lambda_{1}(\alpha )-4 \cos (\alpha ) \lambda_{1}(\alpha )+4
   \lambda_{2}(\alpha )+4 \cos (\alpha ) \lambda_{2}(\alpha )-2 \sin (2 \alpha ) r'(\alpha )+8 \sin (2 \alpha ) t'(\alpha )\\
   &+4t''(\alpha )-4 \cos (2 \alpha ) t''(\alpha )
      \end{split}
\end{equation}

\begin{equation}
ODE3=\lambda_1(\alpha )\left [ \cos (\alpha ) r(\alpha
   )-r(0)+\sin (\alpha ) \tan \left(\frac{\alpha }{2}\right)
   t(\alpha ) \right ] 
\end{equation}

\begin{equation}
ODE4=\lambda_2(\alpha ) \left [ -\cos (\alpha ) r(\alpha
   )+r(\pi )+2 \cos ^2\left(\frac{\alpha }{2}\right) t(\alpha
   ) \right ] 
\end{equation}

We also derive boundary conditions due to symmetry
\begin{equation*}
r'\left ( \frac{\pi }{2}  \right ) =0
\end{equation*}

\begin{equation*}
t'\left ( \frac{\pi }{2}  \right ) =0
\end{equation*}

The boundary values $r(0), t(0)$ are to be determined. 

\end{definition}

Using similar numerical methods as in Section 5, we solve the ODEs in equations (20)-(23) with the boundary conditions and obtain the shape and area. 

\subsection{Results of numerical calculations}
We search for $r(0), t(0)$ from 0.6 to 0.7 and calculate area via numerical calculations (see Appendix 5 for details). The resulting sofa curves consist of 3 curves and 3 line segments. The obtained area of the sofa is close to $2/\pi\approx 0.63662$ and Hammersley’s sofa area $2/\pi+\pi/2\approx 2.20741$.

\section{Asymmetric conditions}

In previous sections, we all assume the shape is symmetry to $\alpha =\pi /2$ to simplify the area calculation. Here in this section, we consider asymmetric conditions, leading to different area calculations from equation (9). Several possible scenarios are depicted in Figure 8.

\begin{figure}[H]
  \centering
  \includegraphics[width=14cm]{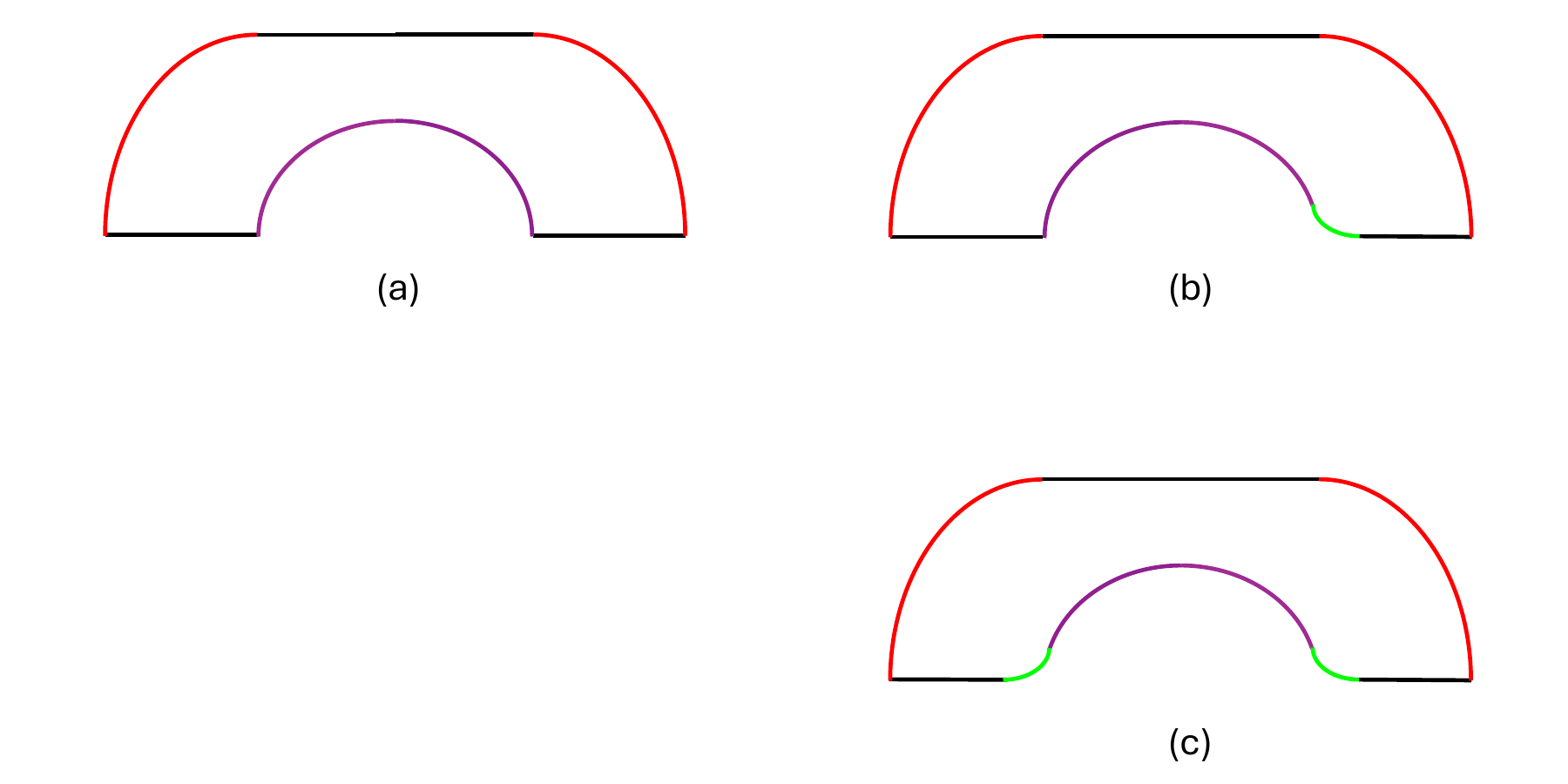}
  \caption{Three possible asymmetric conditions.}
\end{figure}

\subsection{EL equations and boundary conditions}
\begin{definition}[intersection points]

We group the intersection points (which may not exist) on both sides as
\begin{equation*}
\begin{split}
\left \{ A_{x},A_{y} \right \} \bigg|_{\alpha =\alpha_{1p}}=\left \{ E_{ACx},E_{ACy} \right \}\bigg|_{\alpha =\alpha_{2p}}\\
\left \{ A_{x},A_{y} \right \} \bigg|_{\alpha =\alpha_{3p}}=\left \{ E_{ABx},E_{ABy} \right \}\bigg|_{\alpha =\alpha_{4p}}\\
\end{split}
\end{equation*}

Without loss of generality, let 
\begin{equation*}
0<\alpha_{4p}\le \alpha_{1p}\le1\le\alpha_{2p}\le\alpha_{3p}<\pi
\end{equation*}
\end{definition}

\begin{itemize}
\item \textbf{Scenario in Figure 8(a)}
\end{itemize}

\begin{definition}[moving sofa problem-to maximize area]

If there are no intersection points of the trajectory of A and either the envelope of AB or AC, and assuming similar conditions as in equations (5)-(8), the area can be calculated as
\begin{equation*}
\begin{split}
Area=& -r\left ( 0 \right ) +2t\left ( 0 \right ) -r\left ( \pi  \right ) +2t\left ( \pi  \right ) \\
&+\int_{0 }^{\pi} \left ( -E_{ApBpy} \frac{dE_{ApBpx}}{d\alpha } \right ) d\alpha \\
&+\int_{0 }^{\pi} \left ( -E_{ApCpy} \frac{dE_{ApCpx}}{d\alpha } \right ) d\alpha \\
&+\int_{0 }^{\pi}\left ( -A_{y} \frac{dA_{x}}{d\alpha } \right ) d\alpha \\
\end{split}
\end{equation*}

Using a Lagrange multiplier, we have
\begin{equation*}
\begin{split}
\text{maximize} & -r\left ( 0 \right ) +2t\left ( 0 \right ) -r\left ( \pi  \right ) +2t\left ( \pi  \right ) \\
&+\int_{0 }^{\pi} \left ( -E_{ApBpy} \frac{dE_{ApBpx}}{d\alpha } \right ) d\alpha \\
&+\int_{0 }^{\pi} \left ( -E_{ApCpy} \frac{dE_{ApCpx}}{d\alpha } \right ) d\alpha \\
&+\int_{0 }^{\pi}\left ( -A_{y} \frac{dA_{x}}{d\alpha } \right ) d\alpha \\
&+\int_{0}^{\pi} \lambda_{1}(\alpha ) \left[r(\alpha ) \cos (\alpha )+t(\alpha ) \sin (\alpha ) \tan \left(\frac{\alpha
   }{2}\right)-r(0)\right] d\alpha\\
&+\int_{0}^{\pi} \lambda_{2}(\alpha )\left \{ r(\pi )-\left[r(\alpha ) \cos (\alpha )-t(\alpha ) 2 \cos ^2\left(\frac{\alpha
   }{2}\right)\right] \right \}  d\alpha\\
\end{split}
\end{equation*}
\end{definition}

\begin{definition}[EL equations]

Similar to equations (15) and (16), we obtain the associated EL equations (see Appendix 6 for details) as

\begin{equation}
\begin{split}
ODE1=&-\dfrac{1}{2} \cos (\alpha ) \left(-\cos \left(\dfrac{\alpha }{2}\right)+4 \cos (\alpha ) r(\alpha )+\sin \left(\dfrac{\alpha }{2}\right)-2 \cos
   (\alpha ) t(\alpha )-2 \lambda_{1}(\alpha )+2 \lambda_{2}(\alpha )\right.\\
   &\left. +8 \sin (\alpha ) r'(\alpha )-2 \sin (\alpha ) t'(\alpha )-4
   \cos (\alpha ) r''(\alpha )\right)
\end{split}
\end{equation}

\begin{equation}
\begin{split}
ODE2=&\dfrac{1}{4} \left(\cos \left(\dfrac{\alpha }{2}\right)-\cos \left(\dfrac{3 \alpha }{2}\right)+\sin \left(\dfrac{\alpha }{2}\right)+4 r(\alpha ) \sin
   ^2(\alpha )+\sin \left(\dfrac{3 \alpha }{2}\right)-8 \sin ^2(\alpha ) t(\alpha )\right.\\
   &\left.+4 \lambda_{1}(\alpha )-4 \cos (\alpha ) \lambda_{1}(\alpha )+4 \lambda_{2}(\alpha )+4 \cos (\alpha ) \lambda_{2}(\alpha )-2 \sin (2 \alpha ) r'(\alpha )+8 \sin (2 \alpha )
   t'(\alpha )\right.\\
   &\left.+4 t''(\alpha )-4 \cos (2 \alpha ) t''(\alpha )\right)
   \end{split}
\end{equation}

\begin{equation}
ODE3=\lambda_{1}(\alpha ) \left(r(\alpha ) \cos (\alpha )+t(\alpha ) \sin (\alpha ) \tan \left(\dfrac{\alpha
   }{2}\right)-r(0)\right)
\end{equation}

\begin{equation}
ODE4=\lambda_{2}(\alpha ) \left(r(\pi )-\left(r(\alpha ) \cos (\alpha )-2t(\alpha )  \cos ^2\left(\dfrac{\alpha
   }{2}\right)\right)\right)
\end{equation}
\end{definition}

The boundary values $r(0), t(0),r(\pi), t(\pi)$ are to be determined.

\begin{itemize}
\item \textbf{Scenario in Figure 8(b)}
\end{itemize}

\begin{definition}[moving sofa problem-to maximize area]

If there is an intersection point on one side, without loss of generality, we can assume there are intersection points of the trajectory of point A with the envelope of AC but not with the envelope of AB. Under similar assumptions as in equations (5)-(8), the area can be calculated as
\begin{equation*}
\begin{split}
Area=& -r\left ( 0 \right ) +2t\left ( 0 \right ) -r\left ( \pi  \right ) +2t\left ( \pi  \right ) \\
&+\int_{\pi }^{0} \left ( E_{ApBpy} \frac{dE_{ApBpx}}{d\alpha } \right ) d\alpha \\
&+\int_{\pi }^{0} \left ( E_{ApCpy} \frac{dE_{ApCpx}}{d\alpha } \right ) d\alpha \\
&-\int_{\pi }^{\alpha _{1p}}\left ( A_{y} \frac{dA_{x}}{d\alpha } \right ) d\alpha \\
&-\int_{\alpha _{2p}}^{\pi}\left ( E_{ACy} \frac{dE_{ACx}}{d\alpha } \right ) d\alpha \\
\end{split}
\end{equation*}

We apply the Lagrange multiplier so that 
\begin{equation*}
\begin{split}
\text{maximize}&\int_{0}^{\pi} -r\left ( 0 \right ) +2t\left ( 0 \right ) -r\left ( \pi  \right ) +2t\left ( \pi  \right ) \\
&+\int_{\pi }^{0} \left ( E_{ApBpy} \frac{dE_{ApBpx}}{d\alpha } \right ) d\alpha \\
&+\int_{\pi }^{0} \left ( E_{ApCpy} \frac{dE_{ApCpx}}{d\alpha } \right ) d\alpha \\
&-\int_{\pi }^{\alpha _{1p}}\left ( A_{y} \frac{dA_{x}}{d\alpha } \right ) d\alpha \\
&-\int_{\alpha _{2p}}^{\pi}\left ( E_{ACy} \frac{dE_{ACx}}{d\alpha } \right ) d\alpha \\
&+\int_{0}^{\pi}\lambda_{2}\left ( \alpha \right )\left \{ -r\left ( \pi   \right )- \left [ r\left ( \alpha  \right ) \cos \alpha -t\left ( \alpha  \right )  2 \cos ^2\left(\dfrac{\alpha}{2}\right) \right ]   \right \}  d\alpha
\end{split} 
\end{equation*}

\end{definition}

\begin{definition}[EL equations]
Analogously to equations (15) and (16), we obtain the associated EL equations (see Appendix 6 for details) as

\begin{equation}
ODE1=\begin{cases}
  &-\cos \left(\dfrac{\alpha }{2}\right)+4 \cos (\alpha ) r(\alpha )+\sin \left(\dfrac{\alpha }{2}\right)-4 \cos (\alpha ) t(\alpha )+2 \lambda_{2}(\alpha )+8 \sin (\alpha ) r'(\alpha )\\
   &-4 \sin (\alpha ) t'(\alpha )-4 \cos (\alpha ) r''(\alpha ), 0<\alpha \leq \alpha_{1p} \\  \\
   
  &-\cos \left(\dfrac{\alpha }{2}\right)+4 \cos (\alpha ) r(\alpha )+\sin \left(\dfrac{\alpha }{2}\right)-2 \cos (\alpha ) t(\alpha )+2 \lambda_{2}(\alpha )+8 \sin (\alpha ) r'(\alpha )\\
   &-2 \sin (\alpha ) t'(\alpha )-4 \cos (\alpha ) r''(\alpha ), \alpha_{1p}<\alpha \leq \alpha_{2p}\\ \\
   
  &-\cos \left(\dfrac{\alpha }{2}\right)+2 (3 \cos (\alpha )+\cos (2 \alpha )) r(\alpha )+\sin \left(\dfrac{\alpha }{2}\right)-2 (2 \cos (\alpha )+\cos
   (2 \alpha )) t(\alpha )\\
   &+2 \lambda_{2}(\alpha )+12 \sin (\alpha ) r'(\alpha )+3 \sin (2 \alpha ) r'(\alpha )-4 \sin (\alpha ) t'(\alpha
   )-3 \sin (2 \alpha ) t'(\alpha )-r''(\alpha )\\
   &-6 \cos (\alpha ) r''(\alpha )-\cos (2 \alpha ) r''(\alpha )-t''(\alpha )+\cos (2 \alpha )
   t''(\alpha ), \alpha_{2p}<\alpha \leq \pi\\
\end{cases}
\end{equation}

\begin{equation}
ODE2=\begin{cases}
  &\dfrac{1}{2} \cos \left(\dfrac{\alpha }{2}\right) \sin (\alpha )+\dfrac{1}{2} \sin \left(\dfrac{\alpha }{2}\right) \sin (\alpha )+2 r(\alpha ) \sin
   ^2(\alpha )-2 \sin ^2(\alpha ) t(\alpha )+\lambda_{2}(\alpha )\\
   &+\cos (\alpha ) \lambda_{2}(\alpha )-2 \cos (\alpha ) \sin (\alpha
   ) r'(\alpha )+4 \cos (\alpha ) \sin (\alpha ) t'(\alpha )+2 \sin ^2(\alpha ) t''(\alpha ),\\
   & 0<\alpha \leq \alpha_{1p} \\  \\
   
  &\dfrac{1}{2} \cos \left(\dfrac{\alpha }{2}\right) \sin (\alpha )+\dfrac{1}{2} \sin \left(\dfrac{\alpha }{2}\right) \sin (\alpha )+r(\alpha ) \sin
   ^2(\alpha )-2 \sin ^2(\alpha ) t(\alpha )+\text{$\lambda $2}(\alpha )\\
   &+\cos (\alpha ) \lambda_{2}(\alpha )-\cos (\alpha ) \sin (\alpha )
   r'(\alpha )+4 \cos (\alpha ) \sin (\alpha ) t'(\alpha )+2 \sin ^2(\alpha ) t''(\alpha ),\\
   & \alpha_{1p}<\alpha \leq \pi- \alpha_{2p}\\ \\
   
  &\dfrac{1}{2} \left(\cos \left(\dfrac{\alpha }{2}\right) \sin (\alpha )+\sin \left(\dfrac{\alpha }{2}\right) \sin (\alpha )-4 (-1+\cos (\alpha ))
   r(\alpha ) \sin ^2(\alpha )\right.\\
   &\left.+2 (-3+2 \cos (\alpha )) \sin ^2(\alpha ) t(\alpha )+2 \lambda_{2}(\alpha )+2 \cos (\alpha ) \lambda_{2}(\alpha )-\sin (\alpha ) r'(\alpha )\right.\\
   &\left.-4 \cos (\alpha ) \sin (\alpha ) r'(\alpha )+3 \cos (2 \alpha ) \sin (\alpha ) r'(\alpha )-\sin (\alpha) t'(\alpha )-3 \cos (2 \alpha ) \sin (\alpha ) t'(\alpha )\right.\\
   &\left.+6 \sin (2 \alpha ) t'(\alpha )+\sin (\alpha ) \sin (2 \alpha ) r''(\alpha )+6 \sin
   ^2(\alpha ) t''(\alpha )-2 \cos (\alpha ) \sin ^2(\alpha ) t''(\alpha )\right), \\
   &\pi- \alpha_{2p}<\alpha \leq \pi/2\\
\end{cases}
\end{equation}

\begin{equation}
ODE3=\lambda_{2}(\alpha ) \left(r(\pi )-\left(r(\alpha ) \cos (\alpha )-t(\alpha ) 2 \cos ^2\left(\dfrac{\alpha }{2}\right)\right)\right)
\end{equation}
\end{definition}

The boundary values $r(0), t(0),r(\pi), t(\pi)$ are to be determined.

\begin{itemize}
\item \textbf{Scenario in Figure 8(c)}
\end{itemize}

\begin{definition}[moving sofa problem-to maximize area]

If there are intersection points of the trajectory of point A with both the envelopes of AB and AC , we use similar assumptions from equations (5)-(8) to calculate the area as
\begin{equation*}
\begin{split}
Area=& -r\left ( 0 \right ) +2t\left ( 0 \right )  -r\left ( \pi  \right ) +2t\left ( \pi  \right )\\
&+\int_{\pi }^{0} \left ( E_{ApBpy} \frac{dE_{ApBpx}}{d\alpha } \right ) d\alpha \\
&+\int_{\pi }^{0} \left ( E_{ApCpy} \frac{dE_{ApCpx}}{d\alpha } \right ) d\alpha \\
&-\int_{\alpha _{3p}}^{\alpha _{1p}}\left ( A_{y} \frac{dA_{x}}{d\alpha } \right ) d\alpha \\
&-\int_{0}^{\alpha _{4p}}\left ( E_{ABy} \frac{dE_{ABx}}{d\alpha } \right ) d\alpha \\
&-\int_{\alpha _{2p}}^{\pi}\left ( E_{ACy} \frac{dE_{ACx}}{d\alpha } \right ) d\alpha \\
\end{split}
\end{equation*}
\end{definition}

\begin{definition}[EL equations]

Following the same procedure, we derive the associated EL equations (see Appendix 6 for details) as

\begin{equation}
ODE1=\begin{cases}
  &-\cos \left(\dfrac{\alpha }{2}\right)-2 (-3 \cos (\alpha )+\cos (2 \alpha )) r(\alpha )+\sin \left(\dfrac{\alpha }{2}\right)+2 (-3 \cos (\alpha
   )+\cos (2 \alpha )) t(\alpha )\\
   &+12 \sin (\alpha ) r'(\alpha )-3 \sin (2 \alpha ) r'(\alpha )-6 \sin (\alpha ) t'(\alpha )+3 \sin (2 \alpha )
   t'(\alpha )+r''(\alpha )\\
   &-6 \cos (\alpha ) r''(\alpha )+\cos (2 \alpha ) r''(\alpha )+t''(\alpha )-\cos (2 \alpha ) t''(\alpha ), 0<\alpha \leq \alpha_{5p} \\  \\
   
  &-\cos \left(\dfrac{\alpha }{2}\right)+4 \cos (\alpha ) r(\alpha )+\sin \left(\dfrac{\alpha }{2}\right)-4 \cos (\alpha ) t(\alpha )+8 \sin (\alpha )
   r'(\alpha )-4 \sin (\alpha ) t'(\alpha )\\
   &-4 \cos (\alpha ) r''(\alpha ), \alpha_{5p}<\alpha \leq \alpha_{1p}\\ \\
   
  &-\cos \left(\dfrac{\alpha }{2}\right)+4 \cos (\alpha ) r(\alpha )+\sin \left(\dfrac{\alpha }{2}\right)-2 \cos (\alpha ) t(\alpha )+8 \sin (\alpha )
   r'(\alpha )-2 \sin (\alpha ) t'(\alpha )\\
   &-4 \cos (\alpha ) r''(\alpha ), \alpha_{1p}<\alpha \leq \alpha_{2p}\\ \\
   
   &-\cos \left(\dfrac{\alpha }{2}\right)+2 (3 \cos (\alpha )+\cos (2 \alpha )) r(\alpha )+\sin \left(\dfrac{\alpha }{2}\right)-2 (2 \cos (\alpha )+\cos
   (2 \alpha )) t(\alpha )\\
   &+12 \sin (\alpha ) r'(\alpha )+3 \sin (2 \alpha ) r'(\alpha )-4 \sin (\alpha ) t'(\alpha )-3 \sin (2 \alpha ) t'(\alpha
   )-r''(\alpha )\\
   &-6 \cos (\alpha ) r''(\alpha )-\cos (2 \alpha ) r''(\alpha )-t''(\alpha )+\cos (2 \alpha ) t''(\alpha ), \alpha_{2p}<\alpha \leq \alpha_{6p}\\ \\
   
   &-\cos \left(\dfrac{\alpha }{2}\right)+2 (3 \cos (\alpha )+\cos (2 \alpha )) r(\alpha )+\sin \left(\dfrac{\alpha }{2}\right)-2 (3 \cos (\alpha )+\cos
   (2 \alpha )) t(\alpha )\\
   &+12 \sin (\alpha ) r'(\alpha )+3 \sin (2 \alpha ) r'(\alpha )-6 \sin (\alpha ) t'(\alpha )-3 \sin (2 \alpha ) t'(\alpha
   )-r''(\alpha )\\
   &-6 \cos (\alpha ) r''(\alpha )-\cos (2 \alpha ) r''(\alpha )-t''(\alpha )+\cos (2 \alpha ) t''(\alpha ), \alpha_{6p}<\alpha \leq \pi\\
\end{cases}
\end{equation}

\begin{equation}
ODE2=\begin{cases}
  &\cos \left(\dfrac{\alpha }{2}\right)+\sin \left(\dfrac{\alpha }{2}\right)+2 r(\alpha ) (3 \sin (\alpha )+\sin (2 \alpha ))-2 (3 \sin (\alpha )+\sin
   (2 \alpha )) t(\alpha )+r'(\alpha )\\
   &-6 \cos (\alpha ) r'(\alpha )-3 \cos (2 \alpha ) r'(\alpha )+t'(\alpha )+12 \cos (\alpha ) t'(\alpha )+3
   \cos (2 \alpha ) t'(\alpha )\\
   &-\sin (2 \alpha ) r''(\alpha )+6 \sin (\alpha ) t''(\alpha )+\sin (2 \alpha ) t''(\alpha ), 0<\alpha \leq \alpha_{5p} \\  \\
   
  &\cos \left(\dfrac{\alpha }{2}\right)+\sin \left(\dfrac{\alpha }{2}\right)+4 r(\alpha ) \sin (\alpha )-4 \sin (\alpha ) t(\alpha )-4 \cos (\alpha )
   r'(\alpha )+8 \cos (\alpha ) t'(\alpha )\\
   &+4 \sin (\alpha ) t''(\alpha ), \alpha_{5p}<\alpha \leq \alpha_{1p}\\ \\
   
  &\cos \left(\dfrac{\alpha }{2}\right)+\sin \left(\dfrac{\alpha }{2}\right)+2 r(\alpha ) \sin (\alpha )-4 \sin (\alpha ) t(\alpha )-2 \cos (\alpha )
   r'(\alpha )+8 \cos (\alpha ) t'(\alpha )\\
   &+4 \sin (\alpha ) t''(\alpha ), \alpha_{1p}<\alpha \leq \alpha_{2p}\\ \\
   
   &\cos \left(\dfrac{\alpha }{2}\right)+\sin \left(\dfrac{\alpha }{2}\right)+r(\alpha ) (4 \sin (\alpha )-2 \sin (2 \alpha ))+2 (-3 \sin (\alpha )+\sin
   (2 \alpha )) t(\alpha )-r'(\alpha )\\
   &-4 \cos (\alpha ) r'(\alpha )+3 \cos (2 \alpha ) r'(\alpha )-t'(\alpha )+12 \cos (\alpha ) t'(\alpha )-3
   \cos (2 \alpha ) t'(\alpha )\\
   &+\sin (2 \alpha ) r''(\alpha )+6 \sin (\alpha ) t''(\alpha )-\sin (2 \alpha ) t''(\alpha ), \alpha_{2p}<\alpha \leq \alpha_{6p}\\ \\
   
   &\cos \left(\dfrac{\alpha }{2}\right)+\sin \left(\dfrac{\alpha }{2}\right)-2 r(\alpha ) (-3 \sin (\alpha )+\sin (2 \alpha ))+2 (-3 \sin (\alpha
   )+\sin (2 \alpha )) t(\alpha )\\
   &-r'(\alpha )-6 \cos (\alpha ) r'(\alpha )+3 \cos (2 \alpha ) r'(\alpha )-t'(\alpha )+12 \cos (\alpha ) t'(\alpha
   )-3 \cos (2 \alpha ) t'(\alpha )\\
   &+\sin (2 \alpha ) r''(\alpha )+6 \sin (\alpha ) t''(\alpha )-\sin (2 \alpha ) t''(\alpha ), \alpha_{6p}<\alpha \leq \pi\\
\end{cases}
\end{equation}
\end{definition}

The boundary values $r(0), t(0),r(\pi), t(\pi)$ are to be determined.

Utilizing numerical methods similar to those described in Section 5, we search for these boundary values and solve the resulting ODEs to determine the shape and area.

\subsection{Results of numerical calculations}

By searching for $r(0),r(\pi),t(0),t(\pi)$ and performing numerical calculations (see Appendix 6 for details), we obtain the sofa curves shown in Figure 9.

For Figure 8(a), the resulting shape consists of 3 curves and 3 line segments. The calculated results are close to $2/\pi\approx 0.63662$ and Hammersley’s sofa area $2/\pi+\pi/2\approx 2.20741$. For Figure 8(b), the resulting shape consists of 4 curves and 3 line segments, with a sofa area of 2.21952. For Figure 8(c), the resulting shape consists of 5 curves and 3 line segments, yielding a sofa area of 2.21953, consistent with Gerver's results. Figure 10 shows the corresponding functions $r(\alpha)$ and $t(\alpha)$.

\begin{figure}[H]
  \centering
  \includegraphics[width=10cm]{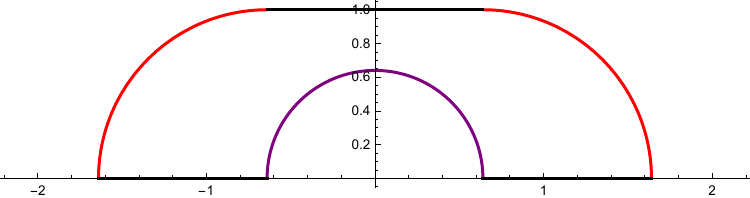}
  \includegraphics[width=10cm]{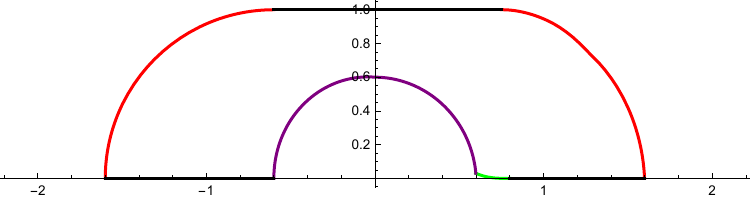}
  \includegraphics[width=10cm]{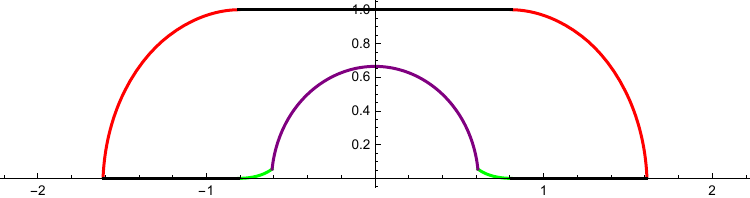}
  \caption{Results of numerical calculations of sofa shape for three asymmetric conditions in Figure 13.}
\end{figure}

\begin{figure}[H]
  \centering
  \includegraphics[width=10cm]{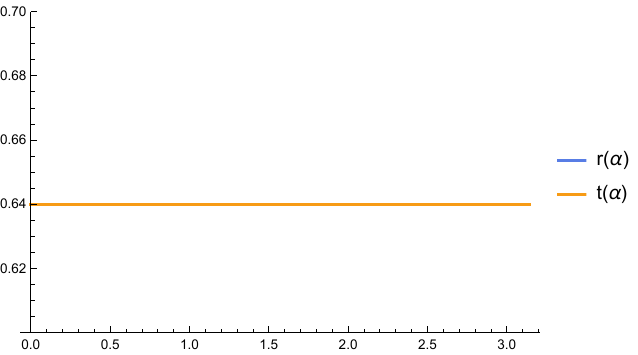}
  \includegraphics[width=10cm]{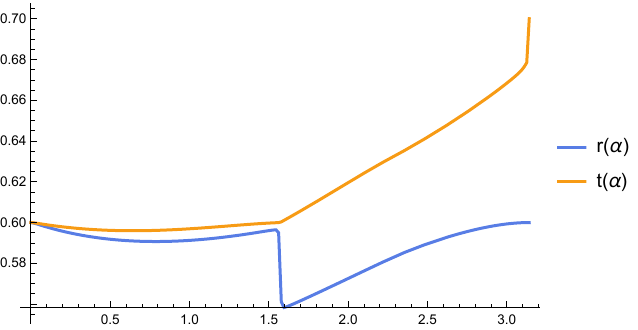}
  \includegraphics[width=10cm]{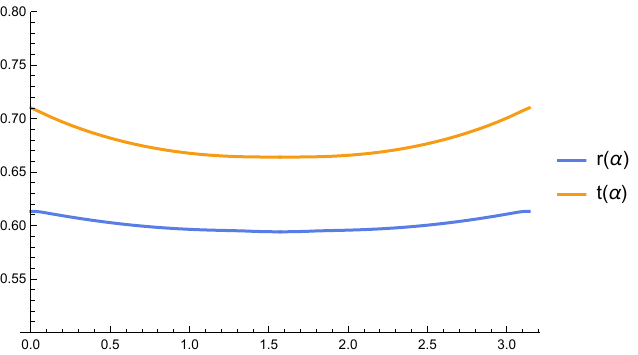}
  \caption{Results of numerical calculations of $r\left ( \alpha  \right ) ,t\left ( \alpha  \right )$ for three asymmetric conditions in Figure 13.}
\end{figure}

\section{Conclusion for moving sofa problem}
In the previous sections, we employed the calculus of variations to solve the moving sofa problem. The shape of the sofa is defined as the envelope curves formed by the walls it must navigate. The sofa's area is formulated as a functional of parametric equations for the curves. By solving the EL equations numerically, we determined the final shape and established the largest area. We verify that Gerver's constructed sofa shape matches our results and EL equations. We also discussed asymmetric conditions, demonstrating the generality and potential of our methods for solving related problems.

\section{Variant problem-Conway’s car}

\subsection{EL equations}
A variant of the problem involves finding a sofa maneuverable around both left and right corners. It is also called Conway’s car problem \cite{Croft2012} and ambidextrous sofa \cite{Romik2018}. Figure 11 illustrates the car shape obtained in the previous study\cite{Romik2018}.

\begin{figure}[H]
  \centering
  \includegraphics[width=9cm]{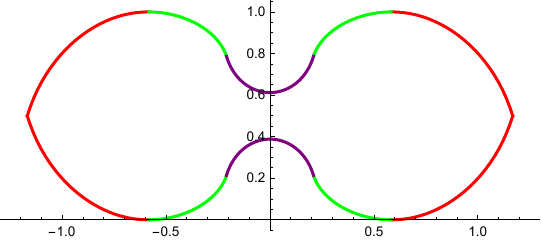}
  \caption{Shape of Conway’s car in the previous study\cite{Romik2018}.}
\end{figure}

\begin{boldremark}
Considering the existence of closed sofa, it is necessary to satisfy
\begin{equation*}
A_{y}\bigg|_{\alpha =\frac{\pi}{2}} \le \frac{1}{2}
\end{equation*}
\end{boldremark}

\begin{definition}[intersection points]
We assume there are two intersection points. 

\begin{equation*}
\left \{ A_{x},A_{y} \right \} \bigg|_{\alpha =\alpha_{1p}}=\left \{ E_{ACx},E_{ACy} \right \} \bigg|_{\alpha =\alpha_{2p}} 
\end{equation*}
\begin{equation*}
\left \{ A_{x},A_{y} \right \} \bigg|_{\alpha =\alpha_{3p}}=\left \{ E_{ABx},E_{ABy} \right \} \bigg|_{\alpha =\alpha_{4p}} 
\end{equation*}
\begin{equation*}
E_{ACy} \bigg|_{\alpha =\alpha_{5p}}=\frac{1}{2} 
\end{equation*}
\begin{equation*}
E_{ABy} \bigg|_{\alpha =\alpha_{6p}}=\frac{1}{2} 
\end{equation*}
\end{definition}

\begin{definition}[area]

Using a similar method as in equation (4), and under similar assumptions as in equations (5)-(8), the area can be expressed as

\begin{equation*}
\begin{split}
Area=&-r\left ( 0 \right ) +2t\left ( 0 \right ) -r\left ( \pi  \right ) +2t\left ( \pi  \right ) \\
&-2\int_{0}^{\alpha _{4p}}\left ( E_{ABy} \frac{dE_{ABx}}{d\alpha } \right ) d\alpha \\
&-2\int_{\alpha _{2p}}^{\pi}\left ( E_{ACy} \frac{dE_{ACx}}{d\alpha } \right ) d\alpha \\
&-2\int_{\alpha _{3p}}^{\alpha _{1p}}\left ( A_{y} \frac{dA_{x}}{d\alpha } \right ) d\alpha \\
&-2\int_{0}^{\alpha _{6p}} \left [ \left ( E_{ApBpy}-\frac{1}{2}\right )  \frac{dE_{ApBpx}}{d\alpha } \right ]   d\alpha \\
&-2\int_{0}^{\alpha _{5p}} \left [ \left ( E_{ApCpy}-\frac{1}{2}\right ) \frac{dE_{ApCpx}}{d\alpha } \right ]   d\alpha \\
\end{split}
\end{equation*}

If assuming symmetry of $\alpha=\pi /2$
\begin{equation*}
\begin{split}
Area=&-2r\left ( 0 \right ) +4t\left ( 0 \right )  \\
&-4\int_{0}^{\pi /2}\left ( E_{ABy} \frac{dE_{ABx}}{d\alpha } \right ) d\alpha \\
&-4\int_{\alpha _{2p}}^{\pi/2}\left ( E_{ACy} \frac{dE_{ACx}}{d\alpha } \right ) d\alpha \\
&+4\int_{\alpha _{1p}}^{\pi/2}\left ( A_{y} \frac{dA_{x}}{d\alpha } \right ) d\alpha \\
&-4\int_{0}^{\pi/2} \left [ \left ( E_{ApBpy}-\frac{1}{2}\right )  \frac{dE_{ApBpx}}{d\alpha } \right ]   d\alpha \\
&-4\int_{\alpha _{5p}}^{\pi/2} \left [ \left ( E_{ApCpy}-\frac{1}{2}\right ) \frac{dE_{ApCpx}}{d\alpha } \right ]   d\alpha \\
\end{split}
\end{equation*}
\end{definition}

Without loss of generality, we may as well assume $0< \alpha_{5p}\le \alpha_{1p}\le 1\le \alpha_{2p}< \dfrac{\pi}{2} $.

\begin{definition}[EL equations]
Similarly, we can construct piecewise functions. To maximize the area functional, we derive the EL equations (See Appendix 7 for details) as

\begin{equation}
ODE1=\begin{cases}
  &1+8 r(\alpha ) \sin \left(\dfrac{3 \alpha }{2}\right)-8 \sin \left(\dfrac{3 \alpha }{2}\right) t(\alpha )+4 \cos \left(\dfrac{\alpha }{2}\right)
   r'(\alpha )-12 \cos \left(\dfrac{3 \alpha }{2}\right) r'(\alpha )\\
   &+4 \cos \left(\dfrac{\alpha }{2}\right) t'(\alpha )+12 \cos \left(\dfrac{3
   \alpha }{2}\right) t'(\alpha )+4 \sin \left(\dfrac{\alpha }{2}\right) r''(\alpha )-4 \sin \left(\dfrac{3 \alpha }{2}\right) r''(\alpha ))\\
   &+4 \sin
   \left(\dfrac{\alpha }{2}\right) t''(\alpha +4 \sin \left(\dfrac{3 \alpha }{2}\right) t''(\alpha ),0<\alpha \leq \alpha_{5p} \\  \\
   
 &-\cos \left(\dfrac{\alpha }{2}\right)-2 (-3 \cos (\alpha )+\cos (2 \alpha )) r(\alpha )+\sin \left(\dfrac{\alpha }{2}\right)+2 (-3 \cos (\alpha
   )+\cos (2 \alpha )) t(\alpha )\\
   &+12 \sin (\alpha ) r'(\alpha )-3 \sin (2 \alpha ) r'(\alpha )-6 \sin (\alpha ) t'(\alpha )+3 \sin (2 \alpha )
   t'(\alpha )+r''(\alpha )\\
   &-6 \cos (\alpha ) r''(\alpha )+\cos (2 \alpha ) r''(\alpha )+t''(\alpha )-\cos (2 \alpha ) t''(\alpha ), \alpha_{5p}<\alpha \leq \alpha_{1p}\\ \\

  &-\cos \left(\dfrac{\alpha }{2}\right)-2 (-3 \cos (\alpha )+\cos (2 \alpha )) r(\alpha )+\sin \left(\dfrac{\alpha }{2}\right)+2 (-2 \cos (\alpha
   )+\cos (2 \alpha )) t(\alpha )\\
   &+12 \sin (\alpha ) r'(\alpha )-3 \sin (2 \alpha ) r'(\alpha )-4 \sin (\alpha ) t'(\alpha )+3 \sin (2 \alpha )
   t'(\alpha )+r''(\alpha )\\
   &-6 \cos (\alpha ) r''(\alpha )+\cos (2 \alpha ) r''(\alpha )+t''(\alpha )-\cos (2 \alpha ) t''(\alpha ), \alpha_{1p}<\alpha \leq \alpha_{2p}\\ \\
   
  &-\cos \left(\dfrac{\alpha }{2}\right)+8 \cos (\alpha ) r(\alpha )+\sin \left(\dfrac{\alpha }{2}\right)-6 \cos (\alpha ) t(\alpha )+16 \sin (\alpha
   ) r'(\alpha )-6 \sin (\alpha ) t'(\alpha )\\
   &-8 \cos (\alpha ) r''(\alpha ),\alpha_{2p}<\alpha \leq \pi/2\\
\end{cases}
\end{equation}

\begin{equation}
ODE2=\begin{cases}
  &\cos \left(\dfrac{\alpha }{2}\right)+4 r(\alpha ) (\sin (\alpha )+\sin (2 \alpha ))-4 (\sin (\alpha )+\sin (2 \alpha )) t(\alpha )+2 r'(\alpha )-4
   \cos (\alpha ) r'(\alpha )\\
   &-6 \cos (2 \alpha ) r'(\alpha )+2 t'(\alpha )+8 \cos (\alpha ) t'(\alpha )+6 \cos (2 \alpha ) t'(\alpha )-2 \sin (2
   \alpha ) r''(\alpha )\\
   &+4 \sin (\alpha ) t''(\alpha )+2 \sin (2 \alpha ) t''(\alpha ), 0<\alpha \leq \alpha_{5p} \\  \\
   
 &\cos \left(\dfrac{\alpha }{2}\right)+\sin \left(\dfrac{\alpha }{2}\right)+2 r(\alpha ) (3 \sin (\alpha )+\sin (2 \alpha ))-2 (3 \sin (\alpha )+\sin
   (2 \alpha )) t(\alpha )+r'(\alpha )\\
   &-6 \cos (\alpha ) r'(\alpha )-3 \cos (2 \alpha ) r'(\alpha )+t'(\alpha )+12 \cos (\alpha ) t'(\alpha )+3
   \cos (2 \alpha ) t'(\alpha )\\
   &-\sin (2 \alpha ) r''(\alpha )+6 \sin (\alpha ) t''(\alpha )+\sin (2 \alpha ) t''(\alpha ), \alpha_{5p}<\alpha \leq \alpha_{1p}\\ \\

  &\cos \left(\dfrac{\alpha }{2}\right)+\sin \left(\dfrac{\alpha }{2}\right)+4 (1+\cos (\alpha )) r(\alpha ) \sin (\alpha )-2 (3 \sin (\alpha )+\sin
   (2 \alpha )) t(\alpha )+r'(\alpha )\\
   &-4 \cos (\alpha ) r'(\alpha )-3 \cos (2 \alpha ) r'(\alpha )+t'(\alpha )+12 \cos (\alpha ) t'(\alpha )+3
   \cos (2 \alpha ) t'(\alpha )\\
   &-\sin (2 \alpha ) r''(\alpha )+6 \sin (\alpha ) t''(\alpha )+\sin (2 \alpha ) t''(\alpha ), \alpha_{1p}<\alpha \leq \alpha_{2p}\\ \\
   
  &\cos \left(\dfrac{\alpha }{2}\right)+\sin \left(\dfrac{\alpha }{2}\right)+6 r(\alpha ) \sin (\alpha )-8 \sin (\alpha ) t(\alpha )-6 \cos (\alpha )
   r'(\alpha )+16 \cos (\alpha ) t'(\alpha )\\
   &+8 \sin (\alpha ) t''(\alpha ),\alpha_{2p}<\alpha \leq \pi/2\\
\end{cases}
\end{equation}
\end{definition}

Similarly, we can also discuss the cases of no intersection point and asymmetry conditions for the car problem, similar to Sections 7 and 8, as illustrated in Figures 12 and 13. These are straightforward for readers.

\begin{figure}[H]
  \centering
  \includegraphics[width=10cm]{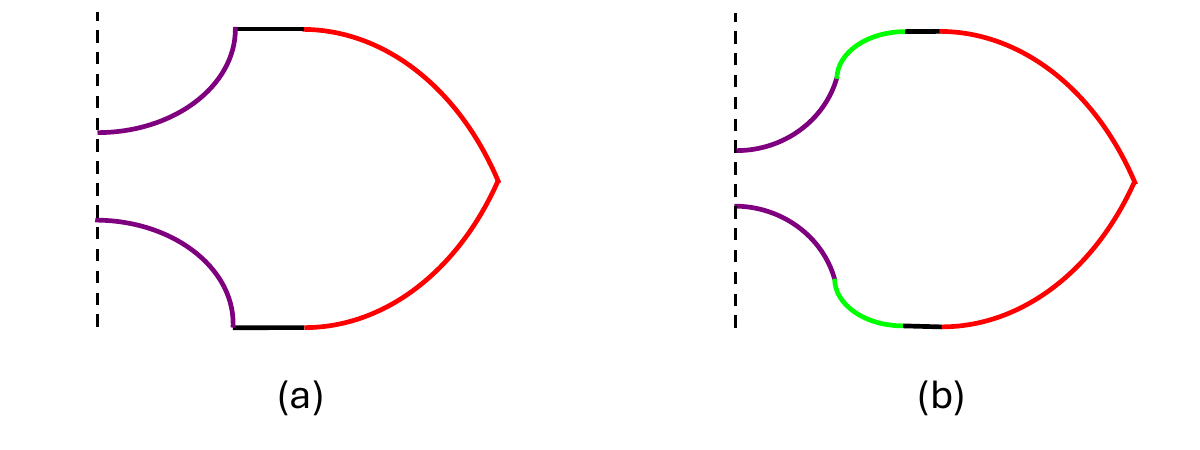}
  \caption{Cases of no intersection point or more intersection points for Conway’s car/ambidextrous sofa.}
\end{figure}

\begin{figure}[H]
  \centering
  \includegraphics[width=12cm]{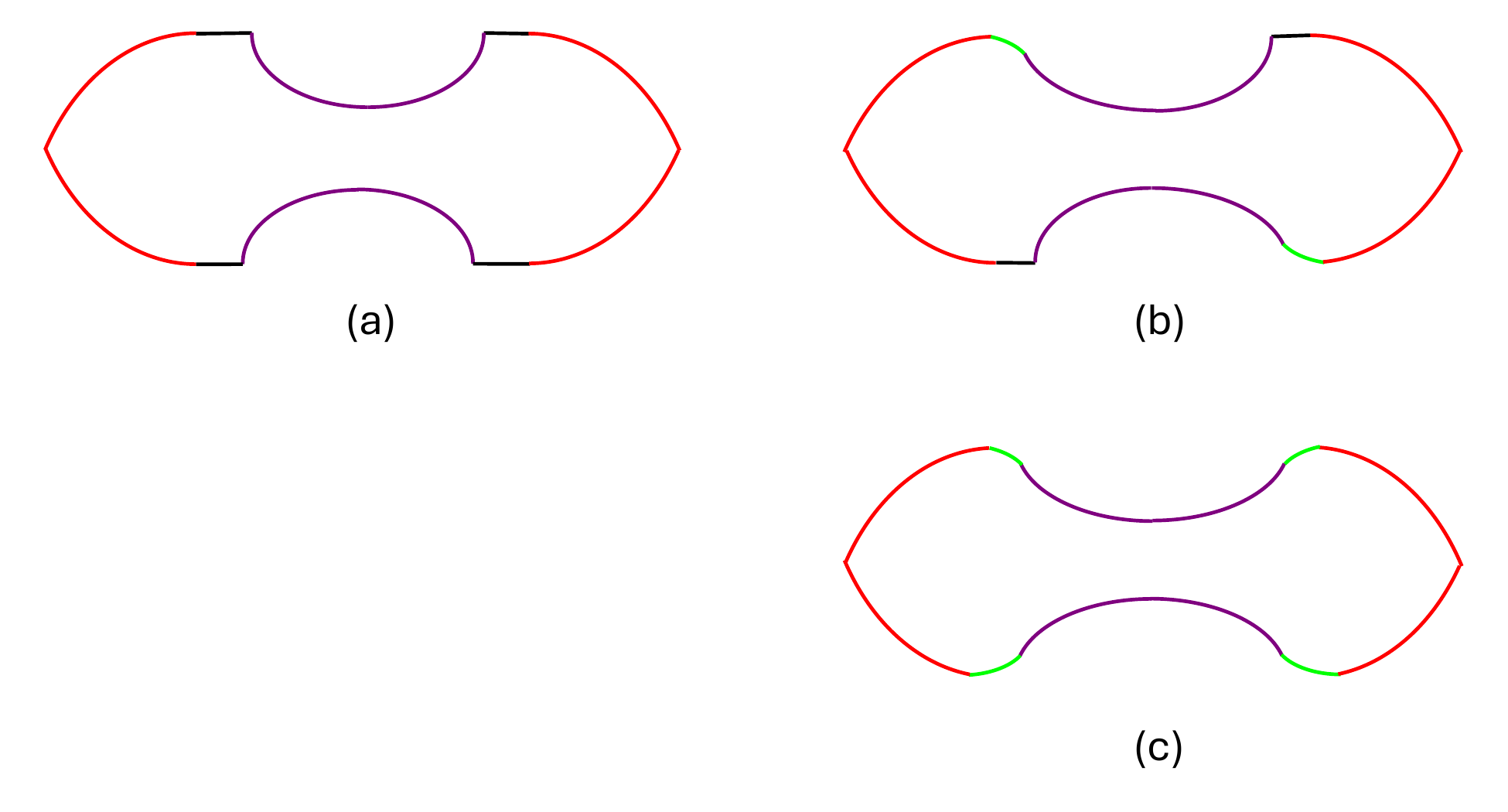}
  \caption{Asymmetry conditions for Conway’s car/ambidextrous sofa.}
\end{figure}

\subsection{Verification of Romik’s results}
\begin{definition}[simplified EL equations for car problem]

As in Section 6, we can verify that the explicit results obtained by Romik \cite{Romik2018} satisfy the EL equations. Specifically, from \cite{Romik2018}, we have $\alpha _{1p}=\alpha _{2p}=\alpha _{5p}=2\beta\approx 0.5793076416$ then the EL equations (33) (34) can be simplified as
\begin{equation}
ODE1=\begin{cases}
  &1+8 r(\alpha ) \sin \left(\dfrac{3 \alpha }{2}\right)-8 \sin \left(\dfrac{3 \alpha }{2}\right) t(\alpha )+4 \cos \left(\dfrac{\alpha }{2}\right)
   r'(\alpha )-12 \cos \left(\dfrac{3 \alpha }{2}\right) r'(\alpha )\\
   &+4 \cos \left(\dfrac{\alpha }{2}\right) t'(\alpha )+12 \cos \left(\dfrac{3
   \alpha }{2}\right) t'(\alpha )+4 \sin \left(\dfrac{\alpha }{2}\right) r''(\alpha )-4 \sin \left(\dfrac{3 \alpha }{2}\right) r''(\alpha )\\
   &+4 \sin
   \left(\dfrac{\alpha }{2}\right) t''(\alpha )+4 \sin \left(\dfrac{3 \alpha }{2}\right) t''(\alpha ),0<\alpha \leq 2\beta\\  
   
 &-\cos \left(\dfrac{\alpha }{2}\right)+8 \cos (\alpha ) r(\alpha )+\sin \left(\dfrac{\alpha }{2}\right)-6 \cos (\alpha ) t(\alpha )+16 \sin (\alpha
   ) r'(\alpha )-6 \sin (\alpha ) t'(\alpha )\\
   &-8 \cos (\alpha ) r''(\alpha ),2\beta<\alpha \leq \pi/2\\
\end{cases}
\end{equation}

\begin{equation}
ODE2=\begin{cases}
  &\cos \left(\dfrac{\alpha }{2}\right)+4 r(\alpha ) (\sin (\alpha )+\sin (2 \alpha ))-4 (\sin (\alpha )+\sin (2 \alpha )) t(\alpha )+2 r'(\alpha
   )-4 \cos (\alpha ) r'(\alpha )\\
   &-6 \cos (2 \alpha ) r'(\alpha )+2 t'(\alpha )+8 \cos (\alpha ) t'(\alpha )+6 \cos (2 \alpha ) t'(\alpha )-2
   \sin (2 \alpha ) r''(\alpha )\\
   &+4 \sin (\alpha ) t''(\alpha )+2 \sin (2 \alpha ) t''(\alpha ),0<\alpha \leq 2\beta\\  
   
 &\cos \left(\dfrac{\alpha }{2}\right)+\sin \left(\dfrac{\alpha }{2}\right)+6 r(\alpha ) \sin (\alpha )-8 \sin (\alpha ) t(\alpha )-6 \cos (\alpha )
   r'(\alpha )+16 \cos (\alpha ) t'(\alpha )\\
   &+8 \sin (\alpha ) t''(\alpha ),2\beta<\alpha \leq \pi/2\\
\end{cases}
\end{equation}
\end{definition}

\begin{definition}[recall Romik's results of car]

Recall the explicit and analytic formulas and symbols from the previous paper \cite{Romik2018} for parametric equation of the trajectory of point A is
\begin{equation*}
\mathbf{x} \left ( \alpha  \right ) =
\begin{cases} 
  \mathbf{x}_{1} \left ( \alpha/2 \right )-\left \{ \dfrac{1}{2} \left(1-e_1+\kappa _{5,1}\right),0 \right \},\text{if } 0\le \alpha< 2\beta \\ 
  \mathbf{x}_{6} \left ( \alpha/2 \right )-\left \{ \dfrac{1}{2} \left(1-e_1+\kappa _{5,1}\right),0 \right \},\text{if } 2\beta\le \alpha\le \pi/2 
\end{cases}
\end{equation*}

The explicit formulas are
\begin{equation*}
\begin{split}
&\mathbf{x}_{1} \left ( \alpha/2 \right )-\left \{ \dfrac{1}{2} \left(1-e_1+\kappa _{5,1}\right),0 \right \}=\\&\left\{-\sin \left(\frac{\alpha }{2}\right) \left(-\dfrac{1}{2}+\sin \left(\frac{\alpha }{2}\right) a_1-\cos \left(\frac{\alpha }{2}\right)
   a_2\right)+\cos \left(\frac{\alpha }{2}\right) \left(-1+\cos \left(\frac{\alpha }{2}\right) a_1+\sin \left(\frac{\alpha }{2}\right)
   a_2\right)\right.\\
&+\kappa _{1,1}+\frac{1}{2} \left(-1+e_1-\kappa _{5,1}\right),\cos \left(\frac{\alpha }{2}\right) \left(-\frac{1}{2}+\sin
   \left(\frac{\alpha }{2}\right) a_1-\cos \left(\frac{\alpha }{2}\right) a_2\right)\\
&\left.+\sin \left(\frac{\alpha }{2}\right) \left(-1+\cos
   \left(\frac{\alpha }{2}\right) a_1+\sin \left(\frac{\alpha }{2}\right) a_2\right)+\kappa _{1,2}\right\}\\
\end{split}
\end{equation*}

\begin{equation*}
\begin{split}
&\mathbf{x}_{6} \left ( \alpha/2 \right )-\left \{ \frac{1}{2} \left(1-e_1+\kappa _{5,1}\right),0 \right \}=\\
&\left\{-\sin \left(\frac{\alpha }{2}\right) \left(-1+\sin \left(\frac{\alpha }{4}\right) f_1-\cos \left(\frac{\alpha }{4}\right) f_2\right)+\cos
   \left(\frac{\alpha }{2}\right) \left(-1+\cos \left(\frac{\alpha }{4}\right) f_1+\sin \left(\frac{\alpha }{4}\right) f_2\right)\right.\\
&+\frac{1}{2}
   \left(-1+e_1-\kappa _{5,1}\right)+\kappa _{6,1},\left.\cos \left(\frac{\alpha }{2}\right) \left(-1+\sin \left(\frac{\alpha }{4}\right) f_1-\cos
   \left(\frac{\alpha }{4}\right) f_2\right)\right.\\
&\left.+\sin \left(\frac{\alpha }{2}\right) \left(-1+\cos \left(\frac{\alpha }{4}\right) f_1+\sin
   \left(\frac{\alpha }{4}\right) f_2\right)+\kappa _{6,2}\right\}\\
\end{split}
\end{equation*}
\end{definition}

\begin{theorem}
Romik’s car results satisfy the EL equations for the necessary condition of maximum area.
\end{theorem}

\begin{proof}
Converting the symbol and coordinate system, we can prove through symbolic calculations that when
\begin{equation*}
\left \{ r\left ( \alpha  \right )\cos \alpha,t\left ( \alpha  \right )\sin \alpha \right \} = \mathbf{x}_{1} \left ( \alpha/2 \right )-\left \{ \frac{1}{2} \left(1-e_1+\kappa _{5,1}\right),0 \right \},\text{if } 0\le \alpha< 2\beta
\end{equation*}
and 
\begin{equation*}
\left \{ r\left ( \alpha  \right )\cos \alpha,t\left ( \alpha  \right )\sin \alpha \right \} = \mathbf{x}_{6} \left ( \alpha/2 \right )-\left \{ \frac{1}{2} \left(1-e_1+\kappa _{5,1}\right),0 \right \},\text{if } 2\beta\le \alpha\le \pi/2
\end{equation*}
$\left \{ r\left ( \alpha  \right ),t\left ( \alpha  \right ) \right \}$ satisfy ODE1 in equation (35) and ODE2 in equation (36) (see Appendix 7 for details). 
\end{proof}

Therefore, we verify that Romik's results of ambidextrous car \cite{Romik2018} satisfy the necessary condition of the maximum.

\section{Discussion}
In this paper, we also addressed most of the open problems raised in previous paper \cite{Romik2018}, such as proof of Gerver’s and Romik’s results as local maxima of the area functional, as well as other maximum shapes for different intersection conditions and asymmetrical conditions.

\subsection{Sufficient conditions for maximum}
The sufficient condition for a maximum is that the second variation is negative semidefinite. We need to calculate the second variation.

\begin{definition}[second variation of area functional]

To find the second variation \( \delta^2 F \), we expand 

$F\left ( r+\epsilon \eta ,r'+\epsilon \eta',r''+\epsilon \eta'',t+\epsilon \xi,t'+\epsilon \xi',t''+\epsilon \xi''  \right )$ up to second order in \( \epsilon \). So that
\begin{equation}
\begin{split}
\delta^2 F = &\int \Bigg[ \eta \frac{\partial^2 F}{\partial r^2} \eta + 2 \eta \frac{\partial^2 F}{\partial r \partial r'} \eta' + 2 \eta \frac{\partial^2 F}{\partial r \partial r''} \eta'' + \eta' \frac{\partial^2 F}{\partial r'^2} \eta' + 2 \eta' \frac{\partial^2 F}{\partial r' \partial r''} \eta'' + \eta'' \frac{\partial^2 F}{\partial r''^2} \eta'' \\
&+ \xi \frac{\partial^2 F}{\partial t^2} \xi + 2 \xi \frac{\partial^2 F}{\partial t \partial t'} \xi' + 2 \xi \frac{\partial^2 F}{\partial t \partial t''} \xi'' + \xi' \frac{\partial^2 F}{\partial t'^2} \xi' + 2 \xi' \frac{\partial^2 F}{\partial t' \partial t''} \xi'' + \xi'' \frac{\partial^2 F}{\partial t''^2} \xi'' \\
&+ 2 \eta \frac{\partial^2 F}{\partial r \partial t} \xi + 2 \eta' \frac{\partial^2 F}{\partial r' \partial t} \xi + 2 \eta'' \frac{\partial^2 F}{\partial r'' \partial t} \xi + 2 \eta \frac{\partial^2 F}{\partial r \partial t'} \xi' + 2 \eta \frac{\partial^2 F}{\partial r \partial t''} \xi''\\
&+ 2 \eta' \frac{\partial^2 F}{\partial r' \partial t'} \xi' + 2 \eta' \frac{\partial^2 F}{\partial r' \partial t''} \xi'' + 2 \eta'' \frac{\partial^2 F}{\partial r'' \partial t'} \xi' + 2 \eta'' \frac{\partial^2 F}{\partial r'' \partial t''} \xi'' \Bigg] d\alpha  
\end{split}
\end{equation}
\end{definition}

\begin{definition}[Hessian matrix]

The sufficient condition for the second variation in equation (37) is related to the negative semidefiniteness of the Hessian matrix as

\begin{equation*}
\mathbf{H} =\begin{bmatrix}
  \dfrac{\partial^2 F}{\partial r^2}         & \dfrac{\partial^2 F}{\partial r\partial r'} & \dfrac{\partial^2 F}{\partial r\partial r''} & \dfrac{\partial^2 F}{\partial r\partial t} & \dfrac{\partial^2 F}{\partial r\partial t'} & \dfrac{\partial^2 F}{\partial r\partial t''}\\
  \dfrac{\partial^2 F}{\partial r'\partial r}& \dfrac{\partial^2 F}{\partial r'^2} & \dfrac{\partial^2 F}{\partial r'\partial r''} & \dfrac{\partial^2 F}{\partial r'\partial t} & \dfrac{\partial^2 F}{\partial r'\partial t'} & \dfrac{\partial^2 F}{\partial r'\partial t''}\\
  \dfrac{\partial^2 F}{\partial r''\partial r}& \dfrac{\partial^2 F}{\partial r''\partial r'} & \dfrac{\partial^2 F}{\partial r''^2} & \dfrac{\partial^2 F}{\partial r''\partial t} & \dfrac{\partial^2 F}{\partial r''\partial t'} & \dfrac{\partial^2 F}{\partial r''\partial t''}\\
  \dfrac{\partial^2 F}{\partial t\partial r} & \dfrac{\partial^2 F}{\partial t\partial r'} & \dfrac{\partial^2 F}{\partial t\partial r''} & \dfrac{\partial^2 F}{\partial t^2} & \dfrac{\partial^2 F}{\partial t\partial t'} & \dfrac{\partial^2 F}{\partial t\partial t''}\\
  \dfrac{\partial^2 F}{\partial t'\partial r}& \dfrac{\partial^2 F}{\partial t'\partial r'} & \dfrac{\partial^2 F}{\partial t'\partial r''} & \dfrac{\partial^2 F}{\partial t'\partial t} & \dfrac{\partial^2 F}{\partial t'^2} & \dfrac{\partial^2 F}{\partial t'\partial t''}\\
  \dfrac{\partial^2 F}{\partial t''\partial r}& \dfrac{\partial^2 F}{\partial t''\partial r'} & \dfrac{\partial^2 F}{\partial t''\partial r''} & \dfrac{\partial^2 F}{\partial t''\partial t} & \dfrac{\partial^2 F}{\partial t''\partial t'} & \dfrac{\partial^2 F}{\partial t''^2}
\end{bmatrix}
\end{equation*}
\end{definition}

All the Hessian Matrices can be found in Appendix 8 for details.

\subsection{More variant problems}
This variational method presented in this paper is general, based on Definitions 3.3-3.5 and 4.5-4.6, and has the potential to solve more variant sofa problems, such as the non-right angle problem \cite{Maruyama1973}. It may also address more general problems as illustrated in Figure 14, such as the rounded hallway, curved hallway, the Piano Mover's Problem \cite{PianoMovers}, and even three-dimensional problems \cite{Croft2012}. We believe that the calculus of variations method can solve more unsolved and open geometric optimization problems involving rotation and translation.

\begin{figure}[H]
  \centering
  \includegraphics[width=3cm]{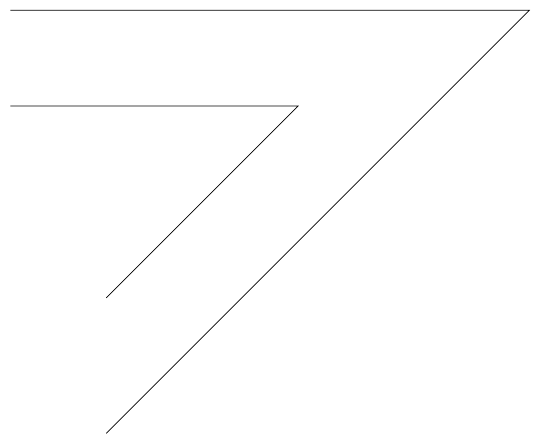}
  \includegraphics[width=3cm]{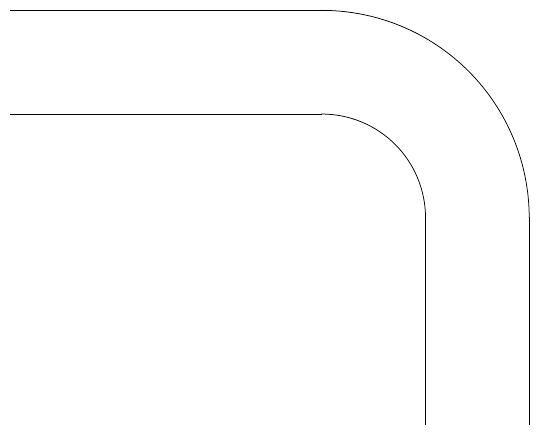}
  \includegraphics[width=3cm]{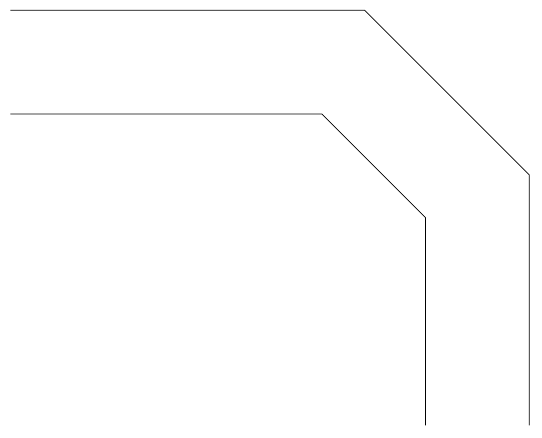}
  \includegraphics[width=4cm]{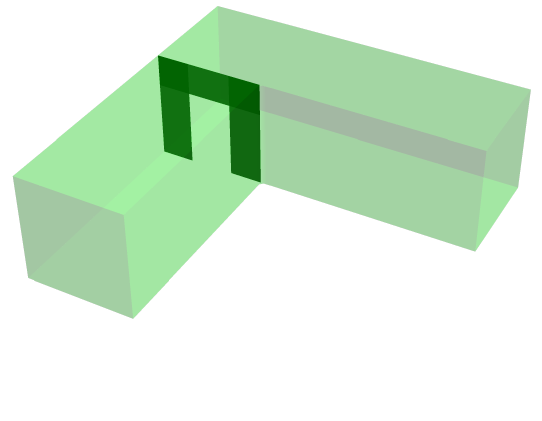}
  \caption{More variant problems with different hallways and 3D problems.}
\end{figure}

\subsection{More complex paths and shapes}
\begin{figure}[H]
  \centering
  \includegraphics[width=12cm]{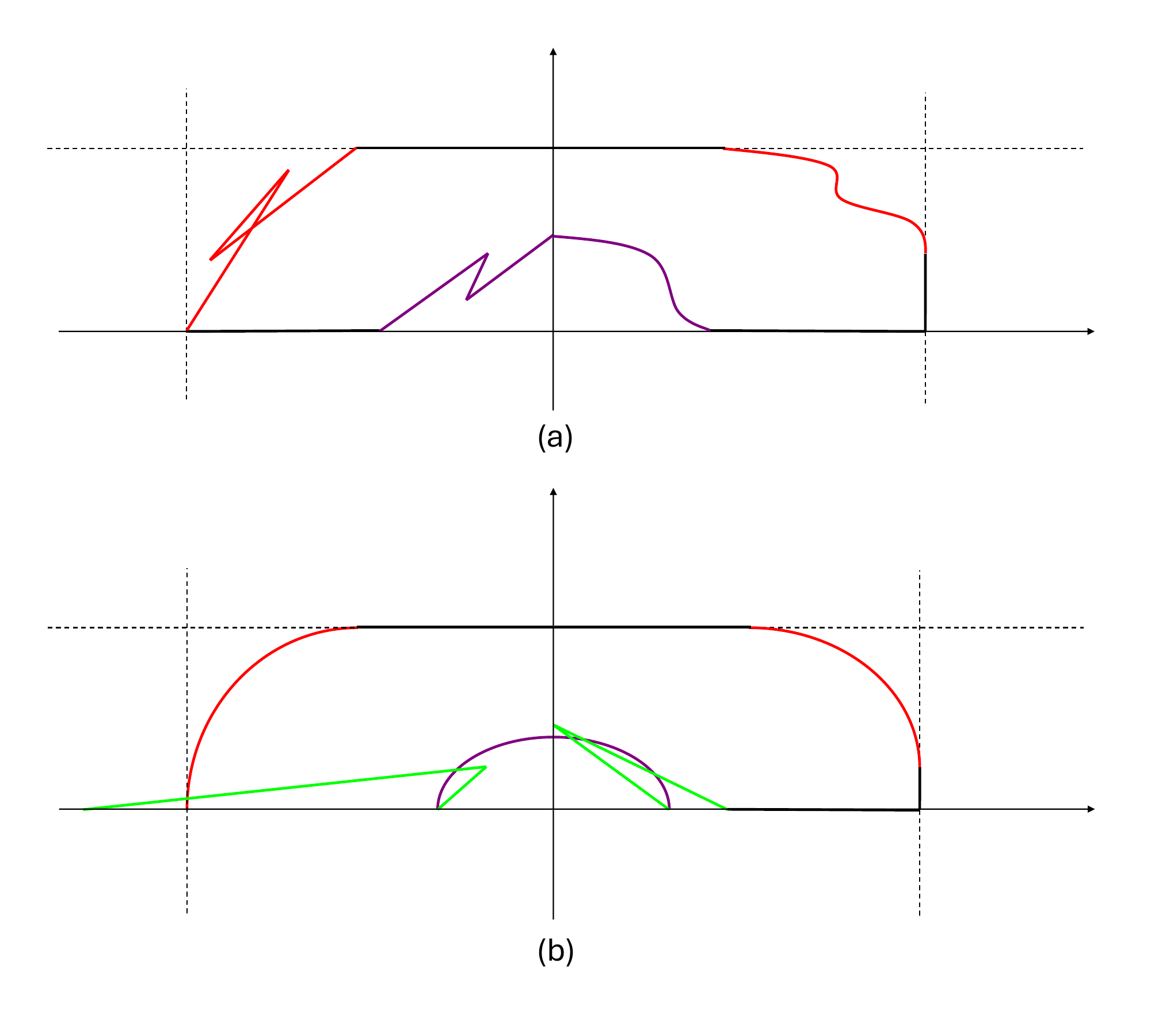}
  \caption{More complex paths and shapes.}
\end{figure}

There are more complex paths and envelopes forming the sofa shapes when the calculus of variations method in this paper cannot work well. These cases require further discussion:
\begin{enumerate}
   \item Self-intersecting paths and envelopes, as shown in Figure 15(a), where the integral calculation area cannot be applied. In this case, the integral in equation (4) needs to be segmented at the intersection point, or further discussion is required to exclude such self-intersection situations from being considered as the maximum area.
   
   \item Similarly, zigzag paths and envelopes, as shown in Figure 15(a). In this case, the integral in equation (4) needs to be segmented at the intersection point, or discussion is needed to exclude such situations from being considered as the maximum area.
       
    \item Non-convex paths and envelopes, as illustrated in Figures 15(a) and 11(b). In this case, some assumptions of boundary constraints and the use of integrals for the area in equation (4) do not hold.
   
   \item Multiple intersection points of the trajectory of point A and the envelope of AC or AB, as shown in Figure 15(b). In this case, the integral in equation (4) needs to be segmented at the multiple intersection points.
   
   \item Intersection of the envelope of ApCp with the envelope of AC outside, as shown in Figure 15(b). In this case, the integral upper and lower limits in equation (4) need to be revised at the intersection point.
   
   \item Part of the envelope of AC or AB below the Y-axis, as shown in Figure 3. In this case, the integral upper and lower limits in equation (4) need to be revised with y = 0.
   
   \item Intersection of the envelope of AC or AB with the X-axis, as shown in Figure 3. In this case, the integral upper and lower limits in equation (4) need to be revised with x = 0.
   
   \begin{boldremark}
The above four cases also complicate the discussion of symmetry and asymmetry conditions in Sections 7 and 8.
\end{boldremark}

   \item Cusp points. In this paper, we assume the differentiability of the function involved when directly applying the classical EL equations, since the functional contains first and second derivatives. This situation arises in problems involving curves that may not be smooth everywhere, such as corners or cusps in the paths or envelopes for the sofa problem. This can be justified by using more sophisticated tools from geometric analysis and geometric measure theory.
\end{enumerate}

%\section{Acknowledgments}

\section{Appendices-Mathematica notebooks}
The appendices provide Mathematica notebooks that contain detailed equation derivations and numerical calculations presented in this paper.
\\
\\
\textbf{Appendix 1} Mathematica notebook for the derivation of EL equations 
\\
\\
\textbf{Appendix 2} Mathematica notebook for searching optimal boundary value t(0)
\\
\\
\textbf{Appendix 3} Mathematica notebook for numerical results of sofa problem
\\
\\
\textbf{Appendix 4} Mathematica notebook for verifying that Gerver’s results satisfy EL equations
\\
\\
\textbf{Appendix 5} Mathematica notebook for the case of no intersection points 
\\
\\
\textbf{Appendix 6} Mathematica notebook for case of asymmetric conditions 
\\
\\
\textbf{Appendix 7} Mathematica notebook for verifying Romik's results for the ambidextrous sofa problem
\\
\\
\textbf{Appendix 8} Mathematica notebook for calculating the Hessian Matrix for all conditions



\section*{Supplementary material}

# Appendix 1

```mathematica
In[ ]:= (*define coordinates of points and envelopes*)

In[ ]:= AA = {r[α] * Cos[α], t[α] * Sin[α]};
       AAp = {r[α] * Cos[α] + √2 * Cos[π/4 + α/2], t[α] * Sin[α] + √2 * Sin[π/4 + α/2]};
       BB = {r[α] * Cos[α] - t[α] * 2 Cos[α/2]^2, 0};
       BBp = {r[α] * Cos[α] - t[α] * 2 Cos[α/2]^2 - 1/Sin[α/2], 0};
       CC = {r[α] * Cos[α] + t[α] * Sin[α] * Tan[α/2], 0};
       CCp = {r[α] * Cos[α] + t[α] * Sin[α] * Tan[α/2] + 1/Cos[α/2], 0};
       EAABBx = 1/2 * (-1 + 2 Cos[α] + Cos[2 α]) r[α] -
           2 Cos[α/2]^2 Cos[α] t[α] + Sin[α] (Cos[α] r'[α] - (1 + Cos[α]) t'[α]);
       EAABBy = -2 Sin[α/2]^2 (r[α] Sin[α] - Sin[α] t[α] - Cos[α] r'[α] + t'[α] + Cos[α] t'[α]);
       EAACCx = r[α] (Cos[α] + Sin[α]^2) -
           2 Cos[α] Sin[α/2]^2 t[α] + Sin[α] (-Cos[α] r'[α] + (-1 + Cos[α]) t'[α]);
       EAACCy = 2 Cos[α/2]^2 (-r[α] Sin[α] + Sin[α] t[α] + Cos[α] r'[α] + t'[α] - Cos[α] t'[α]);
       EAApBBpx = 1/2 * ((-1 + 2 Cos[α] + Cos[2 α]) r[α] - 2 Sin[α/2] -
               4 Cos[α/2]^2 Cos[α] t[α] + Sin[2 α] r'[α] - 2 Sin[α] t'[α] - Sin[2 α] t'[α]);
       EAApBBpy = 1/2 * (2 Cos[α/2] + r[α] (-2 Sin[α] + Sin[2 α]) - 2 (-1 + Cos[α]) Sin[α] t[α] -
               r'[α] + 2 Cos[α] r'[α] - Cos[2 α] r'[α] - t'[α] + Cos[2 α] t'[α]);
       EAApCCpx = 1/2 * (2 Cos[α/2] + 2 r[α] (Cos[α] + Sin[α]^2) +
               (1 - 2 Cos[α] + Cos[2 α]) t[α] - Sin[2 α] r'[α] - 2 Sin[α] t'[α] + Sin[2 α] t'[α]);
       EAApCCpy = 1/2 * (2 Sin[α/2] - 2 (1 + Cos[α]) r[α] Sin[α] + (2 Sin[α] + Sin[2 α]) t[α] +
               r'[α] + 2 Cos[α] r'[α] + Cos[2 α] r'[α] + t'[α] - Cos[2 α] t'[α]);

In[ ]:= (*assume intersection point was p for trajectory of point A and AC envelope
           AA〚1〛/.α→α1p⩵EAACCx/.α→α2p
           AA〚2〛/.α→α1p⩵EAACCy/.α→α2p*)

In[ ]:= (*assume approximate values of α1 and α2 for piecewise function*)

In[ ]:= α1p = 77/1000;

In[ ]:= α2p = 177/100;

In[ ]:= (*define integral functional for sofa area*)
```

```mathematica
In[ ]:= FF = Piecewise[{{0, α < 0},
         {-(EAApBBpx /. α → 0) - EAApBBpy * D[EAApBBpx, α] -
           EAApCCpy * D[EAApCCpx, α]                    - EAABBy * D[EAABBx, α], 0 ≤ α ≤ α1p},
         
         {-(EAApBBpx /. α → 0) - EAApBBpy * D[EAApBBpx, α] - EAApCCpy * D[EAApCCpx, α] +
           AA[[2]] * D[AA[[1]], α] - EAABBy * D[EAABBx, α], α1p < α ≤ 1},
         
         {                      - EAApBBpy * D[EAApBBpx, α] - EAApCCpy * D[EAApCCpx, α] +
           AA[[2]] * D[AA[[1]], α] - EAABBy * D[EAABBx, α], 1 < α ≤ π - α2p},
         
         {                      - EAApBBpy * D[EAApBBpx, α] - EAApCCpy * D[EAApCCpx, α] +
           AA[[2]] * D[AA[[1]], α]                        , π - α2p < α ≤ π / 2},
         
         {0, π / 2 < α}}] // Simplify;

In[ ]:= (*Euler-Lagrange equations*)

In[ ]:= (D[FF, r[α]] - D[D[FF, r'[α]], α] + D[D[FF, r''[α]], {α, 2}]) // Simplify
```

$$Out[\ ]=\begin{cases} -\frac{1}{2}\cos[\alpha]\left(-\cos\left[\frac{\alpha}{2}\right]+4\cos[\alpha]\,r[\alpha]+\sin\left[\frac{\alpha}{2}\right]-\right. \\ \quad\left. 2\cos[\alpha]\,t[\alpha]+8\sin[\alpha]\,r'[\alpha]-2\sin[\alpha]\,t'[\alpha]-4\cos[\alpha]\,r''[\alpha]\right) & -\frac{177}{100}+\pi<\alpha<\frac{\pi}{2} \\ -\frac{1}{2}\cos[\alpha]\left(-\cos\left[\frac{\alpha}{2}\right]-2\,(-3\cos[\alpha]+\cos[2\alpha])\,r[\alpha]+\right. & \frac{77}{1000}<\alpha<-\frac{177}{100}+\pi \\ \quad \sin\left[\frac{\alpha}{2}\right]+2\,(-2\cos[\alpha]+\cos[2\alpha])\,t[\alpha]+12\sin[\alpha]\,r'[\alpha]- \\ \quad 3\sin[2\alpha]\,r'[\alpha]-4\sin[\alpha]\,t'[\alpha]+3\sin[2\alpha]\,t'[\alpha]+r''[\alpha]- \\ \quad\left. 6\cos[\alpha]\,r''[\alpha]+\cos[2\alpha]\,r''[\alpha]+t''[\alpha]-\cos[2\alpha]\,t''[\alpha]\right) \\ -\frac{1}{2}\cos[\alpha]\left(-\cos\left[\frac{\alpha}{2}\right]-2\,(-3\cos[\alpha]+\cos[2\alpha])\,r[\alpha]+\right. & 0<\alpha<\frac{77}{1000} \\ \quad \sin\left[\frac{\alpha}{2}\right]+2\,(-3\cos[\alpha]+\cos[2\alpha])\,t[\alpha]+12\sin[\alpha]\,r'[\alpha]- \\ \quad 3\sin[2\alpha]\,r'[\alpha]-6\sin[\alpha]\,t'[\alpha]+3\sin[2\alpha]\,t'[\alpha]+r''[\alpha]- \\ \quad\left. 6\cos[\alpha]\,r''[\alpha]+\cos[2\alpha]\,r''[\alpha]+t''[\alpha]-\cos[2\alpha]\,t''[\alpha]\right) \\ 0 & 2\alpha>\pi\ ||\ \alpha<0 \\ \text{Indeterminate} & \text{True} \end{cases}$$

*In[ ]:=* `(D[FF, t[α]] - D[D[FF, t'[α]], α] + D[D[FF, t''[α]], {α, 2}]) // Simplify`

*Out[ ]=*

$\begin{cases} \frac{1}{2} \sin[\alpha] \left( \cos\left[\frac{\alpha}{2}\right] + \sin\left[\frac{\alpha}{2}\right] + 2\, r[\alpha] \sin[\alpha] - 4 \sin[\alpha]\, t[\alpha] - 2 \cos[\alpha]\, r'[\alpha] + 8 \cos[\alpha]\, t'[\alpha] + 4 \sin[\alpha]\, t''[\alpha] \right) & -\frac{177}{100} + \pi < \alpha < \frac{\pi}{2} \\[4pt] \frac{1}{2} \sin[\alpha] \big( \cos\left[\frac{\alpha}{2}\right] + \sin\left[\frac{\alpha}{2}\right] + 2\, r[\alpha]\, (3 \sin[\alpha] + \sin[2\alpha]) - 2\, (3 \sin[\alpha] + \sin[2\alpha])\, t[\alpha] + r'[\alpha] - 6 \cos[\alpha]\, r'[\alpha] - 3 \cos[2\alpha]\, r'[\alpha] + t'[\alpha] + 12 \cos[\alpha]\, t'[\alpha] + 3 \cos[2\alpha]\, t'[\alpha] - \sin[2\alpha]\, r''[\alpha] + 6 \sin[\alpha]\, t''[\alpha] + \sin[2\alpha]\, t''[\alpha] \big) & 0 < \alpha \le \frac{77}{1000} \\[4pt] \frac{1}{2} \sin[\alpha] \big( \cos\left[\frac{\alpha}{2}\right] + \sin\left[\frac{\alpha}{2}\right] + 4\, (1 + \cos[\alpha])\, r[\alpha] \sin[\alpha] - 2\, (3 \sin[\alpha] + \sin[2\alpha])\, t[\alpha] + r'[\alpha] - 4 \cos[\alpha]\, r'[\alpha] - 3 \cos[2\alpha]\, r'[\alpha] + t'[\alpha] + 12 \cos[\alpha]\, t'[\alpha] + 3 \cos[2\alpha]\, t'[\alpha] - \sin[2\alpha]\, r''[\alpha] + 6 \sin[\alpha]\, t''[\alpha] + \sin[2\alpha]\, t''[\alpha] \big) & \frac{77}{1000} < \alpha < -\frac{177}{100} + \pi \\[4pt] 0 & 2\alpha > \pi \;||\; \alpha < 0 \\[4pt] \text{Indeterminate} & \text{True} \end{cases}$

*In[ ]:=* `(*use internal package and function to obtain the same Euler-Lagrange equations*)`

*In[ ]:=* **`Needs["VariationalMethods`"]`**

```mathematica
In[ ]:= EulerEquations[FF, {r[α], t[α]}, α] // Simplify
```

Out[ ]=

$$\left\{\alpha<0 \,||\, \left(0<\alpha<\frac{77}{1000} \,\&\&\, \right.\right.$$
$$\cos[\alpha]\left(-\cos\left[\frac{\alpha}{2}\right]-2(-3\cos[\alpha]+\cos[2\alpha])r[\alpha]+\sin\left[\frac{\alpha}{2}\right]+2(-3\cos[\alpha]+\cos[2\alpha])t[\alpha]+\right.$$
$$12\sin[\alpha]r'[\alpha]-3\sin[2\alpha]r'[\alpha]-6\sin[\alpha]t'[\alpha]+3\sin[2\alpha]t'[\alpha]+$$
$$\left.\left.r''[\alpha]-6\cos[\alpha]r''[\alpha]+\cos[2\alpha]r''[\alpha]+t''[\alpha]-\cos[2\alpha]t''[\alpha]\right)==0\right)\,||\,$$
$$\left(1000\alpha>77\,\&\&\,\frac{177}{100}+\alpha<\pi\,\&\&\,\cos[\alpha]\left(-\cos\left[\frac{\alpha}{2}\right]-2(-3\cos[\alpha]+\cos[2\alpha])r[\alpha]+\right.\right.$$
$$\sin\left[\frac{\alpha}{2}\right]+2(-2\cos[\alpha]+\cos[2\alpha])t[\alpha]+12\sin[\alpha]r'[\alpha]-$$
$$3\sin[2\alpha]r'[\alpha]-4\sin[\alpha]t'[\alpha]+3\sin[2\alpha]t'[\alpha]+r''[\alpha]-$$
$$\left.\left.6\cos[\alpha]r''[\alpha]+\cos[2\alpha]r''[\alpha]+t''[\alpha]-\cos[2\alpha]t''[\alpha]\right)==0\right)\,||\,$$
$$\left(\frac{177}{100}+\alpha>\pi\,\&\&\,2\alpha<\pi\,\&\&\,\cos[\alpha]\left(-\cos\left[\frac{\alpha}{2}\right]+4\cos[\alpha]r[\alpha]+\sin\left[\frac{\alpha}{2}\right]-2\cos[\alpha]t[\alpha]+\right.\right.$$
$$\left.\left.8\sin[\alpha]r'[\alpha]-2\sin[\alpha]t'[\alpha]-4\cos[\alpha]r''[\alpha]\right)==0\right)\,||\,2\alpha>\pi,$$

$$\alpha<0\,||\,\left(0<\alpha\leq\frac{77}{1000}\,\&\&\,\sin[\alpha]\left(\cos\left[\frac{\alpha}{2}\right]+\sin\left[\frac{\alpha}{2}\right]+2r[\alpha](3\sin[\alpha]+\sin[2\alpha])-\right.\right.$$
$$2(3\sin[\alpha]+\sin[2\alpha])t[\alpha]+r'[\alpha]-6\cos[\alpha]r'[\alpha]-3\cos[2\alpha]r'[\alpha]+t'[\alpha]+12\cos[\alpha]$$
$$\left.\left.t'[\alpha]+3\cos[2\alpha]t'[\alpha]-\sin[2\alpha]r''[\alpha]+6\sin[\alpha]t''[\alpha]+\sin[2\alpha]t''[\alpha]\right)==0\right)\,||\,$$
$$\left(1000\alpha>77\,\&\&\,\frac{177}{100}+\alpha<\pi\,\&\&\,\sin[\alpha]\left(\cos\left[\frac{\alpha}{2}\right]+\sin\left[\frac{\alpha}{2}\right]+4(1+\cos[\alpha])r[\alpha]\sin[\alpha]-\right.\right.$$
$$2(3\sin[\alpha]+\sin[2\alpha])t[\alpha]+r'[\alpha]-4\cos[\alpha]r'[\alpha]-3\cos[2\alpha]r'[\alpha]+t'[\alpha]+12\cos[\alpha]$$
$$\left.\left.t'[\alpha]+3\cos[2\alpha]t'[\alpha]-\sin[2\alpha]r''[\alpha]+6\sin[\alpha]t''[\alpha]+\sin[2\alpha]t''[\alpha]\right)==0\right)\,||\,$$
$$\left(\frac{177}{100}+\alpha>\pi\,\&\&\,2\alpha<\pi\,\&\&\,\sin[\alpha]\left(\cos\left[\frac{\alpha}{2}\right]+\sin\left[\frac{\alpha}{2}\right]+2r[\alpha]\sin[\alpha]-4\sin[\alpha]t[\alpha]-\right.\right.$$
$$\left.\left.\left.2\cos[\alpha]r'[\alpha]+8\cos[\alpha]t'[\alpha]+4\sin[\alpha]t''[\alpha]\right)==0\right)\,||\,2\alpha>\pi\right\}$$

```mathematica
In[ ]:= (*Boundary condition for α=0*)

In[ ]:= FF0 = -(EAApBBpx /. α → 0) - EAApBBpy * D[EAApBBpx, α] -
         EAApCCpy * D[EAApCCpx, α] - EAABBy * D[EAABBx, α];

In[ ]:= (D[FF0, r'[α], α] - D[D[FF0, r''[α], α], α]) /. α → 0
```

Out[ ]=

$$-\frac{1}{2}+2r[0]-2t[0]-2r''[0]$$

```mathematica
In[ ]:= (D[FF0, t'[α], α] - D[D[FF0, t''[α], α], α]) /. α → 0
```

Out[ ]=

0

```mathematica
In[ ]:= NotebookSave[];
```

# Appendix 2

```mathematica
(*search for optimal boundary value of t[0]*)

ParallelTable[
  (*initial guess of α1p*)
  α1p = 0.078;
  (*initial guess of α2p*)
  α2p = 1.78;
  (*initial value for solving r[α]*)
  ri = 0.61;
  (*numerical search for optimal boundary value t[0]*)
  t0 = tt;
  (*number of discretization points*)
  pts = 2001;
  (*discretization: αα is discrete parameter*)
  αα = Range[0., π/2, (π/2 - 0.)/(pts - 1)];
  (*discretization setting of finite difference method with 4th difference order*)
  dfdα = NDSolve`FiniteDifferenceDerivative[{1}, {αα}, "DifferenceOrder" -> 4][
     "DifferentiationMatrix"];
  d2fdα2 = NDSolve`FiniteDifferenceDerivative[{2}, {αα}, "DifferenceOrder" -> 4][
     "DifferentiationMatrix"];
  Do[
   nL = Length[αα];
   varsr = Table[Unique[r$], {nL}];
   varst = Table[Unique[t$], {nL}];
   (*discrete ODE1*)
   eqns1 =
    (-Cos[αα/2] - 2 (-3 Cos[αα] + Cos[2 αα]) varsr + Sin[αα/2] + 2 (-3 Cos[αα] + Cos[2 αα]) varst +
        12 Sin[αα] dfdα.varsr - 3 Sin[2 αα] dfdα.varsr - 6 Sin[αα] dfdα.varst +
        3 Sin[2 αα] dfdα.varst + d2fdα2.varsr - 6 Cos[αα] d2fdα2.varsr +
        Cos[2 αα] d2fdα2.varsr + d2fdα2.varst - Cos[2 αα] d2fdα2.varst) *
      (Boole[0 ≤ # < α1p] & /@ αα) + (-Cos[αα/2] - 2 (-3 Cos[αα] + Cos[2 αα]) varsr +
        Sin[αα/2] + 2 (-2 Cos[αα] + Cos[2 αα]) varst + 12 Sin[αα] dfdα.varsr -
        3 Sin[2 αα] dfdα.varsr - 4 Sin[αα] dfdα.varst + 3 Sin[2 αα] dfdα.varst +
        d2fdα2.varsr - 6 Cos[αα] d2fdα2.varsr + Cos[2 αα] d2fdα2.varsr +
        d2fdα2.varst - Cos[2 αα] d2fdα2.varst) * (Boole[α1p ≤ # < π - α2p] & /@ αα) +
      (-Cos[αα/2] + 4 Cos[αα] varsr + Sin[αα/2] - 2 Cos[αα] varst + 8 Sin[αα] dfdα.varsr -
        2 Sin[αα] dfdα.varst - 4 Cos[αα] d2fdα2.varsr) * (Boole[π - α2p ≤ # ≤ π/2] & /@ αα);
   (*discrete ODE2*)
   eqns2 =
```

```mathematica
(Cos[αα/2] + Sin[αα/2] + 2 varsr (3 Sin[αα] + Sin[2 αα]) - 2 (3 Sin[αα] + Sin[2 αα]) varst +
    dfdα.varsr - 6 Cos[αα] dfdα.varsr - 3 Cos[2 αα] dfdα.varsr + dfdα.varst +
    12 Cos[αα] dfdα.varst + 3 Cos[2 αα] dfdα.varst - Sin[2 αα] d2fdα2.varsr +
    6 Sin[αα] d2fdα2.varst + Sin[2 αα] d2fdα2.varst) * (Boole[0 ≤ # < α1p] & /@ αα) +
(Cos[αα/2] + Sin[αα/2] + 4 (1 + Cos[αα]) varsr Sin[αα] - 2 (3 Sin[αα] + Sin[2 αα]) varst +
    dfdα.varsr - 4 Cos[αα] dfdα.varsr - 3 Cos[2 αα] dfdα.varsr + dfdα.varst +
    12 Cos[αα] dfdα.varst + 3 Cos[2 αα] dfdα.varst - Sin[2 αα] d2fdα2.varsr +
    6 Sin[αα] d2fdα2.varst + Sin[2 αα] d2fdα2.varst) * (Boole[α1p ≤ # < π - α2p] & /@ αα) +
(Cos[αα/2] + Sin[αα/2] + 2 varsr Sin[αα] - 4 Sin[αα] varst - 2 Cos[αα] dfdα.varsr +
    8 Cos[αα] dfdα.varst + 4 Sin[αα] d2fdα2.varst) * (Boole[π - α2p ≤ # ≤ π/2] & /@ αα);
(*boundary conditions*)
(*-1/2+2 r[0]-2 t[0]-2 r''[0]==0*)
eqns1[[1]] = (-1/2 + 2 varsr[[1]] - 2 varst[[1]] - 2 (d2fdα2.varsr)[[1]]) - 0;
(*r'[π/2]==0*)
eqns1[[nL]] = (dfdα.varsr)[[nL]] - 0;
(* t[0]==t0*)
eqns2[[1]] = varst[[1]] - t0;
(*rt[π/2]==0*)
eqns2[[nL]] = (dfdα.varst)[[nL]] - 0;
(*governing equations and variables*)
govEq = Join[eqns1, eqns2];
allVars = Join[varsr, varst];
(*solve discretized ODEs with initial value r[αα]==ri and t[αα]==t0*)
root = FindRoot[govEq, Thread[{allVars, Flatten[{Table[ri, nL], Table[t0, nL]}]}],
    MaxIterations → 50, AccuracyGoal → Infinity, PrecisionGoal → 8];
sol = allVars /. root;
(*interpolation to obtain r[α] and t[α]*)
{rsoltemp, tsoltemp} = Partition[sol, nL];
rsol = Interpolation[Transpose[{αα, rsoltemp}]];
tsol = Interpolation[Transpose[{αα, tsoltemp}]];
(*obtain curves*)
AA = {rsol[α] * Cos[α], tsol[α] * Sin[α]};
AAp =
    {rsol[α] * Cos[α] + √2 * Cos[π/4 + α/2], tsol[α] * Sin[α] + √2 * Sin[π/4 + α/2]};
BB = {rsol[α] * Cos[α] - tsol[α] * 2 Cos[α/2]^2, 0};
BBp = {rsol[α] * Cos[α] - tsol[α] * 2 Cos[α/2]^2 - 1/Sin[α/2], 0};
CC = {rsol[α] * Cos[α] + tsol[α] * Sin[α] * Tan[α/2], 0};
CCp = {rsol[α] * Cos[α] + tsol[α] * Sin[α] * Tan[α/2] + 1/Cos[α/2], 0};
EAABBx = 1/2 * (-1 + 2 Cos[α] + Cos[2 α]) rsol[α] -
    2 Cos[α/2]^2 Cos[α] tsol[α] + Sin[α] (Cos[α] rsol'[α] - (1 + Cos[α]) tsol'[α]);
EAABBy = -2 Sin[α/2]^2
```

```
                (rsol[α] Sin[α] - Sin[α] tsol[α] - Cos[α] rsol'[α] + tsol'[α] + Cos[α] tsol'[α]);
      EAACCx = rsol[α] (Cos[α] + Sin[α]^2) - 2 Cos[α] Sin[α/2]^2 tsol[α] +
                Sin[α] (-Cos[α] rsol'[α] + (-1 + Cos[α]) tsol'[α]);
      EAACCy = 2 Cos[α/2]^2
                (-rsol[α] Sin[α] + Sin[α] tsol[α] + Cos[α] rsol'[α] + tsol'[α] - Cos[α] tsol'[α]);
      EAApBBpx =
           1/2 * ((-1 + 2 Cos[α] + Cos[2 α]) rsol[α] - 2 Sin[α/2] - 4 Cos[α/2]^2 Cos[α] tsol[α] +
                Sin[2 α] rsol'[α] - 2 Sin[α] tsol'[α] - Sin[2 α] tsol'[α]);
      EAApBBpy =
           1/2 * (2 Cos[α/2] + rsol[α] (-2 Sin[α] + Sin[2 α]) - 2 (-1 + Cos[α]) Sin[α] tsol[α] -
                rsol'[α] + 2 Cos[α] rsol'[α] - Cos[2 α] rsol'[α] - tsol'[α] + Cos[2 α] tsol'[α]);
      EAApCCpx = 1/2 * (2 Cos[α/2] + 2 rsol[α] (Cos[α] + Sin[α]^2) + (1 - 2 Cos[α] + Cos[2 α]) tsol[α] -
                Sin[2 α] rsol'[α] - 2 Sin[α] tsol'[α] + Sin[2 α] tsol'[α]);
      EAApCCpy = 1/2 * (2 Sin[α/2] - 2 (1 + Cos[α]) rsol[α] Sin[α] + (2 Sin[α] + Sin[2 α]) tsol[α] +
                rsol'[α] + 2 Cos[α] rsol'[α] + Cos[2 α] rsol'[α] + tsol'[α] - Cos[2 α] tsol'[α]);
           (*find α1p α2p*)
      ααres = FindRoot[{rsol[αα1] * Cos[αα1] == -EAABBx /. α → (π - αα2),
                tsol[αα1] * Sin[αα1] == EAABBy /. α → (π - αα2)}, {αα1, 0.1}, {αα2, 1.8}];
      α1p = ααres[[1]][[2]];
      α2p = ααres[[2]][[2]];
     , {5}(*iterate for 5 times*)] // Quiet;
     {tt, RealDigits[-2 rsol[0] + 4 tsol[0] +
         2 * NIntegrate[-EAApBBpy * D[EAApBBpx, α] - EAApCCpy * D[EAApCCpx, α], {α, 0, π/2},
           WorkingPrecision → 100, MaxRecursion → 20, AccuracyGoal → 8] + 2 * NIntegrate[
           tsol[α] * Sin[α] * D[rsol[α] * Cos[α], α], {α, α1p, π/2}, WorkingPrecision → 100,
           MaxRecursion → 20, AccuracyGoal → 8] - 2 * NIntegrate[EAABBy * D[EAABBx, α], {α, 0,
             π - α2p}, WorkingPrecision → 100, MaxRecursion → 20, AccuracyGoal → 8] // Quiet},
    {tt, 0.709, 0.711, 0.0001}(*search from 0.709 to 0.711 in 0.0001 intervals*)]
```

Out[ ]=
```
{{0.709, {{2, 2, 1, 9, 5, 2, 7, 9, 1, 6, 9, 5, 5, 1, 9, 4}, 1}},
 {0.7091, {{2, 2, 1, 9, 5, 2, 7, 8, 7, 4, 9, 1, 8, 7, 0, 5}, 1}},
 {0.7092, {{2, 2, 1, 9, 5, 2, 8, 3, 9, 8, 1, 3, 2, 7, 5, 7}, 1}},
 {0.7093, {{2, 2, 1, 9, 5, 2, 9, 4, 1, 0, 2, 3, 3, 6, 8, 2}, 1}},
 {0.7094, {{2, 2, 1, 9, 5, 2, 9, 2, 8, 4, 7, 1, 5, 4, 1, 1}, 1}},
 {0.7095, {{2, 2, 1, 9, 5, 2, 9, 6, 0, 8, 1, 1, 8, 4, 3, 3}, 1}},
 {0.7096, {{2, 2, 1, 9, 5, 2, 9, 9, 5, 9, 7, 6, 1, 7, 9, 6}, 1}},
 {0.7097, {{2, 2, 1, 9, 5, 3, 0, 7, 6, 9, 6, 3, 2, 0, 0, 4}, 1}},
 {0.7098, {{2, 2, 1, 9, 5, 3, 0, 4, 5, 8, 0, 2, 6, 6, 1, 8}, 1}},
 {0.7099, {{2, 2, 1, 9, 5, 3, 0, 6, 2, 6, 3, 6, 5, 9, 1, 1}, 1}},
 {0.71, {{2, 2, 1, 9, 5, 3, 1, 4, 0, 1, 6, 8, 7, 0, 5, 7}, 1}},
 {0.7101, {{2, 2, 1, 9, 5, 3, 0, 9, 1, 9, 8, 8, 5, 8, 3, 3}, 1}},
 {0.7102, {{2, 2, 1, 9, 5, 3, 1, 0, 0, 5, 9, 3, 1, 5, 9, 7}, 1}},
 {0.7103, {{2, 2, 1, 9, 5, 3, 1, 6, 9, 8, 8, 6, 0, 8, 5, 1}, 1}},
 {0.7104, {{2, 2, 1, 9, 5, 3, 1, 0, 4, 3, 7, 0, 5, 4, 3, 3}, 1}},
 {0.7105, {{2, 2, 1, 9, 5, 3, 0, 9, 6, 5, 6, 7, 1, 8, 5, 3}, 1}},
 {0.7106, {{2, 2, 1, 9, 5, 3, 1, 4, 5, 3, 1, 4, 0, 8, 9, 0}, 1}},
 {0.7107, {{2, 2, 1, 9, 5, 3, 0, 7, 2, 6, 4, 1, 0, 7, 9, 0}, 1}},
 {0.7108, {{2, 2, 1, 9, 5, 3, 0, 5, 6, 1, 5, 0, 2, 2, 8, 9}, 1}},
 {0.7109, {{2, 2, 1, 9, 5, 3, 0, 4, 4, 4, 7, 7, 3, 5, 6, 7}, 1}},
 {0.711, {{2, 2, 1, 9, 5, 3, 0, 6, 6, 7, 4, 7, 8, 8, 4, 5}, 1}}}
```

In[ ]:= (*optimal t[0] is close to 0.7103*)

In[ ]:= **NotebookSave[];**

In[ ]:= (*From Romik's paper and code, the optimal value of t[0] is 0.710322
         movingsofas-v1.3
         https://www.math.ucdavis.edu/~romik/data/uploads/software/movingsofas-v1.3.nb*)

```
In[•]:= GerverPrecomputedNumericalConstants = {
        ϕ → 0.039177364790083641863217875242424,
        θ → 0.681301509382724894473855757083,
        κ_{1,1} → -0.210322422072688751416185718488,
        κ_{1,2} → 1/4,
        κ_{2,1} → -0.919179292771593322274696102894,
        κ_{2,2} → 0.472406619750805465181760762512,
        κ_{3,1} → -0.613763229430251668554914291318,
        κ_{3,2} → 0.889626479003221860727043050048,
        κ_{4,1} → -0.308347166088910014835132479741,
        κ_{4,2} → 0.472406619750805465181760762512,
        κ_{5,1} → -1.017204036787814585693642864415,
        κ_{5,2} → 1/4,
        a_1 → 1.210322422072688751416185718490,
        a_2 → -1/4,
        b_1 → -0.527624598026784624160503809373,
        b_2 → 0.920258385160637622893705795012,
        c_1 → 0.626045522848465867552329310386,
        c_2 → -0.944750803946430751678992381254,
        d_1 → 1.313022761424232933776164655190,
        d_2 → -0.525382670414554437202836294305,
        e_1 → 1.210322422072688751416185718490,
        e_2 → 1/4
        };

In[•]:= x_1[t_] := RotationMatrix[t].{a_1 Cos[t] + a_2 Sin[t] - 1, -a_2 Cos[t] + a_1 Sin[t] - 1/2} + {κ_{1,1}, κ_{1,2}};

In[•]:= Limit[((x_1[α/2] - {(1 - e_1 + κ_{5,1})/2, 0})[[2]] / Sin[α]) /.
        GerverPrecomputedNumericalConstants, α → 0]

Out[•]=
        0.710322422072688751416185720

In[•]:= (*after simplification, it is $-\frac{1}{2} + a_1$*)
```

# Appendix 3

```mathematica
(*initial guess of α1p*)
α1p = 0.078;
(*initial guess of α2p*)
α2p = 1.78;
(*numerical search for optimal boundary value t[0]*)
t0 = 0.7103;
(*initial value for solving r[α]*)
ri = 0.61;
(*number of discretization points*)
pts = 2001;
(*discretization: αα is discrete parameter*)
αα = Range[0., π/2, (π/2 - 0.)/(pts - 1)];
(*discretization setting of finite difference method with 4th difference order*)
dfdα = NDSolve`FiniteDifferenceDerivative[{1}, {αα}, "DifferenceOrder" → 4][
   "DifferentiationMatrix"];
d2fdα2 = NDSolve`FiniteDifferenceDerivative[{2}, {αα}, "DifferenceOrder" → 4][
   "DifferentiationMatrix"];
```

```mathematica
(*Find α1p and α2p iteratively*)
Do[
 nL = Length[αα];
  varsr = Table[Unique[r$], {nL}];
  varst = Table[Unique[t$], {nL}];
 (*discrete ODE1-piecewise EL equation*)
  eqns1 =
   (-Cos[αα/2] - 2 (-3 Cos[αα] + Cos[2 αα]) varsr + Sin[αα/2] + 2 (-3 Cos[αα] + Cos[2 αα]) varst +
      12 Sin[αα] dfdα.varsr - 3 Sin[2 αα] dfdα.varsr - 6 Sin[αα] dfdα.varst +
      3 Sin[2 αα] dfdα.varst + d2fdα2.varsr - 6 Cos[αα] d2fdα2.varsr +
      Cos[2 αα] d2fdα2.varsr + d2fdα2.varst - Cos[2 αα] d2fdα2.varst) *
     (Boole[0 ≤ # < α1p] & /@ αα) + (-Cos[αα/2] - 2 (-3 Cos[αα] + Cos[2 αα]) varsr +
      Sin[αα/2] + 2 (-2 Cos[αα] + Cos[2 αα]) varst + 12 Sin[αα] dfdα.varsr -
      3 Sin[2 αα] dfdα.varsr - 4 Sin[αα] dfdα.varst + 3 Sin[2 αα] dfdα.varst +
      d2fdα2.varsr - 6 Cos[αα] d2fdα2.varsr + Cos[2 αα] d2fdα2.varsr +
      d2fdα2.varst - Cos[2 αα] d2fdα2.varst) * (Boole[α1p ≤ # < π - α2p] & /@ αα) +
    (-Cos[αα/2] + 4 Cos[αα] varsr + Sin[αα/2] - 2 Cos[αα] varst + 8 Sin[αα] dfdα.varsr -
      2 Sin[αα] dfdα.varst - 4 Cos[αα] d2fdα2.varsr) * (Boole[π - α2p ≤ # ≤ π/2] & /@ αα);
 (*discrete ODE2-piecewise EL equation*)
  eqns2 =
   (Cos[αα/2] + Sin[αα/2] + 2 varsr (3 Sin[αα] + Sin[2 αα]) - 2 (3 Sin[αα] + Sin[2 αα]) varst +
```

```mathematica
        dfdα.varsr - 6 Cos[αα] dfdα.varsr - 3 Cos[2 αα] dfdα.varsr + dfdα.varst +
        12 Cos[αα] dfdα.varst + 3 Cos[2 αα] dfdα.varst - Sin[2 αα] d2fdα2.varsr +
        6 Sin[αα] d2fdα2.varst + Sin[2 αα] d2fdα2.varst) * (Boole[0 ≤ # < α1p] & /@ αα) +
     (Cos[αα/2] + Sin[αα/2] + 4 (1 + Cos[αα]) varsr Sin[αα] - 2 (3 Sin[αα] + Sin[2 αα]) varst +
        dfdα.varsr - 4 Cos[αα] dfdα.varsr - 3 Cos[2 αα] dfdα.varsr + dfdα.varst +
        12 Cos[αα] dfdα.varst + 3 Cos[2 αα] dfdα.varst - Sin[2 αα] d2fdα2.varsr +
        6 Sin[αα] d2fdα2.varst + Sin[2 αα] d2fdα2.varst) * (Boole[α1p ≤ # < π - α2p] & /@ αα) +
     (Cos[αα/2] + Sin[αα/2] + 2 varsr Sin[αα] - 4 Sin[αα] varst - 2 Cos[αα] dfdα.varsr +
        8 Cos[αα] dfdα.varst + 4 Sin[αα] d2fdα2.varst) * (Boole[π - α2p ≤ # ≤ π/2] & /@ αα);
(*boundary conditions
    -1/2+2 r[0]-2 t[0]-2 r''[0]=0
      r'[π/2]=0
        t[0]=t0
         t'[π/2]=0*)
    eqns1〚1〛 = (-1/2 + 2 varsr〚1〛 - 2 varst〚1〛 - 2 (d2fdα2.varsr)〚1〛) - 0;
eqns1〚nL〛 = (dfdα.varsr)〚nL〛 - 0;
eqns2〚1〛 = varst〚1〛 - t0;
eqns2〚nL〛 = (dfdα.varst)〚nL〛 - 0;
(*all equatoins and variables*)
    govEq = Join[eqns1, eqns2];
    allVars = Join[varsr, varst];
(*solve discretized ODEs*)
    root = FindRoot[govEq, Thread[{allVars, Flatten[{Table[ri, nL], Table[t0, nL]}]}],
       MaxIterations → 50, AccuracyGoal → Infinity, PrecisionGoal → 8];
sol = allVars /. root;
(*interpolation to obtain r[α] and t[α]*)
    {rsoltemp, tsoltemp} = Partition[sol, nL];
    rsol = Interpolation[Transpose[{αα, rsoltemp}]];
    tsol = Interpolation[Transpose[{αα, tsoltemp}]];
(*obtain curves*)
    AA = {rsol[α] * Cos[α], tsol[α] * Sin[α]};
    AAp =
      {rsol[α] * Cos[α] + √2 * Cos[π/4 + α/2], tsol[α] * Sin[α] + √2 * Sin[π/4 + α/2]};
BB = {rsol[α] * Cos[α] - tsol[α] * 2 Cos[α/2]^2, 0};
    BBp = {rsol[α] * Cos[α] - tsol[α] * 2 Cos[α/2]^2 - 1/Sin[α/2], 0};
CC = {rsol[α] * Cos[α] + tsol[α] * Sin[α] * Tan[α/2], 0};
    CCp = {rsol[α] * Cos[α] + tsol[α] * Sin[α] * Tan[α/2] + 1/Cos[α/2], 0};
EAABBx = 1/2 * (-1 + 2 Cos[α] + Cos[2 α]) rsol[α] -
    2 Cos[α/2]^2 Cos[α] tsol[α] + Sin[α] (Cos[α] rsol'[α] - (1 + Cos[α]) tsol'[α]);
EAABBy = -2 Sin[α/2]^2
      (rsol[α] Sin[α] - Sin[α] tsol[α] - Cos[α] rsol'[α] + tsol'[α] + Cos[α] tsol'[α]);
```

```
        EAACCx = rsol[α] (Cos[α] + Sin[α]^2) - 2 Cos[α] Sin[α/2]^2 tsol[α] +
            Sin[α] (-Cos[α] rsol'[α] + (-1 + Cos[α]) tsol'[α]);
        EAACCy = 2 Cos[α/2]^2
            (-rsol[α] Sin[α] + Sin[α] tsol[α] + Cos[α] rsol'[α] + tsol'[α] - Cos[α] tsol'[α]);
        EAApBBpx =
            1/2 * ((-1 + 2 Cos[α] + Cos[2 α]) rsol[α] - 2 Sin[α/2] - 4 Cos[α/2]^2 Cos[α] tsol[α] +
                Sin[2 α] rsol'[α] - 2 Sin[α] tsol'[α] - Sin[2 α] tsol'[α]);
        EAApBBpy =
            1/2 * (2 Cos[α/2] + rsol[α] (-2 Sin[α] + Sin[2 α]) - 2 (-1 + Cos[α]) Sin[α] tsol[α] -
                rsol'[α] + 2 Cos[α] rsol'[α] - Cos[2 α] rsol'[α] - tsol'[α] + Cos[2 α] tsol'[α]);
        EAApCCpx = 1/2 * (2 Cos[α/2] + 2 rsol[α] (Cos[α] + Sin[α]^2) + (1 - 2 Cos[α] + Cos[2 α]) tsol[α] -
                Sin[2 α] rsol'[α] - 2 Sin[α] tsol'[α] + Sin[2 α] tsol'[α]);
        EAApCCpy = 1/2 * (2 Sin[α/2] - 2 (1 + Cos[α]) rsol[α] Sin[α] + (2 Sin[α] + Sin[2 α]) tsol[α] +
                rsol'[α] + 2 Cos[α] rsol'[α] + Cos[2 α] rsol'[α] + tsol'[α] - Cos[2 α] tsol'[α]);
        (*find α1p α2p*)
            ααres = FindRoot[{rsol[αα1] * Cos[αα1] == -EAABBx /. α → (π - αα2),
                tsol[αα1] * Sin[αα1] == EAABBy /. α → (π - αα2)}, {αα1, 0.1}, {αα2, 1.8}];
        (*print iteration difference of α1p α2p*)
            Print[{ααres[[1]][[2]], ααres[[2]][[2]]} - {α1p, α2p}];
            α1p = ααres[[1]][[2]];
        α2p = ααres[[2]][[2]];
        , {5}(*iterate for 5 times*)] // Quiet;

        {0.000464261, -0.00123466}

        {-0.000120681, 0.000489759}

        {0.000060287, -0.000418018}

        {-0.000060287, 0.000418018}

        {0.000060287, -0.000418018}

In[○]:= (*a little fluctuations*)

In[○]:= (*calculate sofa area*)

In[○]:= -2 rsol[0] + 4 tsol[0] + 2 * NIntegrate[-EAApBBpy * D[EAApBBpx, α] - EAApCCpy * D[EAApCCpx, α],
            {α, 0, π/2}, WorkingPrecision → 500, MaxRecursion → 30, AccuracyGoal → 8] +
        2 * NIntegrate[tsol[α] * Sin[α] * D[rsol[α] * Cos[α], α], {α, α1p, π/2},
            WorkingPrecision → 500, MaxRecursion → 30, AccuracyGoal → 8] -
        2 * NIntegrate[EAABBy * D[EAABBx, α], {α, 0, π - α2p},
            WorkingPrecision → 500, MaxRecursion → 30, AccuracyGoal → 8]
Out[○]=
        2.21953

In[○]:= RealDigits[%]
Out[○]=
        {{2, 2, 1, 9, 5, 3, 1, 6, 2, 7, 1, 9, 8, 4, 8, 1}, 1}
```

In[ ]:= (*Small error of 4e-8 from Gerver's 2.21953166 was due to approximate boundary values, discretization, and interpolation*)

In[ ]:= (*obtained α1p α2p*)

In[ ]:= {α1p, α2p}

Out[ ]=
{0.0784039, 1.77884}

In[ ]:= (*display sofa shape*)

In[ ]:= Show[ParametricPlot[{EAApCCpx, EAApCCpy}, {α, 0, π/2},
    PerformanceGoal → "Speed", PlotStyle → {Red}], ParametricPlot[
   {-EAApCCpx, EAApCCpy}, {α, 0, π/2}, PerformanceGoal → "Speed", PlotStyle → {Red}],
  ParametricPlot[{-EAApBBpx, EAApBBpy}, {α, 0, π/2}, PerformanceGoal → "Speed",
   PlotStyle → {Red}], ParametricPlot[{EAApBBpx, EAApBBpy}, {α, 0, π/2},
   PerformanceGoal → "Speed", PlotStyle → {Red}], ParametricPlot[{-EAABBx, EAABBy},
   {α, 0, π - α2p}, PerformanceGoal → "Speed", PlotStyle → {Green}],
  ParametricPlot[{EAABBx, EAABBy}, {α, 0, π - α2p}, PerformanceGoal → "Speed",
   PlotStyle → {Green}], ParametricPlot[{x, 1}, {x, -EAABBx /. α → 0, 0},
   PerformanceGoal → "Speed", PlotStyle → {Black}], ParametricPlot[{x, 1},
   {x, 0, EAABBx /. α → 0}, PerformanceGoal → "Speed", PlotStyle → {Black}],
  ParametricPlot[AA, {α, α1p, π/2}, PerformanceGoal → "Speed", PlotStyle → {Purple}],
  ParametricPlot[{-AA〚1〛, AA〚2〛}, {α, α1p, π/2},
   PerformanceGoal → "Speed", PlotStyle → {Purple}], ParametricPlot[{x, 0},
   {x, -EAABBx /. α → 0, EAApCCpx /. α → 0}, PerformanceGoal → "Speed", PlotStyle → {Black}],
  ParametricPlot[{x, 0}, {x, EAABBx /. α → 0, -EAApCCpx /. α → 0}, PerformanceGoal → "Speed",
   PlotStyle → {Black}], PlotRange → All, AxesOrigin → {0, 0}, ImageSize → 500]

Out[ ]=

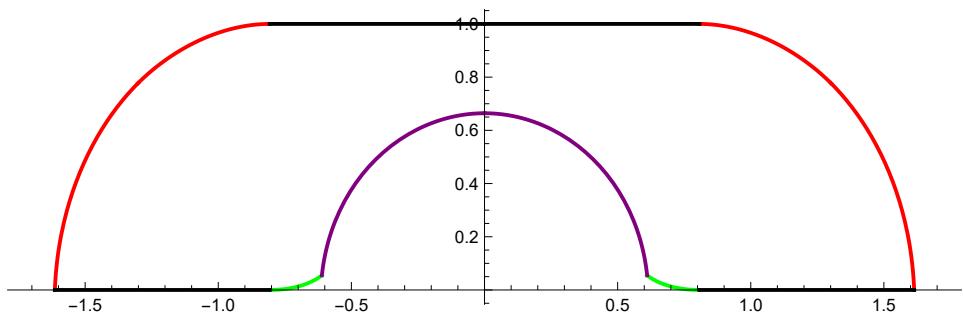

In[ ]:= (*display sofa shape with wall lines and envelop curves*)

```
In[•]:= Show[Table[ListLinePlot[{{AA, BB}, {AA, CC}, {AAp, BBp}, {AAp, CCp}},
          PlotRange → {{-2, 2}, {0, 1.5}}, AspectRatio → Automatic,
          PerformanceGoal → "Speed", PlotStyle → Directive[Black, Thin]], {α, 0, π, π/100}],
        PerformanceGoal → "Speed", ImageSize → 500] // Quiet
```

Out[•]=

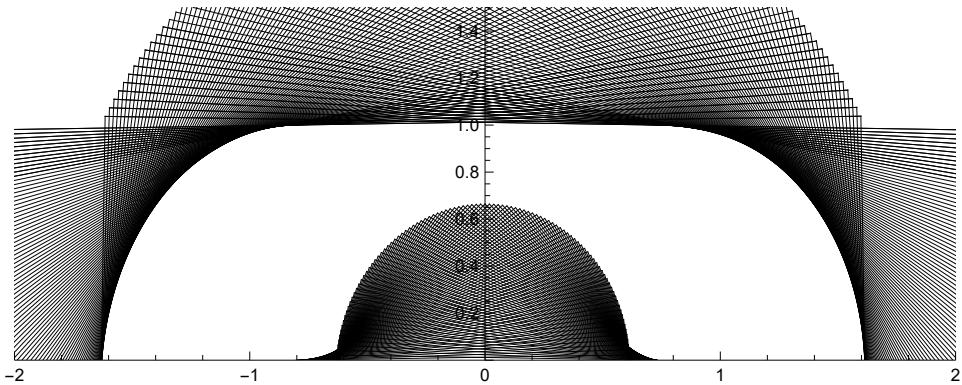

```
In[•]:= (*obtained parametric equations r(α),t(α)*)

In[•]:= Show[Plot[{rsol[α], tsol[α]}, {α, 0, π/2}, PerformanceGoal → "Speed",
         PlotLegends → {"r(α)", "t(α)"}, PlotTheme → "BoldColor", PlotRange → {0.58, 0.72}],
        Plot[{rsol[π - α], tsol[π - α]}, {α, π/2, π}, PerformanceGoal → "Speed",
         PlotLegends → None, PlotTheme → "BoldColor", PlotRange → {0.58, 0.72}],
        PlotRange → {0.58, 0.72}, AspectRatio → 0.5, ImageSize → 500]
```

Out[•]=

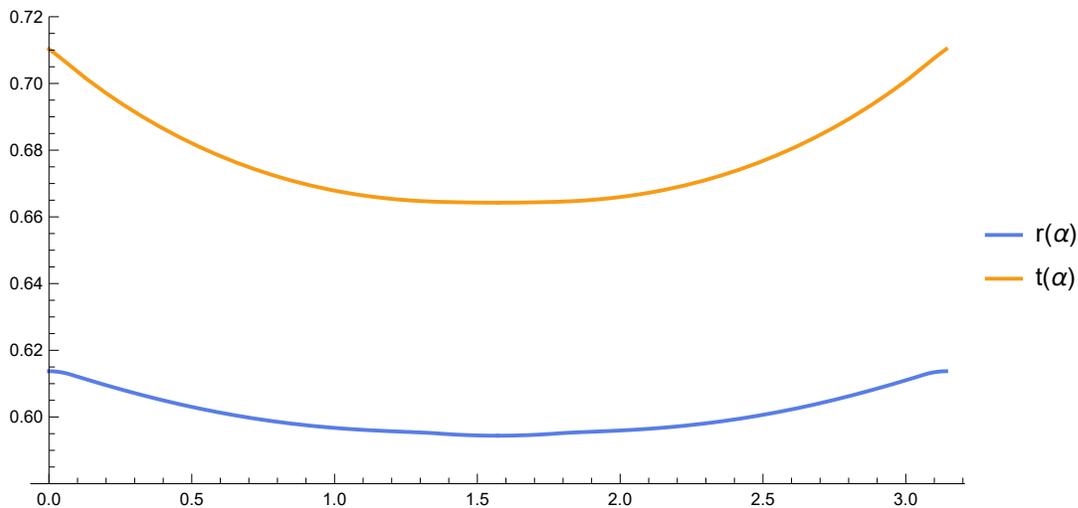

```
In[ ]:= Show[Plot[{rsol'[α], tsol'[α]}, {α, 0, π/2}, PerformanceGoal → "Speed",
    PlotLegends → {"r'(α)", "t'(α)"}, PlotTheme → "BoldColor"],
  Plot[{-rsol'[π - α], -tsol'[π - α]}, {α, π/2, π}, PerformanceGoal → "Speed",
    PlotTheme → "BoldColor"], PlotRange → All, AspectRatio → 0.5, ImageSize → 500]
```

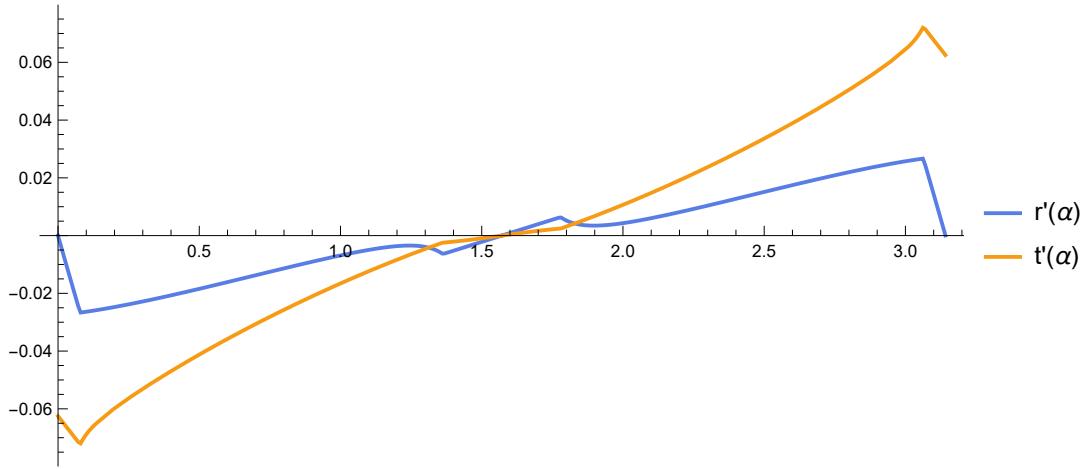

# Appendix 4

```mathematica
In[ ]:= (*Verify Gerver's results from previous paper*)

In[ ]:= (*Romik's explicit results from previous paper and code,
    explicit expression, code from movingsofas-v1.3
      https://www.math.ucdavis.edu/~romik/data/uploads/software/movingsofas-v1.3.nb*)

In[ ]:= X₁[t_] :=
    RotationMatrix[t].{a₁ Cos[t] + a₂ Sin[t] - 1, -a₂ Cos[t] + a₁ Sin[t] - 1/2} + {κ₁,₁, κ₁,₂};

  X₂[t_] := RotationMatrix[t].{-1/4 t² + b₁ t + b₂, 1/2 t - b₁ - 1} + {κ₂,₁, κ₂,₂};

  X₃[t_] := RotationMatrix[t].{c₁ - t, c₂ + t} + {κ₃,₁, κ₃,₂};

In[ ]:=


In[ ]:= (*verify Gerver's and Romik's results satisfy E-L equations ODE1 and ODE2 when 0<α<α1p*)

In[ ]:= (*convert r[α] and t[α] due to different symbols and coordinate systems,
    then substitute r[α] and t[α] into EL equations:
      r[α]→((X₁[α/2]-{(1-e₁+κ₅,₁)/2,0})⟦1⟧)/Cos[α] and  t[α]→
      (X₁[α/2]-{(1-e₁+κ₅,₁)/2,0})⟦2⟧/Sin[α]*)

In[ ]:= {(-Cos[α/2] - 2 (-3 Cos[α] + Cos[2 α]) r[α] + Sin[α/2] + 2 (-3 Cos[α] + Cos[2 α]) t[α] +
      12 Sin[α] r'[α] - 3 Sin[2 α] r'[α] - 6 Sin[α] t'[α] + 3 Sin[2 α] t'[α] +
      r''[α] - 6 Cos[α] r''[α] + Cos[2 α] r''[α] + t''[α] - Cos[2 α] t''[α]),

    (Cos[α/2] + Sin[α/2] + 2 r[α] (3 Sin[α] + Sin[2 α]) - 2 (3 Sin[α] + Sin[2 α]) t[α] +
      r'[α] - 6 Cos[α] r'[α] - 3 Cos[2 α] r'[α] + t'[α] + 12 Cos[α] t'[α] +
      3 Cos[2 α] t'[α] - Sin[2 α] r''[α] + 6 Sin[α] t''[α] + Sin[2 α] t''[α])} /.
    {r[α] → ((X₁[α / 2] - {(1 - e₁ + κ₅,₁) / 2, 0})⟦1⟧) / Cos[α],
     t[α] → (X₁[α / 2] - {(1 - e₁ + κ₅,₁) / 2, 0})⟦2⟧ / Sin[α],
     r'[α] → D[((X₁[α / 2] - {(1 - e₁ + κ₅,₁) / 2, 0})⟦1⟧) / Cos[α], α],
     t'[α] → D[(X₁[α / 2] - {(1 - e₁ + κ₅,₁) / 2, 0})⟦2⟧ / Sin[α], α],
     r''[α] → D[((X₁[α / 2] - {(1 - e₁ + κ₅,₁) / 2, 0})⟦1⟧) / Cos[α], {α, 2}],
     t''[α] → D[(X₁[α / 2] - {(1 - e₁ + κ₅,₁) / 2, 0})⟦2⟧ / Sin[α], {α, 2}]} // Simplify

Out[ ]=
  {0, 0}

  (*{0,0} indicated Gerver's and Romik's results satisfy E-L equations!*)
```

In[ ]:= (*Verify Gerver's and Romik's results satisfy satisfy E-L equations ODE1 and ODE2 when $\alpha 1p<\alpha\leq\pi-\alpha 2p$*)

In[ ]:= $\{(-\cos[\frac{\alpha}{2}] - 2(-3\cos[\alpha] + \cos[2\alpha])r[\alpha] + \sin[\frac{\alpha}{2}] + 2(-2\cos[\alpha] + \cos[2\alpha])t[\alpha] +$
$12\sin[\alpha]r'[\alpha] - 3\sin[2\alpha]r'[\alpha] - 4\sin[\alpha]t'[\alpha] + 3\sin[2\alpha]t'[\alpha] +$
$r''[\alpha] - 6\cos[\alpha]r''[\alpha] + \cos[2\alpha]r''[\alpha] + t''[\alpha] - \cos[2\alpha]t''[\alpha]),$
$(\cos[\frac{\alpha}{2}] + \sin[\frac{\alpha}{2}] + 4(1+\cos[\alpha])r[\alpha]\sin[\alpha] - 2(3\sin[\alpha] + \sin[2\alpha])t[\alpha] +$
$r'[\alpha] - 4\cos[\alpha]r'[\alpha] - 3\cos[2\alpha]r'[\alpha] + t'[\alpha] + 12\cos[\alpha]t'[\alpha] +$
$3\cos[2\alpha]t'[\alpha] - \sin[2\alpha]r''[\alpha] + 6\sin[\alpha]t''[\alpha] + \sin[2\alpha]t''[\alpha])\}$ /.
$\{r[\alpha] \to ((x_2[\alpha/2] - \{(1-e_1+\kappa_{5,1})/2, 0\})[[1]])/\cos[\alpha],$
$t[\alpha] \to (x_2[\alpha/2] - \{(1-e_1+\kappa_{5,1})/2, 0\})[[2]]/\sin[\alpha],$
$r'[\alpha] \to D[((x_2[\alpha/2] - \{(1-e_1+\kappa_{5,1})/2, 0\})[[1]])/\cos[\alpha], \alpha],$
$t'[\alpha] \to D[(x_2[\alpha/2] - \{(1-e_1+\kappa_{5,1})/2, 0\})[[2]]/\sin[\alpha], \alpha],$
$r''[\alpha] \to D[((x_2[\alpha/2] - \{(1-e_1+\kappa_{5,1})/2, 0\})[[1]])/\cos[\alpha], \{\alpha, 2\}],$
$t''[\alpha] \to D[(x_2[\alpha/2] - \{(1-e_1+\kappa_{5,1})/2, 0\})[[2]]/\sin[\alpha], \{\alpha, 2\}]\}$ // Simplify

Out[ ]=
$\{0, 0\}$

In[ ]:= (*Verify Gerver's and Romik's results satisfy E-L equations ODE1 and ODE2 when $\pi-\alpha 2p<\alpha\leq\pi/2$*)

In[ ]:= $\{(-\cos[\frac{\alpha}{2}] + 4\cos[\alpha]r[\alpha] + \sin[\frac{\alpha}{2}] - 2\cos[\alpha]t[\alpha] + 8\sin[\alpha]r'[\alpha] -$
$2\sin[\alpha]t'[\alpha] - 4\cos[\alpha]r''[\alpha]), (\cos[\frac{\alpha}{2}] + \sin[\frac{\alpha}{2}] + 2r[\alpha]\sin[\alpha] -$
$4\sin[\alpha]t[\alpha] - 2\cos[\alpha]r'[\alpha] + 8\cos[\alpha]t'[\alpha] + 4\sin[\alpha]t''[\alpha])\}$ /.
$\{r[\alpha] \to ((x_3[\alpha/2] - \{(1-e_1+\kappa_{5,1})/2, 0\})[[1]])/\cos[\alpha],$
$t[\alpha] \to (x_3[\alpha/2] - \{(1-e_1+\kappa_{5,1})/2, 0\})[[2]]/\sin[\alpha],$
$r'[\alpha] \to D[((x_3[\alpha/2] - \{(1-e_1+\kappa_{5,1})/2, 0\})[[1]])/\cos[\alpha], \alpha],$
$t'[\alpha] \to D[(x_3[\alpha/2] - \{(1-e_1+\kappa_{5,1})/2, 0\})[[2]]/\sin[\alpha], \alpha],$
$r''[\alpha] \to D[((x_3[\alpha/2] - \{(1-e_1+\kappa_{5,1})/2, 0\})[[1]])/\cos[\alpha], \{\alpha, 2\}],$
$t''[\alpha] \to D[(x_3[\alpha/2] - \{(1-e_1+\kappa_{5,1})/2, 0\})[[2]]/\sin[\alpha], \{\alpha, 2\}]\}$ // Simplify

Out[ ]=
$\{0, 0\}$

# Appendix 5

```mathematica
In[ ]:= (*no intersection point of trajectory of A and the AC envelope*)

In[ ]:= (*coordinates of points and envelopes*)

In[ ]:= AA = {r[α] * Cos[α], t[α] * Sin[α]};
       AAp = {r[α] * Cos[α] + √2 * Cos[π/4 + α/2], t[α] * Sin[α] + √2 * Sin[π/4 + α/2]};
       BB = {r[α] * Cos[α] - t[α] * 2 Cos[α/2]^2, 0};
       BBp = {r[α] * Cos[α] - t[α] * 2 Cos[α/2]^2 - 1/Sin[α/2], 0};
       CC = {r[α] * Cos[α] + t[α] * Sin[α] * Tan[α/2], 0};
       CCp = {r[α] * Cos[α] + t[α] * Sin[α] * Tan[α/2] + 1/Cos[α/2], 0};
       EAABBx = 1/2 * (-1 + 2 Cos[α] + Cos[2 α]) r[α] -
           2 Cos[α/2]^2 Cos[α] t[α] + Sin[α] (Cos[α] r'[α] - (1 + Cos[α]) t'[α]);
       EAABBy = -2 Sin[α/2]^2 (r[α] Sin[α] - Sin[α] t[α] - Cos[α] r'[α] + t'[α] + Cos[α] t'[α]);
       EAACCx = r[α] (Cos[α] + Sin[α]^2) -
           2 Cos[α] Sin[α/2]^2 t[α] + Sin[α] (-Cos[α] r'[α] + (-1 + Cos[α]) t'[α]);
       EAACCy = 2 Cos[α/2]^2 (-r[α] Sin[α] + Sin[α] t[α] + Cos[α] r'[α] + t'[α] - Cos[α] t'[α]);
       EAApBBpx = 1/2 * ((-1 + 2 Cos[α] + Cos[2 α]) r[α] - 2 Sin[α/2] -
               4 Cos[α/2]^2 Cos[α] t[α] + Sin[2 α] r'[α] - 2 Sin[α] t'[α] - Sin[2 α] t'[α]);
       EAApBBpy = 1/2 * (2 Cos[α/2] + r[α] (-2 Sin[α] + Sin[2 α]) - 2 (-1 + Cos[α]) Sin[α] t[α] -
               r'[α] + 2 Cos[α] r'[α] - Cos[2 α] r'[α] - t'[α] + Cos[2 α] t'[α]);
       EAApCCpx = 1/2 * (2 Cos[α/2] + 2 r[α] (Cos[α] + Sin[α]^2) +
               (1 - 2 Cos[α] + Cos[2 α]) t[α] - Sin[2 α] r'[α] - 2 Sin[α] t'[α] + Sin[2 α] t'[α]);
       EAApCCpy = 1/2 * (2 Sin[α/2] - 2 (1 + Cos[α]) r[α] Sin[α] + (2 Sin[α] + Sin[2 α]) t[α] +
               r'[α] + 2 Cos[α] r'[α] + Cos[2 α] r'[α] + t'[α] - Cos[2 α] t'[α]);

In[ ]:= (*functional*)

In[ ]:= FF = Piecewise[{{0, α < 0},

           {-r[0] + 2 t[0] - EAApBBpy * D[EAApBBpx, α] -
               EAApCCpy * D[EAApCCpx, α] + AA[[2]] * D[AA[[1]], α], 0 ≤ α ≤ 1},

           {              - EAApBBpy * D[EAApBBpx, α] -
               EAApCCpy * D[EAApCCpx, α] + AA[[2]] * D[AA[[1]], α], 1 < α ≤ π/2},

           {0, π/2 < α}}] // Simplify;

In[ ]:= (*Lagrange multipliers for constraints*)

In[ ]:= HH = FF + λ1[α] * (r[α] * Cos[α] + t[α] * Sin[α] * Tan[α/2] - r0) +
           λ2[α] * (rπ - (r[α] * Cos[α] - t[α] * 2 Cos[α/2]^2));

In[ ]:= (*EL equations*)
```

```
In[•]:= (D[HH, r[α]] - D[D[HH, r'[α]], α] + D[D[HH, r''[α]], {α, 2}]) // Simplify
Out[•]=
        ⎧ Cos[α] (λ1[α] - λ2[α])                                                              2 α > π || α < 0
        ⎪ -½ Cos[α] (-Cos[α/2] + 4 Cos[α] r[α] + Sin[α/2] - 2 Cos[α] t[α] -                   α > 0 && 2 α < π
        ⎨    2 λ1[α] + 2 λ2[α] + 8 Sin[α] r'[α] - 2 Sin[α] t'[α] - 4 Cos[α] r''[α])
        ⎪ Indeterminate                                                                       True

In[•]:= (D[HH, t[α]] - D[D[HH, t'[α]], α] + D[D[HH, t''[α]], {α, 2}]) // Simplify
Out[•]=
        ⎧ λ1[α] - Cos[α] λ1[α] + λ2[α] + Cos[α] λ2[α]                                         2 α > π || α < 0
        ⎪ ½ Cos[α/2] Sin[α] + ½ Sin[α/2] Sin[α] + r[α] Sin[α]² -                              α > 0 && 2 α < π
        ⎨    2 Sin[α]² t[α] + λ1[α] - Cos[α] λ1[α] + λ2[α] + Cos[α] λ2[α] -
        ⎪    Cos[α] Sin[α] r'[α] + 4 Cos[α] Sin[α] t'[α] + 2 Sin[α]² t''[α]
        ⎩ Indeterminate                                                                       True
```

---

```mathematica
In[•]:= (*fine search for optimal boundary values r[0] and t[0]*)

        searchres = ParallelTable[
            (*search for optimal boundary value r[0] and t[0]*)
            r0 = rr; t0 = tt;
            (*number of discretization points*)
            pts = 201;
            αα = Range[0., π/2, (π/2 - 0.)/(pts - 1)];
            dfdα = NDSolve`FiniteDifferenceDerivative[
                {1}, {αα}, "DifferenceOrder" -> 4]["DifferentiationMatrix"];
            d2fdα2 = NDSolve`FiniteDifferenceDerivative[
                {2}, {αα}, "DifferenceOrder" -> 4]["DifferentiationMatrix"];
        nL = Length[αα];
            varsr = Table[Unique[r$], {nL}];
            varst = Table[Unique[t$], {nL}];
            (*Lagrange multipliers λ1 and λ2*)
            varsλ1 = Table[Unique[λ1$], {nL}];
            varsλ2 = Table[Unique[λ2$], {nL}];
        (*discrete ODE*)
            eqns1 = -Cos[αα/2] + 4 Cos[αα] varsr + Sin[αα/2] - 2 Cos[αα] varst - 2 varsλ1 +
                2 varsλ2 + 8 Sin[αα] dfdα.varsr - 2 Sin[αα] dfdα.varst - 4 Cos[αα] d2fdα2.varsr;
            eqns2 = ½ Cos[αα/2] Sin[αα] + ½ Sin[αα/2] Sin[αα] + varsr Sin[αα]² -
                2 Sin[αα]² varst + varsλ1 - Cos[αα] varsλ1 + varsλ2 + Cos[αα] varsλ2 -
                Cos[αα] Sin[αα] dfdα.varsr + 4 Cos[αα] Sin[αα] dfdα.varst + 2 Sin[αα]² d2fdα2.varst;
            eqns3 = varsλ1 * (varsr * Cos[αα] + varst * Sin[αα] * Tan[αα/2] - varsr[[1]]);
        eqns4 = varsλ2 * (-varsr[[1]] - (varsr * Cos[αα] - varst * 2 Cos[αα/2]²));
            (*boundary conditions*)
        eqns1[[1]] = varsr[[1]] - r0;
            eqns1[[nL]] = (dfdα.varsr)[[nL]] - 0;
```

```mathematica
eqns2[[1]] = varst[[1]] - t0;
eqns2[[nL]] = (dfdα.varst)[[nL]] - 0;

govEq = Join[eqns1, eqns2, eqns3, eqns4];
   allVars = Join[varsr, varst, varsλ1, varsλ2];
(*solve discretized ODEs*)root = FindRoot[govEq, Thread[
      {allVars, Flatten[{Table[r0, nL], Table[t0, nL], Table[1, nL], Table[1, nL]}]}],
    MaxIterations → 30, AccuracyGoal → Infinity, PrecisionGoal → 8]; // Quiet;
sol = allVars /. root;
{rsoltemp, tsoltemp, λ1soltemp, λ2soltemp} = Partition[sol, nL];
   (*interpolation*)
   rsol = Interpolation[Transpose[{αα, rsoltemp}]];
   tsol = Interpolation[Transpose[{αα, tsoltemp}]];
   λ1sol = Interpolation[Transpose[{αα, λ1soltemp}]];
   λ2sol = Interpolation[Transpose[{αα, λ2soltemp}]];
(*obtain curves*)
   AA = {rsol[α] * Cos[α], tsol[α] * Sin[α]};
   AAp =
    {rsol[α] * Cos[α] + √2 * Cos[π/4 + α/2], tsol[α] * Sin[α] + √2 * Sin[π/4 + α/2]};
BB = {rsol[α] * Cos[α] - tsol[α] * 2 Cos[α/2]^2, 0};
   BBp = {rsol[α] * Cos[α] - tsol[α] * 2 Cos[α/2]^2 - 1/Sin[α/2], 0};
CC = {rsol[α] * Cos[α] + tsol[α] * Sin[α] * Tan[α/2], 0};
   CCp = {rsol[α] * Cos[α] + tsol[α] * Sin[α] * Tan[α/2] + 1/Cos[α/2], 0};
EAABBx = 1/2 * (-1 + 2 Cos[α] + Cos[2 α]) rsol[α] -
     2 Cos[α/2]^2 Cos[α] tsol[α] + Sin[α] (Cos[α] rsol'[α] - (1 + Cos[α]) tsol'[α]);
EAABBy = -2 Sin[α/2]^2
     (rsol[α] Sin[α] - Sin[α] tsol[α] - Cos[α] rsol'[α] + tsol'[α] + Cos[α] tsol'[α]);
EAACCx = rsol[α] (Cos[α] + Sin[α]^2) - 2 Cos[α] Sin[α/2]^2 tsol[α] +
     Sin[α] (-Cos[α] rsol'[α] + (-1 + Cos[α]) tsol'[α]);
EAACCy = 2 Cos[α/2]^2
     (-rsol[α] Sin[α] + Sin[α] tsol[α] + Cos[α] rsol'[α] + tsol'[α] - Cos[α] tsol'[α]);
EAApBBpx =
    1/2 * ((-1 + 2 Cos[α] + Cos[2 α]) rsol[α] - 2 Sin[α/2] - 4 Cos[α/2]^2 Cos[α] tsol[α] +
       Sin[2 α] rsol'[α] - 2 Sin[α] tsol'[α] - Sin[2 α] tsol'[α]);
EAApBBpy =
    1/2 * (2 Cos[α/2] + rsol[α] (-2 Sin[α] + Sin[2 α]) - 2 (-1 + Cos[α]) Sin[α] tsol[α] -
       rsol'[α] + 2 Cos[α] rsol'[α] - Cos[2 α] rsol'[α] - tsol'[α] + Cos[2 α] tsol'[α]);
EAApCCpx = 1/2 * (2 Cos[α/2] + 2 rsol[α] (Cos[α] + Sin[α]^2) + (1 - 2 Cos[α] + Cos[2 α]) tsol[α] -
       Sin[2 α] rsol'[α] - 2 Sin[α] tsol'[α] + Sin[2 α] tsol'[α]);
EAApCCpy = 1/2 * (2 Sin[α/2] - 2 (1 + Cos[α]) rsol[α] Sin[α] + (2 Sin[α] + Sin[2 α]) tsol[α] +
       rsol'[α] + 2 Cos[α] rsol'[α] + Cos[2 α] rsol'[α] + tsol'[α] - Cos[2 α] tsol'[α]);
   {rr, tt, -2 rsol[0] + 4 tsol[0] +
      2 * NIntegrate[-EAApBBpy * D[EAApBBpx, α] - EAApCCpy * D[EAApCCpx, α],
        {α, 0, π/2}, WorkingPrecision → 100, MaxRecursion → 20, AccuracyGoal → 8] +
      2 * NIntegrate[tsol[α] * Sin[α] * D[rsol[α] * Cos[α], α], {α, 0, π/2},
```

```
          WorkingPrecision -> 100, MaxRecursion -> 20, AccuracyGoal -> 8] //
        Quiet, Norm /@ Partition[govEq /. root, nL],
      Apply[And, # <= 0 & /@ Chop[(varsr * Cos[αα] + varst * Sin[αα] * Tan[αα / 2] - varsr[[1]]) /.
            root, 10^-3]] && Apply[And, # <= 0 & /@ Chop[
          Rest[(-varsr[[1]] - (varsr * Cos[αα] - varst * 2 Cos[αα / 2]^2))] /. root, 10^-3]] //
        Quiet}, {rr, 0.6, 0.7, 0.01}, {tt, 0.6, 0.7, 0.01}
      (*search from 0.6 to 0.7 in 0.01 intervals*)];
    (*output r[0], t[0], residual, area,
    whether satisfy constraints*)(*running may take long time*)

In[ ]:= Select[Partition[Flatten[searchres], 8], #[[8]] == True &]
```

Out[●]= {{0.6, 0.6, 2.20531, 7.02405×10⁻¹¹, 2.85556×10⁻¹¹, 2.89974×10⁻¹⁷, 2.36133×10⁻¹⁶, True},
{0.6, 0.61, 2.21864, 7.83819×10⁻¹¹, 3.48429×10⁻¹¹, 6.73841×10⁻¹⁷, 1.0844×10⁻¹⁶, True},
{0.61, 0.6, 2.19455, 6.78237×10⁻¹¹, 3.23976×10⁻¹¹, 4.10028×10⁻¹⁷, 1.02609×10⁻¹⁶, True},
{0.61, 0.61, 2.2063, 6.33×10⁻¹¹, 2.76269×10⁻¹¹, 3.17786×10⁻¹⁷, 2.39426×10⁻¹⁶, True},
{0.61, 0.62, 2.21919, 6.64875×10⁻¹¹, 3.75848×10⁻¹¹, 5.82144×10⁻¹⁴, 8.82544×10⁻¹⁷, True},
{0.62, 0.6, 2.18573, 6.81378×10⁻¹¹, 2.69295×10⁻¹¹, 9.82966×10⁻¹⁷, 1.01102×10⁻¹⁶, True},
{0.62, 0.61, 2.19519, 7.11227×10⁻¹¹, 2.80777×10⁻¹¹, 1.07296×10⁻¹², 8.3245×10⁻¹⁷, True},
{0.62, 0.62, 2.20698, 6.06517×10⁻¹¹, 2.51886×10⁻¹¹, 3.55536×10⁻¹⁷, 2.4383×10⁻¹⁶, True},
{0.62, 0.63, 2.22032, 7.99066×10⁻¹¹, 3.36933×10⁻¹¹, 1.1156×10⁻¹⁵, 1.40861×10⁻¹⁶, True},
{0.63, 0.6, 2.17694, 8.37392×10⁻¹¹, 4.58727×10⁻¹¹, 2.34064×10⁻¹⁴, 1.6824×10⁻¹², True},
{0.63, 0.61, 2.1861, 8.92692×10⁻¹¹, 4.49359×10⁻¹¹, 6.59015×10⁻¹², 1.24066×10⁻¹⁶, True},
{0.63, 0.62, 2.1956, 8.08539×10⁻¹¹, 3.01089×10⁻¹¹, 1.05924×10⁻¹⁶, 9.68475×10⁻¹⁷, True},
{0.63, 0.63, 2.20735, 7.93504×10⁻¹¹, 3.61679×10⁻¹¹, 4.38741×10⁻¹⁷, 2.39801×10⁻¹⁶, True},
{0.64, 0.62, 2.18614, 8.06305×10⁻¹¹, 3.72309×10⁻¹¹, 1.68576×10⁻¹⁵, 1.12865×10⁻¹⁶, True},
{0.64, 0.63, 2.19565, 8.42259×10⁻¹¹, 3.94143×10⁻¹¹, 3.90342×10⁻¹², 2.12446×10⁻¹⁵, True},
{0.64, 0.64, 2.2074, 7.82523×10⁻¹¹, 3.16466×10⁻¹¹, 4.19533×10⁻¹⁷, 2.37217×10⁻¹⁶, True},
{0.64, 0.65, 2.22059, 7.47214×10⁻¹¹, 3.95693×10⁻¹¹, 8.3052×10⁻¹⁷, 1.04303×10⁻¹⁶, True},
{0.65, 0.63, 2.1859, 8.05378×10⁻¹¹, 4.6477×10⁻¹¹, 6.1571×10⁻¹⁵, 1.03387×10⁻¹⁶, True},
{0.65, 0.64, 2.19539, 9.31907×10⁻¹¹, 3.45003×10⁻¹¹, 1.27535×10⁻¹², 1.07101×10⁻¹⁶, True},
{0.65, 0.65, 2.20713, 8.44036×10⁻¹¹, 3.29749×10⁻¹¹, 4.59373×10⁻¹⁷, 2.34513×10⁻¹⁶, True},
{0.66, 0.64, 2.18533, 7.81207×10⁻¹¹, 4.50975×10⁻¹¹, 9.84332×10⁻¹⁷, 1.05823×10⁻¹⁶, True},
{0.66, 0.65, 2.19482, 8.0644×10⁻¹¹, 3.86282×10⁻¹¹, 1.10685×10⁻¹⁶, 1.13141×10⁻¹⁶, True},
{0.66, 0.66, 2.20656, 6.3298×10⁻¹¹, 2.79202×10⁻¹¹, 4.77973×10⁻¹⁷, 2.38953×10⁻¹⁶, True},
{0.67, 0.62, 2.15634, 7.6201×10⁻¹¹, 3.35499×10⁻¹¹, 8.92505×10⁻¹⁷, 1.32469×10⁻¹⁶, True},
{0.67, 0.65, 2.18427, 8.45975×10⁻¹¹, 4.64582×10⁻¹¹, 1.10368×10⁻¹⁶, 1.41391×10⁻¹⁶, True},
{0.67, 0.66, 2.19392, 7.94298×10⁻¹¹, 3.49735×10⁻¹¹, 1.69556×10⁻¹⁶, 9.56943×10⁻¹⁷, True},
{0.67, 0.67, 2.20567, 6.61425×10⁻¹¹, 2.71405×10⁻¹¹, 5.44969×10⁻¹⁷, 2.41123×10⁻¹⁶, True},
{0.68, 0.66, 2.18325, 7.68797×10⁻¹¹, 4.31148×10⁻¹¹, 1.54381×10⁻¹⁶, 1.07384×10⁻¹⁶, True},
{0.68, 0.67, 2.19274, 7.84117×10⁻¹¹, 4.02052×10⁻¹¹, 1.70831×10⁻¹⁶, 1.13418×10⁻¹⁶, True},
{0.68, 0.68, 2.20446, 6.02925×10⁻¹¹, 2.6976×10⁻¹¹, 5.11972×10⁻¹⁷, 2.37513×10⁻¹⁶, True},
{0.69, 0.66, 2.1721, 6.89528×10⁻¹¹, 3.21429×10⁻¹¹, 1.2822×10⁻¹⁴, 1.14412×10⁻¹⁶, True},
{0.69, 0.67, 2.18227, 8.55566×10⁻¹¹, 4.0313×10⁻¹¹, 9.37952×10⁻¹⁴, 1.29043×10⁻¹⁴, True},
{0.69, 0.68, 2.19121, 7.47129×10⁻¹¹, 3.7053×10⁻¹¹, 1.53615×10⁻¹⁶, 8.62366×10⁻¹⁷, True},
{0.69, 0.69, 2.20294, 5.97984×10⁻¹¹, 2.42536×10⁻¹¹, 5.51179×10⁻¹⁷, 2.35012×10⁻¹⁶, True},
{0.7, 0.68, 2.17988, 7.07795×10⁻¹¹, 2.56251×10⁻¹¹, 1.29718×10⁻¹⁴, 9.17757×10⁻¹⁷, True},
{0.7, 0.69, 2.1894, 7.02606×10⁻¹¹, 3.6752×10⁻¹¹, 1.05705×10⁻¹⁶, 7.72796×10⁻¹⁷, True},
{0.7, 0.7, 2.20111, 6.19456×10⁻¹¹, 1.99005×10⁻¹¹, 6.11734×10⁻¹⁷, 2.34494×10⁻¹⁶, True}}

In[●]:= **NotebookSave[];**

(*although larger results, further valification found {r0,t0}={0.6,0.61},
{0.62,0.63}and{0.64,0.65} cannot satisfy the contraint equations*)

```mathematica
In[ ]:=

In[ ]:= (*results*)

In[ ]:= r0 = 0.64; t0 = 0.64;
pts = 201;
αα = Range[0., π/2, (π/2 - 0.)/(pts - 1)];
dfdα = NDSolve`FiniteDifferenceDerivative[{1}, {αα}, "DifferenceOrder" -> 4][
    "DifferentiationMatrix"];
d2fdα2 = NDSolve`FiniteDifferenceDerivative[{2}, {αα}, "DifferenceOrder" -> 4][
    "DifferentiationMatrix"];
nL = Length[αα];
varsr = Table[Unique[r$], {nL}];
varst = Table[Unique[t$], {nL}];
varsλ1 = Table[Unique[λ1$], {nL}];
varsλ2 = Table[Unique[λ2$], {nL}];
eqns1 = -Cos[αα/2] + 4 Cos[αα] varsr + Sin[αα/2] - 2 Cos[αα] varst - 2 varsλ1 +
    2 varsλ2 + 8 Sin[αα] dfdα.varsr - 2 Sin[αα] dfdα.varst - 4 Cos[αα] d2fdα2.varsr;
eqns2 = 1/2 Cos[αα/2] Sin[αα] + 1/2 Sin[αα/2] Sin[αα] + varsr Sin[αα]^2 -
    2 Sin[αα]^2 varst + varsλ1 - Cos[αα] varsλ1 + varsλ2 + Cos[αα] varsλ2 -
    Cos[αα] Sin[αα] dfdα.varsr + 4 Cos[αα] Sin[αα] dfdα.varst + 2 Sin[αα]^2 d2fdα2.varst;
eqns3 = varsλ1 * (varsr * Cos[αα] + varst * Sin[αα] * Tan[αα/2] - varsr[[1]]);
eqns4 = varsλ2 * (-varsr[[1]] - (varsr * Cos[αα] - varst * 2 Cos[αα/2]^2));

eqns1[[1]] = varsr[[1]] - r0;
eqns1[[nL]] = (dfdα.varsr)[[nL]] - 0;
eqns2[[1]] = varst[[1]] - t0;
eqns2[[nL]] = (dfdα.varst)[[nL]] - 0;

govEq = Join[eqns1, eqns2, eqns3, eqns4];
allVars = Join[varsr, varst, varsλ1, varsλ2];
root = FindRoot[govEq, Thread[
      {allVars, Flatten[{Table[r0, nL], Table[t0, nL], Table[1, nL], Table[1, nL]}]}],
    MaxIterations -> 30, AccuracyGoal -> Infinity, PrecisionGoal -> 8]; // Quiet;
sol = allVars /. root;
Norm /@ Partition[govEq /. root, nL]
(*residual*)

Out[ ]=
    {7.24483 × 10^-11, 4.08916 × 10^-11, 4.29963 × 10^-17, 2.39894 × 10^-16}
```

```mathematica
In[ ]:= {rsoltemp, tsoltemp, λ1soltemp, λ2soltemp} = Partition[sol, nL];
       rsol = Interpolation[Transpose[{αα, rsoltemp}]];
       tsol = Interpolation[Transpose[{αα, tsoltemp}]];
       λ1sol = Interpolation[Transpose[{αα, λ1soltemp}]];
       λ2sol = Interpolation[Transpose[{αα, λ2soltemp}]];
       (*obtain curves*)
       AA = {rsol[α] * Cos[α], tsol[α] * Sin[α]};
       AAp = {rsol[α] * Cos[α] + √2 * Cos[π/4 + α/2], tsol[α] * Sin[α] + √2 * Sin[π/4 + α/2]};
       BB = {rsol[α] * Cos[α] - tsol[α] * 2 Cos[α/2]^2, 0};
       BBp = {rsol[α] * Cos[α] - tsol[α] * 2 Cos[α/2]^2 - 1/Sin[α/2], 0};
       CC = {rsol[α] * Cos[α] + tsol[α] * Sin[α] * Tan[α/2], 0};
       CCp = {rsol[α] * Cos[α] + tsol[α] * Sin[α] * Tan[α/2] + 1/Cos[α/2], 0};
       EAABBx = 1/2 * (-1 + 2 Cos[α] + Cos[2 α]) rsol[α] -
          2 Cos[α/2]^2 Cos[α] tsol[α] + Sin[α] (Cos[α] rsol'[α] - (1 + Cos[α]) tsol'[α]);
       EAABBy = -2 Sin[α/2]^2
          (rsol[α] Sin[α] - Sin[α] tsol[α] - Cos[α] rsol'[α] + tsol'[α] + Cos[α] tsol'[α]);
       EAACCx = rsol[α] (Cos[α] + Sin[α]^2) - 2 Cos[α] Sin[α/2]^2 tsol[α] +
          Sin[α] (-Cos[α] rsol'[α] + (-1 + Cos[α]) tsol'[α]);
       EAACCy = 2 Cos[α/2]^2
          (-rsol[α] Sin[α] + Sin[α] tsol[α] + Cos[α] rsol'[α] + tsol'[α] - Cos[α] tsol'[α]);
       EAApBBpx = 1/2 * ((-1 + 2 Cos[α] + Cos[2 α]) rsol[α] - 2 Sin[α/2] -
             4 Cos[α/2]^2 Cos[α] tsol[α] + Sin[2 α] rsol'[α] - 2 Sin[α] tsol'[α] - Sin[2 α] tsol'[α]);
       EAApBBpy =
         1/2 * (2 Cos[α/2] + rsol[α] (-2 Sin[α] + Sin[2 α]) - 2 (-1 + Cos[α]) Sin[α] tsol[α] -
             rsol'[α] + 2 Cos[α] rsol'[α] - Cos[2 α] rsol'[α] - tsol'[α] + Cos[2 α] tsol'[α]);
       EAApCCpx = 1/2 * (2 Cos[α/2] + 2 rsol[α] (Cos[α] + Sin[α]^2) + (1 - 2 Cos[α] + Cos[2 α]) tsol[α] -
             Sin[2 α] rsol'[α] - 2 Sin[α] tsol'[α] + Sin[2 α] tsol'[α]);
       EAApCCpy = 1/2 * (2 Sin[α/2] - 2 (1 + Cos[α]) rsol[α] Sin[α] + (2 Sin[α] + Sin[2 α]) tsol[α] +
             rsol'[α] + 2 Cos[α] rsol'[α] + Cos[2 α] rsol'[α] + tsol'[α] - Cos[2 α] tsol'[α]);

In[ ]:= -2 rsol[0] + 4 tsol[0] + 2 * NIntegrate[-EAApBBpy * D[EAApBBpx, α] - EAApCCpy * D[EAApCCpx, α],
          {α, 0, π/2}, WorkingPrecision → 100, MaxRecursion → 20, AccuracyGoal → 8] +
        2 * NIntegrate[tsol[α] * Sin[α] * D[rsol[α] * Cos[α], α], {α, 0, π/2},
          WorkingPrecision → 100, MaxRecursion → 20, AccuracyGoal → 8]
Out[ ]=
       2.2074

In[ ]:= RealDigits[%]
Out[ ]=
       {{2, 2, 0, 7, 3, 9, 8, 1, 5, 1, 3, 3, 9, 7, 8, 7}, 1}
```

```
In[ ]:= Show[ParametricPlot[{EAApCCpx, EAApCCpy}, {α, 0, π/2},
     PerformanceGoal → "Speed", PlotStyle → {Red}], ParametricPlot[
    {-EAApCCpx, EAApCCpy}, {α, 0, π/2}, PerformanceGoal → "Speed", PlotStyle → {Red}],
   ParametricPlot[{-EAApBBpx, EAApBBpy}, {α, 0, π/2}, PerformanceGoal → "Speed",
    PlotStyle → {Red}], ParametricPlot[{EAApBBpx, EAApBBpy}, {α, 0, π/2},
    PerformanceGoal → "Speed", PlotStyle → {Red}], ParametricPlot[
    {-EAABBx, EAABBy}, {α, 0, π}, PerformanceGoal → "Speed", PlotStyle → {Green}],
   ParametricPlot[{EAABBx, EAABBy}, {α, 0, π}, PerformanceGoal → "Speed",
    PlotStyle → {Green}], ParametricPlot[{x, 1}, {x, -EAABBx /. α → 0, 0},
    PerformanceGoal → "Speed", PlotStyle → {Black}], ParametricPlot[{x, 1},
    {x, 0, EAABBx /. α → 0}, PerformanceGoal → "Speed", PlotStyle → {Black}],
   ParametricPlot[AA, {α, 0, π/2}, PerformanceGoal → "Speed", PlotStyle → {Purple}],
   ParametricPlot[{-rsol[α] * Cos[α], tsol[α] * Sin[α]},
    {α, 0, π/2}, PerformanceGoal → "Speed", PlotStyle → {Purple}],
   ParametricPlot[{x, 0}, {x, -EAABBx /. α → 0, EAApCCpx /. α → 0},
    PerformanceGoal → "Speed", PlotStyle → {Black}], ParametricPlot[{x, 0},
    {x, EAABBx /. α → 0, -EAApCCpx /. α → 0}, PerformanceGoal → "Speed", PlotStyle → {Black}],
   PlotRange → {{-2.0, 2.0}, {0, 1}}, AxesOrigin → {0, 0}, ImageSize → 500]
```

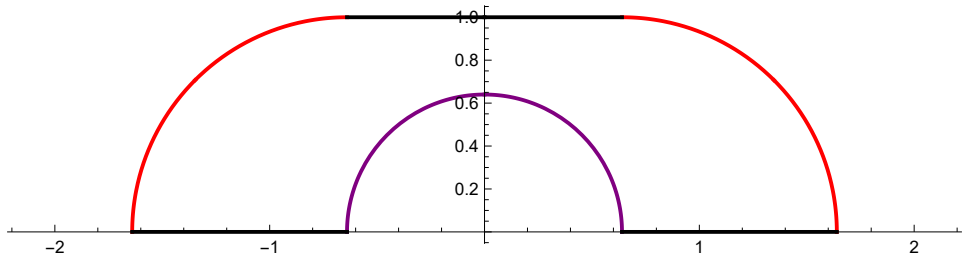

```
In[ ]:=  (*obtained parametric equations r(α),t(α)*)
```

```
In[ ]:= Show[Plot[{rsol[α], tsol[α]}, {α, 0, π/2}, PerformanceGoal → "Speed",
    PlotRange → {0.6, 0.7}, PerformanceGoal → "Speed", PlotLegends → {"r(α)", "t(α)"},
    PlotTheme → "BoldColor"], Plot[{rsol[π - α], tsol[π - α]},
   {α, π/2, π}, PerformanceGoal → "Speed", PlotRange → {0.6, 0.7},
   PerformanceGoal → "Speed", PlotTheme → "BoldColor"], PlotRange → {0.6, 0.7}]
```

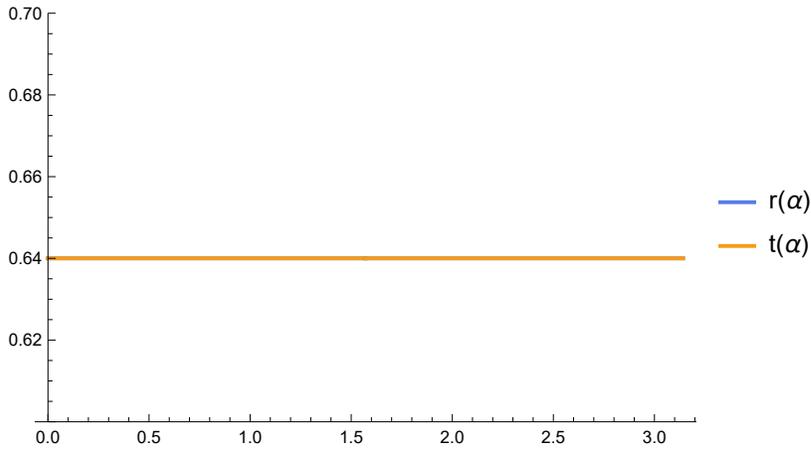

```
In[ ]:= (*test if constraints are satisfied with error less than 10^-3*)

In[ ]:= Apply[And, # ≤ 0 & /@
     Chop[(varsr * Cos[αα] + varst * Sin[αα] * Tan[αα/2] - varsr[[1]]) /. root, 10^-3]] &&
   Apply[And, # ≤ 0 & /@ Chop[Rest[(-varsr[[1]] - (varsr * Cos[αα] - varst * 2 Cos[αα/2]^2))] /.
      root, 10^-3]] // Quiet

Out[ ]= True

In[ ]:= (*results were close to 2/π=0.63661977 and Hammersley's sofa area 2/π+π/2=2.20742*)

In[ ]:= NotebookSave[];
```

# Appendix 6

(*asymmetric conditions*)

(*a*)

(*no intersection point of trajectory of A and the AC envelope and AB envelope*)

(*define coordinates of points and envelopes*)

```mathematica
In[ ]:= AA = {r[α] * Cos[α], t[α] * Sin[α]};
       AAp = {r[α] * Cos[α] + √2 * Cos[π/4 + α/2], t[α] * Sin[α] + √2 * Sin[π/4 + α/2]};
       BB = {r[α] * Cos[α] - t[α] * 2 Cos[α/2]^2, 0};
       BBp = {r[α] * Cos[α] - t[α] * 2 Cos[α/2]^2 - 1/Sin[α/2], 0};
       CC = {r[α] * Cos[α] + t[α] * Sin[α] * Tan[α/2], 0};
       CCp = {r[α] * Cos[α] + t[α] * Sin[α] * Tan[α/2] + 1/Cos[α/2], 0};
       EAABBx = 1/2 * (-1 + 2 Cos[α] + Cos[2 α]) r[α] -
           2 Cos[α/2]^2 Cos[α] t[α] + Sin[α] (Cos[α] r'[α] - (1 + Cos[α]) t'[α]);
       EAABBy = -2 Sin[α/2]^2 (r[α] Sin[α] - Sin[α] t[α] - Cos[α] r'[α] + t'[α] + Cos[α] t'[α]);
       EAACCx = r[α] (Cos[α] + Sin[α]^2) -
           2 Cos[α] Sin[α/2]^2 t[α] + Sin[α] (-Cos[α] r'[α] + (-1 + Cos[α]) t'[α]);
       EAACCy = 2 Cos[α/2]^2 (-r[α] Sin[α] + Sin[α] t[α] + Cos[α] r'[α] + t'[α] - Cos[α] t'[α]);
       EAApBBpx = 1/2 * ((-1 + 2 Cos[α] + Cos[2 α]) r[α] - 2 Sin[α/2] -
               4 Cos[α/2]^2 Cos[α] t[α] + Sin[2 α] r'[α] - 2 Sin[α] t'[α] - Sin[2 α] t'[α]);
       EAApBBpy = 1/2 * (2 Cos[α/2] + r[α] (-2 Sin[α] + Sin[2 α]) - 2 (-1 + Cos[α]) Sin[α] t[α] -
               r'[α] + 2 Cos[α] r'[α] - Cos[2 α] r'[α] - t'[α] + Cos[2 α] t'[α]);
       EAApCCpx = 1/2 * (2 Cos[α/2] + 2 r[α] (Cos[α] + Sin[α]^2) +
               (1 - 2 Cos[α] + Cos[2 α]) t[α] - Sin[2 α] r'[α] - 2 Sin[α] t'[α] + Sin[2 α] t'[α]);
       EAApCCpy = 1/2 * (2 Sin[α/2] - 2 (1 + Cos[α]) r[α] Sin[α] + (2 Sin[α] + Sin[2 α]) t[α] +
               r'[α] + 2 Cos[α] r'[α] + Cos[2 α] r'[α] + t'[α] - Cos[2 α] t'[α]);

In[ ]:= (*define integral functional for sofa area*)

In[ ]:= FF = -r[0] + 2 t[0] - r[π] + 2 t[π] -
           EAApBBpy * D[EAApBBpx, α] - EAApCCpy * D[EAApCCpx, α] + AA[[2]] * D[AA[[1]], α];

In[ ]:= (*Lagrange multipliers for constraints*)

In[ ]:= HH = FF + λ1[α] * (r[α] * Cos[α] + t[α] * Sin[α] * Tan[α/2] - r[0]) +
           λ2[α] * (r[π] - (r[α] * Cos[α] - t[α] * 2 Cos[α/2]^2));

In[ ]:= (*EL equations*)

In[ ]:= (D[HH, r[α]] - D[D[HH, r'[α]], α] + D[D[HH, r''[α]], {α, 2}]) // Simplify
Out[ ]=
       -1/2 Cos[α] (-Cos[α/2] + 4 Cos[α] r[α] + Sin[α/2] - 2 Cos[α] t[α] -
           2 λ1[α] + 2 λ2[α] + 8 Sin[α] r'[α] - 2 Sin[α] t'[α] - 4 Cos[α] r''[α])
```

```
In[ ]:= (D[HH, t[α]] - D[D[HH, t'[α]], α] + D[D[HH, t''[α]], {α, 2}]) // Simplify
```

Out[ ]=
$$\frac{1}{4}\left(\text{Cos}\left[\frac{\alpha}{2}\right] - \text{Cos}\left[\frac{3\alpha}{2}\right] + \text{Sin}\left[\frac{\alpha}{2}\right] + 4\,r[\alpha]\,\text{Sin}[\alpha]^2 + \text{Sin}\left[\frac{3\alpha}{2}\right] - 8\,\text{Sin}[\alpha]^2\,t[\alpha] + 4\,\lambda1[\alpha] - 4\,\text{Cos}[\alpha]\,\lambda1[\alpha] + 4\,\lambda2[\alpha] + 4\,\text{Cos}[\alpha]\,\lambda2[\alpha] - 2\,\text{Sin}[2\alpha]\,r'[\alpha] + 8\,\text{Sin}[2\alpha]\,t'[\alpha] + 4\,t''[\alpha] - 4\,\text{Cos}[2\alpha]\,t''[\alpha]\right)$$

---

```
(*search for boundary values r[0], r[π], t[0], t[π]*)

searchres = ParallelTable[
    r0 = rr0;
    rπ = rrπ;
    t0 = tt0;
    tπ = ttπ;
    pts = 201;
    (*discretization*)
    αα = Range[0., π, (π - 0.) / (pts - 1)];
    dfdα = NDSolve`FiniteDifferenceDerivative[
        {1}, {αα}, "DifferenceOrder" → 4]["DifferentiationMatrix"];
    d2fdα2 = NDSolve`FiniteDifferenceDerivative[
        {2}, {αα}, "DifferenceOrder" → 4]["DifferentiationMatrix"];
nL = Length[αα];
    varsr = Table[Unique[r$], {nL}];
    varst = Table[Unique[t$], {nL}];
    varsλ1 = Table[Unique[λ1$], {nL}];
    varsλ2 = Table[Unique[λ2$], {nL}];
  (*discrete ODE1*)
    eqns1 = -Cos[αα/2] + 4 Cos[αα] varsr + Sin[αα/2] - 2 Cos[αα] varst - 2 varsλ1 +
      2 varsλ2 + 8 Sin[αα] dfdα.varsr - 2 Sin[αα] dfdα.varst - 4 Cos[αα] d2fdα2.varsr;
  (*discrete ODE2*)
    eqns2 = Cos[αα/2] - Cos[3 αα/2] + Sin[αα/2] + 4 varsr Sin[αα]^2 + Sin[3 αα/2] - 8 Sin[αα]^2 varst +
      4 varsλ1 - 4 Cos[αα] varsλ1 + 4 varsλ2 + 4 Cos[αα] varsλ2 - 2 Sin[2 αα] dfdα.varsr +
      8 Sin[2 αα] dfdα.varst + 4 d2fdα2.varst - 4 Cos[2 αα] d2fdα2.varst;
eqns3 = varsλ1 * (varsr * Cos[αα] + varst * Sin[αα] * Tan[αα / 2] - varsr[[1]]);
eqns4 = varsλ2 * (-varsr[[nL]] - (varsr * Cos[αα] - varst * 2 Cos[αα / 2]^2));
  (*boundary conditions*)
    eqns1[[1]] = varsr[[1]] - r0;
```

```mathematica
eqns1[[nL]] = varsr[[nL]] - rπ;
   eqns2[[1]] = varst[[1]] - t0;
eqns2[[nL]] = varst[[nL]] - tπ;
(*governing equations and variables*)
govEq = Join[eqns1, eqns2, eqns3, eqns4];
   allVars = Join[varsr, varst, varsλ1, varsλ2];
(*solve discretized ODEs with initial value r[αα]==r0, t[αα]==t0, λ1=1, λ2=1*)
   root = FindRoot[govEq, Thread[
       {allVars, Flatten[{Table[r0, nL], Table[t0, nL], Table[1, nL], Table[1, nL]}]}],
      MaxIterations → 30, AccuracyGoal → Infinity, PrecisionGoal → 8] // Quiet;
   sol = allVars /. root;
   (*interpolation to obtain r[α] and t[α]*)
   {rsoltemp, tsoltemp, λ1soltemp, λ2soltemp} = Partition[sol, nL];
   rsol = Interpolation[Transpose[{αα, rsoltemp}]];
   tsol = Interpolation[Transpose[{αα, tsoltemp}]];
   λ1sol = Interpolation[Transpose[{αα, λ1soltemp}]];
   λ2sol = Interpolation[Transpose[{αα, λ2soltemp}]];
   (*obtain curves*)
   AA = {rsol[α] * Cos[α], tsol[α] * Sin[α]};
   AAp =
    {rsol[α] * Cos[α] + √2 * Cos[π/4 + α/2], tsol[α] * Sin[α] + √2 * Sin[π/4 + α/2]};
BB = {rsol[α] * Cos[α] - tsol[α] * 2 Cos[α/2]^2, 0};
   BBp = {rsol[α] * Cos[α] - tsol[α] * 2 Cos[α/2]^2 - 1/Sin[α/2], 0};
CC = {rsol[α] * Cos[α] + tsol[α] * Sin[α] * Tan[α/2], 0};
   CCp = {rsol[α] * Cos[α] + tsol[α] * Sin[α] * Tan[α/2] + 1/Cos[α/2], 0};
EAABBx = 1/2 * (-1 + 2 Cos[α] + Cos[2 α]) rsol[α] -
     2 Cos[α/2]^2 Cos[α] tsol[α] + Sin[α] (Cos[α] rsol'[α] - (1 + Cos[α]) tsol'[α]);
EAABBy = -2 Sin[α/2]^2
     (rsol[α] Sin[α] - Sin[α] tsol[α] - Cos[α] rsol'[α] + tsol'[α] + Cos[α] tsol'[α]);
EAACCx = rsol[α] (Cos[α] + Sin[α]^2) - 2 Cos[α] Sin[α/2]^2 tsol[α] +
     Sin[α] (-Cos[α] rsol'[α] + (-1 + Cos[α]) tsol'[α]);
EAACCy = 2 Cos[α/2]^2
     (-rsol[α] Sin[α] + Sin[α] tsol[α] + Cos[α] rsol'[α] + tsol'[α] - Cos[α] tsol'[α]);
EAApBBpx =
    1/2 * ((-1 + 2 Cos[α] + Cos[2 α]) rsol[α] - 2 Sin[α/2] - 4 Cos[α/2]^2 Cos[α] tsol[α] +
       Sin[2 α] rsol'[α] - 2 Sin[α] tsol'[α] - Sin[2 α] tsol'[α]);
EAApBBpy =
    1/2 * (2 Cos[α/2] + rsol[α] (-2 Sin[α] + Sin[2 α]) - 2 (-1 + Cos[α]) Sin[α] tsol[α] -
       rsol'[α] + 2 Cos[α] rsol'[α] - Cos[2 α] rsol'[α] - tsol'[α] + Cos[2 α] tsol'[α]);
EAApCCpx = 1/2 * (2 Cos[α/2] + 2 rsol[α] (Cos[α] + Sin[α]^2) + (1 - 2 Cos[α] + Cos[2 α]) tsol[α] -
       Sin[2 α] rsol'[α] - 2 Sin[α] tsol'[α] + Sin[2 α] tsol'[α]);
EAApCCpy = 1/2 * (2 Sin[α/2] - 2 (1 + Cos[α]) rsol[α] Sin[α] + (2 Sin[α] + Sin[2 α]) tsol[α] +
       rsol'[α] + 2 Cos[α] rsol'[α] + Cos[2 α] rsol'[α] + tsol'[α] - Cos[2 α] tsol'[α]);
   {rr0, rrπ, tt0, ttπ, -rsol[0] + 2 tsol[0] - rsol[π] + 2 tsol[π] +
      NIntegrate[-EAApBBpy * D[EAApBBpx, α] - EAApCCpy * D[EAApCCpx, α],
```

```
              {α, 0, π}, WorkingPrecision → 100, MaxRecursion → 20, AccuracyGoal → 8] +
            NIntegrate[tsol[α] * Sin[α] * D[rsol[α] * Cos[α], α], {α, 0, π},
             WorkingPrecision → 100, MaxRecursion → 20, AccuracyGoal → 8] // Quiet,
         Apply[And, # ≤ 0 & /@ Chop[(varsr * Cos[αα] + varst * Sin[αα] * Tan[αα / 2] - varsr[[1]]) /.
              root, 10^-3]] && Apply[And, # ≤ 0 & /@
             Chop[Rest[(-varsr[[1]] - (varsr * Cos[αα] - varst * Sin[αα] / Tan[αα / 2]))] /. root,
              10^-3]] // Quiet, Norm /@ Partition[govEq /. root, nL]},
        {rr0, 0.6, 0.7, 0.02}, {rrπ, 0.6, 0.7, 0.02}, {tt0, 0.6, 0.7, 0.02},
        {ttπ, 0.6, 0.7, 0.02}];
     (*search from 0.6 to 0.7
      in
      0.02
      intervals*)
     (*running
      may
      take
      very
      long
      time*)
```

```
In[ ]:= Select[Partition[Flatten[searchres], 10],
       #[[6]] == True && #[[7]] < 10^-3 && #[[8]] < 10^-3 && #[[9]] < 10^-3 && #[[10]] < 10^-3 &]

Out[ ]= {{0.6, 0.6, 0.6, 0.6, 2.20531, True, 1.68941×10⁻¹¹, 2.85903×10⁻¹¹,
        6.6008×10⁻¹⁷, 6.58523×10⁻¹⁷}, {0.62, 0.62, 0.62, 0.62, 2.20698,
        True, 1.42688×10⁻¹¹, 2.54382×10⁻¹¹, 5.02515×10⁻¹⁷, 4.00934×10⁻¹⁷},
        {0.64, 0.64, 0.64, 0.64, 2.2074, True, 1.81579×10⁻¹¹, 2.66595×10⁻¹¹,
        6.19592×10⁻¹⁷, 1.08115×10⁻¹⁴}, {0.66, 0.66, 0.66, 0.66, 2.20656,
        True, 1.47005×10⁻¹¹, 2.34574×10⁻¹¹, 5.52483×10⁻¹⁷, 7.09996×10⁻¹⁷},
        {0.68, 0.68, 0.68, 0.68, 2.20446, True, 1.68808×10⁻¹¹, 1.81816×10⁻¹¹,
        5.81064×10⁻¹⁷, 6.68945×10⁻¹⁷}, {0.7, 0.7, 0.7, 0.7, 2.20111, True,
        1.90713×10⁻¹¹, 2.17172×10⁻¹¹, 5.55645×10⁻¹⁷, 6.2887×10⁻¹⁷}}
```


```
     (*Even assuming possible asymmetric, the obtained results were still symmetric*)

In[ ]:= NotebookSave[];
```

```
     (*results*)
```

```mathematica
r0 = 0.64;
rπ = 0.64;
t0 = 0.64;
tπ = 0.64;
pts = 201;
αα = Range[0., π, (π - 0.)/(pts - 1)];
dfdα = NDSolve`FiniteDifferenceDerivative[{1}, {αα}, "DifferenceOrder" -> 4][
    "DifferentiationMatrix"];
d2fdα2 = NDSolve`FiniteDifferenceDerivative[{2}, {αα}, "DifferenceOrder" -> 4][
    "DifferentiationMatrix"];
nL = Length[αα];
varsr = Table[Unique[r$], {nL}];
varst = Table[Unique[t$], {nL}];
varsλ1 = Table[Unique[λ1$], {nL}];
varsλ2 = Table[Unique[λ2$], {nL}];
```

```mathematica
eqns1 = -Cos[αα/2] + 4 Cos[αα] varsr + Sin[αα/2] - 2 Cos[αα] varst - 2 varsλ1 +
    2 varsλ2 + 8 Sin[αα] dfdα.varsr - 2 Sin[αα] dfdα.varst - 4 Cos[αα] d2fdα2.varsr;
eqns2 = Cos[αα/2] - Cos[3 αα/2] + Sin[αα/2] + 4 varsr Sin[αα]^2 + Sin[3 αα/2] - 8 Sin[αα]^2 varst +
    4 varsλ1 - 4 Cos[αα] varsλ1 + 4 varsλ2 + 4 Cos[αα] varsλ2 - 2 Sin[2 αα] dfdα.varsr +
    8 Sin[2 αα] dfdα.varst + 4 d2fdα2.varst - 4 Cos[2 αα] d2fdα2.varst;
eqns3 = varsλ1 * (varsr * Cos[αα] + varst * Sin[αα] * Tan[αα/2] - varsr[[1]]);
eqns4 = varsλ2 * (-varsr[[nL]] - (varsr * Cos[αα] - varst * 2 Cos[αα/2]^2));
```

```mathematica
eqns1[[1]] = varsr[[1]] - r0;
eqns1[[nL]] = varsr[[nL]] - rπ;
eqns2[[1]] = varst[[1]] - t0;
eqns2[[nL]] = varst[[nL]] - tπ;
```

```mathematica
govEq = Join[eqns1, eqns2, eqns3, eqns4];
allVars = Join[varsr, varst, varsλ1, varsλ2];
```

```mathematica
root = FindRoot[govEq,
    Thread[{allVars, Flatten[{Table[r0, nL], Table[t0, nL], Table[1, nL], Table[1, nL]}]}],
    MaxIterations -> 30, AccuracyGoal -> Infinity, PrecisionGoal -> 8];
sol = allVars /. root;
```

> ⋯ FindRoot: The line search decreased the step size to within tolerance specified by AccuracyGoal and PrecisionGoal but was unable to find a sufficient decrease in the merit function. You may need more than MachinePrecision digits of working precision to meet these tolerances.

```mathematica
(*residual*)
```

```mathematica
Norm /@ Partition[govEq /. root, nL]
```

$\{1.83749 \times 10^{-11}, 2.77146 \times 10^{-11}, 4.99483 \times 10^{-17}, 6.77182 \times 10^{-17}\}$

```mathematica
In[•]:= {rsoltemp, tsoltemp, λ1soltemp, λ2soltemp} = Partition[sol, nL];
       rsol = Interpolation[Transpose[{αα, rsoltemp}]];
       tsol = Interpolation[Transpose[{αα, tsoltemp}]];
       λ1sol = Interpolation[Transpose[{αα, λ1soltemp}]];
       λ2sol = Interpolation[Transpose[{αα, λ2soltemp}]];

In[•]:= (*obtain curves*)
       AA = {rsol[α] * Cos[α], tsol[α] * Sin[α]};
       AAp = {rsol[α] * Cos[α] + √2 * Cos[π/4 + α/2], tsol[α] * Sin[α] + √2 * Sin[π/4 + α/2]};
       BB = {rsol[α] * Cos[α] - tsol[α] * 2 Cos[α/2]^2, 0};
       BBp = {rsol[α] * Cos[α] - tsol[α] * 2 Cos[α/2]^2 - 1/Sin[α/2], 0};
       CC = {rsol[α] * Cos[α] + tsol[α] * Sin[α] * Tan[α/2], 0};
       CCp = {rsol[α] * Cos[α] + tsol[α] * Sin[α] * Tan[α/2] + 1/Cos[α/2], 0};
       EAABBx = 1/2 * (-1 + 2 Cos[α] + Cos[2 α]) rsol[α] -
           2 Cos[α/2]^2 Cos[α] tsol[α] + Sin[α] (Cos[α] rsol'[α] - (1 + Cos[α]) tsol'[α]);
       EAABBy = -2 Sin[α/2]^2
           (rsol[α] Sin[α] - Sin[α] tsol[α] - Cos[α] rsol'[α] + tsol'[α] + Cos[α] tsol'[α]);
       EAACCx = rsol[α] (Cos[α] + Sin[α]^2) - 2 Cos[α] Sin[α/2]^2 tsol[α] +
           Sin[α] (-Cos[α] rsol'[α] + (-1 + Cos[α]) tsol'[α]);
       EAACCy = 2 Cos[α/2]^2
           (-rsol[α] Sin[α] + Sin[α] tsol[α] + Cos[α] rsol'[α] + tsol'[α] - Cos[α] tsol'[α]);
       EAApBBpx = 1/2 * ((-1 + 2 Cos[α] + Cos[2 α]) rsol[α] - 2 Sin[α/2] -
               4 Cos[α/2]^2 Cos[α] tsol[α] + Sin[2 α] rsol'[α] - 2 Sin[α] tsol'[α] - Sin[2 α] tsol'[α]);
       EAApBBpy =
           1/2 * (2 Cos[α/2] + rsol[α] (-2 Sin[α] + Sin[2 α]) - 2 (-1 + Cos[α]) Sin[α] tsol[α] -
               rsol'[α] + 2 Cos[α] rsol'[α] - Cos[2 α] rsol'[α] - tsol'[α] + Cos[2 α] tsol'[α]);
       EAApCCpx = 1/2 * (2 Cos[α/2] + 2 rsol[α] (Cos[α] + Sin[α]^2) + (1 - 2 Cos[α] + Cos[2 α]) tsol[α] -
               Sin[2 α] rsol'[α] - 2 Sin[α] tsol'[α] + Sin[2 α] tsol'[α]);
       EAApCCpy = 1/2 * (2 Sin[α/2] - 2 (1 + Cos[α]) rsol[α] Sin[α] + (2 Sin[α] + Sin[2 α]) tsol[α] +
               rsol'[α] + 2 Cos[α] rsol'[α] + Cos[2 α] rsol'[α] + tsol'[α] - Cos[2 α] tsol'[α]);

In[•]:= (*sofa area*)

In[•]:= RealDigits[-rsol[0] + 2 tsol[0] - rsol[π] + 2 tsol[π] +
           NIntegrate[-EAApBBpy * D[EAApBBpx, α] - EAApCCpy * D[EAApCCpx, α],
             {α, 0, π}, WorkingPrecision → 100, MaxRecursion → 20, AccuracyGoal → 8] +
           NIntegrate[tsol[α] * Sin[α] * D[rsol[α] * Cos[α], α], {α, 0, π},
             WorkingPrecision → 100, MaxRecursion → 20, AccuracyGoal → 8]]

Out[•]=
       {{2, 2, 0, 7, 3, 9, 8, 1, 5, 1, 3, 3, 9, 2, 1, 4}, 1}
```

```
Show[ParametricPlot[{EAApCCpx, EAApCCpy}, {α, 0, π}, PerformanceGoal → "Speed",
  PlotStyle → {Red}], ParametricPlot[{EAApBBpx, EAApBBpy}, {α, 0, π},
  PerformanceGoal → "Speed", PlotStyle → {Red}], ParametricPlot[
  {EAACCx, EAACCy}, {α, 0, π}, PerformanceGoal → "Speed", PlotStyle → {Green}],
 ParametricPlot[{EAABBx, EAABBy}, {α, 0, π}, PerformanceGoal → "Speed",
  PlotStyle → {Green}], ParametricPlot[{x, 1}, {x, EAApCCpx /. α → π, 0},
  PerformanceGoal → "Speed", PlotStyle → {Black}], ParametricPlot[{x, 1},
  {x, 0, EAApBBpx /. α → 0}, PerformanceGoal → "Speed", PlotStyle → {Black}],
 ParametricPlot[AA, {α, 0, π}, PerformanceGoal → "Speed", PlotStyle → {Purple}],
 ParametricPlot[{x, 0}, {x, EAABBx /. α → π, EAApBBpx /. α → π},
  PerformanceGoal → "Speed", PlotStyle → {Black}], ParametricPlot[{x, 0},
  {x, EAACCx /. α → 0, EAApCCpx /. α → 0}, PerformanceGoal → "Speed", PlotStyle → {Black}],
 PlotRange → {{-2.0, 2.0}, {0, 1}}, AxesOrigin → {0, 0}, ImageSize → 500]
```

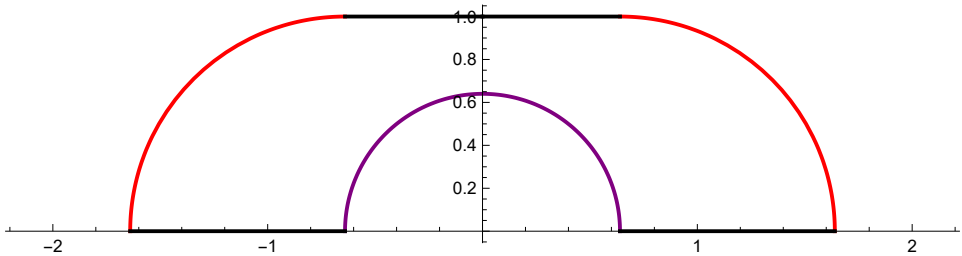

```
Plot[{rsol[α], tsol[α]}, {α, 0, π}, PerformanceGoal → "Speed", PlotRange → {0.6, 0.7},
  PerformanceGoal → "Speed", PlotLegends → {"r(α)", "t(α)"}, PlotTheme → "BoldColor"]
```

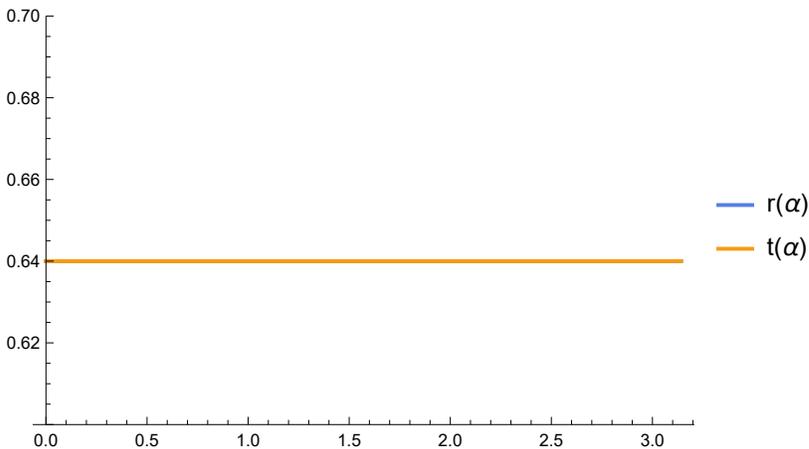

```
(*test if constraints are satisfied with error less than 10^-3*)
```

```
Apply[And, # ≤ 0 & /@
    Chop[(varsr * Cos[αα] + varst * Sin[αα] * Tan[αα / 2] - varsr[[1]]) /. root, 10^-3]] && Apply[
   And, # ≤ 0 & /@ Chop[Rest[(-varsr[[1]] - (varsr * Cos[αα] - varst * Sin[αα] / Tan[αα / 2]))] /.
      root, 10^-3]] // Quiet
```

True

```
In[ ]:= (*more careful search could yield more accurate results*)
In[ ]:= NotebookSave[];
```

(*b*)

In[ ]:= (*intersection point of trajectory of A and AC envelope,
no intersection point of AB envelope*)

In[ ]:= (*define coordinates of points and envelopes*)

In[ ]:= AA = {r[α] * Cos[α], t[α] * Sin[α]};
AAp = {r[α] * Cos[α] + √2 * Cos[π / 4 + α / 2], t[α] * Sin[α] + √2 * Sin[π / 4 + α / 2]};
BB = {r[α] * Cos[α] - t[α] * 2 Cos[α / 2]^2, 0};
BBp = {r[α] * Cos[α] - t[α] * 2 Cos[α / 2]^2 - 1 / Sin[α / 2], 0};
CC = {r[α] * Cos[α] + t[α] * Sin[α] * Tan[α / 2], 0};
CCp = {r[α] * Cos[α] + t[α] * Sin[α] * Tan[α / 2] + 1 / Cos[α / 2], 0};
EAABBx = 1 / 2 * (-1 + 2 Cos[α] + Cos[2 α]) r[α] -
    2 Cos[α / 2]^2 Cos[α] t[α] + Sin[α] (Cos[α] r'[α] - (1 + Cos[α]) t'[α]);
EAABBy = -2 Sin[α / 2]^2 (r[α] Sin[α] - Sin[α] t[α] - Cos[α] r'[α] + t'[α] + Cos[α] t'[α]);
EAACCx = r[α] (Cos[α] + Sin[α]^2) -
    2 Cos[α] Sin[α / 2]^2 t[α] + Sin[α] (-Cos[α] r'[α] + (-1 + Cos[α]) t'[α]);
EAACCy = 2 Cos[α / 2]^2 (-r[α] Sin[α] + Sin[α] t[α] + Cos[α] r'[α] + t'[α] - Cos[α] t'[α]);
EAApBBpx = 1 / 2 * ((-1 + 2 Cos[α] + Cos[2 α]) r[α] - 2 Sin[α / 2] -
        4 Cos[α / 2]^2 Cos[α] t[α] + Sin[2 α] r'[α] - 2 Sin[α] t'[α] - Sin[2 α] t'[α]);
EAApBBpy = 1 / 2 * (2 Cos[α / 2] + r[α] (-2 Sin[α] + Sin[2 α]) - 2 (-1 + Cos[α]) Sin[α] t[α] -
        r'[α] + 2 Cos[α] r'[α] - Cos[2 α] r'[α] - t'[α] + Cos[2 α] t'[α]);
EAApCCpx = 1 / 2 * (2 Cos[α / 2] + 2 r[α] (Cos[α] + Sin[α]^2) +
        (1 - 2 Cos[α] + Cos[2 α]) t[α] - Sin[2 α] r'[α] - 2 Sin[α] t'[α] + Sin[2 α] t'[α]);
EAApCCpy = 1 / 2 * (2 Sin[α / 2] - 2 (1 + Cos[α]) r[α] Sin[α] + (2 Sin[α] + Sin[2 α]) t[α] +
        r'[α] + 2 Cos[α] r'[α] + Cos[2 α] r'[α] + t'[α] - Cos[2 α] t'[α]);

In[ ]:= (*assume intersection point was p for trajectory of point A and AC envelope
    AA〚1〛/.α→α1p⩵EAACCx/.α→α2p
    AA〚2〛/.α→α1p⩵EAACCy/.α→α2p*)

In[ ]:= α1p = 77 / 1000;

In[ ]:= α2p = 177 / 100;

In[ ]:= (*area functional*)

```mathematica
In[ ]:= FF = Piecewise[{{0, α < 0},
         {-EAApBBpy * D[EAApBBpx, α] - EAApCCpy * D[EAApCCpx, α]                                    -
            r[0] + 2 t[0] - r[π] + 2 t[π], 0 ≤ α ≤ α1p},

         {-EAApBBpy * D[EAApBBpx, α] - EAApCCpy * D[EAApCCpx, α] +
            AA[[2]] * D[AA[[1]], α] - r[0] + 2 t[0] - r[π] + 2 t[π], α1p < α ≤ 1},

         {-EAApBBpy * D[EAApBBpx, α] - EAApCCpy * D[EAApCCpx, α] +
            AA[[2]] * D[AA[[1]], α]                                    , 1 < α ≤ α2p},

         {-EAApBBpy * D[EAApBBpx, α] - EAApCCpy * D[EAApCCpx, α] +
            AA[[2]] * D[AA[[1]], α] - EAACCy * D[EAACCx, α]            , α2p < α ≤ π},

         {0, π < α}}] // Simplify;

In[ ]:= (*Lagrange multipliers for constraints*)

In[ ]:= HH = FF + λ2[α] * (r[π] - (r[α] * Cos[α] - t[α] * 2 Cos[α / 2]^2));

In[ ]:= (*EL equations*)

In[ ]:= (D[HH, r[α]] - D[D[HH, r'[α]], α] + D[D[HH, r''[α]], {α, 2}]) // Simplify
```

Out[ ]=

$$\begin{cases} -\cos[\alpha]\,\lambda 2[\alpha] & \alpha > \pi \;||\; \alpha < 0 \\ -\tfrac{1}{2}\cos[\alpha]\,\bigl(-\cos[\tfrac{\alpha}{2}] + 4\cos[\alpha]\,r[\alpha] + \sin[\tfrac{\alpha}{2}] - 2\cos[\alpha]\,t[\alpha] + \\ \quad 2\,\lambda 2[\alpha] + 8\sin[\alpha]\,r'[\alpha] - 2\sin[\alpha]\,t'[\alpha] - 4\cos[\alpha]\,r''[\alpha]\bigr) & \tfrac{77}{1000} < \alpha < \tfrac{177}{100} \\ -\tfrac{1}{2}\cos[\alpha]\,\bigl(-\cos[\tfrac{\alpha}{2}] + 4\cos[\alpha]\,r[\alpha] + \sin[\tfrac{\alpha}{2}] - 4\cos[\alpha]\,t[\alpha] + \\ \quad 2\,\lambda 2[\alpha] + 8\sin[\alpha]\,r'[\alpha] - 4\sin[\alpha]\,t'[\alpha] - 4\cos[\alpha]\,r''[\alpha]\bigr) & 0 < \alpha < \tfrac{77}{1000} \\ -\tfrac{1}{2}\cos[\alpha]\,\bigl(-\cos[\tfrac{\alpha}{2}] + 2\,(3\cos[\alpha] + \cos[2\alpha])\,r[\alpha] + \\ \quad \sin[\tfrac{\alpha}{2}] - 2\,(2\cos[\alpha] + \cos[2\alpha])\,t[\alpha] + 2\,\lambda 2[\alpha] + 12\sin[\alpha]\,r'[\alpha] + \\ \quad 3\sin[2\alpha]\,r'[\alpha] - 4\sin[\alpha]\,t'[\alpha] - 3\sin[2\alpha]\,t'[\alpha] - r''[\alpha] - \\ \quad 6\cos[\alpha]\,r''[\alpha] - \cos[2\alpha]\,r''[\alpha] - t''[\alpha] + \cos[2\alpha]\,t''[\alpha]\bigr) & \tfrac{177}{100} < \alpha < \pi \\ \text{Indeterminate} & \text{True} \end{cases}$$

```
In[ ]:= (D[HH, t[α]] - D[D[HH, t'[α]], α] + D[D[HH, t''[α]], {α, 2}]) // Simplify
Out[ ]=
    ⎧ (1 + Cos[α]) λ2[α]                                                                              α > π || α < 0
    ⎪ 1/2 Cos[α/2] Sin[α] + 1/2 Sin[α/2] Sin[α] + 2 r[α] Sin[α]² - 2 Sin[α]² t[α] + λ2[α] +            0 < α ≤ 77/1000
    ⎪   Cos[α] λ2[α] - 2 Cos[α] Sin[α] r'[α] + 4 Cos[α] Sin[α] t'[α] + 2 Sin[α]² t''[α]
    ⎨ 1/2 Cos[α/2] Sin[α] + 1/2 Sin[α/2] Sin[α] + r[α] Sin[α]² - 2 Sin[α]² t[α] + λ2[α] +              77/1000 < α < 177/100
    ⎪   Cos[α] λ2[α] - Cos[α] Sin[α] r'[α] + 4 Cos[α] Sin[α] t'[α] + 2 Sin[α]² t''[α]
    ⎪ 1/2 (Cos[α/2] Sin[α] + Sin[α/2] Sin[α] - 4 (-1 + Cos[α]) r[α] Sin[α]² +                          177/100 < α < π
    ⎪   2 (-3 + 2 Cos[α]) Sin[α]² t[α] + 2 λ2[α] + 2 Cos[α] λ2[α] -
    ⎪   Sin[α] r'[α] - 4 Cos[α] Sin[α] r'[α] + 3 Cos[2 α] Sin[α] r'[α] -
    ⎪   Sin[α] t'[α] - 3 Cos[2 α] Sin[α] t'[α] + 6 Sin[2 α] t'[α] +
    ⎪   Sin[α] Sin[2 α] r''[α] + 6 Sin[α]² t''[α] - 2 Cos[α] Sin[α]² t''[α])
    ⎩ Indeterminate                                                                                   True
```

---

```
In[ ]:= (*fine search r[0] r[π] t[0] t[π]*)
    searchres = ParallelTable[
    α1p = 0.078;
    α2p = 1.77;
    r0 = rr0;
      rπ = rrπ;
      t0 = tt0;
      tπ = ttπ;
      pts = 301;
      αα = Range[0., π, (π - 0.) / (pts - 1)];
      dfdα = NDSolve`FiniteDifferenceDerivative[
        {1}, {αα}, "DifferenceOrder" → 4]["DifferentiationMatrix"];
      d2fdα2 = NDSolve`FiniteDifferenceDerivative[
        {2}, {αα}, "DifferenceOrder" → 4]["DifferentiationMatrix"];
    nL = Length[αα];
    Do[
    varsr = Table[Unique[r$], {nL}];
      varst = Table[Unique[t$], {nL}];
      
      varsλ2 = Table[Unique[λ2$], {nL}];
      
    (*ODE1*)
    eqns1 = -(1 / 2) Cos[αα] * (-Cos[αα / 2] + 4 Cos[αα] varsr + Sin[αα / 2] -
          4 Cos[αα] varst + 2 varsλ2 + 8 Sin[αα] dfdα.varsr - 4 Sin[αα] dfdα.varst -
```

```
            4 Cos[αα] d2fdα2.varsr) * (Boole[0 ≤ # < α1p] & /@ αα) +
        -(1/2) Cos[αα] * (-Cos[αα/2] + 4 Cos[αα] varsr + Sin[αα/2] -
            2 Cos[αα] varst + 2 varsλ2 + 8 Sin[αα] dfdα.varsr - 2 Sin[αα] dfdα.varst -
            4 Cos[αα] d2fdα2.varsr) * (Boole[α1p ≤ # < α2p] & /@ αα) +
        -(1/2) Cos[αα] * (-Cos[αα/2] + 2 (3 Cos[αα] + Cos[2 αα]) varsr + Sin[αα/2] -
            2 (2 Cos[αα] + Cos[2 αα]) varst + 2 varsλ2 + 12 Sin[αα] dfdα.varsr +
            3 Sin[2 αα] dfdα.varsr - 4 Sin[αα] dfdα.varst - 3 Sin[2 αα] dfdα.varst -
            d2fdα2.varsr - 6 Cos[αα] d2fdα2.varsr - Cos[2 αα] d2fdα2.varsr -
            d2fdα2.varst + Cos[2 αα] d2fdα2.varst) * (Boole[α2p ≤ # ≤ π] & /@ αα);
(*ODE2*)
eqns2 = ((1/2) Cos[αα/2] Sin[αα] + (1/2) Sin[αα/2] Sin[αα] + 2 varsr Sin[αα]^2 -
            2 Sin[αα]^2 varst + varsλ2 + Cos[αα] varsλ2 - 2 Cos[αα] Sin[αα] dfdα.varsr + 4 Cos[αα]
             Sin[αα] dfdα.varst + 2 Sin[αα]^2 d2fdα2.varst) * (Boole[0 ≤ # < α1p] & /@ αα) +
        ((1/2) Cos[αα/2] Sin[αα] + (1/2) Sin[αα/2] Sin[αα] + varsr Sin[αα]^2 -
            2 Sin[αα]^2 varst + varsλ2 + Cos[αα] varsλ2 - Cos[αα] Sin[αα] dfdα.varsr + 4 Cos[αα]
             Sin[αα] dfdα.varst + 2 Sin[αα]^2 d2fdα2.varst) * (Boole[α1p ≤ # < α2p] & /@ αα) +
        (Cos[αα/2] Sin[αα] + Sin[αα/2] Sin[αα] - 4 (-1 + Cos[αα]) varsr Sin[αα]^2 +
            2 (-3 + 2 Cos[αα]) Sin[αα]^2 varst + 2 varsλ2 + 2 Cos[αα] varsλ2 - Sin[αα] dfdα.varsr -
            4 Cos[αα] Sin[αα] dfdα.varsr + 3 Cos[2 αα] Sin[αα] dfdα.varsr -
            Sin[αα] dfdα.varst - 3 Cos[2 αα] Sin[αα] dfdα.varst + 6 Sin[2 αα] dfdα.varst +
            Sin[αα] Sin[2 αα] d2fdα2.varsr + 6 Sin[αα]^2 d2fdα2.varst -
            2 Cos[αα] Sin[αα]^2 d2fdα2.varst) * (Boole[α2p ≤ # ≤ π] & /@ αα);
eqns4 = varsλ2 * (-varsr[[nL]] - (varsr * Cos[αα] - varst * 2 Cos[αα/2]^2));
(*boundary conditions*)

eqns1[[1]] = varsr[[1]] - r0;
    eqns1[[nL]] = varsr[[nL]] - rπ;
eqns2[[1]] = varst[[1]] - t0;
eqns2[[nL]] = varst[[nL]] - tπ;

    govEq = Join[eqns1, eqns2, eqns4];
    allVars = Join[varsr, varst, varsλ2];
(*solve discretized ODEs*) root = FindRoot[govEq,
        Thread[{allVars, Flatten[{Table[r0, nL], Table[t0, nL], Table[0.1, nL]}]}],
        MaxIterations → 20, AccuracyGoal → Infinity, PrecisionGoal → 8] // Quiet;
    sol = allVars /. root;

{rsoltemp, tsoltemp, λ2soltemp} = Partition[sol, nL];

    (*adjust r[π/2] outlier value*)
    rsoltemp[[Ceiling[pts/2]]] =
     (rsoltemp[[Ceiling[pts/2] - 1]] + rsoltemp[[Ceiling[pts/2] + 1]])/2;
    rsol = Interpolation[Transpose[{αα, rsoltemp}]];
    tsol = Interpolation[Transpose[{αα, tsoltemp}]];
    λ2sol = Interpolation[Transpose[{αα, λ2soltemp}]];
(*obtain curves*)
```

```mathematica
      AA = {rsol[α] * Cos[α], tsol[α] * Sin[α]};
      AAp =
       {rsol[α] * Cos[α] + √2 * Cos[π / 4 + α / 2], tsol[α] * Sin[α] + √2 * Sin[π / 4 + α / 2]};
     BB = {rsol[α] * Cos[α] - tsol[α] * 2 Cos[α / 2]^2, 0};
      BBp = {rsol[α] * Cos[α] - tsol[α] * 2 Cos[α / 2]^2 - 1 / Sin[α / 2], 0};
     CC = {rsol[α] * Cos[α] + tsol[α] * Sin[α] * Tan[α / 2], 0};
      CCp = {rsol[α] * Cos[α] + tsol[α] * Sin[α] * Tan[α / 2] + 1 / Cos[α / 2], 0};
     EAABBx = 1 / 2 * (-1 + 2 Cos[α] + Cos[2 α]) rsol[α] -
        2 Cos[α / 2]^2 Cos[α] tsol[α] + Sin[α] (Cos[α] rsol'[α] - (1 + Cos[α]) tsol'[α]);
     EAABBy = -2 Sin[α / 2]^2
        (rsol[α] Sin[α] - Sin[α] tsol[α] - Cos[α] rsol'[α] + tsol'[α] + Cos[α] tsol'[α]);
     EAACCx = rsol[α] (Cos[α] + Sin[α]^2) - 2 Cos[α] Sin[α / 2]^2 tsol[α] +
        Sin[α] (-Cos[α] rsol'[α] + (-1 + Cos[α]) tsol'[α]);
     EAACCy = 2 Cos[α / 2]^2
        (-rsol[α] Sin[α] + Sin[α] tsol[α] + Cos[α] rsol'[α] + tsol'[α] - Cos[α] tsol'[α]);
     EAApBBpx = 1 / 2 * ((-1 + 2 Cos[α] + Cos[2 α]) rsol[α] - 2 Sin[α / 2] - 4 Cos[α / 2]^2 Cos[α]
            tsol[α] + Sin[2 α] rsol'[α] - 2 Sin[α] tsol'[α] - Sin[2 α] tsol'[α]);
     EAApBBpy =
      1 / 2 * (2 Cos[α / 2] + rsol[α] (-2 Sin[α] + Sin[2 α]) - 2 (-1 + Cos[α]) Sin[α] tsol[α] -
          rsol'[α] + 2 Cos[α] rsol'[α] - Cos[2 α] rsol'[α] - tsol'[α] + Cos[2 α] tsol'[α]);
     EAApCCpx = 1 / 2 * (2 Cos[α / 2] + 2 rsol[α] (Cos[α] + Sin[α]^2) + (1 - 2 Cos[α] + Cos[2 α]) tsol[α] -
          Sin[2 α] rsol'[α] - 2 Sin[α] tsol'[α] + Sin[2 α] tsol'[α]);
     EAApCCpy = 1 / 2 * (2 Sin[α / 2] - 2 (1 + Cos[α]) rsol[α] Sin[α] + (2 Sin[α] + Sin[2 α]) tsol[α] +
          rsol'[α] + 2 Cos[α] rsol'[α] + Cos[2 α] rsol'[α] + tsol'[α] - Cos[2 α] tsol'[α]);
     (*find α1p α2p*)
     ααres1 = FindRoot[{rsol[αα1] * Cos[αα1] == (EAACCx /. α → αα2),
          tsol[αα1] * Sin[αα1] == (EAACCy /. α → αα2)}, {αα1, 0.05}, {αα2, 2.2}] // Quiet;
     α1p = ααres1[[1]][[2]];
     α2p = ααres1[[2]][[2]];
     , {6}];
     {rr0, rrπ, tt0, ttπ, -rsol[0] + 2 tsol[0] - rsol[π] + 2 tsol[π] +
         NIntegrate[-EAApBBpy * D[EAApBBpx, α] - EAApCCpy * D[EAApCCpx, α],
          {α, 0, π}, WorkingPrecision → 100, MaxRecursion → 20, AccuracyGoal → 8] +
         NIntegrate[tsol[α] * Sin[α] * D[rsol[α] * Cos[α], α], {α, α1p, π},
          WorkingPrecision → 100, MaxRecursion → 20, AccuracyGoal → 8] +
         NIntegrate[-EAACCy * D[EAACCx, α], {α, α2p, π}, WorkingPrecision → 100,
          MaxRecursion → 20, AccuracyGoal → 8] // Quiet,
        {α1p, α2p}, Norm /@ Partition[govEq /. root, nL]}, {rr0, 0.6, 0.7, 0.05},
       {rrπ, 0.6, 0.7, 0.05}, {tt0, 0.6, 0.7, 0.05}, {ttπ, 0.6, 0.7, 0.05}];
     (*search from 0.6 to 0.7 in
     0.05
     intervals*)
```

In[ ]:= (*running may take very long time*)

```
(*obtain r[0], r[π], t[0], t[π], area, α1p,
 α2p, residual with area>0 and residual<10^-3*)
```

```
Select[Partition[Flatten[searchres], 10],
  3 ≥ #[[5]] ≥ 0 && #[[8]] ≤ 10^-3 && #[[9]] ≤ 10^-3 && #[[10]] ≤ 10^-3 &]
```

Out[ ]=
{{0.6, 0.6, 0.6, 0.6, 2.17662, 0.0504319, 2.25611, 3.19912×10⁻¹¹,
  3.59937×10⁻¹¹, 1.6422×10⁻¹⁶}, {0.6, 0.6, 0.6, 0.65, 2.19808, 0.0505413,
  2.25687, 2.95515×10⁻¹¹, 3.51137×10⁻¹¹, 2.05247×10⁻¹⁶}, {0.6, 0.6, 0.6, 0.7,
  2.21952, 0.0506508, 2.25766, 2.91749×10⁻¹¹, 3.64782×10⁻¹¹, 1.98345×10⁻¹⁶},
 {0.65, 0.65, 0.65, 0.6, 2.15828, 0.0434649, 2.27242, 3.28816×10⁻¹¹,
  3.67609×10⁻¹¹, 1.65433×10⁻¹⁶}, {0.65, 0.65, 0.65, 0.65, 2.17975,
  0.0435685, 2.27303, 3.03736×10⁻¹¹, 3.45379×10⁻¹¹, 1.50472×10⁻¹⁶},
 {0.65, 0.65, 0.65, 0.7, 2.20121, 0.0436722, 2.27368, 3.48702×10⁻¹¹,
  3.81316×10⁻¹¹, 1.62786×10⁻¹⁶}, {0.7, 0.7, 0.7, 0.6, 2.13259,
  0.0370562, 2.29066, 3.9583×10⁻¹¹, 4.7185×10⁻¹¹, 1.47247×10⁻¹⁶},
 {0.7, 0.7, 0.7, 0.65, 2.15408, 0.0371551, 2.2917, 3.6526×10⁻¹¹, 4.7135×10⁻¹¹, 1.7059×10⁻¹⁶},
 {0.7, 0.7, 0.7, 0.7, 2.17544, 0.0375987,
  2.29336, 3.35095×10⁻¹¹, 3.89714×10⁻¹¹, 1.50256×10⁻¹⁶}}

```
NotebookSave[];
```

```
(*results*)
```

```
(*initial guess of α1p α2p*)
```

```
α1p = 0.078;
```

```
α2p = 1.77;
```

```mathematica
In[ ]:= ri = 0.6;
r0 = 0.6;
rπ = 0.6;
t0 = 0.6;
tπ = 0.7;
pts = 301;
αα = Range[0., π, (π - 0.) / (pts - 1)];
dfdα = NDSolve`FiniteDifferenceDerivative[{1}, {αα}, "DifferenceOrder" → 4][
    "DifferentiationMatrix"];
d2fdα2 = NDSolve`FiniteDifferenceDerivative[{2}, {αα}, "DifferenceOrder" → 4][
    "DifferentiationMatrix"];
nL = Length[αα];

In[ ]:= (*Find α1p and α2p iteratively*)
Do[
 varsr = Table[Unique[r$], {nL}];
  varst = Table[Unique[t$], {nL}];

  varsλ2 = Table[Unique[λ2$], {nL}];

 (*ODE1*)
 eqns1 = -1/2 Cos[αα] *
      (-Cos[αα/2] + 4 Cos[αα] varsr + Sin[αα/2] - 4 Cos[αα] varst + 2 varsλ2 + 8 Sin[αα] dfdα.varsr -
        4 Sin[αα] dfdα.varst - 4 Cos[αα] d2fdα2.varsr) * (Boole[0 ≤ # < α1p] & /@ αα) +
    -1/2 Cos[αα] *
      (-Cos[αα/2] + 4 Cos[αα] varsr + Sin[αα/2] - 2 Cos[αα] varst + 2 varsλ2 + 8 Sin[αα] dfdα.varsr -
        2 Sin[αα] dfdα.varst - 4 Cos[αα] d2fdα2.varsr) * (Boole[α1p ≤ # < α2p] & /@ αα) +
    -1/2 Cos[αα] *
      (-Cos[αα/2] + 2 (3 Cos[αα] + Cos[2 αα]) varsr + Sin[αα/2] - 2 (2 Cos[αα] + Cos[2 αα]) varst +
        2 varsλ2 + 12 Sin[αα] dfdα.varsr + 3 Sin[2 αα] dfdα.varsr - 4 Sin[αα] dfdα.varst -
        3 Sin[2 αα] dfdα.varst - d2fdα2.varsr - 6 Cos[αα] d2fdα2.varsr - Cos[2 αα]
          d2fdα2.varsr - d2fdα2.varst + Cos[2 αα] d2fdα2.varst) * (Boole[α2p ≤ # ≤ π] & /@ αα);
 (*ODE2*)
 eqns2 = (1/2 Cos[αα/2] Sin[αα] + 1/2 Sin[αα/2] Sin[αα] + 2 varsr Sin[αα]^2 - 2 Sin[αα]^2 varst + varsλ2 +
        Cos[αα] varsλ2 - 2 Cos[αα] Sin[αα] dfdα.varsr + 4 Cos[αα] Sin[αα] dfdα.varst +
        2 Sin[αα]^2 d2fdα2.varst) * (Boole[0 ≤ # < α1p] & /@ αα) +
```

```mathematica
    (1/2 Cos[αα/2] Sin[αα] + 1/2 Sin[αα/2] Sin[αα] + varsr Sin[αα]^2 - 2 Sin[αα]^2 varst + varsλ2 + 
       Cos[αα] varsλ2 - Cos[αα] Sin[αα] dfdα.varsr + 4 Cos[αα] Sin[αα] dfdα.varst + 
       2 Sin[αα]^2 d2fdα2.varst) * (Boole[α1p ≤ # < α2p] & /@ αα) + 
    (Cos[αα/2] Sin[αα] + Sin[αα/2] Sin[αα] - 4 (-1 + Cos[αα]) varsr Sin[αα]^2 + 
       2 (-3 + 2 Cos[αα]) Sin[αα]^2 varst + 2 varsλ2 + 2 Cos[αα] varsλ2 - Sin[αα] dfdα.varsr - 
       4 Cos[αα] Sin[αα] dfdα.varsr + 3 Cos[2 αα] Sin[αα] dfdα.varsr - 
       Sin[αα] dfdα.varst - 3 Cos[2 αα] Sin[αα] dfdα.varst + 6 Sin[2 αα] dfdα.varst + 
       Sin[αα] Sin[2 αα] d2fdα2.varsr + 6 Sin[αα]^2 d2fdα2.varst - 
       2 Cos[αα] Sin[αα]^2 d2fdα2.varst) * (Boole[α2p ≤ # ≤ π] & /@ αα);
eqns4 = varsλ2 * (-varsr[[nL]] - (varsr * Cos[αα] - varst * 2 Cos[αα/2]^2));
(*boundary conditions*)

eqns1[[1]] = varsr[[1]] - r0;
  eqns1[[nL]] = varsr[[nL]] - rπ;
eqns2[[1]] = varst[[1]] - t0;
eqns2[[nL]] = varst[[nL]] - tπ;

  govEq = Join[eqns1, eqns2, eqns4];
  allVars = Join[varsr, varst, varsλ2];
(*solve discretized ODEs*) root = FindRoot[govEq,
    Thread[{allVars, Flatten[{Table[r0, nL], Table[t0, nL], Table[0.1, nL]}]}],
    MaxIterations → 20, AccuracyGoal → Infinity, PrecisionGoal → 8] // Quiet;
  sol = allVars /. root;
{rsoltemp, tsoltemp, λ2soltemp} = Partition[sol, nL];

  (*adjust r[π/2] outlier value*)
  rsoltemp[[Ceiling[pts / 2]]] = 
   (rsoltemp[[Ceiling[pts / 2] - 1]] + rsoltemp[[Ceiling[pts / 2] + 1]]) / 2;
  rsol = Interpolation[Transpose[{αα, rsoltemp}]];
  tsol = Interpolation[Transpose[{αα, tsoltemp}]];
  λ2sol = Interpolation[Transpose[{αα, λ2soltemp}]];
(*obtain curves*)
  AA = {rsol[α] * Cos[α], tsol[α] * Sin[α]};
  AAp = {rsol[α] * Cos[α] + √2 * Cos[π / 4 + α / 2], tsol[α] * Sin[α] + √2 * Sin[π / 4 + α / 2]};
BB = {rsol[α] * Cos[α] - tsol[α] * 2 Cos[α / 2]^2, 0};
  BBp = {rsol[α] * Cos[α] - tsol[α] * 2 Cos[α / 2]^2 - 1 / Sin[α / 2], 0};
CC = {rsol[α] * Cos[α] + tsol[α] * Sin[α] * Tan[α / 2], 0};
  CCp = {rsol[α] * Cos[α] + tsol[α] * Sin[α] * Tan[α / 2] + 1 / Cos[α / 2], 0};
EAABBx = 1 / 2 * (-1 + 2 Cos[α] + Cos[2 α]) rsol[α] - 
    2 Cos[α / 2]^2 Cos[α] tsol[α] + Sin[α] (Cos[α] rsol'[α] - (1 + Cos[α]) tsol'[α]);
EAABBy = -2 Sin[α / 2]^2
    (rsol[α] Sin[α] - Sin[α] tsol[α] - Cos[α] rsol'[α] + tsol'[α] + Cos[α] tsol'[α]);
```

```
        EAACCx = rsol[α] (Cos[α] + Sin[α]^2) - 2 Cos[α] Sin[α/2]^2 tsol[α] +
          Sin[α] (-Cos[α] rsol'[α] + (-1 + Cos[α]) tsol'[α]);
        EAACCy = 2 Cos[α/2]^2
          (-rsol[α] Sin[α] + Sin[α] tsol[α] + Cos[α] rsol'[α] + tsol'[α] - Cos[α] tsol'[α]);
        EAApBBpx =
          1/2 * ((-1 + 2 Cos[α] + Cos[2 α]) rsol[α] - 2 Sin[α/2] - 4 Cos[α/2]^2 Cos[α] tsol[α] +
            Sin[2 α] rsol'[α] - 2 Sin[α] tsol'[α] - Sin[2 α] tsol'[α]);
        EAApBBpy =
          1/2 * (2 Cos[α/2] + rsol[α] (-2 Sin[α] + Sin[2 α]) - 2 (-1 + Cos[α]) Sin[α] tsol[α] -
            rsol'[α] + 2 Cos[α] rsol'[α] - Cos[2 α] rsol'[α] - tsol'[α] + Cos[2 α] tsol'[α]);
        EAApCCpx = 1/2 * (2 Cos[α/2] + 2 rsol[α] (Cos[α] + Sin[α]^2) + (1 - 2 Cos[α] + Cos[2 α]) tsol[α] -
            Sin[2 α] rsol'[α] - 2 Sin[α] tsol'[α] + Sin[2 α] tsol'[α]);
        EAApCCpy = 1/2 * (2 Sin[α/2] - 2 (1 + Cos[α]) rsol[α] Sin[α] + (2 Sin[α] + Sin[2 α]) tsol[α] +
            rsol'[α] + 2 Cos[α] rsol'[α] + Cos[2 α] rsol'[α] + tsol'[α] - Cos[2 α] tsol'[α]);
        (*find α1p α2p*)
        ααres1 = FindRoot[{rsol[αα1] * Cos[αα1] == (EAACCx /. α → αα2),
            tsol[αα1] * Sin[αα1] == (EAACCy /. α → αα2)}, {αα1, 0.1}, {αα2, 1.8}];
        α1p = ααres1[[1]][[2]];
        α2p = ααres1[[2]][[2]];, {6}(*iterate for 6 times*)];
```

In[ ]:= `{α1p, α2p}`

Out[ ]=
`{0.0506508, 2.25766}`

In[ ]:= 
```
    -rsol[0] + 2 tsol[0] - rsol[π] + 2 tsol[π] +
     NIntegrate[-EAApBBpy * D[EAApBBpx, α] - EAApCCpy * D[EAApCCpx, α],
       {α, 0, π}, WorkingPrecision → 100, MaxRecursion → 20, AccuracyGoal → 8] +
     NIntegrate[tsol[α] * Sin[α] * D[rsol[α] * Cos[α], α], {α, α1p, π}, WorkingPrecision → 100,
       MaxRecursion → 20, AccuracyGoal → 8] + NIntegrate[-EAACCy * D[EAACCx, α],
       {α, α2p, π}, WorkingPrecision → 100, MaxRecursion → 20, AccuracyGoal → 8]
```

Out[ ]=
2.21952

In[ ]:= `RealDigits[%]`

Out[ ]=
`{{2, 2, 1, 9, 5, 2, 3, 3, 7, 6, 3, 6, 3, 9, 9, 9}, 1}`

```
In[ ]:= Show[ParametricPlot[{EAApCCpx, EAApCCpy},
    {α, 0 + 0.0923, π - 0.01(*slightly adjust to eliminate error*)},
    PerformanceGoal → "Speed", PlotStyle → {Red}], ParametricPlot[{EAApBBpx, EAApBBpy},
    {α, 0, π}, PerformanceGoal → "Speed", PlotStyle → {Red}], ParametricPlot[
   {EAACCx, EAACCy}, {α, α2p, π}, PerformanceGoal → "Speed", PlotStyle → {Green}],
  (*ParametricPlot[{EAABBx,EAABBy},{α,0,π},PerformanceGoal→"Speed",
   PlotStyle→{Green}],*)ParametricPlot[{x, 1}, {x, EAApCCpx /. α → π - 0.01, 0},
    PerformanceGoal → "Speed", PlotStyle → {Black}], ParametricPlot[{x, 1},
   {x, 0, EAApBBpx /. α → 0}, PerformanceGoal → "Speed", PlotStyle → {Black}],
  ParametricPlot[AA, {α, α1p, π}, PerformanceGoal → "Speed", PlotStyle → {Purple}],
  ParametricPlot[{x, 0}, {x, EAABBx /. α → π, EAApBBpx /. α → π},
    PerformanceGoal → "Speed", PlotStyle → {Black}], ParametricPlot[{x, 0},
   {x, EAACCx /. α → π, EAApCCpx /. α → 0}, PerformanceGoal → "Speed", PlotStyle → {Black}],
  PlotRange → {{-2.0, 2.0}, {0, 1}}, AxesOrigin → {0, 0}, ImageSize → 500]
```

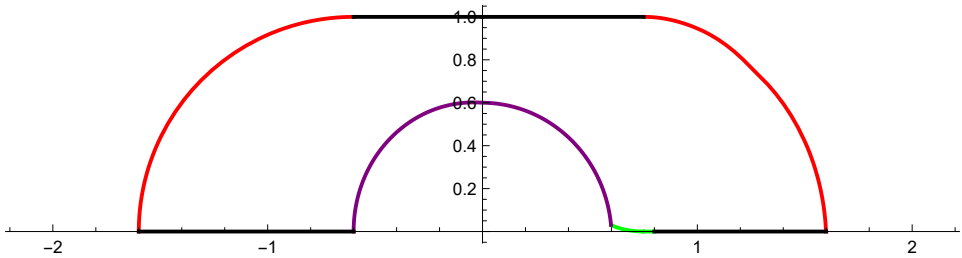

```
In[ ]:= (*obtained parametric equations r(α),t(α)*)

In[ ]:= Plot[{rsol[α], tsol[α]}, {α, 0, π}, PerformanceGoal → "Speed", PlotRange → All,
    PerformanceGoal → "Speed", PlotLegends → {"r(α)", "t(α)"}, PlotTheme → "BoldColor"]
```

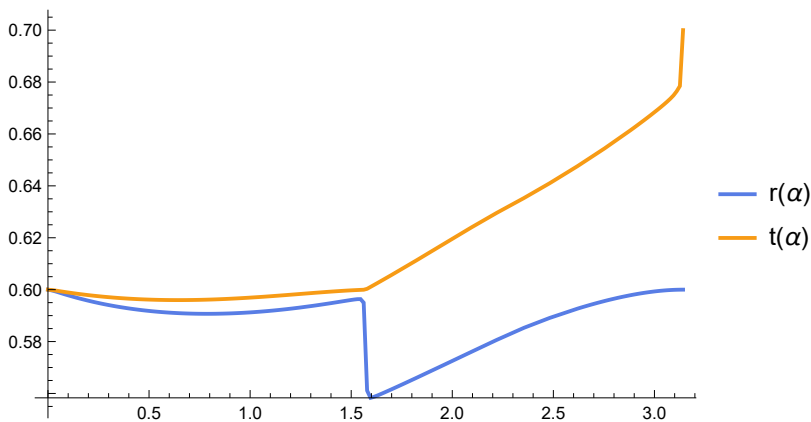

```
In[ ]:= (*more careful search could yield more accurate results*)

In[ ]:= NotebookSave[];
```

```
In[ ]:= (*C*)

In[ ]:= (*intersection point of trajectory of A and AC envelope and AB envelope*)

In[ ]:= AA = {r[α]*Cos[α], t[α]*Sin[α]};
    AAp = {r[α]*Cos[α] + √2 *Cos[π/4+α/2], t[α]*Sin[α] + √2 *Sin[π/4+α/2]};
    BB = {r[α]*Cos[α] - t[α]*2 Cos[α/2]^2, 0};
    BBp = {r[α]*Cos[α] - t[α]*2 Cos[α/2]^2 - 1/Sin[α/2], 0};
    CC = {r[α]*Cos[α] + t[α]*Sin[α]*Tan[α/2], 0};
    CCp = {r[α]*Cos[α] + t[α]*Sin[α]*Tan[α/2] + 1/Cos[α/2], 0};
    EAABBx = 1/2*(-1+2 Cos[α]+Cos[2 α]) r[α] -
       2 Cos[α/2]^2 Cos[α] t[α] + Sin[α] (Cos[α] r'[α] - (1+Cos[α]) t'[α]);
    EAABBy = -2 Sin[α/2]^2 (r[α] Sin[α] - Sin[α] t[α] - Cos[α] r'[α] + t'[α] + Cos[α] t'[α]);
    EAACCx = r[α] (Cos[α]+Sin[α]^2) -
       2 Cos[α] Sin[α/2]^2 t[α] + Sin[α] (-Cos[α] r'[α] + (-1+Cos[α]) t'[α]);
    EAACCy = 2 Cos[α/2]^2 (-r[α] Sin[α] + Sin[α] t[α] + Cos[α] r'[α] + t'[α] - Cos[α] t'[α]);
    EAApBBpx = 1/2*((-1+2 Cos[α]+Cos[2 α]) r[α] - 2 Sin[α/2] -
          4 Cos[α/2]^2 Cos[α] t[α] + Sin[2 α] r'[α] - 2 Sin[α] t'[α] - Sin[2 α] t'[α]);
    EAApBBpy = 1/2*(2 Cos[α/2] + r[α] (-2 Sin[α]+Sin[2 α]) - 2 (-1+Cos[α]) Sin[α] t[α] -
          r'[α] + 2 Cos[α] r'[α] - Cos[2 α] r'[α] - t'[α] + Cos[2 α] t'[α]);
    EAApCCpx = 1/2*(2 Cos[α/2] + 2 r[α] (Cos[α]+Sin[α]^2) +
          (1 - 2 Cos[α]+Cos[2 α]) t[α] - Sin[2 α] r'[α] - 2 Sin[α] t'[α] + Sin[2 α] t'[α]);
    EAApCCpy = 1/2*(2 Sin[α/2] - 2 (1+Cos[α]) r[α] Sin[α] + (2 Sin[α]+Sin[2 α]) t[α] +
          r'[α] + 2 Cos[α] r'[α] + Cos[2 α] r'[α] + t'[α] - Cos[2 α] t'[α]);

In[ ]:= (*assume intersection point was p for trajectory of point A and AC envelope
           AA〚1〛/.α→α1p==EAACCx/.α→α2p
           AA〚2〛/.α→α1p==EAACCy/.α→α2p
          AA〚1〛/.α→α3p==EAABBx/.α→α4p
          AA〚2〛/.α→α3p==EAABBy/.α→α4p*)

In[ ]:= α1p = 77/100;

In[ ]:= α2p = 177/100;

In[ ]:= α3p = 200/100;

In[ ]:= α4p = 30/100;
```

```
In[ ]:= FF = Piecewise[{{0, α < 0},
          {-EAApBBpy * D[EAApBBpx, α] - EAApCCpy * D[EAApCCpx, α]                                    -
            EAABBy * D[EAABBx, α] - r[0] + 2 t[0] - r[π] + 2 t[π], 0 ≤ α ≤ α4p},

          {-EAApBBpy * D[EAApBBpx, α] -
            EAApCCpy * D[EAApCCpx, α]                                                                -
            r[0] + 2 t[0] - r[π] + 2 t[π], α4p ≤ α ≤ α1p},

          {-EAApBBpy * D[EAApBBpx, α] - EAApCCpy * D[EAApCCpx, α] +
            AA〚2〛 * D[AA〚1〛, α]                             - r[0] + 2 t[0] - r[π] + 2 t[π], α1p < α ≤ 1},

          {-EAApBBpy * D[EAApBBpx, α] - EAApCCpy * D[EAApCCpx, α] +
            AA〚2〛 * D[AA〚1〛, α]                                                                      , 1 < α ≤ α2p},

          {-EAApBBpy * D[EAApBBpx, α] - EAApCCpy * D[EAApCCpx, α] + AA〚2〛 * D[AA〚1〛, α] -
            EAACCy * D[EAACCx, α]                       , α2p < α ≤ α3p},

          {-EAApBBpy * D[EAApBBpx, α] - EAApCCpy * D[EAApCCpx, α]                                    -
            EAACCy * D[EAACCx, α]                       , α3p < α ≤ π},

          {0, π < α}}] // Simplify;

In[ ]:= (*EL equations*)

In[ ]:= (D[FF, r[α]] - D[D[FF, r'[α]], α] + D[D[FF, r''[α]], {α, 2}]) // Simplify
```

Out[ ]=

$$\begin{cases}
-\frac{1}{2}\cos[\alpha]\left(-\cos\left[\frac{\alpha}{2}\right]+4\cos[\alpha]r[\alpha]+\sin\left[\frac{\alpha}{2}\right]-\right. & \frac{77}{100}<\alpha<\frac{177}{100} \\
\quad 2\cos[\alpha]t[\alpha]+8\sin[\alpha]r'[\alpha]-2\sin[\alpha]t'[\alpha]-4\cos[\alpha]r''[\alpha]\Big) & \\
-\frac{1}{2}\cos[\alpha]\left(-\cos\left[\frac{\alpha}{2}\right]+4\cos[\alpha]r[\alpha]+\sin\left[\frac{\alpha}{2}\right]-\right. & \frac{3}{10}<\alpha<\frac{77}{100} \\
\quad 4\cos[\alpha]t[\alpha]+8\sin[\alpha]r'[\alpha]-4\sin[\alpha]t'[\alpha]-4\cos[\alpha]r''[\alpha]\Big) & \\
-\frac{1}{2}\cos[\alpha] & \frac{177}{100}<\alpha<2 \\
\quad\left(-\cos\left[\frac{\alpha}{2}\right]+2(3\cos[\alpha]+\cos[2\alpha])r[\alpha]+\sin\left[\frac{\alpha}{2}\right]-2(2\cos[\alpha]+\cos[2\alpha])t[\alpha]+\right. & \\
\quad 12\sin[\alpha]r'[\alpha]+3\sin[2\alpha]r'[\alpha]-4\sin[\alpha]t'[\alpha]-3\sin[2\alpha]t'[\alpha]- & \\
\quad r''[\alpha]-6\cos[\alpha]r''[\alpha]-\cos[2\alpha]r''[\alpha]-t''[\alpha]+\cos[2\alpha]t''[\alpha]\Big) & \\
-\frac{1}{2}\cos[\alpha] & 2<\alpha<\pi \\
\quad\left(-\cos\left[\frac{\alpha}{2}\right]+2(3\cos[\alpha]+\cos[2\alpha])r[\alpha]+\sin\left[\frac{\alpha}{2}\right]-2(3\cos[\alpha]+\cos[2\alpha])t[\alpha]+\right. & \\
\quad 12\sin[\alpha]r'[\alpha]+3\sin[2\alpha]r'[\alpha]-6\sin[\alpha]t'[\alpha]-3\sin[2\alpha]t'[\alpha]- & \\
\quad r''[\alpha]-6\cos[\alpha]r''[\alpha]-\cos[2\alpha]r''[\alpha]-t''[\alpha]+\cos[2\alpha]t''[\alpha]\Big) & \\
-\frac{1}{2}\cos[\alpha]\left(-\cos\left[\frac{\alpha}{2}\right]-2(-3\cos[\alpha]+\cos[2\alpha])r[\alpha]+\right. & 0<\alpha<\frac{3}{10} \\
\quad \sin\left[\frac{\alpha}{2}\right]+2(-3\cos[\alpha]+\cos[2\alpha])t[\alpha]+12\sin[\alpha]r'[\alpha]- & \\
\quad 3\sin[2\alpha]r'[\alpha]-6\sin[\alpha]t'[\alpha]+3\sin[2\alpha]t'[\alpha]+r''[\alpha]- & \\
\quad 6\cos[\alpha]r''[\alpha]+\cos[2\alpha]r''[\alpha]+t''[\alpha]-\cos[2\alpha]t''[\alpha]\Big) & \\
0 & \alpha>\pi\ ||\ \alpha<0 \\
\text{Indeterminate} & \text{True}
\end{cases}$$

```
In[•]:= (D[FF, t[α]] - D[D[FF, t'[α]], α] + D[D[FF, t''[α]], {α, 2}]) // Simplify
```

Out[•]=

$$\begin{cases} \frac{1}{2} \sin[\alpha] \left( \cos\left[\frac{\alpha}{2}\right] + \sin\left[\frac{\alpha}{2}\right] + 4\, r[\alpha] \sin[\alpha] - \right. \\ \qquad \left. 4 \sin[\alpha]\, t[\alpha] - 4 \cos[\alpha]\, r'[\alpha] + 8 \cos[\alpha]\, t'[\alpha] + 4 \sin[\alpha]\, t''[\alpha] \right) & \frac{3}{10} < \alpha \leq \frac{77}{100} \\[4pt] \frac{1}{2} \sin[\alpha] \left( \cos\left[\frac{\alpha}{2}\right] + \sin\left[\frac{\alpha}{2}\right] + 2\, r[\alpha] \sin[\alpha] - \right. \\ \qquad \left. 4 \sin[\alpha]\, t[\alpha] - 2 \cos[\alpha]\, r'[\alpha] + 8 \cos[\alpha]\, t'[\alpha] + 4 \sin[\alpha]\, t''[\alpha] \right) & \frac{77}{100} < \alpha < \frac{177}{100} \\[4pt] \frac{1}{2} \sin[\alpha] & 2 < \alpha < \pi \\[2pt] \left( \cos\left[\frac{\alpha}{2}\right] + \sin\left[\frac{\alpha}{2}\right] - 2\, r[\alpha]\, (-3\sin[\alpha] + \sin[2\alpha]) + 2\, (-3\sin[\alpha] + \sin[2\alpha]) \right. \\ \qquad t[\alpha] - r'[\alpha] - 6 \cos[\alpha]\, r'[\alpha] + 3 \cos[2\alpha]\, r'[\alpha] - t'[\alpha] + 12 \cos[\alpha]\, t'[\alpha] - \\ \qquad \left. 3 \cos[2\alpha]\, t'[\alpha] + \sin[2\alpha]\, r''[\alpha] + 6 \sin[\alpha]\, t''[\alpha] - \sin[2\alpha]\, t''[\alpha] \right) \\[4pt] \frac{1}{2} \sin[\alpha] & \frac{177}{100} < \alpha \leq 2 \\[2pt] \left( \cos\left[\frac{\alpha}{2}\right] + \sin\left[\frac{\alpha}{2}\right] + r[\alpha]\, (4\sin[\alpha] - 2\sin[2\alpha]) + 2\, (-3\sin[\alpha] + \sin[2\alpha])\, t[\alpha] - \right. \\ \qquad r'[\alpha] - 4 \cos[\alpha]\, r'[\alpha] + 3 \cos[2\alpha]\, r'[\alpha] - t'[\alpha] + 12 \cos[\alpha]\, t'[\alpha] - \\ \qquad \left. 3 \cos[2\alpha]\, t'[\alpha] + \sin[2\alpha]\, r''[\alpha] + 6 \sin[\alpha]\, t''[\alpha] - \sin[2\alpha]\, t''[\alpha] \right) \\[4pt] \frac{1}{2} \sin[\alpha] & 0 < \alpha < \frac{3}{10} \\[2pt] \left( \cos\left[\frac{\alpha}{2}\right] + \sin\left[\frac{\alpha}{2}\right] + 2\, r[\alpha]\, (3\sin[\alpha] + \sin[2\alpha]) - 2\, (3\sin[\alpha] + \sin[2\alpha])\, t[\alpha] + \right. \\ \qquad r'[\alpha] - 6 \cos[\alpha]\, r'[\alpha] - 3 \cos[2\alpha]\, r'[\alpha] + t'[\alpha] + 12 \cos[\alpha]\, t'[\alpha] + \\ \qquad \left. 3 \cos[2\alpha]\, t'[\alpha] - \sin[2\alpha]\, r''[\alpha] + 6 \sin[\alpha]\, t''[\alpha] + \sin[2\alpha]\, t''[\alpha] \right) \\[4pt] 0 & \alpha > \pi\ ||\ \alpha < 0 \\ \text{Indeterminate} & \text{True} \end{cases}$$

```
In[•]:= (*Boundary condition*)

In[•]:= FF0 = -EAApBBpy * D[EAApBBpx, α] - EAApCCpy * D[EAApCCpx, α] -
        EAABBy * D[EAABBx, α] - r[0] + 2 t[0] - r[π] + 2 t[π];

In[•]:= (D[FF0, r'[α], α] - D[D[FF0, r''[α], α], α]) /. α → 0
```

Out[•]=

$$-\frac{1}{2} + 2\, r[0] - 2\, t[0] - 2\, r''[0]$$

```
In[•]:= FFπ = -EAApBBpy * D[EAApBBpx, α] - EAApCCpy * D[EAApCCpx, α] - EAACCy * D[EAACCx, α];

In[•]:= (D[FFπ, r'[α], α] - D[D[FFπ, r''[α], α], α]) /. α → π
```

Out[•]=

$$-\frac{1}{2} + 2\, r[\pi] - 2\, t[\pi] - 2\, r''[\pi]$$

```
In[•]:= NotebookSave[];
```

```mathematica
In[ ]:= (*fine search for t[0], t[π]*)
    searchres = ParallelTable[(*iteratively find α1p α2p α3p α4p*)
    (*initial guess of α1p α2p α3p α4p*)
    α1p = 0.077;
    α2p = 1.77;
    α3p = 3.06;
    α4p = 1.37;
    ri = 0.6;
       t0 = tt0;
       tπ = ttπ;
       (*number of discretization points*)
       pts = 302;
    (*αα is discrete parameter*)
       αα = Range[0., π, (π - 0.) / (pts - 1)];
       (*discretization setting*)
    dfdα = NDSolve`FiniteDifferenceDerivative[
         {1}, {αα}, "DifferenceOrder" → 4]["DifferentiationMatrix"];
       d2fdα2 = NDSolve`FiniteDifferenceDerivative[
         {2}, {αα}, "DifferenceOrder" → 4]["DifferentiationMatrix"];
    nL = Length[αα];
    (*Find α1p α2p α3p α4p iteratively*)
       Do[
    varsr = Table[Unique[r$], {nL}];
       varst = Table[Unique[t$], {nL}];

    (*discrete ODE1*)
    eqns1 =
         (-Cos[αα/2] - 2 (-3 Cos[αα] + Cos[2 αα]) varsr + Sin[αα/2] + 2 (-3 Cos[αα] + Cos[2 αα])
              varst + 12 Sin[αα] dfdα.varsr - 3 Sin[2 αα] dfdα.varsr - 6 Sin[αα] dfdα.varst +
              3 Sin[2 αα] dfdα.varst + d2fdα2.varsr - 6 Cos[αα] d2fdα2.varsr + Cos[2 αα] d2fdα2.
                varsr + d2fdα2.varst - Cos[2 αα] d2fdα2.varst) * (Boole[0 ≤ # < α1p] & /@ αα) +

           (-Cos[αα/2] - 2 (-3 Cos[αα] + Cos[2 αα]) varsr +
               Sin[αα/2] + 2 (-2 Cos[αα] + Cos[2 αα]) varst + 12 Sin[αα] dfdα.varsr -
               3 Sin[2 αα] dfdα.varsr - 4 Sin[αα] dfdα.varst + 3 Sin[2 αα] dfdα.varst +
               d2fdα2.varsr - 6 Cos[αα] d2fdα2.varsr + Cos[2 αα] d2fdα2.varsr +
               d2fdα2.varst - Cos[2 αα] d2fdα2.varst) * (Boole[α1p ≤ # < α4p] & /@ αα) +

           (-Cos[αα/2] + 4 Cos[αα] varsr + Sin[αα/2] - 2 Cos[αα] varst + 8 Sin[αα] dfdα.varsr -
               2 Sin[αα] dfdα.varst - 4 Cos[αα] d2fdα2.varsr) * (Boole[α4p ≤ # < α2p] & /@ αα) +
```

$$\left(-\text{Cos}\left[\frac{\alpha\alpha}{2}\right] + 2\,(3\,\text{Cos}[\alpha\alpha] + \text{Cos}[2\,\alpha\alpha])\,\text{varsr} + \text{Sin}\left[\frac{\alpha\alpha}{2}\right] - 2\,(2\,\text{Cos}[\alpha\alpha] + \text{Cos}[2\,\alpha\alpha])\,\text{varst} + \right.$$
$$12\,\text{Sin}[\alpha\alpha]\,\text{dfd}\alpha.\text{varsr} + 3\,\text{Sin}[2\,\alpha\alpha]\,\text{dfd}\alpha.\text{varsr} - 4\,\text{Sin}[\alpha\alpha]\,\text{dfd}\alpha.\text{varst} - 3\,\text{Sin}[2\,\alpha\alpha]$$
$$\text{dfd}\alpha.\text{varst} - \text{d2fd}\alpha2.\text{varsr} - 6\,\text{Cos}[\alpha\alpha]\,\text{d2fd}\alpha2.\text{varsr} - \text{Cos}[2\,\alpha\alpha]\,\text{d2fd}\alpha2.\text{varsr} -$$
$$\left.\text{d2fd}\alpha2.\text{varst} + \text{Cos}[2\,\alpha\alpha]\,\text{d2fd}\alpha2.\text{varst}\right) * (\text{Boole}[\alpha2p \leq \# < \alpha3p]\,\&\,/@\,\alpha\alpha) +$$
$$\left(-\text{Cos}\left[\frac{\alpha\alpha}{2}\right] + 2\,(3\,\text{Cos}[\alpha\alpha] + \text{Cos}[2\,\alpha\alpha])\,\text{varsr} + \text{Sin}\left[\frac{\alpha\alpha}{2}\right] - 2\,(3\,\text{Cos}[\alpha\alpha] + \text{Cos}[2\,\alpha\alpha])\,\text{varst} + \right.$$
$$12\,\text{Sin}[\alpha\alpha]\,\text{dfd}\alpha.\text{varsr} + 3\,\text{Sin}[2\,\alpha\alpha]\,\text{dfd}\alpha.\text{varsr} - 6\,\text{Sin}[\alpha\alpha]\,\text{dfd}\alpha.\text{varst} - 3\,\text{Sin}[2\,\alpha\alpha]$$
$$\text{dfd}\alpha.\text{varst} - \text{d2fd}\alpha2.\text{varsr} - 6\,\text{Cos}[\alpha\alpha]\,\text{d2fd}\alpha2.\text{varsr} - \text{Cos}[2\,\alpha\alpha]\,\text{d2fd}\alpha2.\text{varsr} -$$
$$\left.\text{d2fd}\alpha2.\text{varst} + \text{Cos}[2\,\alpha\alpha]\,\text{d2fd}\alpha2.\text{varst}\right) * (\text{Boole}[\alpha3p \leq \# \leq \pi]\,\&\,/@\,\alpha\alpha);$$

(*discrete ODE2*)
eqns2 =
$$\left(\text{Cos}\left[\frac{\alpha\alpha}{2}\right] + \text{Sin}\left[\frac{\alpha\alpha}{2}\right] + 2\,\text{varsr}\,(3\,\text{Sin}[\alpha\alpha] + \text{Sin}[2\,\alpha\alpha]) - 2\,(3\,\text{Sin}[\alpha\alpha] + \text{Sin}[2\,\alpha\alpha])\,\text{varst} + \right.$$
$$\text{dfd}\alpha.\text{varsr} - 6\,\text{Cos}[\alpha\alpha]\,\text{dfd}\alpha.\text{varsr} - 3\,\text{Cos}[2\,\alpha\alpha]\,\text{dfd}\alpha.\text{varsr} + \text{dfd}\alpha.\text{varst} +$$
$$12\,\text{Cos}[\alpha\alpha]\,\text{dfd}\alpha.\text{varst} + 3\,\text{Cos}[2\,\alpha\alpha]\,\text{dfd}\alpha.\text{varst} - \text{Sin}[2\,\alpha\alpha]\,\text{d2fd}\alpha2.\text{varsr} +$$
$$\left.6\,\text{Sin}[\alpha\alpha]\,\text{d2fd}\alpha2.\text{varst} + \text{Sin}[2\,\alpha\alpha]\,\text{d2fd}\alpha2.\text{varst}\right) * (\text{Boole}[0 \leq \# < \alpha1p]\,\&\,/@\,\alpha\alpha) +$$
$$\left(\text{Cos}\left[\frac{\alpha\alpha}{2}\right] + \text{Sin}\left[\frac{\alpha\alpha}{2}\right] + 4\,(1 + \text{Cos}[\alpha\alpha])\,\text{varsr}\,\text{Sin}[\alpha\alpha] - 2\,(3\,\text{Sin}[\alpha\alpha] + \text{Sin}[2\,\alpha\alpha])\,\text{varst} + \right.$$
$$\text{dfd}\alpha.\text{varsr} - 4\,\text{Cos}[\alpha\alpha]\,\text{dfd}\alpha.\text{varsr} - 3\,\text{Cos}[2\,\alpha\alpha]\,\text{dfd}\alpha.\text{varsr} + \text{dfd}\alpha.\text{varst} +$$
$$12\,\text{Cos}[\alpha\alpha]\,\text{dfd}\alpha.\text{varst} + 3\,\text{Cos}[2\,\alpha\alpha]\,\text{dfd}\alpha.\text{varst} - \text{Sin}[2\,\alpha\alpha]\,\text{d2fd}\alpha2.\text{varsr} +$$
$$\left.6\,\text{Sin}[\alpha\alpha]\,\text{d2fd}\alpha2.\text{varst} + \text{Sin}[2\,\alpha\alpha]\,\text{d2fd}\alpha2.\text{varst}\right) * (\text{Boole}[\alpha1p \leq \# < \alpha4p]\,\&\,/@\,\alpha\alpha) +$$
$$\left(\text{Cos}\left[\frac{\alpha\alpha}{2}\right] + \text{Sin}\left[\frac{\alpha\alpha}{2}\right] + 2\,\text{varsr}\,\text{Sin}[\alpha\alpha] - 4\,\text{Sin}[\alpha\alpha]\,\text{varst} - 2\,\text{Cos}[\alpha\alpha]\,\text{dfd}\alpha.\text{varsr} + \right.$$
$$\left.8\,\text{Cos}[\alpha\alpha]\,\text{dfd}\alpha.\text{varst} + 4\,\text{Sin}[\alpha\alpha]\,\text{d2fd}\alpha2.\text{varst}\right) * (\text{Boole}[\alpha4p \leq \# < \alpha2p]\,\&\,/@\,\alpha\alpha) +$$
$$\left(\text{Cos}\left[\frac{\alpha\alpha}{2}\right] + \text{Sin}\left[\frac{\alpha\alpha}{2}\right] + \text{varsr}\,(4\,\text{Sin}[\alpha\alpha] - 2\,\text{Sin}[2\,\alpha\alpha]) + 2\,(-3\,\text{Sin}[\alpha\alpha] + \text{Sin}[2\,\alpha\alpha])\,\text{varst} - \right.$$
$$\text{dfd}\alpha.\text{varsr} - 4\,\text{Cos}[\alpha\alpha]\,\text{dfd}\alpha.\text{varsr} + 3\,\text{Cos}[2\,\alpha\alpha]\,\text{dfd}\alpha.\text{varsr} - \text{dfd}\alpha.\text{varst} +$$
$$12\,\text{Cos}[\alpha\alpha]\,\text{dfd}\alpha.\text{varst} - 3\,\text{Cos}[2\,\alpha\alpha]\,\text{dfd}\alpha.\text{varst} + \text{Sin}[2\,\alpha\alpha]\,\text{d2fd}\alpha2.\text{varsr} +$$
$$\left.6\,\text{Sin}[\alpha\alpha]\,\text{d2fd}\alpha2.\text{varst} - \text{Sin}[2\,\alpha\alpha]\,\text{d2fd}\alpha2.\text{varst}\right) * (\text{Boole}[\alpha2p \leq \# < \alpha3p]\,\&\,/@\,\alpha\alpha) +$$
$$\left(\text{Cos}\left[\frac{\alpha\alpha}{2}\right] + \text{Sin}\left[\frac{\alpha\alpha}{2}\right] - 2\,\text{varsr}\,(-3\,\text{Sin}[\alpha\alpha] + \text{Sin}[2\,\alpha\alpha]) + 2\,(-3\,\text{Sin}[\alpha\alpha] + \text{Sin}[2\,\alpha\alpha])\,\right.$$
$$\text{varst} - \text{dfd}\alpha.\text{varsr} - 6\,\text{Cos}[\alpha\alpha]\,\text{dfd}\alpha.\text{varsr} + 3\,\text{Cos}[2\,\alpha\alpha]\,\text{dfd}\alpha.\text{varsr} - \text{dfd}\alpha.\text{varst} +$$
$$12\,\text{Cos}[\alpha\alpha]\,\text{dfd}\alpha.\text{varst} - 3\,\text{Cos}[2\,\alpha\alpha]\,\text{dfd}\alpha.\text{varst} + \text{Sin}[2\,\alpha\alpha]\,\text{d2fd}\alpha2.\text{varsr} +$$
$$\left.6\,\text{Sin}[\alpha\alpha]\,\text{d2fd}\alpha2.\text{varst} - \text{Sin}[2\,\alpha\alpha]\,\text{d2fd}\alpha2.\text{varst}\right) * (\text{Boole}[\alpha3p \leq \# \leq \pi]\,\&\,/@\,\alpha\alpha);$$

(*boundary conditions*)
eqns1[[1]] = (-1/2 + 2 varsr[[1]] - 2 varst[[1]] - 2 (d2fdα2.varsr)[[1]]) - 0;
eqns1[[nL]] = (-1/2 + 2 varsr[[nL]] - 2 varst[[nL]] - 2 (d2fdα2.varsr)[[nL]]) - 0;
eqns2[[1]] = varst[[1]] - t0;
eqns2[[nL]] = varst[[nL]] - tπ;
    govEq = Join[eqns1, eqns2];

```mathematica
    allVars = Join[varsr, varst];
(*solve discretized ODEs*)
    root = FindRoot[govEq, Thread[{allVars, Flatten[{Table[ri, nL], Table[t0, nL]}]}],
       MaxIterations → 20, AccuracyGoal → Infinity, PrecisionGoal → 8] // Quiet;
    sol = allVars /. root;
(*interpolation to obtain r[α] and t[α]*)
    {rsoltemp, tsoltemp} = Partition[sol, nL];
    rsol = Interpolation[Transpose[{αα, rsoltemp}]];
    tsol = Interpolation[Transpose[{αα, tsoltemp}]];
    (*obtain curves*)
AA = {rsol[α] * Cos[α], tsol[α] * Sin[α]};
    AAp =
      {rsol[α] * Cos[α] + √2 * Cos[π/4 + α/2], tsol[α] * Sin[α] + √2 * Sin[π/4 + α/2]};
BB = {rsol[α] * Cos[α] - tsol[α] * 2 Cos[α/2]^2, 0};
    BBp = {rsol[α] * Cos[α] - tsol[α] * 2 Cos[α/2]^2 - 1/Sin[α/2], 0};
CC = {rsol[α] * Cos[α] + tsol[α] * Sin[α] * Tan[α/2], 0};
    CCp = {rsol[α] * Cos[α] + tsol[α] * Sin[α] * Tan[α/2] + 1/Cos[α/2], 0};
EAABBx = 1/2 * (-1 + 2 Cos[α] + Cos[2 α]) rsol[α] -
    2 Cos[α/2]^2 Cos[α] tsol[α] + Sin[α] (Cos[α] rsol'[α] - (1 + Cos[α]) tsol'[α]);
EAABBy = -2 Sin[α/2]^2
     (rsol[α] Sin[α] - Sin[α] tsol[α] - Cos[α] rsol'[α] + tsol'[α] + Cos[α] tsol'[α]);
EAACCx = rsol[α] (Cos[α] + Sin[α]^2) - 2 Cos[α] Sin[α/2]^2 tsol[α] +
    Sin[α] (-Cos[α] rsol'[α] + (-1 + Cos[α]) tsol'[α]);
EAACCy = 2 Cos[α/2]^2
     (-rsol[α] Sin[α] + Sin[α] tsol[α] + Cos[α] rsol'[α] + tsol'[α] - Cos[α] tsol'[α]);
EAApBBpx = 1/2 * ((-1 + 2 Cos[α] + Cos[2 α]) rsol[α] - 2 Sin[α/2] - 4 Cos[α/2]^2 Cos[α]
        tsol[α] + Sin[2 α] rsol'[α] - 2 Sin[α] tsol'[α] - Sin[2 α] tsol'[α]);
EAApBBpy =
    1/2 * (2 Cos[α/2] + rsol[α] (-2 Sin[α] + Sin[2 α]) - 2 (-1 + Cos[α]) Sin[α] tsol[α] -
       rsol'[α] + 2 Cos[α] rsol'[α] - Cos[2 α] rsol'[α] - tsol'[α] + Cos[2 α] tsol'[α]);
EAApCCpx = 1/2 * (2 Cos[α/2] + 2 rsol[α] (Cos[α] + Sin[α]^2) + (1 - 2 Cos[α] + Cos[2 α]) tsol[α] -
       Sin[2 α] rsol'[α] - 2 Sin[α] tsol'[α] + Sin[2 α] tsol'[α]);
EAApCCpy = 1/2 * (2 Sin[α/2] - 2 (1 + Cos[α]) rsol[α] Sin[α] + (2 Sin[α] + Sin[2 α]) tsol[α] +
       rsol'[α] + 2 Cos[α] rsol'[α] + Cos[2 α] rsol'[α] + tsol'[α] - Cos[2 α] tsol'[α]);
(*find α1p α2p*)
ααres1 = FindRoot[{rsol[αα1] * Cos[αα1] == (EAACCx /. α → αα2),
       tsol[αα1] * Sin[αα1] == (EAACCy /. α → αα2)}, {αα1, 0.1}, {αα2, 1.8}] // Quiet;
(*find α3p α4p*)
ααres2 = FindRoot[{rsol[αα3] * Cos[αα3] == (EAABBx /. α → αα4),
       tsol[αα3] * Sin[αα3] == (EAABBy /. α → αα4)}, {αα3, 3.0}, {αα4, 1.0}] // Quiet;
    α1p = ααres1[[1]][[2]];
α2p = ααres1[[2]][[2]];
    α3p = ααres2[[1]][[2]];
α4p = ααres2[[2]][[2]];, {5}(*iterate for 5 times*)];
{tt0, ttπ, -rsol[0] + 2 tsol[0] - rsol[π] + 2 tsol[π] +
```

```
            NIntegrate[-EAApBBpy * D[EAApBBpx, α] - EAApCCpy * D[EAApCCpx, α], {α, 0, π},
               WorkingPrecision → 100, MaxRecursion → 20, AccuracyGoal → 8] + NIntegrate[
               tsol[α] * Sin[α] * D[rsol[α] * Cos[α], α], {α, α1p, α3p}, WorkingPrecision → 100,
               MaxRecursion → 20, AccuracyGoal → 8] + NIntegrate[-EAABBy * D[EAABBx, α],
               {α, 0, α4p}, WorkingPrecision → 100, MaxRecursion → 20, AccuracyGoal → 8] +
              NIntegrate[-EAACCy * D[EAACCx, α], {α, α2p, π},
               WorkingPrecision → 100, MaxRecursion → 20, AccuracyGoal → 8] // Quiet,
             {α1p, α2p, α3p, α4p}, Norm /@ Partition[govEq /. root, nL]},
           {tt0, 0.7, 0.72, 0.001}, {ttπ, 0.7, 0.72, 0.001}];
         (*search from 0.7 to 0.72 in 0.001 intervals*)
In[ ]:=  (*running may take long time*)

In[ ]:=  (*obtain t[0], t[π], area, α1p, α2p, α3p, α4p, residual with area>0 and residual<10^-3*)

In[ ]:=  Select[Partition[Flatten[searchres], 9], #[[3]] ≥ 0 && #[[8]] ≤ 10^-3 && #[[9]] ≤ 10^-3 &]
Out[ ]=
         {{0.7, 0.7, 2.2193, 0.0853935, 1.72169, 3.0562, 1.4199, 8.87668 × 10^-11, 1.20715 × 10^-10},
           {0.701, 0.701, 2.21934, 0.0843159, 1.72785,
            3.05728, 1.41374, 8.78872 × 10^-11, 1.08595 × 10^-10},
           {0.703, 0.703, 2.21942, 0.0829828, 1.73919, 3.05861, 1.4024, 8.63498 × 10^-11, 1.05998 × 10^-10},
           {0.704, 0.704, 2.21944, 0.0826942, 1.74337, 3.0589, 1.39823, 8.25423 × 10^-11, 1.08879 × 10^-10},
           {0.705, 0.705, 2.21947, 0.0816284, 1.75079, 3.05996, 1.39081,
            7.93628 × 10^-11, 1.13934 × 10^-10}, {0.706, 0.706, 2.21949, 0.0813535,
            1.75456, 3.06024, 1.38704, 8.55004 × 10^-11, 1.14495 × 10^-10},
           {0.707, 0.707, 2.21951, 0.0802957, 1.76231, 3.0613, 1.37928, 8.85518 × 10^-11, 1.03577 × 10^-10},
           {0.708, 0.708, 2.21952, 0.0800343, 1.76567, 3.06156, 1.37593, 9.3471 × 10^-11, 1.05709 × 10^-10},
           {0.709, 0.709, 2.21953, 0.0789847, 1.77508,
            3.06261, 1.36652, 8.15778 × 10^-11, 1.14061 × 10^-10},
           {0.71, 0.71, 2.21953, 0.0787365, 1.77669, 3.06286, 1.3649, 8.95577 × 10^-11, 1.09105 × 10^-10},
           {0.711, 0.711, 2.21953, 0.0776949, 1.78644, 3.0639, 1.35515, 8.55745 × 10^-11, 1.08817 × 10^-10},
           {0.712, 0.712, 2.21952, 0.0774596, 1.78762,
            3.06413, 1.35397, 8.71775 × 10^-11, 1.11663 × 10^-10},
           {0.713, 0.713, 2.21952, 0.076426, 1.79765, 3.06517, 1.34395, 8.73475 × 10^-11, 1.17649 × 10^-10},
           {0.714, 0.714, 2.2195, 0.0762034, 1.79847, 3.06539, 1.34313, 8.46762 × 10^-11, 1.1145 × 10^-10},
           {0.715, 0.715, 2.21949, 0.0751777, 1.80869,
            3.06641, 1.33291, 8.75324 × 10^-11, 1.17535 × 10^-10},
           {0.716, 0.716, 2.21946, 0.0749675, 1.80922, 3.06663, 1.33238, 8.5241 × 10^-11, 1.18078 × 10^-10},
           {0.717, 0.717, 2.21944, 0.0739497, 1.81957, 3.06764, 1.32202, 7.98793 × 10^-11, 1.1406 × 10^-10},
           {0.718, 0.718, 2.21941, 0.0737517, 1.81987, 3.06784, 1.32172,
            8.67537 × 10^-11, 1.18384 × 10^-10}, {0.719, 0.719, 2.21938, 0.0727417,
            1.83031, 3.06885, 1.31129, 8.53631 × 10^-11, 1.08449 × 10^-10},
           {0.72, 0.72, 2.21934, 0.0725559, 1.83042, 3.06904, 1.31117, 8.79518 × 10^-11, 1.15666 × 10^-10}}
```

```mathematica
(*very close to symmetric results with t[0]=0.710322*)

(*With other different t[0] t[π], it cannot obtain convergent results*)

(*Even assuming possible asymmetric,the obtained results were still symmetric*)

NotebookSave[];
```

```mathematica
(*results*)
```

```mathematica
(*initial guess of α1p α2p α3p α4p*)
α1p = 0.077;
α2p = 1.77;
α3p = 3.06;
α4p = 1.37;
ri = 0.6;
(*obtained optimal t[0] and t[π]*)
t0 = 0.710;
tπ = 0.710;
(*number of discretization points*)
pts = 302;
(*discretization points*)
αα = Range[0., π, (π - 0.) / (pts - 1)];
(*discretization setting*)
dfdα = NDSolve`FiniteDifferenceDerivative[{1}, {αα}, "DifferenceOrder" → 4][
   "DifferentiationMatrix"];
d2fdα2 = NDSolve`FiniteDifferenceDerivative[{2}, {αα}, "DifferenceOrder" → 4][
   "DifferentiationMatrix"];
nL = Length[αα];
(*Find α1p α2p α3p α4p iteratively*)
Do[
 varsr = Table[Unique[r$], {nL}];
  varst = Table[Unique[t$], {nL}];

 (*discrete ODE1*)
 eqns1 =
   (-Cos[αα/2] - 2 (-3 Cos[αα] + Cos[2 αα]) varsr + Sin[αα/2] + 2 (-3 Cos[αα] + Cos[2 αα]) varst +
       12 Sin[αα] dfdα.varsr - 3 Sin[2 αα] dfdα.varsr - 6 Sin[αα] dfdα.varst + 3 Sin[2 αα]
         dfdα.varst + d2fdα2.varsr - 6 Cos[αα] d2fdα2.varsr + Cos[2 αα] d2fdα2.varsr +
```

```mathematica
       d2fdα2.varst - Cos[2 αα] d2fdα2.varst) * (Boole[0 ≤ # < α1p] & /@ αα) +
   (-Cos[αα/2] - 2 (-3 Cos[αα] + Cos[2 αα]) varsr + Sin[αα/2] + 2 (-2 Cos[αα] + Cos[2 αα]) varst +
       12 Sin[αα] dfdα.varsr - 3 Sin[2 αα] dfdα.varsr - 4 Sin[αα] dfdα.varst + 3 Sin[2 αα]
         dfdα.varst + d2fdα2.varsr - 6 Cos[αα] d2fdα2.varsr + Cos[2 αα] d2fdα2.varsr +
       d2fdα2.varst - Cos[2 αα] d2fdα2.varst) * (Boole[α1p ≤ # < α4p] & /@ αα) +
   (-Cos[αα/2] + 4 Cos[αα] varsr + Sin[αα/2] - 2 Cos[αα] varst + 8 Sin[αα] dfdα.varsr -
       2 Sin[αα] dfdα.varst - 4 Cos[αα] d2fdα2.varsr) * (Boole[α4p ≤ # < α2p] & /@ αα) +
   (-Cos[αα/2] + 2 (3 Cos[αα] + Cos[2 αα]) varsr + Sin[αα/2] - 2 (2 Cos[αα] + Cos[2 αα]) varst +
       12 Sin[αα] dfdα.varsr + 3 Sin[2 αα] dfdα.varsr - 4 Sin[αα] dfdα.varst - 3 Sin[2 αα]
         dfdα.varst - d2fdα2.varsr - 6 Cos[αα] d2fdα2.varsr - Cos[2 αα] d2fdα2.varsr -
       d2fdα2.varst + Cos[2 αα] d2fdα2.varst) * (Boole[α2p ≤ # < α3p] & /@ αα) +
   (-Cos[αα/2] + 2 (3 Cos[αα] + Cos[2 αα]) varsr + Sin[αα/2] - 2 (3 Cos[αα] + Cos[2 αα]) varst +
       12 Sin[αα] dfdα.varsr + 3 Sin[2 αα] dfdα.varsr - 6 Sin[αα] dfdα.varst - 3 Sin[2 αα]
         dfdα.varst - d2fdα2.varsr - 6 Cos[αα] d2fdα2.varsr - Cos[2 αα] d2fdα2.varsr -
       d2fdα2.varst + Cos[2 αα] d2fdα2.varst) * (Boole[α3p ≤ # ≤ π] & /@ αα);

(*discrete ODE2*)
eqns2 =
   (Cos[αα/2] + Sin[αα/2] + 2 varsr (3 Sin[αα] + Sin[2 αα]) - 2 (3 Sin[αα] + Sin[2 αα]) varst +
       dfdα.varsr - 6 Cos[αα] dfdα.varsr - 3 Cos[2 αα] dfdα.varsr + dfdα.varst +
       12 Cos[αα] dfdα.varst + 3 Cos[2 αα] dfdα.varst - Sin[2 αα] d2fdα2.varsr +
       6 Sin[αα] d2fdα2.varst + Sin[2 αα] d2fdα2.varst) * (Boole[0 ≤ # < α1p] & /@ αα) +
   (Cos[αα/2] + Sin[αα/2] + 4 (1 + Cos[αα]) varsr Sin[αα] - 2 (3 Sin[αα] + Sin[2 αα]) varst +
       dfdα.varsr - 4 Cos[αα] dfdα.varsr - 3 Cos[2 αα] dfdα.varsr + dfdα.varst +
       12 Cos[αα] dfdα.varst + 3 Cos[2 αα] dfdα.varst - Sin[2 αα] d2fdα2.varsr +
       6 Sin[αα] d2fdα2.varst + Sin[2 αα] d2fdα2.varst) * (Boole[α1p ≤ # < α4p] & /@ αα) +
   (Cos[αα/2] + Sin[αα/2] + 2 varsr Sin[αα] - 4 Sin[αα] varst - 2 Cos[αα] dfdα.varsr +
       8 Cos[αα] dfdα.varst + 4 Sin[αα] d2fdα2.varst) * (Boole[α4p ≤ # < α2p] & /@ αα) +
   (Cos[αα/2] + Sin[αα/2] + varsr (4 Sin[αα] - 2 Sin[2 αα]) + 2 (-3 Sin[αα] + Sin[2 αα]) varst -
       dfdα.varsr - 4 Cos[αα] dfdα.varsr + 3 Cos[2 αα] dfdα.varsr - dfdα.varst +
       12 Cos[αα] dfdα.varst - 3 Cos[2 αα] dfdα.varst + Sin[2 αα] d2fdα2.varsr +
       6 Sin[αα] d2fdα2.varst - Sin[2 αα] d2fdα2.varst) * (Boole[α2p ≤ # < α3p] & /@ αα) +
```

```mathematica
      (Cos[αα/2] + Sin[αα/2] - 2 varsr (-3 Sin[αα] + Sin[2 αα]) + 2 (-3 Sin[αα] + Sin[2 αα]) varst -
          dfdα.varsr - 6 Cos[αα] dfdα.varsr + 3 Cos[2 αα] dfdα.varsr - dfdα.varst +
          12 Cos[αα] dfdα.varst - 3 Cos[2 αα] dfdα.varst + Sin[2 αα] d2fdα2.varsr +
          6 Sin[αα] d2fdα2.varst - Sin[2 αα] d2fdα2.varst) * (Boole[α3p ≤ # ≤ π] & /@ αα);
(*boundary conditions*)
eqns1〚1〛 = (-1 / 2 + 2 varsr〚1〛 - 2 varst〚1〛 - 2 (d2fdα2.varsr)〚1〛) - 0;
eqns1〚nL〛 = (-1 / 2 + 2 varsr〚nL〛 - 2 varst〚nL〛 - 2 (d2fdα2.varsr)〚nL〛) - 0;
eqns2〚1〛 = varst〚1〛 - t0;
eqns2〚nL〛 = varst〚nL〛 - tπ;
   govEq = Join[eqns1, eqns2];
   allVars = Join[varsr, varst];
(*solve discretized ODEs*)
   root = FindRoot[govEq, Thread[{allVars, Flatten[{Table[ri, nL], Table[t0, nL]}]}],
       MaxIterations → 20, AccuracyGoal → Infinity, PrecisionGoal → 8] // Quiet;
   sol = allVars /. root;
{rsoltemp, tsoltemp} = Partition[sol, nL];
   rsol = Interpolation[Transpose[{αα, rsoltemp}]];
   tsol = Interpolation[Transpose[{αα, tsoltemp}]];
   (*obtain curves*)
AA = {rsol[α] * Cos[α], tsol[α] * Sin[α]};
   AAp = {rsol[α] * Cos[α] + √2 * Cos[π / 4 + α / 2], tsol[α] * Sin[α] + √2 * Sin[π / 4 + α / 2]};
BB = {rsol[α] * Cos[α] - tsol[α] * 2 Cos[α / 2]^2, 0};
   BBp = {rsol[α] * Cos[α] - tsol[α] * 2 Cos[α / 2]^2 - 1 / Sin[α / 2], 0};
CC = {rsol[α] * Cos[α] + tsol[α] * Sin[α] * Tan[α / 2], 0};
   CCp = {rsol[α] * Cos[α] + tsol[α] * Sin[α] * Tan[α / 2] + 1 / Cos[α / 2], 0};
EAABBx = 1 / 2 * (-1 + 2 Cos[α] + Cos[2 α]) rsol[α] -
     2 Cos[α / 2]^2 Cos[α] tsol[α] + Sin[α] (Cos[α] rsol'[α] - (1 + Cos[α]) tsol'[α]);
EAABBy = -2 Sin[α / 2]^2
     (rsol[α] Sin[α] - Sin[α] tsol[α] - Cos[α] rsol'[α] + tsol'[α] + Cos[α] tsol'[α]);
EAACCx = rsol[α] (Cos[α] + Sin[α]^2) - 2 Cos[α] Sin[α / 2]^2 tsol[α] +
     Sin[α] (-Cos[α] rsol'[α] + (-1 + Cos[α]) tsol'[α]);
EAACCy = 2 Cos[α / 2]^2
     (-rsol[α] Sin[α] + Sin[α] tsol[α] + Cos[α] rsol'[α] + tsol'[α] - Cos[α] tsol'[α]);
EAApBBpx =
   1 / 2 * ((-1 + 2 Cos[α] + Cos[2 α]) rsol[α] - 2 Sin[α / 2] - 4 Cos[α / 2]^2 Cos[α] tsol[α] +
      Sin[2 α] rsol'[α] - 2 Sin[α] tsol'[α] - Sin[2 α] tsol'[α]);
EAApBBpy =
   1 / 2 * (2 Cos[α / 2] + rsol[α] (-2 Sin[α] + Sin[2 α]) - 2 (-1 + Cos[α]) Sin[α] tsol[α] -
      rsol'[α] + 2 Cos[α] rsol'[α] - Cos[2 α] rsol'[α] - tsol'[α] + Cos[2 α] tsol'[α]);
EAApCCpx = 1 / 2 * (2 Cos[α / 2] + 2 rsol[α] (Cos[α] + Sin[α]^2) + (1 - 2 Cos[α] + Cos[2 α]) tsol[α] -
      Sin[2 α] rsol'[α] - 2 Sin[α] tsol'[α] + Sin[2 α] tsol'[α]);
EAApCCpy = 1 / 2 * (2 Sin[α / 2] - 2 (1 + Cos[α]) rsol[α] Sin[α] + (2 Sin[α] + Sin[2 α]) tsol[α] +
      rsol'[α] + 2 Cos[α] rsol'[α] + Cos[2 α] rsol'[α] + tsol'[α] - Cos[2 α] tsol'[α]);
(*find α1p α2p*)
```

```
αɑres1 = FindRoot[{rsol[ɑɑ1] * Cos[ɑɑ1] == (EAACCx /. α → ɑɑ2),
        tsol[ɑɑ1] * Sin[ɑɑ1] == (EAACCy /. α → ɑɑ2)}, {ɑɑ1, 0.1}, {ɑɑ2, 1.8}] // Quiet;
(*find α3p α4p*)
αɑres2 = FindRoot[{rsol[ɑɑ3] * Cos[ɑɑ3] == (EAABBx /. α → ɑɑ4),
        tsol[ɑɑ3] * Sin[ɑɑ3] == (EAABBy /. α → ɑɑ4)}, {ɑɑ3, 3.0}, {ɑɑ4, 1.0}] // Quiet;
   α1p = αɑres1[[1]][[2]];
α2p = αɑres1[[2]][[2]];
   α3p = αɑres2[[1]][[2]];
α4p = αɑres2[[2]][[2]];, {5}];
```

In[ ]:= (*calculate sofa area*)

In[ ]:= -rsol[0] + 2 tsol[0] - rsol[π] + 2 tsol[π] +
   NIntegrate[-EAApBBpy * D[EAApBBpx, α] - EAApCCpy * D[EAApCCpx, α],
    {α, 0, π}, WorkingPrecision → 100, MaxRecursion → 20, AccuracyGoal → 8] +
   NIntegrate[tsol[α] * Sin[α] * D[rsol[α] * Cos[α], α], {α, α1p, α3p}, WorkingPrecision → 100,
    MaxRecursion → 20, AccuracyGoal → 8] + NIntegrate[-EAABBy * D[EAABBx, α],
    {α, 0, α4p}, WorkingPrecision → 100, MaxRecursion → 20, AccuracyGoal → 8] +
   NIntegrate[-EAACCy * D[EAACCx, α], {α, α2p, π},
    WorkingPrecision → 100, MaxRecursion → 20, AccuracyGoal → 8]

Out[ ]=
2.21953

In[ ]:= **RealDigits[%]**

Out[ ]=
{{2, 2, 1, 9, 5, 2, 6, 9, 2, 0, 9, 4, 1, 7, 0, 9}, 1}

In[ ]:= (*obtained α1p α2p α3p α4p*)

In[ ]:= {α1p, α2p, α3p, α4p}

Out[ ]=
{0.0787365, 1.77669, 3.06286, 1.3649}

In[ ]:= (*display sofa shape*)

```
In[•]:= Show[ParametricPlot[{EAApCCpx, EAApCCpy}, {α, 0, π}, PerformanceGoal → "Speed",
    PlotStyle → {Red}], ParametricPlot[{EAApBBpx, EAApBBpy}, {α, 0, π},
    PerformanceGoal → "Speed", PlotStyle → {Red}], ParametricPlot[
    {EAACCx, EAACCy}, {α, α2p, π}, PerformanceGoal → "Speed", PlotStyle → {Green}],
  ParametricPlot[{EAABBx, EAABBy}, {α, 0, α4p}, PerformanceGoal → "Speed",
    PlotStyle → {Green}], ParametricPlot[{x, 1}, {x, EAApCCpx /. α → π, 0},
    PerformanceGoal → "Speed", PlotStyle → {Black}], ParametricPlot[{x, 1},
    {x, 0, EAApBBpx /. α → 0}, PerformanceGoal → "Speed", PlotStyle → {Black}],
  ParametricPlot[AA, {α, α1p, α3p}, PerformanceGoal → "Speed", PlotStyle → {Purple}],
  ParametricPlot[{x, 0}, {x, EAABBx /. α → 0, EAApBBpx /. α → π},
    PerformanceGoal → "Speed", PlotStyle → {Black}], ParametricPlot[{x, 0},
    {x, EAACCx /. α → π, EAApCCpx /. α → 0}, PerformanceGoal → "Speed", PlotStyle → {Black}],
  PlotRange → {{-2.0, 2.0}, {0, 1}}, AxesOrigin → {0, 0}, ImageSize → 500]
```

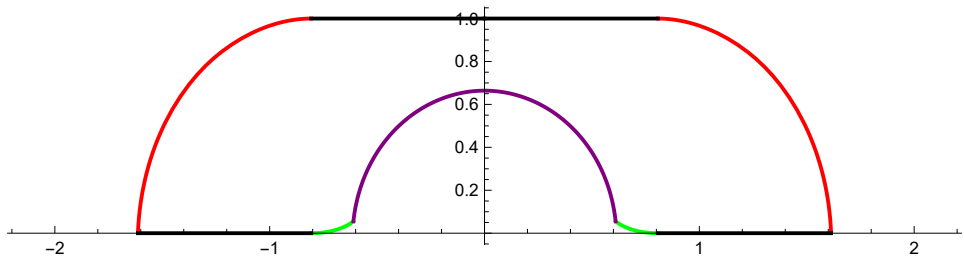

```
In[•]:= (*obtained parametric equations r(α),t(α)*)

In[•]:= Plot[{rsol[α], tsol[α]}, {α, 0, π}, PerformanceGoal → "Speed", PlotRange → {0.5, 0.8},
    PerformanceGoal → "Speed", PlotLegends → {"r(α)", "t(α)"}, PlotTheme → "BoldColor"]
```

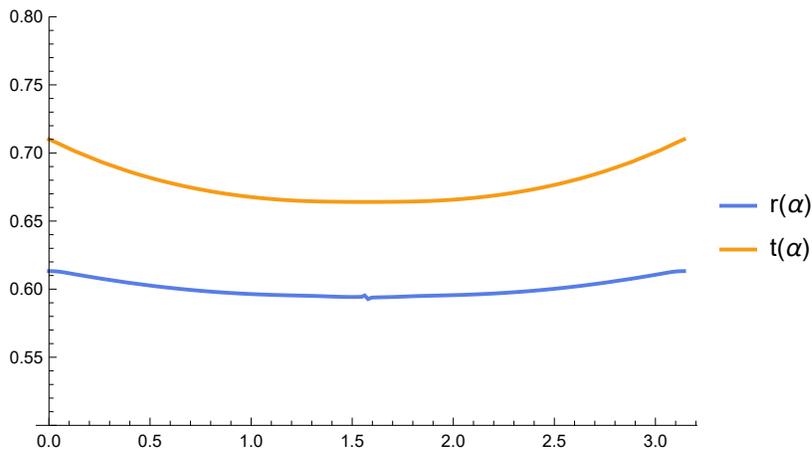

```
In[•]:= (*final obtained r[α] was fluctuate a little at π/2*)

In[•]:= (*more careful search could yield more accurate results*)

In[•]:= NotebookSave[];
```

# Appendix 7

```mathematica
(*Verify Romik's results for ambidextrous car from previous paper*)

AA = {r[α] * Cos[α], t[α] * Sin[α]};
AAp = {r[α] * Cos[α] + √2 * Cos[π/4 + α/2], t[α] * Sin[α] + √2 * Sin[π/4 + α/2]};
BB = {r[α] * Cos[α] - t[α] * 2 Cos[α/2]^2, 0};
BBp = {r[α] * Cos[α] - t[α] * 2 Cos[α/2]^2 - 1/Sin[α/2], 0};
CC = {r[α] * Cos[α] + t[α] * Sin[α] * Tan[α/2], 0};
CCp = {r[α] * Cos[α] + t[α] * Sin[α] * Tan[α/2] + 1/Cos[α/2], 0};
EAABBx = 1/2 * (-1 + 2 Cos[α] + Cos[2 α]) r[α] -
    2 Cos[α/2]^2 Cos[α] t[α] + Sin[α] (Cos[α] r'[α] - (1 + Cos[α]) t'[α]);
EAABBy = -2 Sin[α/2]^2 (r[α] Sin[α] - Sin[α] t[α] - Cos[α] r'[α] + t'[α] + Cos[α] t'[α]);
EAACCx = r[α] (Cos[α] + Sin[α]^2) -
    2 Cos[α] Sin[α/2]^2 t[α] + Sin[α] (-Cos[α] r'[α] + (-1 + Cos[α]) t'[α]);
EAACCy = 2 Cos[α/2]^2 (-r[α] Sin[α] + Sin[α] t[α] + Cos[α] r'[α] + t'[α] - Cos[α] t'[α]);
EAApBBpx = 1/2 * ((-1 + 2 Cos[α] + Cos[2 α]) r[α] - 2 Sin[α/2] -
        4 Cos[α/2]^2 Cos[α] t[α] + Sin[2 α] r'[α] - 2 Sin[α] t'[α] - Sin[2 α] t'[α]);
EAApBBpy = 1/2 * (2 Cos[α/2] + r[α] (-2 Sin[α] + Sin[2 α]) - 2 (-1 + Cos[α]) Sin[α] t[α] -
        r'[α] + 2 Cos[α] r'[α] - Cos[2 α] r'[α] - t'[α] + Cos[2 α] t'[α]);
EAApCCpx = 1/2 * (2 Cos[α/2] + 2 r[α] (Cos[α] + Sin[α]^2) +
        (1 - 2 Cos[α] + Cos[2 α]) t[α] - Sin[2 α] r'[α] - 2 Sin[α] t'[α] + Sin[2 α] t'[α]);
EAApCCpy = 1/2 * (2 Sin[α/2] - 2 (1 + Cos[α]) r[α] Sin[α] + (2 Sin[α] + Sin[2 α]) t[α] +
        r'[α] + 2 Cos[α] r'[α] + Cos[2 α] r'[α] + t'[α] - Cos[2 α] t'[α]);

(*
approximate values for piecewise functions from previous paper
 0.2896538208173209417435216117362
*)

α1p = 5792/10000;

α2p = 5792/10000;

α5p = 5792/10000;

(*functional*)
```

```mathematica
In[ ]:= FF = Piecewise[{{0, α < 0},

        {-r[0] + 2 t[0] - 2 (EAApBBpy - 1 / 2) *
           D[EAApBBpx, α]                                                         -
          2 EAABBy * D[EAABBx, α], 0 ≤ α ≤ α5p},

        {-r[0] + 2 t[0] - 2 (EAApBBpy - 1 / 2) * D[EAApBBpx, α] -
          2 (EAApCCpy - 1 / 2) * D[EAApCCpx, α]                    -
          2 EAABBy * D[EAABBx, α]                        , α5p < α ≤ α1p},

        {-r[0] + 2 t[0] - 2 (EAApBBpy - 1 / 2) * D[EAApBBpx, α] -
          2 (EAApCCpy - 1 / 2) * D[EAApCCpx, α] + 2 AA[[2]] * D[AA[[1]], α] -
          2 EAABBy * D[EAABBx, α]                        , α1p < α ≤ α2p},

        {-r[0] + 2 t[0] - 2 (EAApBBpy - 1 / 2) * D[EAApBBpx, α] -
          2 (EAApCCpy - 1 / 2) * D[EAApCCpx, α] + 2 AA[[2]] * D[AA[[1]], α] -
          2 EAABBy * D[EAABBx, α] - 2 EAACCy * D[EAACCx, α], α2p < α ≤ 1},

        {            - 2 (EAApBBpy - 1 / 2) * D[EAApBBpx, α] -
          2 (EAApCCpy - 1 / 2) * D[EAApCCpx, α] + 2 AA[[2]] * D[AA[[1]], α] -
          2 EAABBy * D[EAABBx, α] - 2 EAACCy * D[EAACCx, α], 1 < α ≤ π / 2},

        {0, π / 2 < α}}] // Simplify;

In[ ]:=


In[ ]:= (*EL equations*)

In[ ]:= (D[FF, r[α]] - D[D[FF, r'[α]], α] + D[D[FF, r''[α]], {α, 2}]) // Simplify
Out[ ]=
   ⎡ -Cos[α] (-Cos[α/2] + 8 Cos[α] r[α] + Sin[α/2] -                              362/625 < α < π/2
   ⎢   6 Cos[α] t[α] + 16 Sin[α] r'[α] - 6 Sin[α] t'[α] - 8 Cos[α] r''[α])
   ⎢ 1/2 (Sin[α/2] - Sin[3α/2])                                                    0 < α < 362/625
   ⎢   (1 + 8 r[α] Sin[3α/2] - 8 Sin[3α/2] t[α] + 4 Cos[α/2] r'[α] - 12 Cos[3α/2] r'[α] +
   ⎢    4 Cos[α/2] t'[α] + 12 Cos[3α/2] t'[α] + 4 Sin[α/2] r''[α] -
   ⎢    4 Sin[3α/2] r''[α] + 4 Sin[α/2] t''[α] + 4 Sin[3α/2] t''[α])
   ⎢ 0                                                                              2α > π || α < 0
   ⎣ Indeterminate                                                                  True
```

```
In[•]:= (D[FF, t[α]] - D[D[FF, t'[α]], α] + D[D[FF, t''[α]], {α, 2}]) // Simplify
```

Out[•]=

$$\begin{cases} \sin[\alpha]\left(\cos\left[\frac{\alpha}{2}\right] + \sin\left[\frac{\alpha}{2}\right] + 6\,r[\alpha]\,\sin[\alpha] - \right. \\ \quad 8\,\sin[\alpha]\,t[\alpha] - 6\,\cos[\alpha]\,r'[\alpha] + 16\,\cos[\alpha]\,t'[\alpha] + 8\,\sin[\alpha]\,t''[\alpha]\big) & \frac{362}{625} < \alpha < \frac{\pi}{2} \\ \sin[\alpha]\left(\cos\left[\frac{\alpha}{2}\right] + 4\,r[\alpha]\,(\sin[\alpha] + \sin[2\,\alpha]) - 4\,(\sin[\alpha] + \sin[2\,\alpha])\,t[\alpha] + \right. \\ \quad 2\,r'[\alpha] - 4\,\cos[\alpha]\,r'[\alpha] - 6\,\cos[2\,\alpha]\,r'[\alpha] + 2\,t'[\alpha] + 8\,\cos[\alpha]\,t'[\alpha] + \\ \quad 6\,\cos[2\,\alpha]\,t'[\alpha] - 2\,\sin[2\,\alpha]\,r''[\alpha] + 4\,\sin[\alpha]\,t''[\alpha] + 2\,\sin[2\,\alpha]\,t''[\alpha]\big) & 0 < \alpha < \frac{362}{625} \\ 0 & 2\,\alpha > \pi \;||\; \alpha < 0 \\ \text{Indeterminate} & \text{True} \end{cases}$$

```
In[•]:=
```

```
In[•]:= (*Romik's results-ambidextrous car*)

In[•]:= (*explicit expression, code from movingsofas-v1.3
        https://www.math.ucdavis.edu/~romik/data/uploads/software/movingsofas-v1.3.nb*)

In[•]:= (*define functions and parameters*)

In[•]:= x₁[t_] :=
        RotationMatrix[t].{a₁ Cos[t] + a₂ Sin[t] - 1, -a₂ Cos[t] + a₁ Sin[t] - 1/2} + {κ₁,₁, κ₁,₂};

       x₂[t_] := RotationMatrix[t].{-1/4 t² + b₁ t + b₂, 1/2 t - b₁ - 1} + {κ₂,₁, κ₂,₂};

       x₃[t_] := RotationMatrix[t].{c₁ - t, c₂ + t} + {κ₃,₁, κ₃,₂};

       x₄[t_] := RotationMatrix[t].{-1/2 t + d₁ - 1, -1/4 t² + d₁ t + d₂} + {κ₄,₁, κ₄,₂};

       x₅[t_] :=
         RotationMatrix[t].{e₁ Cos[t] + e₂ Sin[t] - 1/2, -e₂ Cos[t] + e₁ Sin[t] - 1} + {κ₅,₁, κ₅,₂};

       x₆[t_] := RotationMatrix[t].
            {f₁ Cos[t/2] + f₂ Sin[t/2] - 1, -f₂ Cos[t/2] + f₁ Sin[t/2] - 1} + {κ₆,₁, κ₆,₂};

In[•]:= μ[t_] := {Cos[t], Sin[t]};
       ν[t_] := {-Sin[t], Cos[t]};
```

```mathematica
In[*]:= AmbidextrousAngleβ = 0.289653820817320941743521611736;
       AmbidextrousPrecomputedNumericalConstants = {
         β → 0.289653820817320941743521611736,
         κ_{1,1} → 0.124712637587267758739932415305,
         κ_{1,2} → 1/2,
         κ_{6,1} → -0.167049816550309655013423446260,
         κ_{6,2} → 1/2,
         κ_{5,1} → -0.458812270687887068766779307825,
         κ_{5,2} → 1/2,
         a_1 → 0.875287362412732241260067584695,
         a_2 → 0,
         f_1 → 1.202938908156911389070222800034,
         f_2 → -0.498273610464875672029397859080,
         e_1 → 0.875287362412732241260067584695,
         e_2 → 0
       };

In[*]:= A_j_[t_] := Collect[Simplify[x_j[t] + μ[t] + (x_j'[t].μ[t]) ν[t]], {Cos[t], Sin[t]}];
       B_j_[t_] := Collect[Simplify[x_j[t] + (x_j'[t].μ[t]) ν[t]], {Cos[t], Sin[t]}];
       C_j_[t_] := Collect[Simplify[x_j[t] + ν[t] - (x_j'[t].ν[t]) μ[t]], {Cos[t], Sin[t]}];
       D_j_[t_] := Collect[Simplify[x_j[t] - (x_j'[t].ν[t]) μ[t]], {Cos[t], Sin[t]}];

In[*]:= Ambix[t_] := ⎧ x_1[t]   0 ≤ t < β
                    ⎨ x_6[t]   β ≤ t < π/2 - β
                    ⎩ x_5[t]   π/2 - β ≤ t ≤ π/2

       AmbiA[t_] := ⎧ A_1[t]   0 ≤ t < β
                    ⎨ A_6[t]   β ≤ t < π/2 - β
                    ⎩ A_5[t]   π/2 - β ≤ t ≤ π/2

       AmbiB[t_] := ⎧ B_1[t]   0 ≤ t < β
                    ⎨ B_6[t]   β ≤ t < π/2 - β
                    ⎩ B_5[t]   π/2 - β ≤ t ≤ π/2

       AmbiC[t_] := ⎧ C_1[t]   0 ≤ t < β
                    ⎨ C_6[t]   β ≤ t ≤ π/2 - β
                    ⎩ C_5[t]   π/2 - β ≤ t ≤ π/2

       AmbiD[t_] := ⎧ D_1[t]   0 ≤ t < β
                    ⎨ D_6[t]   β ≤ t < π/2 - β
                    ⎩ D_5[t]   π/2 - β ≤ t ≤ π/2

In[*]:= (*show curves*)
```

```mathematica
In[ ]:= Show[ListLinePlot[Table[Ambix[t] /. AmbidextrousPrecomputedNumericalConstants,
      {t, AmbidextrousAngleβ, π/2 - AmbidextrousAngleβ, 0.01}],
     PerformanceGoal → "Speed", AspectRatio → Automatic, PlotStyle → Purple],
   ListLinePlot[
    Table[AmbiA[t] /. AmbidextrousPrecomputedNumericalConstants, {t, 0, π/2, 0.01}],
    PerformanceGoal → "Speed", AspectRatio → Automatic, PlotStyle → Red],
   ListLinePlot[Table[AmbiB[t] /. AmbidextrousPrecomputedNumericalConstants,
     {t, AmbidextrousAngleβ, π/2, 0.01}], PerformanceGoal → "Speed",
    AspectRatio → Automatic, PlotStyle → Green], ListLinePlot[
    Table[AmbiC[t] /. AmbidextrousPrecomputedNumericalConstants, {t, 0, π/2, 0.01}],
    PerformanceGoal → "Speed", AspectRatio → Automatic, PlotStyle → Red],
   ListLinePlot[Table[AmbiD[t] /. AmbidextrousPrecomputedNumericalConstants,
     {t, 0, π/2 - AmbidextrousAngleβ, 0.01}], PerformanceGoal → "Speed",
    AspectRatio → Automatic, PlotStyle → Green], ListLinePlot[
    (# * {1, -1} + {0, 1}) & /@ Table[Ambix[t] /. AmbidextrousPrecomputedNumericalConstants,
      {t, AmbidextrousAngleβ, π/2 - AmbidextrousAngleβ, 0.01}], PerformanceGoal → "Speed",
    AspectRatio → Automatic, PlotStyle → Purple, PlotStyle → Purple],
   ListLinePlot[(# * {1, -1} + {0, 1}) & /@
     Table[AmbiA[t] /. AmbidextrousPrecomputedNumericalConstants, {t, 0, π/2, 0.01}],
    PerformanceGoal → "Speed", AspectRatio → Automatic, PlotStyle → Red], ListLinePlot[
    (# * {1, -1} + {0, 1}) & /@ Table[AmbiB[t] /. AmbidextrousPrecomputedNumericalConstants,
     {t, AmbidextrousAngleβ, π/2, 0.01}], PerformanceGoal → "Speed",
    AspectRatio → Automatic, PlotStyle → Green], ListLinePlot[(# * {1, -1} + {0, 1}) & /@
     Table[AmbiC[t] /. AmbidextrousPrecomputedNumericalConstants, {t, 0, π/2, 0.01}],
    PerformanceGoal → "Speed", AspectRatio → Automatic, PlotStyle → Red], ListLinePlot[
    (# * {1, -1} + {0, 1}) & /@ Table[AmbiD[t] /. AmbidextrousPrecomputedNumericalConstants,
     {t, 0, π/2 - AmbidextrousAngleβ, 0.01}], PerformanceGoal → "Speed",
    AspectRatio → Automatic, PlotStyle → Green], PlotRange → All,
   PerformanceGoal → "Speed", AspectRatio → Automatic, AxesOrigin → {0, 0}]
```

Out[ ]= 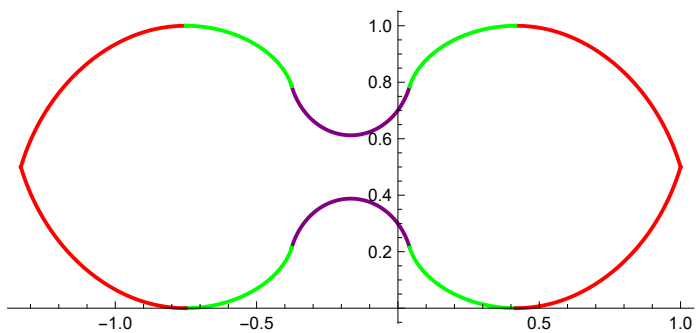

In[ ]:=

In[ ]:= (*Verify r[α] and t[α] satisfy E-L equations*)

```mathematica
In[•]:= {1/2 (Sin[α/2] - Sin[3α/2]) (1 + 8 r[α] Sin[3α/2] - 8 Sin[3α/2] t[α] +
          4 Cos[α/2] r'[α] - 12 Cos[3α/2] r'[α] + 4 Cos[α/2] t'[α] + 12 Cos[3α/2] t'[α] +
          4 Sin[α/2] r''[α] - 4 Sin[3α/2] r''[α] + 4 Sin[α/2] t''[α] + 4 Sin[3α/2] t''[α]),
    Sin[α] (Cos[α/2] + 4 r[α] (Sin[α] + Sin[2 α]) - 4 (Sin[α] + Sin[2 α]) t[α] +
          2 r'[α] - 4 Cos[α] r'[α] - 6 Cos[2 α] r'[α] + 2 t'[α] + 8 Cos[α] t'[α] +
          6 Cos[2 α] t'[α] - 2 Sin[2 α] r''[α] + 4 Sin[α] t''[α] + 2 Sin[2 α] t''[α])} /.
    {r[α] → ((x_1[α/2] - {(1 - e_1 + κ_{5,1})/2, 0})[[1]])/Cos[α],
     t[α] → (x_1[α/2] - {(1 - e_1 + κ_{5,1})/2, 0})[[2]]/Sin[α],
     r'[α] → D[((x_1[α/2] - {(1 - e_1 + κ_{5,1})/2, 0})[[1]])/Cos[α], α],
     t'[α] → D[(x_1[α/2] - {(1 - e_1 + κ_{5,1})/2, 0})[[2]]/Sin[α], α],
     r''[α] → D[((x_1[α/2] - {(1 - e_1 + κ_{5,1})/2, 0})[[1]])/Cos[α], {α, 2}],
     t''[α] → D[(x_1[α/2] - {(1 - e_1 + κ_{5,1})/2, 0})[[2]]/Sin[α], {α, 2}]} // Simplify

Out[•]= {0, 0}

In[•]:= {-Cos[α] (-Cos[α/2] + 8 Cos[α] r[α] + Sin[α/2] -
          6 Cos[α] t[α] + 16 Sin[α] r'[α] - 6 Sin[α] t'[α] - 8 Cos[α] r''[α]),
    Sin[α] (Cos[α/2] + Sin[α/2] + 6 r[α] Sin[α] - 8 Sin[α] t[α] - 6 Cos[α] r'[α] +
          16 Cos[α] t'[α] + 8 Sin[α] t''[α])} /.
    {r[α] → ((x_6[α/2] - {(1 - e_1 + κ_{5,1})/2, 0})[[1]])/Cos[α],
     t[α] → (x_6[α/2] - {(1 - e_1 + κ_{5,1})/2, 0})[[2]]/Sin[α],
     r'[α] → D[((x_6[α/2] - {(1 - e_1 + κ_{5,1})/2, 0})[[1]])/Cos[α], α],
     t'[α] → D[(x_6[α/2] - {(1 - e_1 + κ_{5,1})/2, 0})[[2]]/Sin[α], α],
     r''[α] → D[((x_6[α/2] - {(1 - e_1 + κ_{5,1})/2, 0})[[1]])/Cos[α], {α, 2}],
     t''[α] → D[(x_6[α/2] - {(1 - e_1 + κ_{5,1})/2, 0})[[2]]/Sin[α], {α, 2}]} // Simplify

Out[•]= {0, 0}
```

**Appendix 8**

(*Hessian Matrix for symmetric conditions*)

(*symmetric condition a*)

AA = {r[α] * Cos[α], t[α] * Sin[α]};

AAp = {r[α] * Cos[α] + √2 * Cos[π/4 + α/2], t[α] * Sin[α] + √2 * Sin[π/4 + α/2]};

BB = {r'[α] * Cos[α] - t[α] * Sin[α] * 2 Cos[α/2]^2, 0};

BBp = {r'[α] * Cos[α] - t[α] * 2 Cos[α/2]^2 - 1/Sin[α/2], 0};

CC = {r'[α] * Cos[α] + t[α] * Sin[α] * Tan[α/2], 0};

CCp = {r'[α] * Cos[α] + t[α] * Sin[α] * Tan[α/2] + 1/Cos[α/2], 0};

EAABBx = 1/2 * (-1 + 2 Cos[α] + Cos[2 α]) r[α] - 2 Cos[α/2]^2 Cos[α] t[α] + Sin[α] (Cos[α] r'[α] - (1 + Cos[α]) t'[α]);

EAABBy = -2 Sin[α/2]^2 (r[α] Sin[α] - Sin[α] t[α] - Cos[α] r'[α] + t'[α] + Cos[α] t'[α]);

EAACCx = r[α] (Cos[α] + Sin[α]^2) - 2 Cos[α] Sin[α/2]^2 t[α] + Sin[α] (-Cos[α] r'[α] + (-1 + Cos[α]) t'[α]);

EAACCy = 2 Cos[α/2]^2 (-r[α] Sin[α] + Sin[α] t[α] + Cos[α] r'[α] + t'[α] - Cos[α] t'[α]);

EAApBBpx =
1/2 * ((-1 + 2 Cos[α] + Cos[2 α]) r[α] - 2 Sin[α/2] - 4 Cos[α/2]^2 Cos[α] t[α] + Sin[2 α] r'[α] - 2 Sin[α] t'[α] - Sin[2 α] t'[α]);

EAApBBpy = 1/2 * (2 Cos[α/2] + r[α] (-2 Sin[α] + Sin[2 α]) -
2 (-1 + Cos[α]) Sin[α] t[α] - r'[α] + 2 Cos[α] r'[α] - Cos[2 α] r'[α] - t'[α] + Cos[2 α] t'[α]);

EAApCCpx =
1/2 * (2 Cos[α/2] + 2 r[α] (Cos[α] + Sin[α]^2) + (1 - 2 Cos[α] + Cos[2 α]) t[α] - Sin[2 α] r'[α] - 2 Sin[α] t'[α] + Sin[2 α] t'[α]);

EAApCCpy = 1/2 * (2 Sin[α/2] - 2 (1 + Cos[α]) r[α] Sin[α] +
(2 Sin[α] + Sin[2 α]) t[α] + r'[α] + 2 Cos[α] r'[α] + Cos[2 α] r'[α] + t'[α] - Cos[2 α] t'[α]);

FF = -EAApBBpy * D[EAApBBpx, α] - EAApCCpy * D[EAApCCpx, α] + AA[[2]] * D[AA[[1]], α];

HH = ResourceFunction["HessianMatrix"][FF, {r[α], t[α], r'[α], t'[α], r''[α], t''[α]}] // Simplify // MatrixForm

_trixForm=_

| 4 Cos[2α] Sin[α]² | -((1 + 4 Cos[2α]) Sin[α]²) | ¼ (6 Sin[2α] - 5 Sin[4α]) | -10 Cos[α] Sin[α]³ | -2 Cos[α]² Sin[α]² |
| -((1 + 4 Cos[2α]) Sin[α]²) | 4 Cos[2α] Sin[α]² | -Sin[2α] + 5/4 Sin[4α] | 10 Cos[α] Sin[α]³ | 2 Cos[α]² Sin[α]² |
| ¼ (6 Sin[2α] - 5 Sin[4α]) | -Sin[2α] + 5/4 Sin[4α] | -3 Sin[2α]² | 2 (1 + 3 Cos[2α]) Sin[α]² | 2 Cos[α]³ Sin[α] | 2 |
| -10 Cos[α] Sin[α]³ | 10 Cos[α] Sin[α]³ | 2 (1 + 3 Cos[2α]) Sin[α]² | -2 (1 + 3 Cos[2α]) Sin[α]² | 2 Cos[α] Sin[α]³ | -; |
| -2 Cos[α]² Sin[α]² | 2 Cos[α]² Sin[α]² | 2 Cos[α]³ Sin[α] | 2 Cos[α] Sin[α]³ | 0 |
| -2 Sin[α]⁴ | 2 Sin[α]⁴ | | | 0 |

$$\begin{pmatrix} -2\sin[\alpha]^4 \\ 2\sin[\alpha]^4 \\ \cos[\alpha]\sin[\alpha]^3 \\ 2\cos[\alpha]\sin[\alpha]^3 \\ 0 \\ 0 \end{pmatrix}$$

(*Symmetric condition b*)

FF1 = -r[0] + 2 t[0] - EAApBBpy * D[EAApBBpx, α] - EAApCCpy * D[EAApCCpx, α]   - EAABBy * D[EAABBx, α];

HH1 = ResourceFunction["HessianMatrix"][FF1, {r[α], t[α], r'[α], t'[α], r''[α], t''[α]}] // Simplify // MatrixForm

$\begin{pmatrix} -2(\cos[\alpha]-3\cos[2\alpha])\sin[\alpha]^2 & 2(\cos[\alpha]-3\cos[2\alpha])\sin[\alpha]^2 & \frac{1}{8}(-10\sin[\alpha]+18\sin[2\alpha]+6 \\ 2(\cos[\alpha]-3\cos[2\alpha])\sin[\alpha]^2 & -2(\cos[\alpha]-3\cos[2\alpha])\sin[\alpha]^2 & \frac{1}{4}(2-3\cos[\alpha]-6\cos[2\alpha]+ \\ \frac{1}{8}(-10\sin[\alpha]+18\sin[2\alpha]+6\sin[3\alpha]-15\sin[4\alpha]) & \frac{1}{4}(2-3\cos[\alpha]-6\cos[2\alpha]+15\cos[3\alpha])\sin[\alpha] & 2(2-9\cos[\alpha])\cos[ \\ -15\cos[\alpha]\sin[\alpha]^3 & 15\cos[\alpha]\sin[\alpha]^3 & 3(1+3\cos[2\alpha]) \\ -\cos[\alpha](-1+3\cos[\alpha])\sin[\alpha]^2 & \cos[\alpha](-1+3\cos[\alpha])\sin[\alpha]^2 & \cos[\alpha]^2(-1+3\cos[ \\ -3\sin[\alpha]^4 & 3\sin[\alpha]^4 & 3\cos[\alpha]\sin \end{pmatrix}$

FF2 = -r[0] + 2 t[0] - EAApBBpy * D[EAApBBpx, α] - EAApCCpy * D[EAApCCpx, α] + AA[[2]] * D[AA[[1]], α] - EAABBy * D[EAABBx, α];

HH2 = ResourceFunction["HessianMatrix"][FF1, {r[α], t[α], r'[α], t'[α], r''[α], t''[α]}] // Simplify // MatrixForm

$\begin{pmatrix} -2(\cos[\alpha]-3\cos[2\alpha])\sin[\alpha]^2 & 2(\cos[\alpha]-3\cos[2\alpha])\sin[\alpha]^2 & \frac{1}{8}(-10\sin[\alpha]+18\sin[2\alpha]+6 \\ 2(\cos[\alpha]-3\cos[2\alpha])\sin[\alpha]^2 & -2(\cos[\alpha]-3\cos[2\alpha])\sin[\alpha]^2 & \frac{1}{4}(2-3\cos[\alpha]-6\cos[2\alpha]+ \\ \frac{1}{8}(-10\sin[\alpha]+18\sin[2\alpha]+6\sin[3\alpha]-15\sin[4\alpha]) & \frac{1}{4}(2-3\cos[\alpha]-6\cos[2\alpha]+15\cos[3\alpha])\sin[\alpha] & 2(2-9\cos[\alpha])\cos[ \\ -15\cos[\alpha]\sin[\alpha]^3 & 15\cos[\alpha]\sin[\alpha]^3 & 3(1+3\cos[2\alpha]) \\ -\cos[\alpha](-1+3\cos[\alpha])\sin[\alpha]^2 & \cos[\alpha](-1+3\cos[\alpha])\sin[\alpha]^2 & \cos[\alpha]^2(-1+3\cos[ \\ -3\sin[\alpha]^4 & 3\sin[\alpha]^4 & 3\cos[\alpha]\sin \end{pmatrix}$

FF3 = -EAApBBpy * D[EAApBBpx, α] - EAApCCpy * D[EAApCCpx, α] + AA[[2]] * D[AA[[1]], α] - EAABBy * D[EAABBx, α];

$$\begin{pmatrix} \cdots \operatorname{Sin}[3\alpha] - 15 \operatorname{Sin}[4\alpha]) & -15 \operatorname{Cos}[\alpha] \operatorname{Sin}[\alpha]^3 & -\operatorname{Cos}[\alpha] (-1 + 3 \operatorname{Cos}[\alpha]) \operatorname{Sin}[\alpha]^2 & -3 \operatorname{Sin}[\alpha]^4 \\ 15 \operatorname{Cos}[3\alpha]) \operatorname{Sin}[\alpha] & 15 \operatorname{Cos}[\alpha] \operatorname{Sin}[\alpha]^3 & \operatorname{Cos}[\alpha] (-1 + 3 \operatorname{Cos}[\alpha]) \operatorname{Sin}[\alpha]^2 & 3 \operatorname{Sin}[\alpha]^4 \\ [\alpha] \operatorname{Sin}[\alpha]^2 & 3 (1 + 3 \operatorname{Cos}[2\alpha]) \operatorname{Sin}[\alpha]^2 & \operatorname{Cos}[\alpha]^2 (-1 + 3 \operatorname{Cos}[\alpha]) \operatorname{Sin}[\alpha] & 3 \operatorname{Cos}[\alpha] \operatorname{Sin}[\alpha]^3 \\ \operatorname{Sin}[\alpha]^2 & -((3 + 4 \operatorname{Cos}[\alpha] + 9 \operatorname{Cos}[2\alpha]) \operatorname{Sin}[\alpha]^2) & 3 \operatorname{Cos}[\alpha] \operatorname{Sin}[\alpha]^3 & -((1 + 3 \operatorname{Cos}[\alpha]) \operatorname{Sin}[\alpha]^3) \\ \alpha]) \operatorname{Sin}[\alpha] & 3 \operatorname{Cos}[\alpha] \operatorname{Sin}[\alpha]^3 & 0 & 0 \\ [\alpha]^3 & -((1 + 3 \operatorname{Cos}[\alpha]) \operatorname{Sin}[\alpha]^3) & 0 & 0 \end{pmatrix}$$

$$\begin{pmatrix} \cdots \operatorname{Sin}[3\alpha] - 15 \operatorname{Sin}[4\alpha]) & -15 \operatorname{Cos}[\alpha] \operatorname{Sin}[\alpha]^3 & -\operatorname{Cos}[\alpha] (-1 + 3 \operatorname{Cos}[\alpha]) \operatorname{Sin}[\alpha]^2 & -3 \operatorname{Sin}[\alpha]^4 \\ 15 \operatorname{Cos}[3\alpha]) \operatorname{Sin}[\alpha] & 15 \operatorname{Cos}[\alpha] \operatorname{Sin}[\alpha]^3 & \operatorname{Cos}[\alpha] (-1 + 3 \operatorname{Cos}[\alpha]) \operatorname{Sin}[\alpha]^2 & 3 \operatorname{Sin}[\alpha]^4 \\ [\alpha] \operatorname{Sin}[\alpha]^2 & 3 (1 + 3 \operatorname{Cos}[2\alpha]) \operatorname{Sin}[\alpha]^2 & \operatorname{Cos}[\alpha]^2 (-1 + 3 \operatorname{Cos}[\alpha]) \operatorname{Sin}[\alpha] & 3 \operatorname{Cos}[\alpha] \operatorname{Sin}[\alpha]^3 \\ \operatorname{Sin}[\alpha]^2 & -((3 + 4 \operatorname{Cos}[\alpha] + 9 \operatorname{Cos}[2\alpha]) \operatorname{Sin}[\alpha]^2) & 3 \operatorname{Cos}[\alpha] \operatorname{Sin}[\alpha]^3 & -((1 + 3 \operatorname{Cos}[\alpha]) \operatorname{Sin}[\alpha]^3) \\ \alpha]) \operatorname{Sin}[\alpha] & 3 \operatorname{Cos}[\alpha] \operatorname{Sin}[\alpha]^3 & 0 & 0 \\ [\alpha]^3 & -((1 + 3 \operatorname{Cos}[\alpha]) \operatorname{Sin}[\alpha]^3) & 0 & 0 \end{pmatrix}$$

```mathematica
In[•]:= HH3 = ResourceFunction["HessianMatrix"][FF1, {r[α], t[α], r'[α], t'[α], r''[α], t''[α], r'''[α], t'''[α]}] // Simplify // MatrixForm
```

Out[•]//MatrixForm=

$$\begin{pmatrix} -2(\cos[\alpha] - 3\cos[2\alpha])\sin[\alpha]^2 & 2(\cos[\alpha] - 3\cos[2\alpha])\sin[\alpha]^2 & \frac{1}{8}(-10\sin[\alpha] + 18\sin[2\alpha] + 6\cdots \\ 2(\cos[\alpha] - 3\cos[2\alpha])\sin[\alpha]^2 & -2(\cos[\alpha] - 3\cos[2\alpha])\sin[\alpha]^2 & \frac{1}{4}(2 - 3\cos[\alpha] - 6\cos[2\alpha]) + \\ \frac{1}{8}(-10\sin[\alpha] + 18\sin[2\alpha] + 6\sin[3\alpha] - 15\sin[4\alpha]) & \frac{1}{4}(2 - 3\cos[\alpha] - 6\cos[2\alpha] + 15\cos[3\alpha])\sin[\alpha] & 2(2 - 9\cos[\alpha])\cos[\\ -15\cos[\alpha]\sin[\alpha]^3 & 15\cos[\alpha]\sin[\alpha]^3 & 3(1 + 3\cos[2\alpha]) \\ -\cos[\alpha](-1 + 3\cos[\alpha])\sin[\alpha]^2 & \cos[\alpha](-1 + 3\cos[\alpha])\sin[\alpha]^2 & \cos[\alpha]^2(-1 + 3\cos[\\ -3\sin[\alpha]^4 & 3\sin[\alpha]^4 & 3\cos[\alpha]\sin \end{pmatrix}$$

```mathematica
In[•]:= FF4 = -EAApBBpy*D[EAApBBpx, α] - EAApCCpy*D[EAApCCpx, α] + AA[[2]]*D[AA[[1]], α];

In[•]:= HH4 = ResourceFunction["HessianMatrix"][FF1, {r[α], t[α], r'[α], t'[α], r''[α], t''[α]}] // Simplify // MatrixForm
```

Out[•]//MatrixForm=

$$\begin{pmatrix} -2(\cos[\alpha] - 3\cos[2\alpha])\sin[\alpha]^2 & 2(\cos[\alpha] - 3\cos[2\alpha])\sin[\alpha]^2 & \frac{1}{8}(-10\sin[\alpha] + 18\sin[2\alpha] + 6\cdots \\ 2(\cos[\alpha] - 3\cos[2\alpha])\sin[\alpha]^2 & -2(\cos[\alpha] - 3\cos[2\alpha])\sin[\alpha]^2 & \frac{1}{4}(2 - 3\cos[\alpha] - 6\cos[2\alpha]) + \\ \frac{1}{8}(-10\sin[\alpha] + 18\sin[2\alpha] + 6\sin[3\alpha] - 15\sin[4\alpha]) & \frac{1}{4}(2 - 3\cos[\alpha] - 6\cos[2\alpha] + 15\cos[3\alpha])\sin[\alpha] & 2(2 - 9\cos[\alpha])\cos[\\ -15\cos[\alpha]\sin[\alpha]^3 & 15\cos[\alpha]\sin[\alpha]^3 & 3(1 + 3\cos[2\alpha]) \\ -\cos[\alpha](-1 + 3\cos[\alpha])\sin[\alpha]^2 & \cos[\alpha](-1 + 3\cos[\alpha])\sin[\alpha]^2 & \cos[\alpha]^2(-1 + 3\cos[\\ -3\sin[\alpha]^4 & 3\sin[\alpha]^4 & 3\cos[\alpha]\sin \end{pmatrix}$$

```mathematica
In[•]:= (*Hessian Matrix for asymmetric conditions*)

In[•]:= (*asymmetric condition a*)

In[•]:= FF = -r[0] + 2t[0] - r[π] + 2t[π] - EAApBBpy*D[EAApBBpx, α] - EAApCCpy*D[EAApCCpx, α] + AA[[2]]*D[AA[[1]], α];
```


$$\begin{pmatrix} \sin[3\alpha] - 15\sin[4\alpha] & -15\cos[\alpha]\sin[\alpha]^3 & -\cos[\alpha](-1+3\cos[\alpha])\sin[\alpha]^2 & -3\sin[\alpha]^4 \\ 15\cos[3\alpha]\sin[\alpha] & 15\cos[\alpha]\sin[\alpha]^3 & \cos[\alpha](-1+3\cos[\alpha])\sin[\alpha]^2 & 3\sin[\alpha]^4 \\ \alpha]\sin[\alpha]^2 & 3(1+3\cos[2\alpha])\sin[\alpha]^2 & \cos[\alpha]^2(-1+3\cos[\alpha])\sin[\alpha] & 3\cos[\alpha]\sin[\alpha]^3 \\ \sin[\alpha]^2 & -((3+4\cos[\alpha]+9\cos[2\alpha])\sin[\alpha]^2) & 3\cos[\alpha]\sin[\alpha]^3 & 0 \\ \alpha])\sin[\alpha] & 3\cos[\alpha]\sin[\alpha]^3 & 0 & 0 \\ \alpha]^3 & -((1+3\cos[\alpha])\sin[\alpha]^3) & 0 & -((1+3\cos[\alpha])\sin[\alpha]^3) \end{pmatrix}$$

$$\begin{pmatrix} \sin[3\alpha] - 15\sin[4\alpha] & -15\cos[\alpha]\sin[\alpha]^3 & -\cos[\alpha](-1+3\cos[\alpha])\sin[\alpha]^2 & -3\sin[\alpha]^4 \\ 15\cos[3\alpha]\sin[\alpha] & 15\cos[\alpha]\sin[\alpha]^3 & \cos[\alpha](-1+3\cos[\alpha])\sin[\alpha]^2 & 3\sin[\alpha]^4 \\ \alpha]\sin[\alpha]^2 & 3(1+3\cos[2\alpha])\sin[\alpha]^2 & \cos[\alpha]^2(-1+3\cos[\alpha])\sin[\alpha] & 3\cos[\alpha]\sin[\alpha]^3 \\ \sin[\alpha]^2 & -((3+4\cos[\alpha]+9\cos[2\alpha])\sin[\alpha]^2) & 3\cos[\alpha]\sin[\alpha]^3 & 0 \\ \alpha])\sin[\alpha] & 3\cos[\alpha]\sin[\alpha]^3 & 0 & 0 \\ \alpha]^3 & -((1+3\cos[\alpha])\sin[\alpha]^3) & 0 & -((1+3\cos[\alpha])\sin[\alpha]^3) \end{pmatrix}$$

**HH = ResourceFunction["HessianMatrix"][FF, {r[α], t[α], r'[α], t'[α], r''[α], t''[α]}] // Simplify // MatrixForm**

$$\begin{pmatrix} 4\cos[2\alpha]\sin[\alpha]^2 & -((1+4\cos[2\alpha])\sin[\alpha]^2) & \frac{1}{4}(6\sin[2\alpha]-5\sin[4\alpha]) & -10\cos[\alpha]\sin[\alpha]^3 & -2\cos[\alpha]^2\sin[\alpha]^2 \\ -((1+4\cos[2\alpha])\sin[\alpha]^2) & 4\cos[2\alpha]\sin[\alpha]^2 & -\sin[2\alpha]+\frac{5}{4}\sin[4\alpha] & 10\cos[\alpha]\sin[\alpha]^3 & 2\cos[\alpha]^2\sin[\alpha]^2 \\ \frac{1}{4}(6\sin[2\alpha]-5\sin[4\alpha]) & -\sin[2\alpha]+\frac{5}{4}\sin[4\alpha] & -3\sin[2\alpha]^2 & 2(1+3\cos[2\alpha])\sin[\alpha]^2 & 2\cos[\alpha]^3\sin[\alpha] \\ -10\cos[\alpha]\sin[\alpha]^3 & 10\cos[\alpha]\sin[\alpha]^3 & 2(1+3\cos[2\alpha])\sin[\alpha]^2 & -2(1+3\cos[2\alpha])\sin[\alpha]^2 & 2\cos[\alpha]^3\sin[\alpha] \\ -2\cos[\alpha]^2\sin[\alpha]^2 & 2\cos[\alpha]^2\sin[\alpha]^2 & 2\cos[\alpha]^3\sin[\alpha] & 2\cos[\alpha]\sin[\alpha]^3 & 0 \\ -2\sin[\alpha]^4 & 2\sin[\alpha]^4 & 2\cos[\alpha]\sin[\alpha]^3 & -2\cos[\alpha]\sin[\alpha]^3 & 0 \end{pmatrix}$$

(*asymmetric condition b*)

**FF1 = -EAApBBpy * D[EAApBBpx, α] - EAApCCpy * D[EAApCCpx, α]      -r[0] + 2 t[0] - r[π] + 2 t[π];**

**HH1 = ResourceFunction["HessianMatrix"][FF1, {r[α], t[α], r'[α], t'[α], r''[α], t''[α]}] // Simplify // MatrixForm**

$$\begin{pmatrix} 4\cos[2\alpha]\sin[\alpha]^2 & -((1+4\cos[2\alpha])\sin[\alpha]^2) & \frac{1}{4}(6\sin[2\alpha]-5\sin[4\alpha]) & -10\cos[\alpha]\sin[\alpha]^3 & -2\cos[\alpha]^2\sin[\alpha]^2 \\ -4\cos[2\alpha]\sin[\alpha]^2 & 4\cos[2\alpha]\sin[\alpha]^2 & -\sin[2\alpha]+\frac{5}{4}\sin[4\alpha] & 10\cos[\alpha]\sin[\alpha]^3 & 2\cos[\alpha]^2\sin[\alpha]^2 \\ \frac{1}{4}(6\sin[2\alpha]-5\sin[4\alpha]) & \frac{1}{4}(-6\sin[2\alpha]+5\sin[4\alpha]) & -3\sin[2\alpha]^2 & 2(1+3\cos[2\alpha])\sin[\alpha]^2 & 2\cos[\alpha]^3\sin[\alpha] \\ -10\cos[\alpha]\sin[\alpha]^3 & 10\cos[\alpha]\sin[\alpha]^3 & 2(1+3\cos[2\alpha])\sin[\alpha]^2 & -2(1+3\cos[2\alpha])\sin[\alpha]^2 & 2\cos[\alpha]^3\sin[\alpha] \\ -2\cos[\alpha]^2\sin[\alpha]^2 & 2\cos[\alpha]^2\sin[\alpha]^2 & 2\cos[\alpha]^3\sin[\alpha] & 2\cos[\alpha]\sin[\alpha]^3 & 0 \\ -2\sin[\alpha]^4 & 2\sin[\alpha]^4 & 2\cos[\alpha]\sin[\alpha]^3 & -2\cos[\alpha]\sin[\alpha]^3 & 0 \end{pmatrix}$$

**FF2 = -EAApBBpy * D[EAApBBpx, α] - EAApCCpy * D[EAApCCpx, α] + AA[[2]] * D[AA[[1]], α] - r[0] + 2 t[0] - r[π] + 2 t[π];**

**HH2 = ResourceFunction["HessianMatrix"][FF2, {r[α], t[α], r'[α], t'[α], r''[α], t''[α]}] // Simplify // MatrixForm**

$$\begin{pmatrix} 4\cos[2\alpha]\sin[\alpha]^2 & -((1+4\cos[2\alpha])\sin[\alpha]^2) & \frac{1}{4}(6\sin[2\alpha]-5\sin[4\alpha]) & -10\cos[\alpha]\sin[\alpha]^3 & -2\cos[\alpha]^2\sin[\alpha]^2 \\ & 4\cos[2\alpha]\sin[\alpha]^2 & -\sin[2\alpha]+\frac{5}{4}\sin[4\alpha] & 10\cos[\alpha]\sin[\alpha]^3 & 2\cos[\alpha]^2\sin[\alpha]^2 \\ \frac{1}{4}(6\sin[2\alpha]-5\sin[4\alpha]) & -\sin[2\alpha]+\frac{5}{4}\sin[4\alpha] & -3\sin[2\alpha]^2 & 2(1+3\cos[2\alpha])\sin[\alpha]^2 & 2\cos[\alpha]^3\sin[\alpha] \\ -10\cos[\alpha]\sin[\alpha]^3 & 10\cos[\alpha]\sin[\alpha]^3 & 2(1+3\cos[2\alpha])\sin[\alpha]^2 & -2(1+3\cos[2\alpha])\sin[\alpha]^2 & 2\cos[\alpha]^3\sin[\alpha] \\ -2\cos[\alpha]^2\sin[\alpha]^2 & 2\cos[\alpha]^2\sin[\alpha]^2 & 2\cos[\alpha]^3\sin[\alpha] & 2\cos[\alpha]\sin[\alpha]^3 & 0 \\ -2\sin[\alpha]^4 & 2\sin[\alpha]^4 & 2\cos[\alpha]\sin[\alpha]^3 & -2\cos[\alpha]\sin[\alpha]^3 & 0 \end{pmatrix}$$

$$\begin{pmatrix} -2\sin[\alpha]^4 \\ 2\sin[\alpha]^4 \\ \cos[\alpha]\sin[\alpha]^3 \\ 2\cos[\alpha]\sin[\alpha]^3 \\ 0 \\ 0 \end{pmatrix}$$

$$\begin{pmatrix} -2\sin[\alpha]^4 \\ 2\sin[\alpha]^4 \\ 2\cos[\alpha]\sin[\alpha]^3 \\ 2\cos[\alpha]\sin[\alpha]^3 \\ 0 \\ 0 \end{pmatrix}$$

$$\begin{pmatrix} -2\sin[\alpha]^4 \\ 2\sin[\alpha]^4 \\ \cos[\alpha]\sin[\alpha]^3 \\ 2\cos[\alpha]\sin[\alpha]^3 \\ 0 \\ 0 \end{pmatrix}$$

**FF3** = **-EAApBBpy * D[EAApBBpx, α] - EAApCCpy * D[EAApCCpx, α] + AA[[2]] * D[EAApCCpx, α], {r[α], t[α], r'[α], t'[α], r''[α], t''[α]}] // Simplify // MatrixForm**

**HH3** = **ResourceFunction["HessianMatrix"][FF3, {r[α], t[α], r'[α], t'[α], r''[α], t''[α]}] // Simplify // MatrixForm**

$$\begin{pmatrix} 4\cos[2\alpha]\sin[\alpha]^2 & -((1+4\cos[2\alpha])\sin[\alpha]^2) & \tfrac{1}{4}(6\sin[2\alpha]-5\sin[4\alpha]) & -10\cos[\alpha]\sin[\alpha]^3 & -2\cos[\alpha]^2\sin[\alpha]^2 \\ -((1+4\cos[2\alpha])\sin[\alpha]^2) & 4\cos[2\alpha]\sin[\alpha]^2 & -\sin[2\alpha]+\tfrac{5}{4}\sin[4\alpha] & 10\cos[\alpha]\sin[\alpha]^3 & 2\cos[\alpha]^2\sin[\alpha]^2 & 2 \\ \tfrac{1}{4}(6\sin[2\alpha]-5\sin[4\alpha]) & -\sin[2\alpha]+\tfrac{5}{4}\sin[4\alpha] & -3\sin[2\alpha]^2 & 2(1+3\cos[2\alpha])\sin[\alpha]^2 & 2\cos[\alpha]^3\sin[\alpha] \\ -10\cos[\alpha]\sin[\alpha]^3 & 10\cos[\alpha]\sin[\alpha]^3 & 2(1+3\cos[2\alpha])\sin[\alpha]^2 & -2(1+3\cos[2\alpha])\sin[\alpha]^2 & 2\cos[\alpha]\sin[\alpha]^3 & 0 \\ -2\cos[\alpha]^2\sin[\alpha]^2 & 2\cos[\alpha]^2\sin[\alpha]^2 & 2\cos[\alpha]^3\sin[\alpha] & 2\cos[\alpha]\sin[\alpha]^3 & 0 \\ -2\sin[\alpha]^4 \end{pmatrix}$$

**FF4** = **-EAApBBpy * D[EAApBBpx, α] - EAApCCpy * D[EAApCCpx, α] + AA[[2]] * D[AA[[1]], α] - EAACCy * D[EAACCx, α];**

**HH4** = **ResourceFunction["HessianMatrix"][FF4, {r[α], t[α], r'[α], t'[α], r''[α], t''[α]}] // Simplify // MatrixForm**

$$\begin{pmatrix} 2(\cos[\alpha]+3\cos[2\alpha])\sin[\alpha]^2 & -((1+2\cos[\alpha]+6\cos[2\alpha])\sin[\alpha]^2) & \tfrac{1}{8}(10\sin[\alpha]+18\sin[2\alpha]-6\sin \\ -((1+2\cos[\alpha]+6\cos[2\alpha])\sin[\alpha]^2) & 2(\cos[\alpha]+3\cos[2\alpha])\sin[\alpha]^2 & \tfrac{1}{4}(-2+\cos[\alpha]+6\cos[2\alpha]+15 \\ \tfrac{1}{8}(10\sin[\alpha]+18\sin[2\alpha]-6\sin[3\alpha]-15\sin[4\alpha]) & \tfrac{1}{4}(-2+\cos[\alpha]+6\cos[2\alpha]+15\cos[3\alpha])\sin[\alpha] & -2\cos[\alpha](2+9\cos[\alpha] \\ -15\cos[\alpha]\sin[\alpha]^3 & 15\cos[\alpha]\sin[\alpha]^3 & 3(1+3\cos[2\alpha])\text{Si} \\ -\cos[\alpha](1+3\cos[\alpha])\sin[\alpha]^2 & \cos[\alpha](1+3\cos[\alpha])\sin[\alpha]^2 & \cos[\alpha]^2(1+3\cos[\alpha]) \\ -3\sin[\alpha]^4 & 3\sin[\alpha]^4 & 3\cos[\alpha]\sin[\alpha] \end{pmatrix}$$

(*asymmetric condition c*)

**FF1** =

**-EAApBBpy * D[EAApBBpx, α] - EAApCCpy * D[EAApCCpx, α]**        **-EAABBy * D[EAABBx, α] - r[0] + 2 t[0] - r[π] + 2 t[π];**

$$-2 \sin[\alpha]^4$$
$$2 \sin[\alpha]^4$$
$$\cos[\alpha] \sin[\alpha]^3$$
$$2 \cos[\alpha] \sin[\alpha]^3$$
$$0$$
$$0$$

$$[3\alpha] - 15 \sin[4\alpha]$$
$$\cos[3\alpha]) \sin[\alpha]$$
$$) \sin[\alpha]^2$$
$$\sin[\alpha]$$
$$) \sin[\alpha]$$
$$]^3$$

$$-15 \cos[\alpha] \sin[\alpha]^3$$
$$15 \cos[\alpha] \sin[\alpha]^3$$
$$3 (1 + 3 \cos[2\alpha]) \sin[\alpha]^2$$
$$(-3 + 4 \cos[\alpha] - 9 \cos[2\alpha]) \sin[\alpha]^2$$
$$3 \cos[\alpha] \sin[\alpha]^3$$
$$-((-1 + 3 \cos[\alpha]) \sin[\alpha]^3)$$

$$-\cos[\alpha] (1 + 3 \cos[\alpha]) \sin[\alpha]^2$$
$$\cos[\alpha] (1 + 3 \cos[\alpha]) \sin[\alpha]^2$$
$$\cos[\alpha]^2 (1 + 3 \cos[\alpha]) \sin[\alpha]$$
$$3 \cos[\alpha] \sin[\alpha]^3$$
$$0$$
$$0$$

$$-3 \sin[\alpha]^4$$
$$3 \sin[\alpha]^4$$
$$3 \cos[\alpha] \sin[\alpha]^3$$
$$-((-1 + 3 \cos[\alpha]) \sin[\alpha]^3)$$
$$0$$
$$0$$

**HH1 = ResourceFunction["HessianMatrix"][FF1, {r[α], t[α], r'[α], t'[α], r''[α], t''[α]}] // Simplify // MatrixForm**

$$\begin{pmatrix} -2(\cos[\alpha]-3\cos[2\alpha])\sin[\alpha]^2 & 2(\cos[\alpha]-3\cos[2\alpha])\sin[\alpha]^2 & \frac{1}{8}(-10\sin[\alpha]+18\sin[2\alpha]+6\ldots \\ 2(\cos[\alpha]-3\cos[2\alpha])\sin[\alpha]^2 & -2(\cos[\alpha]-3\cos[2\alpha])\sin[\alpha]^2 & \frac{1}{4}(2-3\cos[\alpha]-6\cos[2\alpha]+ \\ \frac{1}{8}(-10\sin[\alpha]+18\sin[2\alpha]+6\sin[3\alpha]-15\sin[4\alpha]) & \frac{1}{4}(2-3\cos[\alpha]-6\cos[2\alpha]+15\cos[3\alpha])\sin[\alpha] & -2\cos[\alpha](-2+9\cos \\ -15\cos[\alpha]\sin[\alpha]^3 & 15\cos[\alpha]\sin[\alpha]^3 & 3(1+3\cos[2\alpha]) \\ -\cos[\alpha](-1+3\cos[\alpha])\sin[\alpha]^2 & \cos[\alpha](-1+3\cos[\alpha])\sin[\alpha]^2 & \cos[\alpha]^2(-1+3\cos[\ldots \\ -3\sin[\alpha]^4 & 3\sin[\alpha]^4 & 3\cos[\alpha]\sin \end{pmatrix}$$

**FF2 =**

**− EAApBBpy * D[EAApBBpx, α] − EAApCCpy * D[EAApCCpx, α]    − r[0] + 2 t[0] − r[π] + 2 t[π];**

**HH2 = ResourceFunction["HessianMatrix"][FF2, {r[α], t[α], r'[α], t'[α], r''[α], t''[α]}] // Simplify // MatrixForm**

$$\begin{pmatrix} 4\cos[2\alpha]\sin[\alpha]^2 & -4\cos[2\alpha]\sin[\alpha]^2 & \frac{1}{4}(6\sin[2\alpha]-5\sin[4\alpha]) & -10\cos[\alpha]\sin[\alpha]^3 & -2\cos[\alpha]^2\sin[\alpha]^2 \\ -4\cos[2\alpha]\sin[\alpha]^2 & 4\cos[2\alpha]\sin[\alpha]^2 & \frac{1}{4}(-6\sin[2\alpha]+5\sin[4\alpha]) & 10\cos[\alpha]\sin[\alpha]^3 & 2\cos[\alpha]^2\sin[\alpha]^2 \\ \frac{1}{4}(6\sin[2\alpha]-5\sin[4\alpha]) & \frac{1}{4}(-6\sin[2\alpha]+5\sin[4\alpha]) & -3\sin[2\alpha]^2 & 2(1+3\cos[2\alpha])\sin[\alpha]^2 & 2\cos[\alpha]^3\sin[\alpha] \quad 2 \\ -10\cos[\alpha]\sin[\alpha]^3 & 10\cos[\alpha]\sin[\alpha]^3 & 2(1+3\cos[2\alpha])\sin[\alpha]^2 & -2(1+3\cos[2\alpha])\sin[\alpha]^2 & 2\cos[\alpha]\sin[\alpha]^3 \quad - \\ -2\cos[\alpha]^2\sin[\alpha]^2 & 2\cos[\alpha]^2\sin[\alpha]^2 & 2\cos[\alpha]^3\sin[\alpha] & 2\cos[\alpha]\sin[\alpha]^3 & 0 \\ -2\sin[\alpha]^4 & 2\sin[\alpha]^4 & 2\cos[\alpha]\sin[\alpha]^3 & -2\cos[\alpha]\sin[\alpha]^3 & 0 \end{pmatrix}$$

**FF3 =**

**− EAApBBpy * D[EAApBBpx, α] − EAApCCpy * D[EAApCCpx, α] + AA[[2]] * D[AA[[1]], α]    − r[0] + 2 t[0] − r[π] + 2 t[π];**

**HH3 = ResourceFunction["HessianMatrix"][FF3, {r[α], t[α], r'[α], t'[α], r''[α], t''[α]}] // Simplify // MatrixForm**

$$\begin{pmatrix} 4\cos[2\alpha]\sin[\alpha]^2 & -((1+4\cos[2\alpha])\sin[\alpha]^2) & \frac{1}{4}(6\sin[2\alpha]-5\sin[4\alpha]) & -10\cos[\alpha]\sin[\alpha]^3 & -2\cos[\alpha]^2\sin[\alpha]^2 \\ -((1+4\cos[2\alpha])\sin[\alpha]^2) & 4\cos[2\alpha]\sin[\alpha]^2 & -\sin[2\alpha]+\frac{5}{4}\sin[4\alpha] & 10\cos[\alpha]\sin[\alpha]^3 & 2\cos[\alpha]^2\sin[\alpha]^2 \\ \frac{1}{4}(6\sin[2\alpha]-5\sin[4\alpha]) & -\sin[2\alpha]+\frac{5}{4}\sin[4\alpha] & -3\sin[2\alpha]^2 & 2(1+3\cos[2\alpha])\sin[\alpha]^2 & 2\cos[\alpha]^3\sin[\alpha] \quad 2 \\ -10\cos[\alpha]\sin[\alpha]^3 & 10\cos[\alpha]\sin[\alpha]^3 & 2(1+3\cos[2\alpha])\sin[\alpha]^2 & -2(1+3\cos[2\alpha])\sin[\alpha]^2 & 2\cos[\alpha]\sin[\alpha]^3 \quad -; \\ -2\cos[\alpha]^2\sin[\alpha]^2 & 2\cos[\alpha]^2\sin[\alpha]^2 & 2\cos[\alpha]^3\sin[\alpha] & 2\cos[\alpha]\sin[\alpha]^3 & 0 \\ -2\sin[\alpha]^4 & 2\sin[\alpha]^4 & 2\cos[\alpha]\sin[\alpha]^3 & -2\cos[\alpha]\sin[\alpha]^3 & 0 \end{pmatrix}$$

$$\begin{pmatrix}
\sin[3\alpha] - 15\sin[4\alpha] & -15\cos[\alpha]\sin[\alpha]^3 & -\cos[\alpha](-1+3\cos[\alpha])\sin[\alpha]^2 & -3\sin[\alpha]^4 \\
15\cos[3\alpha]\sin[\alpha] & 15\cos[\alpha]\sin[\alpha]^3 & \cos[\alpha](-1+3\cos[\alpha])\sin[\alpha]^2 & 3\sin[\alpha]^4 \\
[\alpha])\sin[\alpha]^2 & 3(1+3\cos[2\alpha])\sin[\alpha]^2 & \cos[\alpha]^2(-1+3\cos[\alpha])\sin[\alpha] & 3\cos[\alpha]\sin[\alpha]^3 \\
\sin[\alpha]^2 & -((3+4\cos[\alpha]+9\cos[2\alpha])\sin[\alpha]^2) & 3\cos[\alpha]\sin[\alpha]^3 & -((1+3\cos[\alpha])\sin[\alpha]^3) \\
\alpha])\sin[\alpha] & 3\cos[\alpha]\sin[\alpha]^3 & 0 & 0 \\
[\alpha]^3 & -((1+3\cos[\alpha])\sin[\alpha]^3) & 0 & 0
\end{pmatrix}$$

$$\begin{pmatrix}
-2\sin[\alpha]^4 \\
2\sin[\alpha]^4 \\
\cos[\alpha]\sin[\alpha]^3 \\
2\cos[\alpha]\sin[\alpha]^3 \\
0 \\
0
\end{pmatrix}$$

$$\begin{pmatrix}
-2\sin[\alpha]^4 \\
2\sin[\alpha]^4 \\
\cos[\alpha]\sin[\alpha]^3 \\
2\cos[\alpha]\sin[\alpha]^3 \\
0 \\
0
\end{pmatrix}$$

**FF4 = -EAApBBpy * D[EAApBBpx, α] - EAApCCpy * D[EAApCCpx, α] + AA[[2]] * D[EAApCCpx, α, {r[α], t[α], r'[α], t'[α], r''[α], t''[α]}] // Simplify // MatrixForm;**

**HH4 = ResourceFunction["HessianMatrix"][FF4, {r[α], t[α], r'[α], t'[α], r''[α], t''[α]}] // Simplify // MatrixForm**

$$\begin{pmatrix} 4\cos[2\alpha]\sin[\alpha]^2 & -((1+4\cos[2\alpha])\sin[\alpha]^2) & \tfrac{1}{4}(6\sin[2\alpha]-5\sin[4\alpha]) & -10\cos[\alpha]\sin[\alpha]^3 & -2\cos[\alpha]^2\sin[\alpha]^2 \\ -((1+4\cos[2\alpha])\sin[\alpha]^2) & 4\cos[2\alpha]\sin[\alpha]^2 & -\sin[2\alpha]+\tfrac{5}{4}\sin[4\alpha] & 10\cos[\alpha]\sin[\alpha]^3 & 2\cos[\alpha]^2\sin[\alpha]^2 \\ \tfrac{1}{4}(6\sin[2\alpha]-5\sin[4\alpha]) & -\sin[2\alpha]+\tfrac{5}{4}\sin[4\alpha] & -3\sin[2\alpha]^2 & 2(1+3\cos[2\alpha])\sin[\alpha]^2 & 2\cos[\alpha]^3\sin[\alpha] \\ -10\cos[\alpha]\sin[\alpha]^3 & 10\cos[\alpha]\sin[\alpha]^3 & 2(1+3\cos[2\alpha])\sin[\alpha]^2 & -2(1+3\cos[2\alpha])\sin[\alpha]^2 & 2\cos[\alpha]^3\sin[\alpha] \\ -2\cos[\alpha]^2\sin[\alpha]^2 & 2\cos[\alpha]^2\sin[\alpha]^2 & 2\cos[\alpha]^3\sin[\alpha] & 2\cos[\alpha]\sin[\alpha]^3 & 0 \\ -2\sin[\alpha]^4 & & & & \end{pmatrix}$$

**FF5 = -EAApBBpy * D[EAApBBpx, α] - EAApCCpy * D[EAApCCpx, α] + AA[[2]] * D[AA[[1]], α] - EAACCy * D[EAACCx, α];**

**HH5 = ResourceFunction["HessianMatrix"][FF5, {r[α], t[α], r'[α], t'[α], r''[α], t''[α]}] // Simplify // MatrixForm**

$$\begin{pmatrix} 2(\cos[\alpha]+3\cos[2\alpha])\sin[\alpha]^2 & -((1+2\cos[\alpha]+6\cos[2\alpha])\sin[\alpha]^2) & -((1+2\cos[\alpha]+6\cos[2\alpha])\sin[\alpha]^2) & \tfrac{1}{8}(10\sin[\alpha]+18\sin[2\alpha]-6\sin[\alpha] \\ -((1+2\cos[\alpha]+6\cos[2\alpha])\sin[\alpha]^2) & 2(\cos[\alpha]+3\cos[2\alpha])\sin[\alpha]^2 & 2(\cos[\alpha]+3\cos[2\alpha])\sin[\alpha]^2 & \tfrac{1}{4}(-2+\cos[\alpha]+6\cos[2\alpha]+15 \\ \tfrac{1}{8}(10\sin[\alpha]+18\sin[2\alpha]-6\sin[3\alpha]-15\sin[4\alpha]) & \tfrac{1}{4}(-2+\cos[\alpha]+6\cos[2\alpha]+15\cos[3\alpha])\sin[\alpha] & -2\cos[\alpha](2+9\cos[\alpha] \\ -15\cos[\alpha]\sin[\alpha]^3 & 15\cos[\alpha]\sin[\alpha]^3 & 3(1+3\cos[2\alpha])\sin \\ -\cos[\alpha](1+3\cos[\alpha])\sin[\alpha]^2 & \cos[\alpha](1+3\cos[\alpha])\sin[\alpha]^2 & \cos[\alpha]^2(1+3\cos[\alpha] \\ -3\sin[\alpha]^4 & 3\sin[\alpha]^4 & 3\cos[\alpha]\sin[\alpha] \end{pmatrix}$$

**FF6 = -EAApBBpy * D[EAApBBpx, α] - EAApCCpy * D[EAApCCpx, α] - EAACCy * D[EAACCx, α];**

**HH6 = ResourceFunction["HessianMatrix"][FF6, {r[α], t[α], r'[α], t'[α], r''[α], t''[α]}] // Simplify // MatrixForm**

$$\begin{pmatrix} 2(\cos[\alpha]+3\cos[2\alpha])\sin[\alpha]^2 & -2(\cos[\alpha]+3\cos[2\alpha])\sin[\alpha]^2 & \tfrac{1}{8}(10\sin[\alpha]+18\sin[2\alpha]-6\sin[\alpha] \\ -2(\cos[\alpha]+3\cos[2\alpha])\sin[\alpha]^2 & 2(\cos[\alpha]+3\cos[2\alpha])\sin[\alpha]^2 & \tfrac{1}{4}(-2-3\cos[\alpha]+6\cos[2\alpha])+ \\ \tfrac{1}{8}(10\sin[\alpha]+18\sin[2\alpha]-6\sin[3\alpha]-15\sin[4\alpha]) & \tfrac{1}{4}(-2-3\cos[\alpha]+6\cos[2\alpha]+15\cos[3\alpha])\sin[\alpha] & -2\cos[\alpha](2+9\cos[\alpha]) \\ -15\cos[\alpha]\sin[\alpha]^3 & 15\cos[\alpha]\sin[\alpha]^3 & 3(1+3\cos[2\alpha]) \\ -\cos[\alpha](1+3\cos[\alpha])\sin[\alpha]^2 & \cos[\alpha](1+3\cos[\alpha])\sin[\alpha]^2 & \cos[\alpha]^2(1+3\cos[\alpha] \\ -3\sin[\alpha]^4 & 3\sin[\alpha]^4 & 3\cos[\alpha]\sin[\alpha] \end{pmatrix}$$

$$-2\operatorname{Sin}[\alpha]^4$$
$$2\operatorname{Sin}[\alpha]^4$$
$$\operatorname{Cos}[\alpha]\operatorname{Sin}[\alpha]^3$$
$$2\operatorname{Cos}[\alpha]\operatorname{Sin}[\alpha]^3$$
$$0$$
$$0$$

$$[3\alpha]-15\operatorname{Sin}[4\alpha]$$
$$\operatorname{Cos}[3\alpha])\operatorname{Sin}[\alpha]$$
$$)\operatorname{Sin}[\alpha]^2$$
$$\operatorname{in}[\alpha]^2$$
$$)\operatorname{Sin}[\alpha]$$
$$]^3$$

$$-15\operatorname{Cos}[\alpha]\operatorname{Sin}[\alpha]^3$$
$$15\operatorname{Cos}[\alpha]\operatorname{Sin}[\alpha]^3$$
$$3(1+3\operatorname{Cos}[2\alpha])\operatorname{Sin}[\alpha]^2$$
$$(-3+4\operatorname{Cos}[\alpha]-9\operatorname{Cos}[2\alpha])\operatorname{Sin}[\alpha]^2$$
$$3\operatorname{Cos}[\alpha]\operatorname{Sin}[\alpha]^3$$
$$-((-1+3\operatorname{Cos}[\alpha])\operatorname{Sin}[\alpha]^3)$$

$$-\operatorname{Cos}[\alpha](1+3\operatorname{Cos}[\alpha])\operatorname{Sin}[\alpha]^2$$
$$\operatorname{Cos}[\alpha](1+3\operatorname{Cos}[\alpha])\operatorname{Sin}[\alpha]^2$$
$$\operatorname{Cos}[\alpha]^2(1+3\operatorname{Cos}[\alpha])\operatorname{Sin}[\alpha]$$
$$3\operatorname{Cos}[\alpha]\operatorname{Sin}[\alpha]^3$$
$$0$$
$$0$$

$$-3\operatorname{Sin}[\alpha]^4$$
$$3\operatorname{Sin}[\alpha]^4$$
$$3\operatorname{Cos}[\alpha]\operatorname{Sin}[\alpha]^3$$
$$-((-1+3\operatorname{Cos}[\alpha])\operatorname{Sin}[\alpha]^3)$$
$$0$$
$$0$$

$$\operatorname{in}[3\alpha]-15\operatorname{Sin}[4\alpha])$$
$$15\operatorname{Cos}[3\alpha])\operatorname{Sin}[\alpha]$$
$$\alpha])\operatorname{Sin}[\alpha]^2$$
$$\operatorname{Sin}[\alpha]^2$$
$$])\operatorname{Sin}[\alpha]$$
$$\alpha]^3$$

$$-15\operatorname{Cos}[\alpha]\operatorname{Sin}[\alpha]^3$$
$$15\operatorname{Cos}[\alpha]\operatorname{Sin}[\alpha]^3$$
$$3(1+3\operatorname{Cos}[2\alpha])\operatorname{Sin}[\alpha]^2$$
$$(-3+4\operatorname{Cos}[\alpha]-9\operatorname{Cos}[2\alpha])\operatorname{Sin}[\alpha]^2$$
$$3\operatorname{Cos}[\alpha]\operatorname{Sin}[\alpha]^3$$
$$-((-1+3\operatorname{Cos}[\alpha])\operatorname{Sin}[\alpha]^3)$$

$$-\operatorname{Cos}[\alpha](1+3\operatorname{Cos}[\alpha])\operatorname{Sin}[\alpha]^2$$
$$\operatorname{Cos}[\alpha](1+3\operatorname{Cos}[\alpha])\operatorname{Sin}[\alpha]^2$$
$$\operatorname{Cos}[\alpha]^2(1+3\operatorname{Cos}[\alpha])\operatorname{Sin}[\alpha]$$
$$3\operatorname{Cos}[\alpha]\operatorname{Sin}[\alpha]^3$$
$$0$$
$$0$$

$$-3\operatorname{Sin}[\alpha]^4$$
$$3\operatorname{Sin}[\alpha]^4$$
$$3\operatorname{Cos}[\alpha]\operatorname{Sin}[\alpha]^3$$
$$-((-1+3\operatorname{Cos}[\alpha])\operatorname{Sin}[\alpha]^3)$$
$$0$$
$$0$$


\end{document}